\documentclass[english]{amsart}
\usepackage[english]{babel}

\usepackage[utf8]{inputenc}
\usepackage[T1]{fontenc}
\usepackage{amsmath}
\usepackage{amsthm}
\usepackage{mathtools}
\usepackage{amsmath}
\usepackage{amssymb}
\usepackage{mathrsfs}
\usepackage{multirow}
\usepackage{stmaryrd}
\SetSymbolFont{stmry}{bold}{U}{stmry}{m}{n}
\usepackage{array}
\usepackage{hyperref}
\usepackage{thmtools}
\usepackage{esint}
\usepackage[all]{xy}
\usepackage{pgf,tikz}
\usepackage{mathrsfs}
\usetikzlibrary{arrows}
\usepackage{enumitem}
\usepackage{xcolor}
\usepackage{bm}
\usepackage{indentfirst}

\newcommand{\A}{\mathrm{A}}
\newcommand{\B}{\mathrm{B}}
\newcommand{\dd}{\mathrm{d}}
\newcommand{\C}{\mathrm{C}}
\newcommand{\D}{\mathrm{D}}
\newcommand{\E}{\mathrm{E}}
\newcommand{\F}{\mathrm{F}}
\let\H\undefined
\newcommand{\H}{\mathrm{H}}
\let\L\undefined
\newcommand{\L}{\mathrm{L}}
\newcommand{\loc}{\mathrm{loc}}

\newcommand{\N}{\mathbb{N}}

\newcommand{\R}{\mathbb{R}}
\newcommand{\T}{\mathbb{T}}
\newcommand{\Z}{\mathbb{Z}}

\newcommand{\CCC}{\mathscr{C}}

\newcommand{\EEE}{\mathscr{E}}

\newcommand{\ddv}{\dd{v}}

\newcommand{\eps}{\varepsilon}
\newcommand{\ueps}{u_{\varepsilon}}
\newcommand{\feps}{f_{\varepsilon}}
\newcommand{\jeps}{j_{\varepsilon}}
\newcommand{\rhoeps}{\rho_{\varepsilon}}

\newcommand{\chieps}{\chi_\eps}

\newcommand{\uegd}{u_{\eps,\gamma,\sigma}}
\newcommand{\fegd}{f_{\eps,\gamma,\sigma}}
\newcommand{\jegd}{j_{\eps,\gamma,\sigma}}
\newcommand{\rhoegd}{\rho_{\eps,\gamma,\sigma}}
\newcommand{\Fegd}{F_{\eps,\gamma,\sigma}}

\newcommand{\egd}{{\eps,\gamma,\sigma}}
\newcommand{\ed}{{\eps,\sigma}}

\newcommand{\Hed}{\mathscr{H}_{\eps,\sigma}}
\newcommand{\He}{\mathscr{H}_{\eps}}
\newcommand{\ued}{u_{\eps,\sigma}}
\newcommand{\fed}{f_{\eps,\sigma}}
\newcommand{\jed}{j_{\eps,\sigma}}
\newcommand{\rhoed}{\rho_{\eps,\sigma}}
\newcommand{\Fed}{F_{\eps,\sigma}}
\newcommand{\wed}{w_{\eps,\sigma}}

\let\div\undefined
\DeclareMathOperator{\div}{div}

\DeclareMathOperator{\Id}{Id}

\newcommand{\na}{\nabla}

\newcommand{\norme}[1]{\left\Vert #1\right\Vert}

\renewcommand{\tilde}{\widetilde}

\let\1\undefined
\newcommand{\1}{\mathbf{1}}

\newcommand{\pa}{\partial}

\definecolor{gris}{RGB}{175,175,175}

\newtheorem{theorem}{Theorem}[section]
\newtheorem{definition}[theorem]{Definition}
\newtheorem{proposition}[theorem]{Proposition}
\newtheorem{corollary}[theorem]{Corollary}
\newtheorem{lemma}[theorem]{Lemma}
\newtheorem{remark}[theorem]{Remark}

\newtheorem{assumption}{Assumption}

\declaretheoremstyle[
    spaceabove=-6pt, 
    spacebelow=6pt, 
    headfont=\normalfont\bfseries, 
    bodyfont = \normalfont,
    postheadspace=1em, 
    qed=$\Box$, 
    headpunct={$\rhd$}
]{stylepreuve} 
\declaretheorem[name={}, style=stylepreuve, unnumbered]{preuve}

\numberwithin{equation}{section}

\usepackage{mathdots}

\makeindex

\title{On hydrodynamic limits of the Vlasov-Navier-Stokes system}
\author{Daniel Han-Kwan}
\address{Centre de Math\'ematiques Laurent Schwartz (UMR 7640), Ecole Polytechnique, Institut Polytechnique de Paris, 91128 Palaiseau Cedex, France} \email{daniel.han-kwan@polytechnique.edu}
\author{David Michel}
\address{Sorbonne Université and Université de Paris, CNRS, Laboratoire Jacques-Louis Lions (LJLL), F-75005 Paris, France}
\email{david.michel@ljll.math.upmc.fr}

\begin{document}

\everymath{\displaystyle}

\setlist[itemize]{label=\textbullet}

\begin{abstract}
    We introduce a framework to justify hydrodynamic limits of the Vlasov-Navier-Stokes system. We specifically study high friction regimes, which take into account the fact that particles of the dispersed phase are  light (resp. small) compared to the fluid part, and lead to the derivation of Transport-Navier-Stokes (resp. Inhomogeneous Navier-Stokes) systems.
\end{abstract}

\maketitle

\tableofcontents

\section{Introduction}

This work is concerned with the mathematical description of the motion of a dispersed phase of small particles (e.g. a spray or an aerosol) flowing in a surrounding incompressible homogeneous fluid. In the class of the so-called \textit{fluid-kinetic} models (see \cite{wil}, \cite{oro_phd}), the cloud of particles (resp. the fluid) is described by its distribution function $f$ in the phase space (resp. its velocity $u$ and its pressure $p$). We assume that the particles do not directly act on one another and that their interaction with the fluid can be represented by a drag acceleration of the particles and, conversely, a drag force applied to the fluid, called the Brinkman force.  This leads to the Vlasov-Navier-Stokes system:
\begin{equation}
\label{VNS}
\left\{
\begin{aligned}
    &\partial_tu+(u\cdot\na_{x})u-\Delta_{x}u+\na_{x}p=j_f-\rho_fu,\\
    &\div_{x}u=0,\\
    &\partial_tf+v\cdot\na_{x}f+\div_{v}\left[f(u-v)\right]=0,
\end{aligned}
\right.
\end{equation}
where
\begin{align*}
    \rho_f(t,x)&:=\int_{\R^3}f(t,x,v)\ddv, \\
    j_f(t,x)&:=\int_{\R^3}vf(t,x,v)\ddv.
\end{align*}

In this work, we study~\eqref{VNS} for $(x,v) \in \T^3 \times \R^3$, where $\T^3 := \R^3/(2\pi\Z)^3$. The flat torus $\T^3$ is equipped with the normalized Lebesgue measure, so that $\text{Leb}(\T^3) = 1$. The equations~\eqref{VNS} are endowed with the initial conditions
\begin{align*}
&u|_{t=0}(x) = u^0(x), \\
&f|_{t=0}(x,v) = f^0(x,v).
\end{align*}
The Vlasov-Navier-Stokes system can be (at least formally) derived as a \emph{mean-field} limit of a N-particle system interacting with a fluid. For some (partial) rigorous results, we refer to \cite{all,des-gol-ric,hil,hil-mou-sue,car-hil,giu-hof,hofe-cmp,mec}, though the full justification is still an outstanding open problem.

In another direction, \cite{ber-des-gol-ric} proposed a program (similar to the Bardos-Golse-Levermore program \cite{bar-gol-lev} for the hydrodynamic limits of the Boltzmann equation) to derive \eqref{VNS} from multiphase Boltzmann equations describing a gas mixture. Again, the  complete rigorous justification is an open problem.

We also refer the reader to \cite{bou-gra-lor-mou} for a discussion of several physical extensions and applications of the model.

The Vlasov-Navier-Stokes system is endowed with a remarkable energy--dissipation structure, on which most of the recent mathematical theories build. Introducing the energy
\[
    \E(t) =   \frac{1}{2}\int_{\T^3\times\R^3}|v|^2 f(t,x,v)\dd x\dd v +   \frac{1}{2}\int_{\T^3}|u(t,x)|^2\dd x,
\]
and the dissipation
\[
    \D(t)=\int_{\T^3}|\na_xu (t,x)|^2\dd x+\int_{\T^3 \times \R^3} |v-u(t,x)|^2 f(t,x,v) \dd x\dd v,
\]
we have (at least formally) the identity
\begin{equation}
\label{energy-dissipation-intro}
\frac{\dd}{\dd t} \E + \D =0.
\end{equation}

In this work we specifically aim at studying {\bf high friction} regimes of the Vlasov-Navier-Stokes system and at justifying its approximation by reduced purely hydrodynamic equations of Navier-Stokes type.

As a concrete illustration, one regime we consider corresponds to
\begin{equation*}
    \left\{
\begin{aligned}
    &\partial_tu+ (u\cdot\na_{x})u-\Delta_{x}u+\na_{x}p= \frac{1}{\eps}\left(j_f-\rho_fu\right),\\
  &\div_{x}u=0,\\
   &\partial_tf+v\cdot\na_{x}f+ \frac{1}{\eps} \div_{v}\left[f( u-v)\right]=0, \\
    &  \rho_f(t,x)=\int_{\R^3}f(t,x,v)\ddv,  \quad
  j_f(t,x)=\int_{\R^3}vf(t,x,v)\ddv,
\end{aligned}
\right.
\end{equation*}
for $\eps \ll 1$, which, at least formally, leads to the {\bf Inhomogeneous Incompressible Navier-Stokes} equations
\begin{equation*}
    \left\{
\begin{aligned}
 &\partial_t \rho + \div_x (\rho u) =0, \\
 &\partial_t ( (1+\rho) u) + \na_x ((1+\rho) u\otimes u) -\Delta_x u +  \nabla_x p =0, \\
&\div_x u = 0.
\end{aligned}
\right.
\end{equation*}
We also study other high friction regimes, which we introduce in the next section. A goal of this work is to develop robust methods to deal with all limits in a general framework.

\subsection{High friction scalings of the Vlasov-Navier-Stokes system}
\label{sec-scalings}

In this work we study two classes of high friction limits for the Vlasov-Navier-Stokes system, referred to as the {\bf light} particles  and the {\bf fine} particles regimes.
To explain their physical relevancy, following \cite{gou-jab-vas04a} (see also the work of Caflisch and Papanicolaou \cite{caf-pap} in which several scalings are discussed in view of experimental parameter values), we write \eqref{VNS} in the following dimensionless form: 
\begin{equation*}
    \left\{
\begin{aligned}
    &\partial_tu+ (u\cdot\na_{x})u-\F \Delta_{x}u+\na_{x}p=\C \left(j_f-\rho_fu\right),\\
   &\div_{x}u=0,\\
   &\partial_tf+ \A v\cdot\na_{x}f+\B \div_{v}\left[f\left(\frac{1}{\A}u-v\right)\right]=0, \\
   &\rho_f(t,x)=\int_{\R^3}f(t,x,v)\ddv, \quad j_f(t,x)={\A} \int_{\R^3}vf(t,x,v)\ddv,
\end{aligned}
\right.
\end{equation*}
with
\begin{itemize}

\item $\A= \frac{\sqrt{\theta}}{U}$, where  $\sqrt{\theta}$ (resp. $U$) is the typical thermal velocity of the particles (resp. velocity of the fluid),  and $ U= \frac{L}{T}$, where $L$ (resp. $T$) stands for the characteristic observation
length (resp. observation time);

\item $\B= \frac{T}{\tau}$, where $\tau$ is  the so-called Stokes relaxation time, given by the formula
$$
\tau = \frac{2 a^2 \rho_{\text{part}}}{9 \mu},
$$
in which $a$ (resp. $\rho_{\text{part}}$) stands for the renormalized\footnote{Very loosely speaking, the Vlasov-Navier-Stokes system comes from the $N\to +\infty$ limit of a $N$-particle system, in which the radius $r$ of a particle also tends to $0$ as N$ \to +\infty$, namely $r\sim \frac{1}{N} a$.} radius (resp. the mass density) of a particle and $\mu$ is the diffusivity constant of the fluid;

\item $\C =\frac{T}{\tau} \frac{\rho_{\text{part}}}{\rho_{\text{fluid}}}= \B\frac{\rho_{\text{part}}}{\rho_{\text{fluid}}} $, where  $\rho_{\text{fluid}}$ is for the density of the fluid;

\item $\F = \frac{2}{9}\left( \frac{a}{L}\right)^2 \frac{T}{\tau}  \frac{\rho_{\text{part}}}{\rho_{\text{fluid}}}= \frac{2}{9} \left( \frac{a}{L}\right)^2 \C$.

\end{itemize}

We shall study in this work the  following classes of regimes.

\subsubsection{\emph{Light} and \emph{light and fast} particle regimes} The {\bf light particle} regime corresponds to
\begin{equation*}
\A= 1, \quad \B= \frac{1}{\eps}, \quad \C=1, \quad \F=1,
\end{equation*}
leading to the system
\begin{equation*}
    \left\{
\begin{aligned}
    &\partial_tu+ (u\cdot\na_{x})u-\Delta_{x}u+\na_{x}p= j_f-\rho_fu,\\
  &\div_{x}u=0,\\
    &\partial_tf+  v\cdot\na_{x}f+ \frac{1}{\eps} \div_{v}\left[f(u-v)\right]=0, \\
   &\rho_f(t,x)=\int_{\R^3}f(t,x,v)\ddv,  \quad
     j_f(t,x)= \int_{\R^3}vf(t,x,v)\ddv.
\end{aligned}
\right.
\end{equation*}
This corresponds to the physical situation where:
\begin{itemize}
\item  $U\sim \sqrt{\theta}$, i.e. the particles and fluid have comparable velocities,
\item $\tau \ll T$, i.e. the Stokes relaxation time is small compared to the observation time, \item ${\rho_{\text{part}}}\ll {\rho_{\text{fluid}}}$, which means that the particles are light compared to the fluid. Asymptotically they become inertialess.
\end{itemize}

A variant of this scaling corresponds to  the  {\bf light and fast particle} regime, for which
\begin{equation*}
\A= \frac{1}{\eps^\alpha}, \quad \B= \frac{1}{\eps}, \quad \C=1, \quad \F=1,
\end{equation*}
where $\alpha>0$, which leads to the system
\begin{equation*}
    \left\{
\begin{aligned}
   &\partial_tu+ (u\cdot\na_{x})u-\Delta_{x}u+\na_{x}p= j_f-\rho_fu,\\
 &\div_{x}u=0,\\
   &\partial_tf+ \frac{1}{\eps^\alpha} v\cdot\na_{x}f+ \frac{1}{\eps} \div_{v}\left[f(\eps^\alpha u-v)\right]=0, \\
   &      \rho_f(t,x)=\int_{\R^3}f(t,x,v)\ddv,  \quad
    j_f(t,x)= \frac{1}{\eps^\alpha} \int_{\R^3}vf(t,x,v)\ddv.
\end{aligned}
\right.
\end{equation*}
Compared to the previous regime, the only difference is that
 $U\ll \sqrt{\theta}$, i.e. the velocity of the particles is large compared to that of the fluid. Otherwise, as previously, the Stokes relaxation time is  small compared to the observation time and the particles are light compared to the fluid.

\subsubsection{Fine particle regime} The {\bf fine particle} regime corresponds to
\begin{equation*}
\A= 1, \quad \B= \frac{1}{\eps}, \quad \C=\frac{1}{\eps}, \quad \F=1,
\end{equation*}
leading to the system
\begin{equation*}
    \left\{
\begin{aligned}
    &\partial_tu+ (u\cdot\na_{x})u-\Delta_{x}u+\na_{x}p= \frac{1}{\eps}\left(j_f-\rho_fu\right),\\
  &\div_{x}u=0,\\
   &\partial_tf+v\cdot\na_{x}f+ \frac{1}{\eps} \div_{v}\left[f( u-v)\right]=0, \\
    &  \rho_f(t,x)=\int_{\R^3}f(t,x,v)\ddv,  \quad
  j_f(t,x)=\int_{\R^3}vf(t,x,v)\ddv.
\end{aligned}
\right.
\end{equation*}
This corresponds to the physical situation where:
\begin{itemize}
\item  $U\sim \sqrt{\theta}$, i.e. the particles and fluid have comparable velocities,
\item $\tau \ll T$, i.e. the Stokes settling time is small compared to the observation time, 
\item ${\rho_{\text{part}}}\sim {\rho_{\text{fluid}}}$, which means that the particles and fluid have comparable mass densities,
\item $a \ll L$, i.e. the size of the particles is small compared to the observation length.
\end{itemize}

\subsubsection{General scalings for the equations} In summary, we shall study the following scalings for the Vlasov-Navier-Stokes system
\begin{equation}
    \label{VNS-general}
     \index{f@$f_{\eps,\gamma,\sigma}(t,x,v)$: kinetic distribution function of the dispersed phase}
    \index{r@$\rho_{\eps,\gamma,\sigma}(t,x)$: density of the dispersed phase}
        \index{j@$j_{\eps,\gamma,\sigma}(t,x)$: momentum of the dispersed phase}
        \index{u@$u_{\eps,\gamma,\sigma}(t,x)$: fluid velocity}
    \left\{
\begin{aligned}
&\partial_t u_{\eps,\gamma,\sigma}+ (u_{\eps,\gamma,\sigma}\cdot\na_{x})u_{\eps,\gamma,\sigma}-\Delta_{x}u_{\eps,\gamma,\sigma}+\na_{x}p_{\eps,\gamma,\sigma}= \frac{1}{\gamma}\left(j_{f_{\eps,\gamma,\sigma}}-\rho_{f_{\eps,\gamma,\sigma}}u_{\eps,\gamma,\sigma}\right),\\
&\div_{x}u_{\eps,\gamma,\sigma}=0,\\
  &\partial_tf_{\eps,\gamma,\sigma}+ \frac{1}{\sigma}  v\cdot\na_{x}f_{\eps,\gamma,\sigma}+ \frac{1}{\eps} \div_{v}\left[f_{\eps,\gamma,\sigma}( \sigma u_{\eps,\gamma,\sigma}-v)\right]=0, \\
 &\rho_{f_{\eps,\gamma,\sigma}}(t,x)=\int_{\R^3}f_{\eps,\gamma,\sigma}(t,x,v)\ddv,  \quad
j_{f_{\eps,\gamma,\sigma}}(t,x)=\frac{1}{\sigma} \int_{\R^3}vf_{\eps,\gamma,\sigma}(t,x,v)\ddv,
\end{aligned}
\right.
\end{equation}
in which we consider
\begin{itemize}
\item $\sigma=1, \gamma=1$, for the light particle regime,
\item $\sigma= \eps^\alpha, \gamma =1$ with $\alpha>0$, for the light and fast particle regime,
\item $\sigma = 1, \gamma =\eps$, for the fine particle regime.
\end{itemize}
The scaled energy and dissipation functionals read:
\begin{align}
\label{energy-dissipation-general1}
 &\E_{\eps,\gamma,\sigma}(t)
                =   \frac{\eps}{ \sigma^2 \gamma} \frac{1}{2}\int_{\T^3\times\R^3}|v|^2f_{\eps,\gamma,\sigma} (t,x,v)\dd x\dd v
                +   \frac{1}{2}\int_{\T^3}|u_{\eps,\gamma,\sigma}(t,x)|^2\dd x,
               \\
\label{energy-dissipation-general2}
                &\D_{\eps,\gamma,\sigma}(t)=\int_{\T^3}|\na_x u_{\eps,\gamma,\sigma}(t,x)|^2\dd x+ \frac{1}{\gamma} \int_{\T^3 \times \R^3}\left|\frac{v}{\sigma}-u_{\eps,\gamma,\sigma}\right|^2 f_{\eps,\gamma,\sigma}(t,x,v)\dd x\dd v.
\end{align}
The Brinkman force is given by
\[
    \Fegd:=\frac{1}{\gamma}(j_{\fegd}-\rho_{\fegd}\uegd).
\]
In the light and fast particle regime, we shall restrict to the range of parameters $\alpha \in (0,1/2]$. Dealing with $\alpha>1/2$ would be possible, but by a view of the scaled energy $\E_{\eps,\gamma,\sigma}$, it would systematically require a well-preparedness assumption in order to ensure that it is initially bounded.
Remark that the parameter $\alpha=1/2$ seems particularly natural as the energies of the kinetic and the fluid part have then the same order in $\E_{\eps,\gamma,\sigma}$.

\subsection{Formal derivation of the limits}
\label{sec-formal}

We explain in this section how to obtain the formal high friction limits for the previous regimes.

\subsubsection{\emph{Light} and \emph{light and fast} particle regimes}

We start with the light particle limit of \eqref{VNS-general}, that corresponds to $\gamma=1,\,\sigma=1$. The conservation of mass and momentum read
\begin{equation*}
    \left\{
\begin{aligned}
&\partial_t \rho_{f_\eps} + \div_x j_{f_\eps}=0, \\
&\partial_t j_{f_\eps} + \div_x \left( \int_{\R^3} v \otimes v f_\eps \, \ddv\right) = \frac{1}{\eps} \left( \rho_{f_\eps} u_\eps- j_{f_\eps} \right).
\end{aligned}
\right.
\end{equation*}
We deduce that we must have
$$
\rho_{f_\eps} u_\eps- j_{f_\eps} \xrightharpoonup[\eps \to 0]{} 0.
$$
Assuming the following convergences
$$
(\rho_{f_\eps}, j_{f_\eps}, u_\eps)\xrightharpoonup[\eps \to 0]{} (\rho, j , u),
$$
we thus formally\footnote{assuming in particular that products also pass to the limit} get
$$
j = \rho u, 
$$
and $(\rho, u)$ has to satisfy the following system:
\begin{equation}
    \label{TNS}
    \left\{
\begin{aligned}
&\partial_t \rho + \div_x (\rho u)=0, \\
&\partial_tu+ u\cdot\na_{x} u-\Delta_{x}u+\na_{x}p=0, \\
&\div_x u =0,
\end{aligned}
\right.
\end{equation}
which corresponds to a transport equation driven by a velocity field satisfying the incompressible Navier-Stokes equation, and which we refer to as {\bf Transport-Navier-Stokes}. 

\begin{remark}
Contrary to the other regimes studied in this paper, the form of the scaled dissipation \eqref{energy-dissipation-general2} does not straightforwardly imply  that the distribution function $f_\eps$ weakly converges to a Dirac distribution in velocity. However, it will be a consequence of the upcoming analysis that 
$$
f_\eps \xrightharpoonup[\eps\to 0]{} \rho \otimes \delta_{v=u}.
$$

\end{remark}

\bigskip

The light and fast particle limit of \eqref{VNS-general}, that corresponds to $\gamma=1,\,\sigma=\eps^{\alpha}$, can be formally analyzed in a similar fashion, except that the expression of the scaled energy and dissipation \eqref{energy-dissipation-general1}-\eqref{energy-dissipation-general2} coupled with the energy--dissipation identity \eqref{energy-dissipation-intro} yields more information in this regime.  Assume the following convergences
$$
(\rho_{f_\eps},  u_\eps)\xrightharpoonup[\eps \to 0]{} (\rho, u).
$$
By the scaled energy--dissipation identity \eqref{energy-dissipation-intro}, we have
$$
\int_0^{+\infty} \int_{\T^3 \times \R^3} f_\eps | v - \eps^{\alpha} u_\eps|^2 \, \ \dd v \dd x \dd s \xrightarrow[\eps\to 0]{}  0
$$
at rate $\eps^{2\alpha}$, and therefore we deduce that formally
\begin{equation}
\label{introdirac-light}
f_\eps (t,x,v)\xrightharpoonup[\eps \to 0]{} \rho (t,x) \otimes \delta_{v=0}.
\end{equation}
The conservation of mass and momentum read in this regime
\begin{equation*}
    \left\{
\begin{aligned}
&\partial_t \rho_{f_\eps} + \div_x j_{f_\eps}=0, \\
&\partial_t j_{f_\eps} +  \frac{1}{\eps^{\alpha}}  \div_x \left( \int_{\R^3} v \otimes v f_\eps \, \ddv\right) = \frac{1}{\eps} \left( \rho_{f_\eps} u_\eps- j_{f_\eps} \right).
\end{aligned}
\right.
\end{equation*}
By~\eqref{introdirac-light}, we expect that 
$$
 \left( \int_{\R^3} v \otimes v f_{\eps} \, \ddv\right) \xrightharpoonup[\eps \to 0]{}  0,
$$
so that 
$$
\rho_{f_\eps} u_\eps- j_{f_\eps} \xrightharpoonup[\eps \to 0]{} 0.
$$
As a conclusion, we again formally obtain that $(\rho,u)$ has to satisfy the Transport-Navier-Stokes system~\eqref{TNS}.

\subsubsection{Fine particle regime}

Finally, we consider the fine particle regime for~\eqref{VNS-general} that is $\gamma=\eps,\,\sigma = 1$. 

Assume the following convergences
$$
(\rho_{f_\eps},  u_\eps)\xrightharpoonup[\eps \to 0]{} (\rho, u).
$$ 
In view of the scaled dissipation \eqref{energy-dissipation-general2} in this regime, 
we expect
\[
    \int_0^{+\infty}\int_{\T^3\times\R^3}\feps|v-\ueps|^2\dd v\dd x\dd s\xrightarrow[\eps\to0]{}0
\]
and thus,
\begin{equation}
\label{introdirac-fine}
f_\eps \xrightharpoonup[\eps \to 0]{} \rho \otimes \delta_{v=u}.
\end{equation}
The conservation of mass and momentum read in this regime
\begin{equation*}
    \left\{
\begin{aligned}
&\partial_t \rho_{f_\eps} +  \div_x j_{f_\eps}=0, \\
&\partial_t j_{f_\eps} + \div_x \left( \int_{\R^3} v \otimes v f_\eps \, \dd v\right)= -F_\eps.
\end{aligned}
\right.
\end{equation*}
where we have set
$$
F_\eps = \frac{1}{\eps} (\rho_{f_\eps}u_\eps - j_{f_\eps}).
$$
By \eqref{introdirac-fine}, we get
$$
\partial_t \rho + \div (\rho u) =0
$$
As we also formally have by~\eqref{introdirac-fine}
$$\int_{\R^3} v \otimes v f_\eps \, \dd v \xrightharpoonup[\eps \to 0]{} \rho  u\otimes u,$$
assuming that 
$$
F_\eps\xrightharpoonup[\eps \to 0]{}  F,
$$
we must formally have
$$
\partial_t (\rho u) + \na_x (\rho u\otimes u)  = - F.
$$
Passing to the limit in the Navier-Stokes equations, we obtain
$$
\partial_t u + \na_x (\rho u\otimes u) +  \nabla_x p =F,  \qquad \div_x u =0,
$$
hence $(\rho,u)$ must satisfy the system
\begin{equation}
    \label{inNS}
    \left\{
\begin{aligned}
 &\partial_t \rho + \div_x (\rho u) =0, \\
 &\partial_t ( (1+\rho) u) + \na_x ((1+\rho) u\otimes u) -\Delta_x u +  \nabla_x p =0, \\
&\div_x u = 0,
\end{aligned}
\right.
\end{equation}
which corresponds to the classical {\bf Inhomogeneous Incompressible Navier-Stokes} equations (in short, referred to as the Inhomogeneous Navier-Stokes system in the rest of the paper).
The Cauchy problem for these equations, for global weak solutions \emph{\`a la Leray} or strong solutions \emph{\`a la Fujita-Kato} has been the object of many research works. We may refer, among many others, to  the monographs of Antontsev, Kazhikhov and Monakhov \cite{ant-kaz-mon} and Lions \cite{LionsI} and   \cite{sim,dan,abi-gui-zha-cpam,abi-gui-zha,pai-zha-zha,pou,dan-muc} (and references therein).

\begin{remark}
Note that while modelling additional fragmentation phenomena,~\cite{ben-des-mou} also formally derived a similar system of equations. We shall come back to this type of model in Section~\ref{sec-mix}.
\end{remark}

\subsection{Main results of this work}

To conclude this introduction, let us state the main results of this work. The system~\eqref{VNS-general} is endowed with initial conditions
\[
\uegd|_{t=0}= \uegd^0, \quad \fegd|_{t=0}= \fegd^0.
\]
that may depend on $\eps$.
We start by presenting the assumptions about the initial data that we will use throughout this paper.

\begin{assumption}\label{hypGeneral}
    Let $(\uegd^0)_{\eps>0}\subset\H^1(\T^3)$ and $(\fegd^0)_{\eps>0}\subset\L^1\cap\L^\infty(\R^3\times\T^3)$ such that
    \begin{itemize}
        \item there exists $q>4$ such that for any $\eps>0$,
            \[
                (1+|v|^q)\fegd^0\in\L^1(\R^3;\L^\infty(\T^3))\cap \L^\infty(\T^3\times\R^3);
            \]
        \item there exist $r \in (2,3)$
        and $p\in\Big(3,\frac{3(2+r)}{4}\Big]$ such that and $\uegd^0\in \B^{1,3}_2(\T^3) \cap \B^{1-2/r,3}_r(\T^3) \cap \B^{s,p}_p(\T^3)$, where $s=2-2/p$;
        \item there exists $M>1$ such that, for any $\eps>0$,
            \[
                \norme{\uegd^0}_{\H^1\cap\B^{s,p}_p(\T^3)}+\norme{\fegd^0|v|^q}_{\L^1(\R^3;\L^\infty(\T^3))}+\sup_{(x,v)\in\T^3\times\R^3}(1+|v|^q)\fegd^0(x,v)\le M.
            \]
    \end{itemize}
\end{assumption}

The definition of the Besov spaces $\B^{s,q}_r(\T^3)$ is recalled in Definition~\ref{def-besov} of the Appendix. 

\begin{remark}
In the \emph{light} and \emph{light and fast} particle regimes, it is actually possible to lower down the required regularity on the initial fluid velocity, namely by only asking $\uegd^0\in\B^{\tilde{s},p}_p(\T^3)$, with $\tilde{s}= 2- 2\frac{p-1}{p}$ instead of $\B^{s,p}_p(\T^3)$.
\end{remark}

\begin{remark}
    Note that the parameters $p,q$ in Assumption~\ref{hypGeneral} verify $p<q$.
\end{remark}

We will also need the standard Fujita-Kato type \cite{fuj-kat} smallness assumption for the initial data of the Navier-Stokes equations: 

\begin{assumption}\label{hypSmallDataNS}
    The initial data $(\uegd^0)_{\eps>0}$ satisfy
    \[
        \forall\eps>0,\qquad\norme{\uegd^0}_{\dot\H^{\frac{1}{2}}(\T^3)}\le C^*/2,
    \]
    where $C_*$ is the universal constant given by Theorem~\ref{propTpsLong}.
\end{assumption}
The definition of the homogeneous Sobolev space $\dot\H^{\frac{1}{2}}(\T^3)$ is recalled in Definition~\ref{def-sobolev} of the Appendix.

To simplify the upcoming statements, we shall assume here that there exist $(u^0, \rho^0) \in \L^2(\T^3) \times \L^\infty(\T^3)$ such that 
\begin{equation*}
\norme{\uegd^0 - u^0}_{\L^2(\T^3)} \xrightarrow[\eps \to0]{} 0, \quad W_1\left(\rho_{\fegd^0}, \rho^0\right) \xrightarrow[\eps \to0]{} 0,
\end{equation*}
which can always be ensured up to an extraction by standard compactness arguments. Here $W_1$ stands for the Wasserstein-1 distance whose definition is recalled in Definition~\ref{def-W1} of the Appendix.

To state the convergence results, let us distinguish between the light, light and fast, and fine particle regimes.

\subsubsection{Light particle regime} In the light particle regime ($(\gamma,\sigma)=(1,1)$ in~\eqref{VNS-general}), we derive the Transport-Navier-Stokes system~\eqref{TNS}. This is achieved within three sets of assumptions:
\begin{itemize}
\item Under the sole Assumption~\ref{hypGeneral}, a situation we refer to as the \emph{general} case, one can find a time $T>0$ below which $u_{\eps,1,1}(t)$ converges towards $u(t)$ in $\L^2(\T^3)$, while $f_{\eps,1,1} \xrightharpoonup[\eps\to0]{} \rho\otimes\delta_{v= u}$ \emph{up to an integration in time}, where $(\rho,u)$ satisfies~\eqref{TNS} with initial conditions $(\rho^0, u^0)$.\\
The fact that the convergence of $f_{\eps,1,1}$ is not pointwise and holds only after integration in time comes from the fact that no such convergence is assumed for $t=0$. It thus formally means that there is an initial layer which turns out to be integrable and negligible when $\eps \to 0$.
\item With additional assumptions, the so-called \emph{mildly well-prepared} case, which loosely speaking consists in assuming that the initial fluid velocity is close enough to its average,
we are able to show that the aforementioned convergence results hold for all times $T>0$ 
\item Finally, under \emph{well-prepared} assumptions, which in particular ensure the initial convergence $f_{\eps,1,1}^0 \xrightharpoonup[\eps\to0]{} \rho^0\otimes\delta_{v=u^0}$, we
also obtain  pointwise convergence for $f_{\eps,1,1}$.

\end{itemize}

\begin{theorem}[Light particle regimes]
\label{thm1} Let $(u_{\eps,1,1},f_{\eps,1,1})$ be a global weak solution associated to the initial condition $(u^0_{\eps,1,1},f^0_{\eps,1,1})$. We have the following convergence results.

  \noindent {\bf 1. General case.}  Under Assumption~\ref{hypGeneral}, there exists $T>0$ such that 
  \begin{equation}
  \label{conv-1}
    \int_0^T W_1(f_{\eps,1,1}(t),\rho(t)\otimes\delta_{v= u(t)})\dd t \xrightarrow[\eps\to0]{}0
  \end{equation}
    and for all $t \in [0,T]$,
    \begin{equation}
    \label{conv-2}
     \norme{u_{\eps,1,1}(t)-u(t)}_{\L^2(\T^3)} \xrightarrow[\eps\to0]{}0,
    \end{equation}
   where $(\rho,u)$ satisfies the Transport-Navier-Stokes system \eqref{TNS}.
    
     \noindent {\bf 2. Mildly well-prepared case.} Under the additional Assumption~\ref{hypSmallDataNS}, there exists $\eta>0$ small enough such that if
    \[
        \norme{u_{\eps,1,1}^0-\langle u_{\eps,1,1}^0\rangle}_{\L^2(\T^3)}\le\eta,
    \]
    the convergences~\eqref{conv-1} and~\eqref{conv-2} hold for any $T>0$.
   
   \noindent {\bf 3. Well-prepared case.}  
    Finally, if we further assume that
    \[
        \int_{\T^3\times\R^3}|v-u_{\eps,1,1}^0(x)|f_{\eps,1,1}^0(x,v)\,\dd x\dd v\xrightarrow[\eps\to0]{}0,
    \]
    then, for all $t>0$,
    \begin{equation}
    \label{conv-3}
        W_1(f_{\eps,1,1}(t),\rho(t)\otimes\delta_{v=u(t)}) \xrightarrow[\eps\to0]{}0.
    \end{equation}
    All convergence results are quantitative with respect to $T$ and $\eps$.
    
\end{theorem}

As a matter of fact, we shall provide more general and more precise statements. We refer to Section~\ref{sec-smallT} for case {\bf 1} and to Section~\ref{SectionAsymptotic-12} for cases {\bf 2} and {\bf 3}.

\subsubsection{Light and fast particle regime} In the light and fast particle regime  ($(\gamma,\sigma)=(1,\eps^\alpha)$ in~\eqref{VNS-general}), we also derive the Transport-Navier-Stokes system~\eqref{TNS}.

We obtain similar results to the light particle regime, and the previous formal discussion mostly still prevails here.
The main difference with the light particle regime is as follows. As its denomination suggests, this regime allows for the description of \emph{fast} particles, which means in particular that at initial time, the  Brinkman force may not be uniformly bounded with respect to $\eps$ (mind the $1/\sigma$ in the definition of $\jegd$).
However our results show that because of the high friction, the order of magnitude of the particle velocities becomes instantaneously of order $1$. Again, there is therefore in general an initial layer whose contribution vanishes when integrated in time as $\eps \to 0$. 

The case $\alpha=1/2$ turns out to be more singular than the cases $\alpha<1/2$, and requires more assumptions on the sequence of initial data. It is convenient to distinguish between the two.

\begin{theorem}[Light and fast particle regimes]
\label{thm2}
Let $(u_{\eps,1,\eps^{\alpha}},f_{\eps,1,\eps^{\alpha}})$ be a global weak solution associated to the initial condition $(u^0_{\eps,1,\eps^{\alpha}},f^0_{\eps,1,\eps^{\alpha}})$. We have the following convergence results.

\bigskip

\noindent {\bf A.} \underline{The case $\alpha<1/2$}.

\medskip

  \noindent {\bf 1. General case.} Under Assumption~\ref{hypGeneral}, there exists $T>0$ such that
  \begin{equation}
  \label{conv-1-LF}
    \int_0^TW_1(f_{\eps,1,\eps^{\alpha}}(t),\rho(t)\otimes\delta_{v=0})\dd t \xrightarrow[\eps\to0]{}0
  \end{equation}
    and for all $t \in [0,T]$,
    \begin{equation}
    \label{conv-2-LF}
     \norme{u_{\eps,1,\eps^{\alpha}}(t)-u(t)}_{\L^2(\T^3)} \xrightarrow[\eps\to0]{}0,
    \end{equation}
   where $(\rho,u)$ satisfies the Transport-Navier-Stokes system \eqref{TNS}.
    
     \noindent {\bf 2. Mildly well-prepared case.}  Under the additional Assumption~\ref{hypSmallDataNS}, there exists $\eta>0$ small enough such that if
\[
          \norme{u_{\eps,1,\eps^{\alpha}}^0-\langle u_{\eps,1,\eps^{\alpha}}^0\rangle}_{\L^2(\T^3)}\le\eta,
\]
then
    the convergences~\eqref{conv-1-LF} and~\eqref{conv-2-LF} hold for any $T>0$.
   
   \noindent {\bf 3. Well-prepared case.}  
    Finally, if we further assume that
    \[
        \int_{\T^3\times\R^3}|v|f_{\eps,1,\eps^{\alpha}}^0(x,v)\,\dd x\dd v\xrightarrow[\eps\to0]{}0,
    \]
    then, for all $t>0$,
    \begin{equation}
    \label{conv-3-LF}
        W_1(f_{\eps,1,\eps^{\alpha}}(t),\rho(t)\otimes\delta_{v=0}) \xrightarrow[\eps\to0]{}0.
    \end{equation}
Furthermore all convergence results are quantitative with respect to $T$ and $\eps$.

\bigskip

\noindent {\bf B.} \underline{The case $\alpha=1/2$}.  

\medskip

The previous results still hold, with the following additional assumptions, for the general case and the mildly well-prepared case: 
\[
\int_{\T^3\times\R^3}|v|^2f_{\eps,1,\eps^{1/2}}^0(x,v)\dd x\dd v \le \eta,
\]
for $\eta>0$ small enough. 
Moreover the  results are quantitative under the assumption
\[
\int_{\T^3 \times \R^3} f^0_{\eps,1,\eps^{1/2}} {|v|}\dd x \dd v \xrightarrow[\eps\to0]{} 0.
\]

\end{theorem}

Again, we emphasize that we shall provide  more general and more precise statements. It turns out that we can prove Theorem~\ref{thm2} together with Theorem~\ref{thm1} in a completely unified manner, however under the additional well-prepared assumption
 \begin{equation}\label{hyp-well-prep-2Intro}
        \int_{\R^3\times\R^3}|v|^pf_{\eps,1,\eps^{\alpha}}^0(x,v)\dd x\dd v\lesssim\eps^{\alpha p-1+\kappa},
    \end{equation}
    for some $\kappa\in(0,1)$, which is significant only for $\alpha>1/p$.
 We therefore also refer to Section~\ref{sec-smallT} for case {\bf 1} and to Section~\ref{SectionAsymptotic-12} for cases {\bf 2} and {\bf 3}. Theorem~\ref{thm2} with~\eqref{hyp-well-prep-2Intro} removed requires a refinement of the method and is proved in Sections~\ref{sec-LF-nohyp} and~\ref{sec-conv12-re}.  

\subsubsection{Fine particle regime} In the fine particle regime ($(\gamma,\sigma)= (\eps,1)$ in~\eqref{VNS-general}) we derive the Inhomogeneous Navier-Stokes system~\eqref{inNS}. This proves to be much more singular than the previous two regimes. As a consequence, more assumptions on the initial data are required to justify the convergence results. In particular, we are no longer able to state any result with the sole Assumption~\ref{hypGeneral} and need to always consider well-prepared initial data. Moreover, a smallness assumption on the kinetic initial distribution function is required.

\begin{theorem}[Fine particle regimes]
\label{thm3}
Let $(u_{\eps,\eps,1},f_{\eps,\eps,1})$ be a global weak solution associated to the initial condition $(u^0_{\eps,\eps,1},f^0_{\eps,\eps,1})$. We have the following convergence results.

  \noindent {\bf 1. Mildly well-prepared case.}   Under Assumption~\ref{hypGeneral}, there exist $\eps_0$, $\eta>0$, $M'>0$,  such that, if for all $\eps \in (0,\eps_0)$,
  \[
   \norme{f_{\eps,\eps,1}^0}_{\L^1(\R^3;\L^\infty(\T^3))} \leq \eta
  \]
  and 
  \[
   \int_{\T^3\times\R^3}\frac{|v-u_{\eps,\eps,1}^0(x)|^p}{\eps^{p-1}}f_{\eps,\eps,1}^0(x,v)\dd x\dd v\le M',
  \]
  then
  there exists $T>0$ such that for all $t \in [0,T]$,
  \begin{equation}
  \label{conv-1-Fine}
  W_1(f_{\eps,\eps,1}(t),\rho(t)\otimes\delta_{v= u(t)})\xrightarrow[\eps\to0]{}0,
  \end{equation}
    and for all $t \in [0,T]$,
    \begin{equation}
    \label{conv-2-Fine}
     \norme{u_{\eps,1,1}(t)-u(t)}_{\L^2(\T^3)} \xrightarrow[\eps\to0]{}0,
    \end{equation}
   where $(\rho,u)$ satisfies the Inhomogeneous Navier-Stokes system \eqref{inNS}.
    
     \noindent {\bf 2. Well-prepared case.} Under the additional Assumption~\ref{hypSmallDataNS}  there exists $\eta'>0$ small enough such that if
    \begin{multline*}
     \int_{\T^3\times\R^3}\left|v-\frac{\langle j_{\eps,1,\eps}^0\rangle}{\langle \rho_{\eps,\eps,1}^0\rangle}\right|^2 f_{\eps,\eps,1}^0\,\dd x\dd v
        +\int_{\T^3}|u_{\eps,\eps,1}^0-\langle u_{\eps,\eps,1}^0\rangle|^2\dd x \\
        + \left|\frac{\langle j_{\eps,\eps,1}^0\rangle}{\langle \rho_{\eps,\eps,1}^0\rangle}-\langle u_{\eps,\eps,1}^0\rangle\right|^2
        \leq \eta,
    \end{multline*}
     the convergences~\eqref{conv-1-Fine} and~\eqref{conv-2-Fine} hold for any $T>0$.
   
  If we further assume that the solution $(\rho,u)$ to~\eqref{inNS} is smooth enough, 
the convergences~\eqref{conv-1-Fine} and~\eqref{conv-2-Fine} are quantitative with respect to $T$ and $\eps$.

\end{theorem}

We have used the following notations for the two first initial moments in the statement of the theorem:
$$
\rho_{\eps,\eps,1}^0= \int_{\R^3} f_{\eps,\eps,1}^0(\cdot,v) \dd v, \quad j_{\eps,\eps,1}^0= \int_{\R^3} v f_{\eps,\eps,1}^0(\cdot,v) \dd v.
$$

Like for the two other regimes, we shall provide more complete statements. Let us refer to Section~\ref{sec-smallT-fine} for case {\bf 1} and to Sections~\ref{sec-first-fine} and \ref{sec-relat-fine} for case {\bf 2}.


We end the statement of the results with a remark that concerns all regimes.

\begin{remark} The assumptions in Theorems~\ref{thm1}, \ref{thm2} and \ref{thm3} are somewhat comparable to those of \cite{hank-mou-moy} for the long time behavior of the Vlasov-Navier-Stokes system \eqref{VNS}, except that we require in Assumption~\ref{hypGeneral} extra bounds in  Besov spaces. The reason is as follows. The analysis in \eqref{VNS} relies on the instantaneous smoothing effect of the heat semi-group on $\T^3$ and uses rough bounds (similar to those of Section~\ref{sec-roughfirst}) close to $t=0$; we cannot argue similarly for the study of high friction regimes as these rough bounds blow up when $\eps \to 0$, see  Section~\ref{sec-roughfirst}.

\end{remark}


\subsection{Acknowledgements}  Partial support by the grant ANR-19-CE40-0004 is acknowledged.

\section{State of the art and methodology of this work}

\subsection{State of the art}

Let us start by reviewing the mathematical literature that is relevant to this work.

\noindent {\bf Cauchy problem for Vlasov-Navier-Stokes.} 
The Cauchy problem was first studied for related reduced models in the pioneering works of Anoshchenko and Boutet de Monvel-Berthier \cite{ano-bou} and Hamdache \cite{ham}.
Existence of global weak solutions to the Vlasov-Navier-Stokes system on $\T^3$ was then established in the seminal work of Boudin, Desvillettes, Grandmont and Moussa \cite{bou-des-gra-mou}. This was later extended to the case of  bounded domains of $\R^3$ in \cite{wang-yu} and to  time-dependent bounded domains of $\R^3$ in \cite{bou-gra-mou}. Lately, in \cite{bou-mic-mou}, the authors obtain similar results for an augmented model in which hygroscopic effects are taken into account and lead to considering the size and temperature variation of the particles in the Vlasov equation. In dimension two, the uniqueness of weak solutions in the torus or in the whole space is obtained in \cite{hank-miot-mou-moy}. 
The existence of smooth solutions has also been studied, for instance \cite{cho-kwo} provides a local well-posedness result in dimension three. 

\bigskip


\noindent {\bf Large time behavior for Vlasov-Navier-Stokes.} 
Large time behavior for a related reduced model was first tackled in the pioneering work of Jabin \cite{jabin}.
Concerning the Vlasov-Navier-Stokes system on $\T^3$, a conditional result, describing the large time behavior of solutions, was obtained by Choi and Kwon in \cite{cho-kwo}; loosely speaking a global bound on the moment $\rho_f$ needs to be \emph{a priori} assumed. More recently, this restriction has been removed in \cite{hank-mou-moy} for initial conditions close to equilibrium, in a framework \emph{\`a la} Fujita-Kato. This was recently extended to the whole space case in \cite{hank} and to bounded domains with absorption boundary conditions in \cite{ert-han-mou}. As we shall soon see, this question is closely linked to the limits we consider in this paper.
In the Vlasov-Navier-Stokes system, concentration effects in velocity are at play, which eventually lead to a \emph{monokinetic} behavior for the kinetic distribution function, that is to say, convergence to a Dirac mass in velocity. This is precisely this behavior that is described in \cite{hank-mou-moy,hank}.

In a somewhat different direction, \cite{gla-hank-mou} provides the existence and stability of regular equilibria for the system in a domain with partly absorbing boundary conditions and injection of particles and fluid.

\bigskip

\noindent {\bf High friction limits for Vlasov-(Fokker-Planck)-(Navier)-Stokes.}  
High friction limits for fluid-kinetic systems have been the topic of several recent works, most of the time with an additional Fokker-Planck operator in the kinetic equation (that accounts for the effect of Brownian motion), resulting in the Vlasov-Fokker-Planck equation
$$
 \partial_tf_{\eps,\gamma,\sigma}+ \frac{1}{\sigma}  v\cdot\na_{x}f_{\eps,\gamma,\sigma}+ \frac{1}{\eps} \div_{v}\left[f_{\eps,\gamma,\sigma}( \sigma u_{\eps,\gamma,\sigma}-v)+\na_v\fegd\right]=0, 
$$
instead of the Vlasov equation in~\eqref{VNS-general}.

\noindent\emph {Light particle regimes.}
Light particle regimes have been studied by Jabin in \cite{jab} for a Vlasov equation where the fluid velocity in \ref{VNS-general} is computed by means of a convolution operator with a smooth kernel applied to a moment of the distribution function; a general framework has been put forward but it does not apply to the full Vlasov-Navier-Stokes system.  A similar problem is dealt with in \cite{gou-pou} where the fluid velocity is a given random vector field. A toy model for \eqref{VNS-general} is tackled in dimension $d=1$ in \cite{gou}; the analysis strongly relies on the one-dimensional framework (notably, Sobolev embeddings are much more favorable) and cannot be directly adapted to the three-dimensional case. 
More recently, the light particle regime was studied by H\"ofer in \cite{hofe} for the Vlasov-Stokes system on the whole space in the presence of a gravity force (a limit referred to as the inertialess limit there).

The light and fast particle regime was studied in a seminal work by Goudon, Jabin and Vasseur \cite{gou-jab-vas04a} in the Vlasov-Fokker-Planck case for the critical parameter $\alpha=1/2$.
They were able to justify the asymptotics in dimension $d=2$, leading to the Kramer-Smoluchowski equations, thanks to global entropy bounds related to the Fokker-Planck operator. 


\noindent\emph {Fine particle regimes.} In another seminal work \cite{gou-jab-vas04b}, Goudon, Jabin and Vasseur consider the fine particle regime in the whole space, also for the Vlasov-Fokker-Planck-Navier-Stokes system. They relied on a \emph{relative entropy} method (see Section~\ref{sec-classicalmethods} below for a short review) to derive the Inhomogeneous Navier-Stokes equations, see also \cite{mel-vas08}.


\bigskip

\noindent{\bf Other hydrodynamic limits with monokinetic behavior.} 
It is important to note that in the Vlasov-Fokker-Planck cases studied in \cite{gou-jab-vas04a,gou-jab-vas04b}, the distribution function converges (at least formally) to a \emph{local Maxwellian}, precisely because of the Fokker-Planck operator. Without it, as for the large time asymptotics, we expect a monokinetic behavior in the limit, that is to say convergence towards a Dirac mass in velocity.
To conclude this review of the literature, let us briefly discuss other hydrodynamic limits for Vlasov-type equations in which Dirac distributions in velocity also show up. In a very influential work \cite{bre}, Brenier studies the quasineutral limit (i.e. the small Debye length regime) of the Vlasov-Poisson system. He shows that if the sequence of initial data converges in a certain sense to a Dirac distribution in velocity (with mass one) -- an assumption one can refer to as \emph{well-prepared} initial data -- then the solution also displays a monokinetic behavior with a velocity that satisfies the incompressible Euler equations.
More recently, Kang and Vasseur \cite{kan-vas} and Kang and Figalli \cite{fig-kan} have studied high friction limits for the so-called kinetic Cucker-Smale equation (see also the recent work \cite{car-cho}). They show that for well-prepared initial data, the solution asymptotically converges towards a Dirac mass in velocity, with a velocity satisfying pressureless Euler-type equations. See also the related recent works on hydrodynamic limits of Fritz-Nagumo kinetic equations \cite{cre-fay-fil,cre}.


\subsection{Strategy}

Before describing the strategy at stake in this work, we provide an overview of the classical methods that are available to tackle hydrodynamic limits for kinetic equations.

\subsubsection{A review of classical methods}
\label{sec-classicalmethods}

There are mainly three classes of methods that have been devised to study such questions. 

\bigskip

\noindent {\bf Weak compactness methods.} The broad principle of this class of methods consists in finding uniform (with respect to the small parameter $\eps$) bounds to obtain the weak compactness of certain relevant quantities in the equations. 
Such bounds often come from conservation laws of the system (e.g. conservation of mass, energy, etc.). From the mathematical point of view, this usually also involves finding some (strong) compactness in order to overcome possible oscillations and eventually pass to the limit in nonlinear terms.

With regard to the Navier-Stokes limit of the Boltzmann equation, a program based on this method was set up by Bardos, Golse and Levermore \cite{bar-gol-lev,bar-gol-lev2} and a complete justification has been obtained, see \cite{gol-sai,gol-sai2,ars} and references therein. Let us also refer to the recent monograph of Ars\'enio and Saint-Raymond \cite{ars-sai} for generalizations in several directions.

With regard to high friction limits, this is the strategy applied in the 1D toy model of \cite{gou}, and in  \cite{gou-jab-vas04a} in the Vlasov-Fokker-Planck case for the light and fast particles  regime for $\alpha=1/2$.

\bigskip

\noindent {\bf Modulated energy/relative entropy methods.} Following the pioneering works of Dafermos \cite{daf} and Yau \cite{Yau} and more specifically Brenier \cite{bre} and Golse \cite{bou-gol-pul} for kinetic equations, this class of methods consists in \emph{modulating} a well chosen functional (often the energy or entropy) with solutions to the target limit equations and showing that this new functional is vanishing as $\eps \to 0$.
One must also ensure that the modulated functional indeed allows to quantify the convergence to the limit.

This method was notably used to treat the incompressible Euler limit of the Boltzmann equation, see \cite{sai,sai2}. We can also again refer to the recent \cite{ars-sai}.

For what concerns high friction limits, this is the strategy applied in \cite{gou-jab-vas04b} in the Vlasov-Fokker-Planck case for the fine particle regime. 
This is also the main strategy at play for the aforementioned hydrodynamic limits for Vlasov-Poisson \cite{bre} and Cucker-Smale \cite{kan-vas,fig-kan,car-cho}. 

\bigskip

\noindent  {\bf Higher regularity/Strong compactness methods.} This class of methods consists in building strong solutions enjoying uniform regularity with respect to the small parameter $\eps$. This has proven efficient for various hydrodynamic limits for the Boltzmann equation,
whether with an approach based on an (Hilbert or Chapman-Enskog) expansion with respect to $\eps$, see e.g. \cite{nis,dem-esp-leb}, or in a regime close to equilibrium \cite{bar-uka}. The latter, which most often requires a fine understanding of the spectral properties of the linearized Boltzmann operator, has witnessed several recent developments, see e.g. \cite{bri,bri-mer-mou,jia-xu-zha,gal-tri,alo-lod-tri} and references therein. 


\bigskip


\subsubsection{Strategy of this work}

The strategy that we implement can be seen as a mix between all aforementioned approaches. It is important to note that an argument solely based on propagation of higher regularity for all quantities is likely to fail, as we expect to obtain Dirac masses in velocity in the limit for the distribution function: as a result, higher regularity with respect to the velocity variable cannot be expected.


As usual in singular limits problems, the key is to obtain appropriate uniform estimates. Here, we aim at controlling the local density $\rho_{f_\eps}$ 
in  $\L^\infty(0,T; \L^\infty(\T^3))$.  Such a control does not straightforwardly follow from the natural conservation laws of the system. 
Inspired by the analysis of \cite{hank-mou-moy} (which, we recall, concerns the long time behavior of solutions to~\eqref{VNS}), such bounds for $\rho_{f_\eps}$ actually follow from a control on the fluid velocity field of the form
\begin{equation}
\label{eq-intro-nau}
\| \na_x u_\eps\|_{\L^1(0,T; \L^\infty(\T^3))} \ll 1,
\end{equation}
thanks to a natural straightening change of variables in velocity in the style of \cite{bar-deg}. This paves the way to the short time results of Theorems~\ref{thm1}, \ref{thm2} and \ref{thm3}, using weak compactness techniques to pass to the limit. For the large time results, thanks to the work of Choi and Kwon \cite{cho-kwo}, it turns out that a uniform control of $\rho_{f_\eps}$ in $\L^\infty(0,T; \L^\infty(\T^3))$ entails a uniform exponential decay in time for the so-called \emph{modulated energy}, a functional encoding the large time monokinetic behavior of the distribution function. We  emphasize that this is not sufficient to directly show any convergence as $\eps \to 0$. However it provides some control, namely integrability in time and smallness, to eventually deduce uniform bounds such as \eqref{eq-intro-nau} for arbitrarily large times.   


We shall set up a bootstrap argument to prove~\eqref{eq-intro-nau}.
To this end, we use higher order estimates for $u_\eps$ which are obtained thanks to maximal regularity estimates for the Stokes equation. These estimates are then interpolated with the pointwise in time $\L^2(\T^3)$ bounds stemming from the modulated energy dissipation.
The main difficulty lies in obtaining appropriate uniform in $\eps$ estimates of the Brinkman force
$$F_{\eps,\gamma,\sigma}=\frac{1}{\gamma}\left(j_{f_{\eps,\gamma,\sigma}}-\rho_{f_{\eps,\gamma,\sigma}}u_{\eps,\gamma,\sigma}\right)$$
in $\L^p(0,T; \L^p(\T^3))$ for large enough values of $p$. To this end, a process of \emph{desingularization} (with respect to $\eps$) of the Brinkman force is required.  

As we shall see below, the \emph{fine particle regime} is significantly more singular than the light particles regimes.
To tackle the former, we rely on some ideas from \cite{hank}.
As already mentioned, the work \cite{hank} deals with the large time behavior of small data solutions to the Vlasov-Navier-Stokes system posed on $\R^3 \times \R^3$. The general strategy follows that of the torus case \cite{hank-mou-moy}, but the slow decay (polynomial on the whole space, vs. exponential on the torus) of the solution to the Stokes equation makes (among other things) a finer understanding of the structure of the Brinkman force compulsory. This has led to a family of identities for a notion of \emph{higher} dissipation (see \cite[Lemma 4.2]{hank}) allowing for a better decay, which we shall also rely on.

Let us finally mention that in the fine particle regime, we also take advantage of a \emph{relative entropy} functional, which is useful to obtain quantitative convergence estimates.

\begin{remark}
The difficulties in order to apply the weak compactness strategy to study the fine particle regime are thoroughly discussed in \cite[Chapitre 2, Section 2.3.3]{mou-hdr}.
\end{remark}

\subsection{Outline of the proofs and organization of the paper}

The aim of this section is to provide a detailed overview of the proofs, give a flavor of the analysis and explain, at the same time, the organization of this paper. 
As we already stated, the study of high friction limits is closely related to the long time behavior problem, and the general strategy we follow is inspired from \cite{hank-mou-moy}.
In this presentation, we shall mostly focus on \emph{large time} results (i.e. for arbitrary times $T>0$) which require some well-prepared initial conditions. We recall that such assumptions can also be dispensed with, at the expense of restricting to \emph{short time} results.

\bigskip

\noindent {\bf Preliminaries.} The beginning of the analysis of all three regimes is common and consists, loosely speaking, in a careful extension of the preliminaries of \cite{hank-mou-moy} (plus some additional technical results), which tracks down the dependence with respect to the small parameter $\eps$. This step is performed in Section~\ref{sec-prelim}.


Let us introduce the main objects of interest in this work. We begin by the the characteristic curves associated to the Vlasov equation
\begin{equation}\label{charac-intro}
    \left\{\begin{aligned}
        \dot{X}_{\egd}(s;t,x,v)&=\frac{1}{\sigma}V_{\egd}(s;t,x,v),\\
        \dot{V}_{\egd}(s;t,x,v)&=\frac{1}{\eps}\left(\sigma u_\egd\left(s,X_{\egd}(s;t,x,v)\right)-V_{\egd}(s;t,x,v)\right),\\
        X_{\egd}(t;t,x,v)&=x,\\
        V_{\egd}(t;t,x,v)&=v,
    \end{aligned}\right.
\end{equation}
which are particularly useful since we have the classical representation formula 
\begin{equation}
\label{characmethod-intro}
f_\egd(t,x,v) = e^{\frac{3t}{\eps}} f^0\left({X}_{\egd}(0;t,x,v), {V}_{\egd}(0;t,x,v)\right).
\end{equation}
We claim that the key is to understand how to obtain \emph{uniform} (with respect to $\eps$) estimates for the moments $\rho_\egd  = \int_{\R^3} f_\egd \dd v$ and $j_\egd = \frac{1}{\sigma} \int_{\R^3} f_\egd v\, \dd v$ of the form
$$
    \| \rho_\egd \|_{\L^\infty(0,T; \L^\infty(\T^3))} +   \| j_\egd \|_{\L^\infty(0,T; \L^\infty(\T^3))} \lesssim 1.
$$

Controlling $\| \rho_\egd \|_{\L^\infty(0,T; \L^\infty(\T^3))}$ proves interesting because of the following facts, that stem from the study of the long time behavior for the Vlasov-Navier-Stokes system.

Building on the scaled energy~\eqref{energy-dissipation-general1} and following Choi and Kwon \cite{cho-kwo}, we define the \emph{modulated energy} as the following functional
\begin{align*}
    \EEE_{\egd}(t)=&
    \frac{\eps}{2\gamma}\int_{\T^3\times\R^3}\left|\frac{v}{\sigma}-\frac{\langle j_\egd(t)\rangle}{\langle\rho_\egd(t)\rangle}\right|^2f_\egd(t,x,v)\,\dd x\dd v\\
    &+\frac{1}{2}\int_{\T^3}|u_\egd(t,x)-\langle u_\egd(t)\rangle|^2\dd x\\
    &+\frac{\eps \langle \rho_\egd(t) \rangle}{2(\gamma+\eps \langle \rho_\egd(t) \rangle)}\left|\frac{\langle j_\egd(t)\rangle}{\langle\rho_\egd(t) \rangle}-\langle u_\egd(t)\rangle\right|^2,
\end{align*}
where $\langle g \rangle$ stands for the mean value of $g$ on the torus $\T^3$. Using the energy--dissipation identity~\eqref{energy-dissipation-intro} and the conservation laws of the system, we obtain for almost all $t\ge s \geq 0$ (including $s=0$),
\begin{equation}
    \label{mod-ener-intro}
    \EEE_{\egd}(t)+\int_s^t\D_{\egd}(\tau)\dd\tau\le\EEE_{\egd}(s).
\end{equation}
As in \cite{cho-kwo}, a remarkable consequence of this identity is a conditional exponential decay of the modulated energy $\EEE_\egd$. Let $T \in \R^+ \cup \{+\infty\}$. \emph{Assuming} a control on the density $\rho_\egd$  of the form
\begin{equation}
\label{rhoeps-intro}
\| \rho_\egd\|_{\L^\infty(0,T ; \L^\infty(\T^3))}\leq \C_0,\quad \forall \eps \in (0,1),
\end{equation}
one can prove (see Lemma~\ref{lemmeDE})
    \begin{equation}\label{DE-intro}
        \forall t\in[0,T],\qquad \EEE_{\egd}(t)\le C e^{-\lambda_{\egd} t}\EEE_{\egd}(0),
    \end{equation}
    where $\lambda_\egd>0$ is bounded from below by a positive constant and $C>0$ depends only on $\lambda_\egd$, \emph{independently} of $\eps$. As a result,  ensuring the global control 
    $$
\| \rho_\egd\|_{\L^\infty(0,+\infty; \L^\infty(\T^3))}\leq \C_0,\quad \forall \eps \in (0,1)
$$
is the key to understanding the long time behavior of $\EEE_\egd$; it turns out that this is then sufficient to describe the long time behavior of the solutions to the Vlasov-Navier-Stokes system. This fact is at the heart of the strategy of \cite{hank-mou-moy}.
In the context of high friction limits, such a control of $\EEE_\egd$ will provide \emph{uniform} integrability in time properties, as well as \emph{smallness} since $\EEE_{\egd}(0)$ can be chosen as small as desired.
    
The strategy to control the moments directly comes from \cite{hank-mou-moy}. 
\emph{Assuming} a control on the fluid velocity $u_\egd$ of the form
 \begin{equation}\label{hyp-na-u-intro}
    \norme{\na_x u_\egd}_{\L^1(0,T;\L^\infty(\T^3))} \leq 1/30,
\end{equation}
one proves (see Lemma \ref{changeVarV}) that for  all $x\in\T^3$, the map
\begin{equation}
\label{def-gamma-intro}
        \Gamma_{\egd}^{t,x}:v\mapsto V_{\egd}(0;t,x,v)
\end{equation}
    is a $\CCC^1$-diffeomorphism from $\R^3$ to itself and satisfies
    \[
        \forall v\in\R^3,\qquad \det\D_v\Gamma_{\egd}^{t,x}(v)\ge\frac{e^{\frac{3t}{\eps}}}{2}.
    \]
    We can then use the change of variables in velocity  $w:=\Gamma_{\egd}^{t,x}(v)$ in the formula for the density $\rho_\egd$:
  \begin{align*}
\rho_\egd (t,x)&= e^{\frac{3t}{\eps}} \int_{\R^3} f^0({X}_{\egd}(0;t,x,v), {V}_{\egd}(0;t,x,v)) \dd v \\
&= e^{\frac{3t}{\eps}} \int_{\R^3} f^0(\tilde{X}_{\egd}^{t,x,v}(0), w)  |\det\D_v\Gamma_{\egd}^{t,x}|^{-1}(w)\dd w,
  \end{align*}
with
$$
  \tilde{X}_{\egd}^{t,x,v}(s):= {X}_{\egd}(s;t,x,[ \Gamma_{\egd}^{t,x}]^{-1}(w)),
$$
 which, using Assumption~\ref{hypGeneral}, yields a control of the form~\eqref{rhoeps-intro}.  Under the same assumption~\eqref{hyp-na-u-intro}, one also proves that the map $x \mapsto   \tilde{X}_{\egd}^{t,x,v}(s)$ is a $\CCC^1$-diffeomorphism from $\T^3$ to itself with Jacobian bounded from below by $1/2$.
    
    Therefore, the goal is now to obtain the bound \eqref{hyp-na-u-intro}. To this end, one looks for higher order estimates for $\uegd$, which rely on the parabolic nature of the Navier-Stokes equation. We shall use $\L^p \L^p$ (for $p>3$) maximal estimates for the Stokes operator on $\T^3$ to get
       \begin{multline*}
       \label{LpLp-intro}
        \norme{\partial_tu_\egd}_{\L^p((0,T)\times\T^3)}+\norme{\Delta_xu_\egd}_{\L^p((0,T)\times\T^3)}\\
        \lesssim \norme{F_\egd}_{\L^p((0,T)\times\T^3)}+\norme{(u_\egd\cdot\na_x)u_\egd}_{\L^p((0,T)\times\T^3)}+\norme{\uegd^0}_{\B^{s,p}_p(\T^3)},
    \end{multline*}
    in which $F_\egd= j_\egd- \rho_\egd u_\egd$ stands for the Brinkman force and $s=2-2/p$.
   To control the  term  due to the nonlinearity $u_\egd\cdot\na_x u_\egd$, one has to deal with strong enough solutions to the Navier-Stokes equation. Namely, we require a $\L^\infty\H^1$--$\L^2\dot\H^2$ control of the form
   \begin{equation}\label{ENSH1-intro}
    \norme{u_\egd(t)}_{\H^1(\T^3)}^2+\int_0^t\norme{\Delta_xu_\egd(s)}_{\L^2(\T^3)}^2\dd s\lesssim \norme{u_\egd^0}_{\H^1(\T^3)}^2+\norme{\F_\egd}_{\L^2((0,t)\times\T^3)}^2.
\end{equation}
To enforce this property, we consider a Fujita-Kato type  framework, which consists in ensuring that
\begin{equation}\label{hypFujitaKato-Intro}
    \int_0^t\norme{e^{s\Delta}\uegd^0}_{\dot\H^1(\T^3)}^4\dd s+\int_0^t\norme{F_\egd(s)}_{\dot{\H}^{-\frac{1}{2}}(\T^3)}^2\dd s \le C^*,
\end{equation}
for some universal constant $C^*\in(0,1)$. We can achieve this with
\begin{itemize}
    \item a smallness assumption in $\dot\H^{\frac{1}{2}}(\T^3)$ for the initial fluid velocity,
    \item or by considering only small times.
\end{itemize}
In either case, we will still need to make sure that $\norme{F_\egd}_{\L^2((0,t)\times\T^3)}\ll 1$.
The idea is to set up a bootstrap argument. 
For a given $\eps>0$, we begin by saying that $T>0$ is a \emph{strong existence time} if
\begin{equation}
    \label{strong-intro}
    \norme{\na_x u_\egd}_{\L^1(0,T;\L^\infty(\T^3))}\leq 1/30 \text{ and }\eqref{hypFujitaKato-Intro}\text{ holds}.
\end{equation}
 Because of the previous considerations, for  a given strong existence time $T>0$,  on the interval $[0,T]$, one has a uniform control on $\| \rho_\egd\|_{\L^\infty(0,T ; \L^\infty(\T^3))}$ which in turn allows to prove the uniform decay for $\EEE_\egd$. 


After proving the existence of a strong existence time (that may at this stage depend on $\eps$), we finally define
\[
    T^*_\eps=\sup\{T>0,\,T\text{ is a strong existence time}\}.
\]
Following these preliminaries, the main part of the analysis is to prove that $T^*_\eps= +\infty$ under the smallness assumption for the initial fluid velocity. We also show that we can dispense with this assumption and still find a strong existence time that does not depend on $\eps$.

The \emph{light} and \emph{light and fast} particle regimes can be treated, at first, in a unified manner, based on this $\L^p \L^p$ strategy. It however requires a well-posedness assumption for treating the cases $\alpha>1/p$. To remove it, a refined approach based on $\L^{\frac{p}{p-1}} \L^p$ parabolic estimates is subsequently developed (see below).

A similar $\L^p \L^p$  strategy could be adapted to close the bootstrap argument in the fine particle regime. Nevertheless, we choose to develop a different argument which has two advantages:
\begin{itemize}
\item it directly yields pointwise in time convergence for the full distribution function towards the limit Dirac mass;
\item it could provide some insight as to how to deal with this limit in the whole space $\R^3$ instead of $\T^3$.
\end{itemize} 

\bigskip

\noindent {\bf The bootstrap argument in the \emph{light} and \emph{light and fast} particle regimes.}  Section~\ref{SectionLandLFParticle} is dedicated to both the light particle and light and fast particle regimes. To ease readability, we drop the parameter $\gamma=1$ in the following formulas.
We shall focus on the proof that $T^\star_\eps<+\infty$ under suitable assumptions on the initial data. We assume by contradiction that $T^\star_\eps =+\infty$ and aim at proving that for all $T<T^\star_\eps$, the condition~\eqref{strong-intro} can be improved, which by a continuity argument, will lead to a contradiction with the definition of $T^\star_\eps$.

Let us first discuss the $\L^2$ estimate for the Brinkman force $\Fed$. The smallness property one needs to enforce follows from the dissipation $\D_{\ed}$.
  By the Cauchy-Schwarz inequality, we indeed have
$$
\| \Fed\|_{\L^2 (0,T; \L^2(\T^3))} \leq \| \rhoed \|^{1/2}_{\L^\infty(0,T;\L^\infty(\T^3))} \left(\int_0^T \int_{\T^3 \times \R^3} \fed |v-\ued|^2 \, \dd v \dd x \dd t\right)^{1/2},
$$
and we use that since $T<T^\star_\eps$, we do have $ \| \rhoed \|_{\L^\infty(0,T;\L^\infty(\T^3))} \lesssim 1$. By the modulated energy--dissipation identity~\eqref{mod-ener-intro}, for all $T>0$,
$$
\int_0^T \int_{\T^3 \times \R^3} \fed |v-\ued|^2 \, \dd v \dd x \dd t \leq  \EEE_{\ed}(0),
$$
which yields the required bound for $\| \Fed \|_{\L^2 (0,T; \L^2(\T^3))} $, choosing  $\EEE_{\ed}(0)$ small enough.

 The heart of the proof is to obtain a refined estimate for the Brinkman force, which stems from a kind of desingularization of its expression with respect to the small parameter $\eps$. Writing
$$
\Fed(t,x) = \int_{\R^d} \fed (t,x,v) \left(\frac{v}{\sigma} - \ued(t,x)\right) \dd v,
$$
the representation formula~\eqref{characmethod-intro}
and the change of variables in velocity  $w:=\Gamma_{\ed}^{t,x}(v)$ (defined in~\eqref{def-gamma-intro})
yield
\begin{multline*}
    \Fed(t,x)=\\
    e^{\frac{3t}{\eps}}\int_{\R^3}\fed^0\left(\tilde{X}_{\ed}^{t,x,w}(0),w\right)\left(\frac{1}{\sigma}[\Gamma_{\ed}^{t,x}]^{-1}(w)-\ued(t,x)\right)|\det\nabla_{w}[\Gamma_{\ed}^{t,x}]^{-1}(w)|\,\dd w.
\end{multline*}
The key is the following identity, obtained from the equation of characteristics~\eqref{charac-intro} and a integration by parts in time, taking advantage of the fast effect of the friction term, embodied by the presence of a  factor $1/\eps$ in the exponential functions in time:
\begin{align*}
    \frac{1}{\sigma}[\Gamma_{\ed}^{t,x}]^{-1}(w)
    -\ued(t,x)&= e^{-\frac{t}{\eps}}\left(\frac{w}{\sigma}-\ued^0\left(\tilde{X}_{\ed}^{t,x,w}(0)\right)\right)\\
    -\int_0^t&e^{\frac{s-t}{\eps}}\partial_s\ued\left(s,\tilde{X}_{\ed}^{t,x,w}(s)\right)\dd s\\
    -\int_0^t&e^{\frac{s-t}{\eps}}V_{\ed}\left(s;t,x,[\Gamma_{\ed}^{t,x}]^{-1}(w)\right)\cdot\na_{x}u_\ed \left(s,\tilde{X}_{\ed}^{t,x,w}(s)\right)\dd s.
\end{align*}
We therefore end up with
    \begin{equation}
\label{decomp-F-intro}
        \Fed (t,x)= \Fed^0+\Fed^{dt}+\Fed^{dx},
    \end{equation}
    where
    \begin{align*}
          \Fed^0&=e^{\frac{-t}{\eps}}\int_{\R^3}\fed^0\left(\tilde{X}_{\ed}^{t,x,w}(0),w\right)\left[\frac{w}{\sigma}-\ued^0\left(\tilde{X}_{\ed}^{t,x,w}(0)\right)\right]\dd w, \\
          \Fed^{dt}&=\int_{\R^3}\int_0^te^{\frac{s-t}{\eps}}\fed^0\left(\tilde{X}_{\ed}^{t,x,w}(0),w\right)\partial_s\ued\left(s,\tilde{X}_{\ed}^{t,x,w}(s)\right)\dd s\dd w, \\
          \Fed^{dx}&=
        \int_{\R^3}\int_0^te^{\frac{s-t}{\eps}}\fed^0\left(\tilde{X}_{\ed}^{t,x,w}(0),w\right)\\
           &\qquad\qquad\qquad V_{\ed}\left(s;t,x,[\Gamma_{\ed}^{t,x}]^{-1}(w)\right)\cdot\na_{x}\ued\left(s,\tilde{X}_{\ed}^{t,x,w}(s)\right)\dd s\dd w.
    \end{align*}
    Let $p>3$ be the regularity index of Assumption~\ref{hypGeneral}.
Building on this decomposition, one proves (see Lemmas~\ref{lemmaF0-12}, \ref{lemmaFdt-12} and~\ref{lemmaFdx-12}) that
\begin{equation}
\label{estimFed-intro}
\begin{aligned}
\|  \Fed\|_{\L^p((0,T) \times \T^3)} \| 
    &\lesssim  \frac{\eps^{\frac{1}{p}}}{\sigma}\norme{|v|^p\fed^0}_{\L^1(\T^3\times\R^3)}^{\frac{1}{p}}+ \eps+\eps \EEE_{\ed}(0)^{\frac{1}{2}}  \\
    &+  \eps \norme{\partial_t\ued}_{\L^p((0,T)\times\T^3)} + \eps\norme{\Delta_x\ued}_{\L^p((0,T)\times\T^3)}.
\end{aligned}
\end{equation}
Combining this with the aforementioned $\L^p$ maximal estimate and taking $\eps$ small enough to absorb some terms of the right-hand side, one finally obtains
\begin{equation}\label{estimPara-12-intro}
    \norme{\partial_t\ued }_{\L^p((0,T)\times\T^3)}+\norme{\Delta_x\ued}_{\L^p((0,T)\times\T^3)}\lesssim 1+\frac{\eps^{\frac{1}{p}}}{\sigma}\norme{|v|^p\fed^0}_{\L^1(\T^3\times\R^3)}^{\frac{1}{p}}.
\end{equation}
We will therefore need an additional \emph{well-prepared} assumption on the initial distribution function $\fed^0$ depending on $\sigma=\eps^{\alpha}$ in the light and fast regime, namely
\begin{equation}
\label{additional-intro}
\frac{\eps^{\frac{1}{p}}}{\sigma}\norme{|v|^p\fed^0}_{\L^1(\T^3\times\R^3)}^{\frac{1}{p}} \lesssim 1,
\end{equation}
in order to make sure that the second term is bounded with respect to $\eps$. As $p$ may be taken arbitrarily close to $3$, this is significant only for the parameter range $\alpha>1/3$.
In particular, note that this is not necessary in the light regime (that corresponds to $\alpha=0$). A refined strategy, based on anisotropic $\L^{\frac{p}{p-1}} \L^p$ parabolic estimates will later allow us to dispense with the additional assumption~\eqref{additional-intro} to conclude the bootstrap argument; it is exposed in Section~\ref{sec-betterLF} (see also below for a short summary). 

To conclude, one may interpolate this higher order estimate with the $\L^2$ norm of $\ued-\langle \ued\rangle$ thanks to the Gagliardo-Nirenberg inequality:
   $$
   \norme{\na_x\ued}_{\L^\infty(\T^3)} \lesssim \| \ued-\langle \ued\rangle\|_{\L^2(\T^3)}^{1-\beta_p} \norme{\Delta_x \ued }_{\L^p(\T^3)}^{\beta_p},
   $$
   for some $\beta_p \in (0,1)$. Using the exponential decay of $\EEE_{\ed}(t)$ (recall~\eqref{DE-intro}), we therefore get
    \begin{multline*}
    \int_0^T\norme{\na_x\ued(t)}_{\L^\infty(\T^3)}\dd t\lesssim\EEE_{\ed}(0)^{\frac{1-\beta_p}{2}}\norme{\Delta_x\ued}_{\L^p((0,T)\times\T^3)}^{\beta_p}
    \lesssim \EEE_{\ed}(0)^{\frac{1-\beta_p}{2}},
\end{multline*}
which can be made as small as desired by picking $\EEE_{\ed}(0)$ small enough. This allows to conclude the bootstrap argument and conclude that necessarily, $T^\star_\eps =+\infty$.
    
Section~\ref{SectionAsymptotic-12} is then dedicated to the proofs of convergence. We first assume that $(\rhoed^0,\ued^0)$ weakly converges to $(\rho^0, u^0)$. Now that we have proven that $T^\star_\eps=+\infty$ for any small $\eps>0$, we can combine the energy--dissipation identity~\eqref{energy-dissipation-intro} and the bound~\eqref{estimPara-12-intro} to deduce from the Aubin-Lions lemma the strong convergence of $(\ued)$ in $\L^2((0,T)\times\T^3)$, for any $T>0$. The control~\eqref{estimPara-12-intro} can also be used to prove the vanishing of $\Fed$ in $\L^p((0,T)\times\T^3)$ thanks to~\eqref{estimFed-intro}, up to a slightly stronger additional assumption in the light and fast regime. This suffices to take the limit in the Navier-Stokes equations and the conservation of mass, and thus obtain the convergence of $(\rhoed,\ued)$ to the solution to the Transport-Navier-Stokes system~\eqref{TNS}.
    
If we further assume the strong convergence of $(\ued^0)$ to $u^0$, then the structure of the Navier-Stokes (resp. transport) equations satisfied by ($\ued$) and $u$ (resp. $(\rhoed)$ and $\rho$) provide a quantitative result for the pointwise convergence of $(\ued)$ in $\L^2(\T^3)$ (resp. $(\rhoed)$ in the Wasserstein-1 metric). The quantitative convergence result for $(\fed)$, up to an integration in time, is then straightforward.

In order to obtain pointwise convergence for $(\fed)$, we note that for any test function $\psi$, and almost every $t\in(0,T)$, we have
\[
    \left|\langle\fed(t)-\rho(t)\otimes\delta_{v=\sigma u(t)},\psi\rangle\right|\le\norme{\na_v\psi}_{\L^\infty(\T^3\times\R^3)}\int_{\T^3\times\R^3}\fed(t)|v-\sigma\ued(t)|\dd x\dd v,
\]
which is reminiscent of the expression of the Brinkman force. Indeed, using the decomposition described above, we derive a quantitative convergence result in the case of well-prepared initial data.

In Section~\ref{sec-smallT}, we explain how to dispense with the assumptions of smallness for the initial fluid velocity and modulated energy and prove short time convergence results. The initial bounds on the initial data in Assumption~\ref{hypGeneral} are in fact sufficient to ensure~\eqref{hypFujitaKato-Intro} holds for every $\eps\in(0,1)$, provided the time $t$ is less than some small time $T_M>0$ that is independent of $\eps$. We then apply the same bootstrap argument as above and prove that, up to reducing the value of $T_M$ (regardless of $\eps$), every time $T\in[0,T_M]$ is a strong existence time. The proof of convergence is then identical to that exposed in the previous paragraphs.

\bigskip

\noindent {\bf The bootstrap argument in the fine particle regimes.}  In Section~\ref{SectionFineParticle}, we tackle the most singular high friction regime studied in this work, that is the fine particle regime. To ease readability, we drop the parameters $\sigma=1$ and $\gamma=\eps$ in this presentation. As already explained, the difficulty lies in the fact that the Brinkman force in this regime, namely
$$
F_\eps(t,x) = \frac{1}{\eps} \int_{\R^d} f_\eps (t,x,v) (v - u_\eps(t,x)) \dd v,
$$
is singular with respect to $\eps$.
By opposition to the light particle case, thanks to the  modulated energy--dissipation identity and assuming the global control 
$$ \| \rhoeps \|_{\L^\infty(0,+\infty \times \T^3)} \leq \C_0,$$ one only obtains the bound
\[ 
    \norme{F_\eps}_{\L^2(\R_+\times\T^3)}\lesssim\frac{1}{\sqrt{\eps}} \EEE_\eps(0),
\]
which in general blows up as $\eps \to 0$, and is therefore not satisfactory in view of the $\L^\infty\H^1$--$\L^2\dot\H^2$ estimate for the fluid velocity that we aim at establishing.  Therefore, the approach developed for the light particle regimes has to be modified.
 
We propose another desingularization procedure of the Brinkman force, that somehow reflects the fine algebraic structure of the Vlasov-Navier-Stokes system.
Inspired by \cite{hank}, the idea is to consider the \emph{higher dissipation} functional
\[
    \D_{\eps}^{(r)}(t):=\int_{\T^3\times\R^3}\feps\frac{|v-\ueps(t,x)|^r}{\eps^r}\dd x\dd v.
\]
for $r\geq 2$. Such an object is useful in view of studying the Brinkman force, because of the  estimate 
\[
    \norme{F_\eps}_{\L^r((0,T)\times\T^3)}\le\norme{\feps^0}_{\L^1(\R^3;\L^\infty(\T^3))}^{1-\frac{1}{r}}\left(\int_0^T\D_{\eps}^{(r)}(t)\,\dd t\right)^{\frac{1}{r}},
\]
see Lemma~\ref{LienFD-3}. The key is an identity satisfied by $\int_0^t  \D_{\eps}^{(r)}(t) \dd t$, similar to the ones introduced in~\cite{hank}, which is provided in Lemma~\ref{fine-key}.
This results in the following estimate:
 \begin{multline*}
    \int_0^t\D_\eps^{(r)}(s)\dd t+\int_{\T^3\times\R^3}\feps(t)\frac{|v-\ueps(t)|^r}{\eps^{r-1}}\dd x\dd v\\
    \lesssim\norme{\feps^0}_{\L^1(\R^3;\L^\infty(\T^3))}\norme{\partial_t\ueps}_{\L^r((0,T)\times\T^3)}^r\\
    +\norme{|\na_x\ueps|m_r^{\frac{1}{r}}}_{\L^r((0,T)\times\T^3)}^r
    +\int_{\T^3\times\R^3}\feps^0\frac{|v-\ueps^0|^r}{\eps^{r-1}}\dd x\dd v,
\end{multline*}
that we use for $r=p$ (where $p$ is the regularity index of Assumption~\ref{hypGeneral}) to get, assuming the smallness condition
 $$
 \norme{\feps^0}_{\L^1(\R^3;\L^\infty(\T^3))}\ll 1,
 $$
the well-preparedness assumption
\[
    \int_{\T^3\times\R^3}\feps^0\frac{|v-\ueps^0|^r}{\eps^{r-1}}\dd x\dd v \lesssim 1,
\]
 and imposing $\eps$ small enough,
\[
    \norme{\partial_t u_\eps }_{\L^p((0,T)\times\T^3)}+\norme{\Delta_x u_\eps}_{\L^p((0,T)\times\T^3)}\lesssim 1.
\]
We point out that the smallness condition on $\feps^0$ is used only in this argument. To prove the bound~\eqref{hyp-na-u-intro}, we also require the smallness of the initial modulated energy. Finally, by a similar analysis for $r=2$, one obtains a relevant control of the Brinkman force in $\L^2((0,T) \times \T^3)$, which allows to conclude the bootstrap argument in this case.

Section~\ref{sec-first-fine} is dedicated to a non-quantitative convergence result. We obtain the strong convergence of $(\ueps)$ and the weak convergence of $(\rhoeps)$ as in the \emph{light} and \emph{light and fast} particle regimes. Furthermore, the energy--dissipation inequality yields the weak convergence of $(\jeps)$. Because of the singularity of the Brinkman force, this does not imply that $(F_\eps)$ vanishes at $\eps\to0$. But thanks to the bounds derived for the bootstrap argument, we still know that $(F_\eps)$ converges weakly. This suffices to take the limit $\eps\to0$ in the conservations of mass and momentum as well as the Navier-Stokes equations and thus obtain the convergence of $(\rhoeps,\ueps)$ to the solution to the Inhomogeneous incompressible Navier-Stokes system. As in the \emph{light} and \emph{light and fast} particle regimes, a convergence result for $(\feps)$, up to an integration in time, then follows.

If we further assume the strong convergence of $(\ueps^0)$ to $u^0$ and some regularity on the limit solution $u$, then we can obtain a pointwise convergence result for $(\rhoeps)$ and, thanks to the well-preparedness assumption and the key estimate on the higher order dissipation, a pointwise convergence result for $(\feps)$. Note that, at this stage, we do not have quantitative estimates on $\ueps-u$, which prevents us from getting a quantitative result for $(\rhoeps)$ and $(\feps)$ as well.

This indeed appears more difficult and is the object of Section~\ref{sec-relat-fine}. We introduce yet another functional, the so-called \emph{relative entropy}
\begin{equation*}
    \begin{aligned}
        \mathscr{H}_\eps(t) &:=\frac{1}{2}\int_{\T^3\times\R^3}\feps(t,x,v)|v-u(t,x)|^2\dd x\dd v
                    +   \frac{1}{2}\int_{\T^3}|\ueps(t,x)-u(t,x)|^2\dd x.
    \end{aligned}
\end{equation*} 
As stated above, such functionals have already proven useful in many related high friction limits, see e.g. \cite{gou-jab-vas04b,fig-kan}.
A side goal of this section is to explain why it cannot be \emph{directly} used in the fine particle regime that we study. Concretely, it turns out that a bound on the density
\begin{equation}
\label{eq-condi-rho-intro}
    \| \rhoeps \|_{\L^\infty(0,T \times \T^3)} \leq \C_0,\quad \forall \eps \in (0,1),
\end{equation}
which is exactly the same as that required in the other parts of the analysis, is required to complete the analysis. It appears that such a bound cannot be enforced with a proof solely based on relative entropy. But obtaining this bound has been precisely the object of the previous analysis.
 
By explicit computations (see Lemma~\ref{lem-evolution-HH}) one gets 
\begin{equation*}
    \begin{aligned}
        \mathscr{H}_\eps(t) &+ \int_0^t \int_{\T^3}|\na_x(\ueps - u )|^2\dd x \dd s + \frac{1}{\eps} \int_0^t \int_{\T^3 \times \R^3}|v-\ueps|^2\feps \dd x\dd v \dd s \\
        &\leq \Hed(0) + \int_0^t  \sum_{j=1}^{4} I_j(s)  \dd s.
    \end{aligned}
\end{equation*}
with
\begin{align*}
    &I_1 := - \int_{\T^3 \times \R^3}  f_\eps (v-u) \otimes (v-u) : \na_x u\, \dd x \dd v, \\
    &I_2 := - \int_{\R^3}  (\ueps - u ) \otimes (\ueps -u) : \na_x u \,\dd x, \\
    &I_3 :=  \int_{\T^3 \times \R^3}  \feps (v-\ueps) \cdot G\,  \dd x \dd v , \\
    &I_4 :=  \int_{\T^3} (\rhoeps- \rho) (\ueps -u ) \cdot G\, \dd x,
\end{align*}
where $G = \frac{\na_x p - \Delta_x u}{1+\rho}$. The terms $I_1$ and $I_2$ are the easiest to handle. To study $I_3$ one may use either the energy--dissipation inequality or, better, the dissipation $\D_\eps^{(2)}$ to obtain the refined inequality
\begin{multline*}
    \int_0^T|I_3(t)|\dd t\le\norme{\rhoeps}_{\L^\infty((0,T)\times\T^3)}^{\frac{1}{2}}\left(\eps^2\int_0^T\D_\eps^{(2)}(t)\dd t\right)^{\frac{1}{2}}\norme{G}_{\L^2((0,T)\times\T^3)}
    \lesssim \eps.
\end{multline*}
The control of $I_4$ relies on the following bound, which is valid under a regularity assumption on $(\rho,u)$ and on the condition that~\eqref{eq-condi-rho-intro} holds:
    \begin{multline}
        \norme{\rhoeps(t)-\rho(t)}_{\dot{\H}^{-1}(\T^3)}
        \lesssim\norme{\rhoeps^0-\rho^0}_{\dot{\H}^{-1}(\T^3)}\\
            +\eps\int_0^t\norme{F_\eps(s)}_{\L^2(\T^3)}\dd s+\int_0^t\norme{\ueps(s)-u(s)}_{\L^2(\T^3)}\dd s.
    \end{multline}
Gathering all pieces together, we finally obtain the estimate
    \[
        \mathscr{H}_\eps(t)\lesssim \norme{\ueps^0-u^0}_{\L^2(\T^3)}^2
            +\norme{\rhoeps^0-\rho^0}_{\dot\H^{-1}(\T^3)}^2+\eps,
    \]
    from which quantitative convergence bounds follow.
    
As in the \emph{light} and \emph{light and fast} particle regimes, we eventually explain in Section~\ref{sec-smallT-fine} how to dispense with the smallness assumptions on the initial fluid velocity and modulated energy, at the expense of a short time constraint.
 
 \bigskip
 
 \noindent {\bf Improvements in the \emph{light} and \emph{light and fast} particle regimes.} 
 In Section~\ref{sec-betterLF}, we come back to the \emph{light} and \emph{light and fast} particle regimes. 
 The main goal is to remove the unnecessary assumption~\eqref{additional-intro} to conclude the bootstrap argument. To this end, as already briefly said, the idea is to rather use a $\L^{\frac{p}{p-1}}(0,T;\L^p(\T^3))$ parabolic estimate instead of the $\L^p((0,T)\times \T^3)$ estimate of Section~\ref{SectionLandLFParticle}.

In Section~\ref{sec-LF-nohyp}, we provide the key $\L^{\frac{p}{p-1}}(0,T;\L^p(\T^3))$ estimate on the Brinkman force, replacing the former~\eqref{estimFed-intro} by
\begin{multline*}
        \norme{\Fed}_{\L^{\frac{p}{p-1}}(0,T;\L^p(\T^3))} \lesssim M^{\mu_p}\frac{\eps^{\frac{1}{p}}}{\sigma^{\frac{p}{p-1}}}+ \eps M^{\mu_p}\EEE_{\ed}(0)^{\frac{1}{2}} \\
        + \eps M\norme{\partial_t\ued}_{\L^{\frac{p}{p-1}}(0,T;\L^p(\T^3))}+ \eps\norme{\Delta_x\ued}_{\L^{\frac{p}{p-1}}(0,T;\L^p(\T^3))},
\end{multline*}
which is based on the same desingularization as in~\eqref{decomp-F-intro}. We then explain how to conclude the bootstrap argument and to adapt the convergence results in Section~\ref{sec-conv12-re}.

In a short Section~\ref{SectionRetour-12}, we explain how the new objects introduced to study the fine particle regime, namely the higher dissipation and relative entropy functionals, can be used to provide alternative proofs for the well-prepared results in the \emph{light} and \emph{light and fast} particle regimes.

\bigskip

\noindent {\bf Further developments.} 
In Section~\ref{sec-further}, we conclude the paper by providing some further applications of the methods developed in this work. First, in Section~\ref{sec-dim2}, we briefly explain how our results can extend in the two-dimensional case; loosely speaking, as the Cauchy problem for the Navier-Stokes equation becomes sub-critical, the assumption of small $\dot\H^{\frac{1}{2}}(\T^3)$ norm for the initial velocity field can be dropped. Section~\ref{sec-mix} considers an extension of the Vlasov-Navier-Stokes system~\eqref{VNS} for mixtures of dispersed phases; the model which we introduce allows fragmentation phenomena, which is inspired from~\cite{ben-des-mou}. We derive two-fluid Inhomogeneous Navier-Stokes systems in the limit.
Section~\ref{sec-BNS} discusses the derivation of Boussinesq type equations for the equations set on the whole space and in the presence of an additional (constant) gravity force; we explain how our techniques yield a short time result.
Finally Section~\ref{sec-open}  discusses some open problems that are related to this work.


\section{Preliminary results}
\label{sec-prelim}

Let us start the analysis by introducing in a unified manner the first common arguments for all the regimes studied in this work.
We recall that we aim at studying 
\begin{equation}
    \label{VNS-general-re}
    \left\{
\begin{aligned}
&\partial_t u_{\eps,\gamma,\sigma}+ (u_{\eps,\gamma,\sigma}\cdot\na_{x})u_{\eps,\gamma,\sigma}-\Delta_{x}u_{\eps,\gamma,\sigma}+\na_{x}p_{\eps,\gamma,\sigma}= \frac{1}{\gamma}\left(j_{f_{\eps,\gamma,\sigma}}-\rho_{f_{\eps,\gamma,\sigma}}u_{\eps,\gamma,\sigma}\right),\\
&\div_{x}u_{\eps,\gamma,\sigma}=0,\\
  &\partial_tf_{\eps,\gamma,\sigma}+ \frac{1}{\sigma}  v\cdot\na_{x}f_{\eps,\gamma,\sigma}+ \frac{1}{\eps} \div_{v}\left[f_{\eps,\gamma,\sigma}( \sigma u_{\eps,\gamma,\sigma}-v)\right]=0, \\
 &\rho_{f_{\eps,\gamma,\sigma}}(t,x)=\int_{\R^3}f_{\eps,\gamma,\sigma}(t,x,v)\ddv,  \quad
j_{f_{\eps,\gamma,\sigma}}(t,x)=\frac{1}{\sigma} \int_{\R^3}vf_{\eps,\gamma,\sigma}(t,x,v)\ddv.
\end{aligned}
\right.
\end{equation}
Throughout this section, we fix $\eps>0$ and consider
\[
    (\gamma,\sigma)\in\{(1,\eps^\alpha),(\eps,1)\},
\]
for the light ($\alpha=0$), light and fast ($\alpha \in (0,1/2]$), and fine particle regimes, respectively. This part of the paper is organized as follows:
\begin{itemize}
\item Section~\ref{sec-weak} is a short reminder about the notion of weak solutions to the Vlasov-Navier-Stokes system.
\item In Sections~\ref{sec-CL}, \ref{sec-higherNS} and \ref{sec-roughfirst}, we gather several rather standard estimates for the Vlasov equation with friction and Navier-Stokes equations .
\item Section~\ref{sec-modenergy} introduces a key object, the \emph{modulated energy} and proves its exponential decay under the assumption of an $\L^\infty$ bound for the kinetic density.
\item  In Section~\ref{SectionChangeVar}, we study straightening changes of variables in velocity (resp. in space)  which will be used at multiple times to obtain uniform bounds for moments and the Brinkman force. They are in particular the key to the aforementioned $\L^\infty$ bound for the kinetic density.
We also introduce the notion of \emph{strong existence times} which are, loosely speaking, times for which the changes of variables are admissible.
\item In Section~\ref{sec-conv}, we gather some estimates for the convective term in the Navier-Stokes equations.
\item  Finally, we initialize the bootstrap argument in Section~\ref{SubsectionInitialization}.
\end{itemize}

\subsection{Weak solutions}
\label{sec-weak}

We begin by presenting the notion of weak solutions of the Vlasov-Navier-Stokes system and introducing useful notations.

\begin{definition}
    We define the kinetic energy of the system~\eqref{VNS-general}, for every $t\ge0$, by
    \[
        \E_{\egd}(t)=\frac{1}{2}\norme{\uegd(t)}_{\L^2(\T^3)}^2+\frac{\eps}{2\gamma\sigma^2}\int_{\T^3\times\R^3}|v|^2\fegd(t,x,v)\dd x\dd v
        \index{$\E_{\egd}(t)$: energy}
    \]
    and the dissipation by
    \[
        \D_{\egd}(t)=\norme{\na_x\uegd(t)}_{\L^2(\T^3)}^2+\frac{1}{\gamma} \int_{\T^3\times\R^3}\left|\frac{v}{\sigma}-\uegd(t,x)\right|^2\fegd(t,x,v)\dd x\dd v.
        \index{$\D_{\egd}(t)$: dissipation}
    \]
\end{definition}

As already said in the introduction, one can formally check that, for every $t\ge0$, the following identity holds:
\[
    \frac{\dd}{\dd t}\E_{\egd}(t)+\D_{\egd}(t)=0.
\]
This identity plays a crucial role in the analysis of the Vlasov-Navier-Stokes system, be it for the proof of the existence of solutions (\cite{bou-des-gra-mou}, \cite{bou-gra-mou}, \cite{bou-mic-mou}), the long time behavior (\cite{cho-kwo}, \cite{hank}, \cite{hank-mou-moy}) or the asymptotic analysis we conduct in this paper. As such, we shall only consider solutions that verify an inequality version of this energy--dissipation identity. 

\begin{definition}\label{defWeakSol}
    Assume
    \begin{enumerate}
        \item $\uegd^0\in\L^2_{\div}(\T^3)=\{U\in\L^2(\T^3),\,\div_xU=0\}$,
        \item $0\le\fegd^0\in\L^1\cap\L^\infty(\T^3\times\R^3)$,
        \item $(x,v)\mapsto\fegd^0(x,v)|v|^2\in\L^1(\T^3\times\R^3)$,
    \end{enumerate}
    then a global weak solution of the Vlasov-Navier-Stokes system~\eqref{VNS-general} with initial condition $(\uegd^0,\fegd^0)$ is a pair $(\uegd,\fegd)$ such that
    \begin{itemize}
        \item the distribution function $\fegd\in\L^\infty_{\loc}(\R_+;\L^1\cap\L^\infty(\T^3\times\R^3))$ is a renormalized solution of the Vlasov equation,
        \item the fluid velocity $\uegd\in\L^\infty_{\loc}(\R_+;\L^2(\T^3))\cap\L^2_{\loc}(\R_+;\H^1(\T^3))$ is a Leray solution of the Navier-Stokes equations,
        \item $\jegd-\rhoegd u\in\L^2_{\loc}(\R_+;\H^{-1}(\T^3))$,
        \item for almost all $t\ge s\ge0$ (including $s=0$), the energy--dissipation inequality holds:
            \begin{equation}\label{EstimationEnergieEq}
                \E_{\egd}(t)+\int_s^t\D_{\egd}(\tau)\dd \tau\le\E_{\egd}(s).
            \end{equation}
    \end{itemize}
\end{definition}

Such global weak solutions are built in~\cite{bou-des-gra-mou}. Note that Assumption~\ref{hypGeneral} ensures that the requirements of Definition~\ref{defWeakSol} are met. In the rest of this work, we consider a weak solution $(\uegd,\fegd)$ in the sense of Definition~\ref{defWeakSol}.

The energy--dissipation estimate~\eqref{EstimationEnergieEq} has a possible consequence on the estimate of the Brinkman force which we recall is defined as
$$
\Fegd := \frac{1}{\gamma} \left( \jegd- \rhoegd\uegd\right).
\index{$\Fegd$: Brinkman force}
$$
This is the object of the next remark.


\begin{remark}\label{RemarkBorneUnifF}

   Let us assume the global uniform bound 
   $$
   \norme{\rhoegd}_{\L^{\infty}(\R_+; \L^\infty(\T^3))} <+\infty,
   $$
   Then by the Cauchy-Schwarz inequality, we obtain that
   \begin{align*}
   &\qquad  \norme{\Fegd}_{\L^2(\R_+;\L^2(\T^3))}^2 \\
   &\leq \norme{\rhoegd}_{\L^{\infty}(\R_+; \L^\infty(\T^3))} \frac{1}{\gamma^2}\int_0^{+\infty}\int_{\T^3\times\R^3}\left|\frac{v}{\sigma}-\uegd(t,x)\right|^2\fegd(t,x,v)\dd x\dd v \dd t\\
   &\leq  \frac{1}{\gamma}  \norme{\rhoegd}_{\L^{\infty}(\R_+; \L^\infty(\T^3))} \int_0^{+\infty} \D_{\egd}(t) \dd t \\
   &\leq \frac{1}{\gamma}  \norme{\rhoegd}_{\L^{\infty}(\R_+; \L^\infty(\T^3))} \E_{\egd}(0),
 \end{align*}
 thanks to the energy--dissipation estimate~\eqref{EstimationEnergieEq}. 
  In the \emph{light} and \emph{light and fast} particle regimes, this provides a uniform bound for $\Fegd$ in $\L^2(\R_+\times\T^3)$; in sharp contrast, in the fine particle regime it only yields
    \[
        \norme{F_{\egd}}_{\L^2(\R_+\times\T^3)}\lesssim\frac{1}{\sqrt{\eps}},
    \]
    which may blow up as $\eps\to0$ and is thus not satisfactory. This is a first indication that the fine particle regime is more singular.
\end{remark}



\subsection{Conservation laws for the Vlasov equation}
\label{sec-CL}

In this paragraph, we focus on the Vlasov equation and state basic conservation laws and estimates that can be justified thanks to the DiPerna-Lions theory \cite{dip-lio}, as described in the following remark.

\begin{remark}\label{remarkDiPernaLions}
    As in \cite[Remark 3.1]{hank-mou-moy}, we shall use the DiPerna-Lions theory~\cite{dip-lio} repeatedly without writing down the argument explicitly. The \emph{implicit} argument will always be:
    \begin{itemize}
        \item consider a sequence of regularized initial data $(f_n^0)_{n\in\N}$, an approximating sequence of fluid velocities $(u_n)_{n\in\N}$, and the associated smooth solutions $(f_n)_{n\in\N}$ of the Vlasov equation;
        \item prove the desired estimate for $f_n$;
        \item pass to the limit using the strong stability property of renormalized solutions \cite{dip-lio}.
    \end{itemize}
    
    Eventually, we will prove that the velocity field $\uegd$ is smooth enough to apply the classical Cauchy-Lipschitz theorem, and therefore the DiPerna-Lions theory will not actually be needed.
\end{remark}

Let us first define the notations associated to the characteristics of the equation. 

\begin{definition}\label{defCharac}
    Let $t\ge0$, $x\in\T^3$ and $v\in\R^3$. We denote by 
    \[
        (X_{\egd}(\cdot;t,x,v),V_{\egd}(\cdot;t,x,v))
    \]
      \index{X@$(X_{\egd}(\cdot;t,x,v),V_{\egd}(\cdot;t,x,v))$: characteristic curves associated to the Vlasov equation}
    the solution to the system of differential equations
    \begin{equation}\label{characV}
        \left\{\begin{aligned}
            \dot{X}_{\egd}(s;t,x,v)&=\frac{1}{\sigma}V_{\egd}(s;t,x,v),\\
            \dot{V}_{\egd}(s;t,x,v)&=\frac{1}{\eps}\left(\sigma\uegd(s,X_{\egd}(s;t,x,v))-V_{\egd}(s;t,x,v)\right),\\
            X_{\egd}(t;t,x,v)&=x,\\
            V_{\egd}(t;t,x,v)&=v.
        \end{aligned}\right.
    \end{equation}
\end{definition}

We can now state the following conservation laws of the Vlasov equation.

\begin{lemma}\label{VlasovConservation}
    Under Assumption \ref{hypGeneral}, for all $\eps>0$, for almost any $t\ge0$,
    \[
        \fegd(t,x,v)\ge0\qquad a.e.\,(x,v)\in\T^3\times\R^3,
    \]
    and
    \[
        \int_{\T^3\times\R^3}\fegd(t,x,v)\dd x\dd v=\int_{\T^3\times\R^3}\fegd^0(x,v)\dd x\dd v.
    \]
    Furthermore, the equations for local conservation of mass and momentum read
    \begin{equation}\label{conservationMass}
        \partial_t\rhoegd+\div_x\jegd=0,
    \end{equation}
    and
    \begin{equation}\label{conservationMomentum}
        \partial_t\jegd+\frac{1}{\sigma^2}\div_x\left(\int_{\R^3}v\otimes v\,\fegd\dd x\dd v\right)=-\frac{\gamma}{\eps}\Fegd.
    \end{equation}
\end{lemma}

\begin{preuve}
    Let us assume that both $\uegd$ and $\fegd$ are smooth functions, the general case follows by the DiPerna-Lions theory (recall Remark~\ref{remarkDiPernaLions}). Then, thanks to the method of characteristics, we have, for $t\ge0$ and $(x,v)\in\T^3\times\R^3$,
    \[
        \fegd(t,x,v)=e^{\frac{3t}{\eps}}\fegd^0(X_{\egd}(0;t,x,v),V_{\egd}(0;t,x,v)),
    \]
    hence the nonnegativity of $\fegd$.
    
    The second statement stems from the integration of the Vlasov equation
    \begin{equation}
    \label{V1}
        \partial_tf_{\eps,\gamma,\sigma}+ \frac{1}{\sigma}  v\cdot\na_{x}f_{\eps,\gamma,\sigma}+ \frac{1}{\eps} \div_{v}\left[f_{\eps,\gamma,\sigma}( \sigma u_{\eps,\gamma,\sigma}-v)\right]=0
    \end{equation}
    over $(0,t)\times\T^3\times\R^3$, for any $t\ge0$, while the conservation of mass results from integrating only over $\R^3$ and the conservation of momentum comes from multiplying~\eqref{V1} by $v/\sigma$ and integrating over $\R^3$.
\end{preuve}

The maximum principle for Vlasov equations yields the following inequality.

\begin{lemma}\label{principe-max}
    Under Assumption~\ref{hypGeneral}, we have
    \[
        \forall t\ge0,\qquad\norme{\fegd}_{\L^\infty(\T^3\times\R^3)}\le e^{\frac{3t}{\eps}}\norme{\fegd^0}_{\L^\infty(\T^3\times\R^3)}.
    \]
\end{lemma}

\begin{preuve}
    This is a straightforward consequence of the method of characteristics.
\end{preuve}

\begin{remark}
    The blow-up as $\eps$ goes to 0 in the right-hand side of this inequality is consistent with the expected convergence of $(\fegd)$ to a Dirac mass, as presented in Section~\ref{sec-formal}.
\end{remark}

Finally, 
we shall need the following conservation property of the total momentum.
\begin{lemma}\label{conservationVNSmoy}
    Under Assumption~\ref{hypGeneral}, for every $t\ge0$,
    \[
        \left\langle\frac{\eps}{\gamma}\jegd(t)+\uegd(t)\right\rangle=\left\langle\frac{\eps}{\gamma}\jegd^0+\uegd^0\right\rangle.
    \]
\end{lemma}

\begin{preuve}
    Let $t\ge0$. By integrating the Navier-Stokes equation and~\eqref{conservationMomentum} on $(0,t)\times\T^3$, we obtain
    \[
        \langle\uegd(t)\rangle=\langle\Fegd(t)\rangle\qquad\text{and}\qquad\langle\jegd(t)\rangle=\frac{-\gamma}{\eps}\langle\Fegd(t)\rangle,
    \]
    hence the result.
\end{preuve}


\subsection{Higher regularity estimate for the Navier-Stokes equations}
\label{sec-higherNS}

In this section, we present a higher regularity estimate for the Navier-Stokes equations. We also introduce some notations that we use to abbreviate the results given in Appendix \ref{RappelsNS}.

According to Theorem~\ref{propTpsLong}, the estimate
\begin{equation}\label{ENSH1}
    \norme{\uegd(t)}_{\H^1(\T^3)}^2+\int_0^t\norme{\Delta_x\uegd(s)}_{\L^2(\T^3)}^2\dd s\lesssim\Psi_{\egd,0},
\end{equation}
where
\[
\index{P@$\Psi_{\egd,0}$: control of the $\H^1$ energy estimate}
    \Psi_{\egd,0}=\norme{\uegd^0}_{\H^1(\T^3)}^2+\norme{\Fegd}_{\L^2((0,t)\times\T^3)}^2.
\] 
holds for $t\ge0$ if, for some universal constant $C^*>0$,
\begin{equation}\label{hypNSSup}
    \int_0^t\norme{e^{s\Delta}\uegd^0}_{\dot\H^1(\T^3)}^4\dd s+\int_0^t\norme{\Fegd(s)}_{\dot\H^{-1/2}(\T^3)}^2\dd s\le C^*
\end{equation}
is satisfied.

\begin{remark}
The  subscript $0$ in $\Psi_{\egd,0}$ may be misleading as in the definition of  $\Psi_{\egd,0}$, the term $\norme{\Fegd}_{\L^2((0,t)\times\T^3)}^2$ depends on time $t$. However we shall always consider so-called strong existence times (see Definition~\ref{StrongExistenceTime}) for which 
$\norme{\Fegd}_{\L^2((0,t)\times\T^3)}^2 \leq C$,
for some fixed $C>0$, so that the notation makes sense.
\end{remark}

Thanks to the Sobolev embeddings $\H^1(\T^3)\hookrightarrow\L^r(\T^3)$ ($r\le6$), we draw the following first consequence of the higher regularity estimate~\eqref{ENSH1}.

\begin{corollary}\label{EstimNormLp}
    Under Assumption~\ref{hypGeneral} and if $T>0$ is such that~\eqref{ENSH1} holds on $[0,T]$, then, for every $2\le r\le 6$,
    \[
        \norme{\uegd}_{\L^\infty(0,T;\L^r(\T^3))}\lesssim\Psi_{\egd,0}^{\frac{1}{2}}.
    \]
\end{corollary}

Now we investigate when~\eqref{hypNSSup} is satisfied. To this end, we start with the following preliminary result which corresponds to \cite[Lemma 4.3]{hank-mou-moy}.
\begin{lemma}
\label{lem-ouf}
  Under Assumption~\ref{hypGeneral}, for all $T>0$, ${\Fegd} \in \L^2(0,T; \L^{3/2} (\T^3))$.
\end{lemma}

\begin{preuve}
    We follow~\cite[Lemmas 4.2 and 4.3]{hank-mou-moy} and only provide the proof for the sake of completeness. First, note that, for any $T>0$, thanks to Hölder's inequality and the Sobolev embedding $\H^1(\T^3)\hookrightarrow\L^6(\T^3)$, we have
    \[
        \norme{\rhoegd\uegd}_{\L^2(0,T;\L^{\frac{3}{2}}(\T^3))}\lesssim\norme{\rhoegd}_{\L^\infty(0,T;\L^2(\T^3))}^{\frac{1}{2}}\norme{\uegd}_{\L^2(0,T;\H^1(\T^3))},
    \]
    so we only need to prove that 
    $$\rhoegd\in\L^\infty(0,T;\L^2(\T^3)) \text{   and   } \jegd\in\L^\infty(0,T;\L^{\frac{3}{2}}(\T^3)).$$
     Recall the following interpolation estimate: for any $0\le\ell\le k$ and any nonnegative $g\in\L^\infty(\T^3\times\R^3)$,
    \begin{equation}\label{interpolationDistrib}
        \norme{m_{\ell}g}_{\L^{\frac{k+3}{\ell+3}}(\T^3)}\lesssim(M_kg)^{\frac{\ell+3}{k+3}}\norme{g}_{\L^\infty(\T^3\times\R^3)}^{\frac{k-\ell}{k+3}},
    \end{equation}
    where
    \[
        m_\ell g=\int_{\R^3}|v|^{\ell}g(\cdot,v)\dd v\qquad\text{and}\qquad M_kg=\int_{\T^3\times\R^3}|v|^kg(x,v)\dd x\dd v.
    \]
    Applying~\eqref{interpolationDistrib} for $\fegd$ with $(\ell,k)=(0,3)$ and $(\ell,k)=(1,3)$ yields, for almost every $t\in(0,T)$,
    \[
        \norme{\rhoed(t)}_{\L^2(\T^3)}=\norme{m_0\fegd(t)}_{\L^2(\T^3)}\lesssim (M_3\fegd)^{\frac{1}{2}}\norme{\fegd(t)}_{\L^\infty(\T^3\times\R^3)}^{\frac{1}{2}},
    \]
    and
    \begin{align*}
        \norme{\jegd(t)}_{\L^{\frac{3}{2}}(\T^3)}&\le\sigma^{-1}\norme{m_1\fegd(t)}_{\L^{\frac{3}{2}}(\T^3)}\\
        &\lesssim\sigma^{-1}(M_3\fegd(t))^{\frac{2}{3}}\norme{\fegd(t)}_{\L^\infty(\T^3\times\R^3)}^{\frac{1}{3}}.
    \end{align*}
    Thanks to Lemma~\ref{principe-max}, we know that $\fegd\in\L^\infty(0,T;\L^\infty(\T^3\times\R^3))$ so that we shall focus on proving that $M_3\fegd\in\L^\infty((0,T))$.
    
    Multiplying the Vlasov equation by $|v|^3$ and integrating over $\T^3\times\R^3$, we get
    \begin{equation}\label{eq-M3-Vlasov}
        \frac{\dd}{\dd t}M_3\fegd(t)+\frac{3}{\eps}M_3\fegd(t)\le\frac{3\sigma}{\eps}\int_{\T^3}|\uegd(t)|\,|m_2\fegd(t)|\dd x.
    \end{equation}
    We apply~\eqref{interpolationDistrib} with $(\ell,k)=(2,3)$ and Lemma~\ref{principe-max} to obtain, for almost every $t\ge0$,
    \begin{align*}
        \norme{m_2\fegd(t)}_{\L^{\frac{6}{5}}(\T^3)}&\lesssim\left(M_3\fegd(t)\right)^{\frac{5}{6}}\norme{\fegd(t)}_{\L^\infty(\T^3\times\R^3)}^{\frac{1}{6}}\\
        &\lesssim e^{\frac{t}{2\eps}}\norme{\fegd^0}_{\L^\infty(\T^3\times\R^3)}^{\frac{1}{6}}\left(M_3\fegd(t)\right)^{\frac{5}{6}}.
    \end{align*}
    Therefore, thanks to Hölder's inequality
    \begin{multline*}
        \int_{\T^3}|\uegd(t)|\,|m_2\fegd(t)|\dd x \\ \lesssim e^{\frac{t}{2\eps}}\norme{\fegd^0}_{\L^\infty(\T^3\times\R^3)}^{\frac{1}{6}}\left(M_3\fegd(t)\right)^{\frac{5}{6}}\norme{\uegd(t)}_{\L^6(\T^3)}.
    \end{multline*}
    Injecting this in~\eqref{eq-M3-Vlasov} yields, using once again the embedding $\H^1(\T^3)\hookrightarrow\L^6(\T^3)$,
    \[
        \frac{\dd}{\dd t}\left(M_3\fegd(t)\right)^{\frac{1}{6}}+\frac{1}{2\eps}\left(M_3\fegd(t)\right)^{\frac{1}{6}}\lesssim e^{\frac{t}{2\eps}}\norme{\fegd^0}_{\L^\infty(\T^3\times\R^3)}^{\frac{1}{6}}\norme{\uegd}_{\H^1(\T^3)},
    \]
    from which we get
    \[
        \frac{\dd}{\dd t}\left(e^{\frac{t}{2\eps}}\left(M_3\fegd\right)^{\frac{1}{6}}\right)\lesssim e^{\frac{t}{\eps}}\norme{\fegd^0}_{\L^\infty(\T^3\times\R^3)}\norme{\uegd(t)}_{\H^1(\T^3)}.
    \]
    Assumption~\ref{hypGeneral} implies that $M_3\fegd^0<\infty$, so that, thanks to the Cauchy-Schwarz inequality and the energy--dissipation estimate~\eqref{EstimationEnergieEq},
    \[
        M_3\fegd(t)\lesssim M_3\fegd^0+(t^3+t^6)e^{\frac{3t}{\eps}}\norme{\fegd^0}_{\L^\infty(\T^3\times\R^3)}\E_{\egd}(0)^3,
    \]
    which concludes the proof.
\end{preuve}


In the following lemma, we highlight two cases for which~\eqref{hypNSSup} is verified.

\begin{lemma}\label{EstimationNSsupGeneral}
    Under Assumption~\ref{hypGeneral},
    \begin{enumerate}
        \item there exists $T_{M,\eps}>0$, depending on $M$ and $\eps$, such that \eqref{hypNSSup} holds for almost every $t\in[0,T_{M,\eps}]$.
        \item under Assumption~\ref{hypSmallDataNS} and if $T>0$ verifies
            \[
                \norme{\Fegd}_{\L^2((0,T)\times\T^3)}\le C^*/2,
            \]
            then, \eqref{hypNSSup} holds for almost every $t\in[0,T]$.
    \end{enumerate}
\end{lemma}

\begin{preuve}
Thanks to Sobolev interpolation and Lemma~\ref{lemmaHeat},
    \begin{align*}
        \int_0^T\norme{e^{t\Delta}\uegd^0}_{\dot\H^1(\T^3)}^4\dd t&\lesssim\left(\sup_{t\in(0,T)}\norme{e^{t\Delta}\uegd^0}_{\dot\H^{\frac{1}{2}}(\T^3)}^2\right)\int_0^T\norme{e^{t\Delta}\uegd^0}_{\dot\H^{\frac{3}{2}}(\T^3)}^2\dd t\\
        &\lesssim\norme{\uegd^0}_{\dot\H^{\frac{1}{2}}(\T^3)}^2\int_0^T\norme{e^{t\Delta}\uegd^0}_{\dot\H^{\frac{3}{2}}(\T^3)}^2\dd t,
    \end{align*}
    and
    \[
        \int_0^T\norme{e^{t\Delta}\uegd^0}_{\dot\H^{\frac{3}{2}}(\T^3)}^2\dd t\le\norme{\uegd^0}_{\dot\H^{\frac{1}{2}}(\T^3)}^2.
    \]
    \begin{itemize}
        \item Under Assumption~\ref{hypSmallDataNS}, since
        $\int_0^t\norme{\Fegd(s)}_{\dot{\H}^{-\frac{1}{2}}(\T^3)}^2\dd s \leq\norme{\Fegd}_{\L^2((0,T)\times\T^3)}$,
        the conclusion is straightforward.
        \item If we only consider Assumption~\ref{hypGeneral}, then $(\norme{\uegd^0}_{\dot\H^{1}(\T^3)})_{\eps>0}$ is uniformly bounded 
        and therefore, for every $t\in\R_+$,
        \[
            \norme{e^{t\Delta}\uegd^0}_{\dot\H^1(\T^3)}\le M,
        \]
        from which we derive the existence of $T_M>0$ such that
        \[
            \int_0^{T_M}\norme{e^{t\Delta}\uegd^0}_{\dot\H^1(\T^3)}^4\le\frac{C^*}{2}.
        \]
        Similarly, thanks to Sobolev's embedding and Lemma~\ref{lem-ouf}, there exists $T_{M,\eps}>0$ such that
        \[
           \int_0^{T_{M,\eps}}\norme{\Fegd(t)}_{\dot{\H}^{-\frac{1}{2}}(\T^3)}^2\dd t \lesssim \int_0^{T_{M,\eps}}\norme{\Fegd(t)}_{\L^{\frac{3}{2}}(\T^3)}^2\dd t\le\frac{C^*}{2},
        \]
     which concludes the proof of the lemma.
        
    \end{itemize}
\end{preuve}

The purpose of the next lemma is to ensure that  there exists $T_{M,\eps}>0$ such that ${\Fegd} \in \L^2((0,T_{M,\eps})\times\T^3)$, in order to be able to apply~\eqref{ENSH1} at least for short times.




\begin{lemma}
\label{lem-FL2local}
  Under Assumption~\ref{hypGeneral}, there exists $T_{M,\eps}>0$ such that ${\Fegd} \in \L^2((0,T_{M,\eps})\times\T^3)$.
\end{lemma}

\begin{preuve}

Following, once again, the proof of~\cite[Lemma 4.2]{hank-mou-moy}, we obtain, for almost every $t\in\R_+$,
\[
    M_5\fegd(t)\lesssim M_5\fegd^0+e^{\frac{3t}{\eps}}\norme{\fegd^0}_{\L^\infty(\T^3\times\R^3)}\left(\int_0^t\norme{\uegd(s)}_{\L^8(\T^3)}\dd s\right)^8.
\]
Thanks to Lemma~\ref{EstimationNSsupGeneral} and Theorem~\ref{propTpsLong}, there exists $T_{M,\eps}>0$ such that $\uegd\in\L^2(0,T_{M,\eps};\H^{\frac{3}{2}}(\T^3))$, which in turn implies that $\uegd\in\L^2(0,T_{M,\eps};\L^r(\T^3))$ for any $r\in[1,\infty)$. Then, thanks to the interpolation estimate~\eqref{interpolationDistrib} with $(\ell,k)=(1,5)$,
\[
    \norme{\jegd(t)}_{\L^2(\T^3)}\lesssim\sigma^{-1}\left(M_5\fegd(t)\right)^{\frac{1}{2}}\norme{\fegd(t)}_{\L^\infty(\T^3\times\R^3)}^{\frac{1}{2}}<\infty,
\]
and $\jegd\in\L^2(0,T_{M,\eps};\L^2(\T^3))$. Similarly, using~\eqref{interpolationDistrib} with $(\ell,k)=(0,5)$,
\[
    \norme{\rhoegd(t)}_{\L^{\frac{4}{3}}(\T^3)}\lesssim\left(M_5\fegd(t)\right)^{\frac{3}{8}}\norme{\fegd(t)}_{\L^\infty(\T^3\times\R^3)}^{\frac{5}{8}}<\infty
\]
and we conclude with Hölder's inequality that $\rhoegd\uegd\in\L^2((0,T_{M,\eps})\times\T^3)$.
\end{preuve}



\subsection{Rough bounds on the first moments}
\label{sec-roughfirst}

In order to have a first rough estimate on the Brinkman force, we will use the following bounds on the moments $\rhoegd$ and $\jegd$.

\begin{lemma}\label{EstimateRhoJLinfini}
    Under Assumption~\ref{hypGeneral}, there exist a time $T_{\egd}>0$ and a continuous function $\varphi_{\egd}^{\mathrm{force}}$, increasing with respect to both its variables, and such that for every $T\in[0,T_{\egd}]$,
    \[
        \norme{\rhoegd}_{\L^\infty((0,T)\times\T^3)}\le\varphi_{\egd}^{\mathrm{force}}(T,M),
    \]
    \[
        \norme{\jegd}_{\L^\infty((0,T)\times\T^3)}\le\varphi_{\egd}^{\mathrm{force}}(T,M),
    \]
    and
    \[
        \norme{\Fegd}_{\L^p((0,T)\times\T^3)}\le\varphi_{\egd}^{\mathrm{force}}(T,M)
    \]
\end{lemma}

\begin{remark}
Note that the functions $\varphi_{\egd,0}^{\mathrm{force}}$ obtained in the proof blow up as $\eps \to 0$.
\end{remark}


\begin{preuve}
    By the method of characteristics, for any $t\ge0$ and $x\in\T^3$, we have
    \[
        \rhoegd(t,x)=e^{\frac{3t}{\eps}}\int_{\R^3}\fegd^0(X_{\egd}(0;t,x,v),V_{\egd}(0;t,x,v)\dd v.
    \]
    Using Assumption~\ref{hypGeneral}, we get
    \begin{equation}\label{roughboundRho}
        |\rhoegd(t,x)|\le Me^{\frac{3t}{\eps}}\int_{\R^3}\frac{\dd v}{1+|V_{\egd}(0;t,x,v)|^q},
    \end{equation}
    with $q>3$. Therefore, we need a lower bound for $|V_{\egd}(0;t,x,v)|$. For $t\ge0$ and $(x,v)\in\T^3\times\R^3$, we can integrate the equation of characteristics~\eqref{characV} to find that
    \[
        V_{\egd}(0;t,x,v)=e^{\frac{t}{\eps}}v-\frac{\sigma}{\eps}\int_0^te^{\frac{\tau}{\eps}}\uegd(\tau,X_{\egd}(\tau;t,x,v))\dd\tau.
    \]
    Furthermore, thanks to Lemma~\ref{EstimationNSsupGeneral}, there exists $T_{\egd}>0$ such that~\eqref{hypNSSup} and therefore~\eqref{ENSH1} and Lemma~\ref{lem-FL2local} hold. This yields
    \[
        \norme{\Delta_x\uegd}_{\L^2((0,T_{\egd})\times\T^3)}\lesssim \Psi_{\egd,0}^{\frac{1}{2}}\lesssim 1.
    \]
    Therefore, thanks to the Gagliardo-Nirenberg inequality (Theorem~\ref{th-Gagliardo-Nirenberg}),
    \[
        \norme{\uegd(t)}_{\L^\infty(\T^3)}\lesssim\norme{\uegd(t)}_{\L^2(\T^3)}^{\frac{4}{7}}\norme{\Delta_x\uegd(t)}_{\L^2(\T^3)}^{\frac{3}{7}}+\norme{\uegd(t)}_{\L^2(\T^3)},
    \]
    and using Hölder's inequality and the energy--dissipation estimate \eqref{EstimationEnergieEq}, we get, for every $T\in[0,T_{\egd}]$,
    \begin{align*}
        &\norme{\uegd}_{\L^1(0,T;\L^\infty(\T^3))}\\
        &\qquad\lesssim T^{\frac{11}{14}}\norme{\uegd}_{\L^\infty(0,T;\L^2(\T^3))}^{\frac{4}{7}}\norme{\Delta_x\uegd}_{\L^2((0,T)\times\T^3)}^{\frac{3}{7}}+T\norme{\uegd}_{\L^\infty(0,T;\L^2(\T^3))}\\
        &\qquad \lesssim T^{\frac{11}{14}}\E_{\egd}(0)^{\frac{2}{7}}+T\,\E_{\egd}(0)^{\frac{2}{7}}.
    \end{align*}
    
    Thus, we infer from~\eqref{roughboundRho} and Assumption~\ref{hypGeneral} that there exists a continuous function $\psi_{\egd}$ that is increasing with respect to both its variables and such that
    \[
        |V(0;t,x,v)|\ge|v|-\frac{\sigma e^{\frac{T}{\eps}}}{\eps}\norme{\ueps}_{\L^1(0,T;\L^\infty(\T^3))}\ge|v|-\psi_{\egd}(T,M).
    \]
    Then, since $q>3$,
    \begin{multline*}
        \rhoegd(t,x)\lesssim Me^{\frac{3t}{\eps}}\left(\int_{|v|\le2\psi_{\egd}(T,M)}\dd v+\int_{|v|\ge2\psi_{\egd}(T,M)}\frac{\dd v}{1+|v|^q}\right)\\
        \lesssim \varphi_{\egd}(T,M)
    \end{multline*}
    for some continuous function $\varphi_{\egd}$ that is increasing with respect to both its variables.
    
    Similarly, by the method of characteristics,
    \[
        |\jegd(t,x)|\le \sigma^{-1}Me^{\frac{3t}{\eps}}\int_{\R^3}\frac{|v|\dd v}{1+|V(0;t,x,v)|^q}
    \]
    and we can proceed as above to obtain the stated result since $q>4$ in Assumption~\ref{hypGeneral}.
    
    We can then conclude, since $\Fegd=(\jegd-\rhoegd\uegd)/\gamma$, by applying Lemma~\ref{EstimNormLp}.
\end{preuve}


\subsection{Modulated energy dissipation}
\label{sec-modenergy}

Following \cite{cho-kwo} and \cite{hank-mou-moy}, we introduce a \emph{modulated} energy for which, under certain conditions, we can prove an exponential decay.

\begin{definition}
    We define the modulated energy, for every $t\ge0$, by
    \begin{multline*}
        \EEE_{\egd}(t)=
        \frac{\eps}{2\gamma}\int_{\T^3\times\R^3}\left|\frac{v}{\sigma}-\frac{\langle\jegd(t)\rangle}{\langle\rhoegd\rangle}\right|^2\fegd(t,x,v)\,\dd x\dd v\\
        +\frac{1}{2}\int_{\T^3}|\uegd(t,x)-\langle\uegd(t)\rangle|^2\dd x\\
        +\frac{\eps \langle\rhoegd\rangle}{2(\gamma+\eps \langle\rhoegd\rangle)}\left|\frac{\langle\jegd(t)\rangle}{\langle\rhoegd\rangle}-\langle\uegd(t)\rangle\right|^2.
    \end{multline*}
        \index{$\EEE_{\egd}(t)$: modulated energy}
\end{definition}

\begin{remark}
    Recall that, thanks to Lemma~\ref{VlasovConservation}, $\langle\rhoegd(t)\rangle$ does not depend on time, which justifies the notation in the functional above.
\end{remark}

The modulated energy satisfies the same type of energy--dissipation estimate as the energy.

\begin{lemma}\label{ModulatedEnergyDissipation}
    Under Assumption~\ref{hypGeneral}, for every $\eps>0$, for almost every $0\le s\le t$ (including $s=0$),
    \[
        \EEE_{\egd}(t)+\int_s^t\D_{\egd}(\tau)\dd\tau\le\EEE_{\egd}(s).
    \]
\end{lemma}

\begin{preuve}
    We expand the quadratic terms and find that for all $t\ge0$
    \begin{align*}
        \EEE_{\egd}(t)&=\E_{\egd}(t)-\frac{\eps}{2\gamma}\frac{|\langle\jegd(t)\rangle|^2}{\langle\rhoegd\rangle}-\frac{|\uegd(t)|^2}{2}\\
            &\quad+\frac{\eps}{2(\gamma+\eps \langle\rhoegd\rangle)}\frac{|\langle\jegd(t)\rangle|^2}{\langle\rhoegd\rangle}+\frac{\eps \langle\rhoegd\rangle}{2(\gamma+\eps \langle\rhoegd\rangle)}|\langle\uegd(t)\rangle|^2\\
            &\quad-\frac{\eps}{\gamma+\eps \langle\rhoegd\rangle}\langle\jegd(t)\rangle\cdot\langle\uegd(t)\rangle\\
        &=\E_{\egd}(t)-\frac{\gamma}{2(\gamma+\eps \langle\rhoegd\rangle)}\left|\frac{\eps}{\gamma}\langle\jegd(t)\rangle+\langle\uegd(t)\rangle\right|^2.
    \end{align*}
    We conclude thanks to the energy--dissipation estimate \eqref{EstimationEnergieEq} and Lemma~\ref{conservationVNSmoy}.

\end{preuve}

As shown in \cite{cho-kwo} and \cite{hank-mou-moy}, under an assumption on the first moment of $\fegd$, we have an exponential decay of the modulated energy.

\begin{lemma}\label{lemmeDE}
    Under Assumption \ref{hypGeneral}, if $T>0$ is such that $\rhoegd\in\L^\infty((0,T)\times\T^3)$, then 
    \begin{equation}\label{DE}
        \forall t\in[0,T],\qquad \EEE_{\egd}(t)\le C e^{-\lambda_{\egd} t}\EEE_{\egd}(0),
    \end{equation}
    with
    \[
        \lambda_{\egd}=\min\left(\frac{c_P}{\eps\left(c_P+4\norme{\rhoegd}_{\L^\infty((0,T)\times\T^3)}\right)},\frac{c_P}{2}\right),
        \index{L@$\lambda_{\egd}$: exponential decay rate}
    \]
    where $c_P$ is a universal constant and $C>0$ depends only on $\lambda_{\egd}$.
\end{lemma}

\begin{preuve}
    We adapt the proof of \cite[Lemma 3.4]{hank-mou-moy}. All we need to prove is that
    \[
        \forall t\in[0,T],\qquad\D_{\egd}(t)\ge\lambda_{\egd}\EEE_{\egd}(t),
    \]
    since the exponential decay~\eqref{DE} then follows from Lemma~\ref{ModulatedEnergyDissipation} and~\cite[Lemma 9.3]{hank-mou-moy}. For the sake of readability, we shall not write the time, position or velocity variables in the rest of the proof.
    
    By the Poincaré-Wirtinger inequality, there exists a universal constant $c_P>0$ such that
    \begin{equation}\label{preuveDEeq0}
        \D_{\egd}\ge\frac{1}{2}\int_{\T^3\times\R^3}\left|\frac{v}{\sigma}-\uegd\right|^2\fegd \, \dd x \dd v\\
            +\frac{c_P}{2}\norme{\uegd-\langle\uegd\rangle}_{\L^2(\T^3)}^2.
    \end{equation}
    Furthermore, expanding the quadratic factor yields
    \begin{multline*}
        \int_{\T^3\times\R^3}\left|\frac{v}{\sigma}-\uegd\right|^2\fegd \, \dd x \dd v\\
        =\int_{\T^3\times\R^3}\left|\frac{v}{\sigma}-\langle\uegd\rangle\right|^2\fegd \, \dd x \dd v
            +\int_{\T^3}\rhoegd|\uegd-\langle\uegd\rangle|^2 \, \dd x \\
            +2\int_{\T^3\times\R^3}\left(\frac{v}{\sigma}-\langle\uegd\rangle\right)\cdot\left(\langle\uegd\rangle-\uegd\right)\fegd \, \dd x \dd v.
    \end{multline*}
    For any $a\in(0,1)$, we have
    \begin{multline*}
        2\int_{\T^3\times\R^3}\left(\frac{v}{\sigma}-\langle\uegd\rangle\right)\cdot\left(\langle\uegd\rangle-\uegd\right)\fegd \, \dd x \dd v\\
        \ge-a\int_{\T^3\times\R^3}\left|\frac{v}{\sigma}-\langle\uegd\rangle\right|^2\fegd \, \dd x \dd v
            -\frac{1}{a}\int_{\T^3}\rhoegd|\uegd-\langle\uegd\rangle|^2 \, \dd x ,
    \end{multline*}
    so that
    \begin{multline}\label{preuveDEeq1}
        \int_{\T^3\times\R^3}\left|\frac{v}{\sigma}-\uegd\right|^2\fegd \, \dd x \dd v
        \ge(1-a)\int_{\T^3\times\R^3}\left|\frac{v}{\sigma}-\langle\uegd\rangle\right|^2\fegd \, \dd x \dd v\\
            -\left(\frac{1}{a}-1\right)\int_{\T^3}\rhoegd|\uegd-\langle\uegd\rangle|^2 \, \dd x.
    \end{multline}
    
    In order to deal with the first term, we write
    \begin{multline*}
        \left|\frac{v}{\sigma}-\langle\uegd\rangle\right|^2
        =\left|\frac{\langle\jegd\rangle}{\langle\rhoegd\rangle}-\langle\uegd\rangle\right|^2 \\
            +2\left(\frac{v}{\sigma}-\frac{\langle\jegd\rangle}{\langle\rhoegd\rangle}\right)\cdot\left(\frac{\langle\jegd\rangle}{\langle\rhoegd\rangle}-\langle\uegd\rangle\right)
            +\left|\frac{v}{\sigma}-\frac{\langle\jegd\rangle}{\langle\rhoegd\rangle}\right|^2
    \end{multline*}
    and obtain, after integrating against $\fegd$,
    \begin{align*}
        &\int_{\T^3\times\R^3}\left|\frac{v}{\sigma}-\langle\uegd\rangle\right|^2\fegd \, \dd x \dd v\\
        &\qquad \qquad \qquad =\langle\rhoegd\rangle\left|\frac{\langle\jegd\rangle}{\langle\rhoegd\rangle}-\langle\uegd\rangle\right|^2
            +\int_{\T^3\times\R^3}\left|\frac{v}{\sigma}-\frac{\langle\jegd\rangle}{\langle\rhoegd\rangle}\right|^2\fegd \, \dd x \dd v\\
         &\qquad \qquad \qquad\ge\frac{2\gamma}{\eps}\left(\EEE_{\egd}-\frac{1}{2}\norme{\uegd-\langle\uegd\rangle}_{\L^2(\T^3)}^2\right).
    \end{align*}
    Therefore, the inequality~\eqref{preuveDEeq1} becomes
    \begin{multline*}
        \int_{\T^3\times\R^3}\left|\frac{v}{\sigma}-\uegd\right|^2\fegd \, \dd x \dd v \\
        \ge\frac{2\gamma(1-a)}{\eps}\left(\EEE_{\egd}-\frac{1}{2}\norme{\uegd-\langle\uegd\rangle}_{\L^2(\T^3)}^2\right)\\
        -\left(\frac{1}{a}-1\right)\norme{\rhoegd}_{\L^\infty((0,T)\times\T^3)}\norme{\uegd-\langle\uegd\rangle}_{\L^2(\T^3)}^2.
    \end{multline*}
    We set
    \[
        a=\frac{4\norme{\rhoegd}_{\L^\infty((0,T)\times\T^3}}{c_p+4\norme{\rhoegd}_{\L^\infty((0,T)\times\T^3)}}\in(0,1)
    \]
    so that
    \[
        \left(\frac{1}{a}-1\right)\norme{\rhoegd}_{\L^\infty((0,T)\times\T^3)}=\frac{c_P}{4}.
    \]
    Injecting this into~\eqref{preuveDEeq0} yields
    \begin{align*}
        \D_{\egd}
        &\ge\frac{c_P}{\eps\left(c_P+4\norme{\rhoegd}_{\L^\infty((0,T)\times\T^3)}\right)}\left(\EEE_{\egd}-\frac{1}{2}\norme{\uegd-\langle\uegd\rangle}_{\L^2(\T^3)}^2\right)\\
            &\quad +\frac{c_P}{2}\times\frac{1}{2}\norme{\uegd-\langle\uegd\rangle}_{\L^2(\T^3)}^2\\
        &\ge\min\left(\frac{c_P}{\eps\left(c_P+4\norme{\rhoegd}_{\L^\infty((0,T)\times\T^3)}\right)},\frac{c_P}{2}\right)\EEE_{\egd},
    \end{align*}
    hence the result.
\end{preuve}

We will heavily rely on this exponential decay of the modulated energy, through the following estimates.

\begin{lemma}\label{estimNabla}
    Under Assumption~\ref{hypGeneral}, for every $r>3$, if $T>0$ is such that $\rhoegd\in\L^\infty((0,T)\times\T^3)$ and $\Delta_x\uegd\in\L^r((0,T)\times\T^3)$, then for any $t\in[0,T]$,
    \[
        \norme{\na_x\uegd(t)}_{\L^r(\T^3)}\lesssim\EEE_{\egd}(0)^{\frac{1-\alpha_r}{2}}e^{-\frac{(1-\alpha_r)\lambda_{\egd}}{2}t}\norme{\Delta_x\uegd(t)}_{\L^r(\T^3)}^{\alpha_r},
    \]
    \[
        \norme{\na_x\uegd(t)}_{\L^\infty(\T^3)}\lesssim\EEE_{\egd}(0)^{\frac{1-\beta_r}{2}}e^{-\frac{(1-\beta_r)\lambda_{\egd}}{2}t}\norme{\Delta_x\uegd(t)}_{\L^r(\T^3)}^{\beta_r},
    \]
    where $\alpha_r,\beta_r$ are defined in Corollary~\ref{Gagliardo-Nirenberg}.
\end{lemma}

\begin{preuve}
    This is a direct consequence of the Gagliardo-Nirenberg inequality (Corollary~\ref{Gagliardo-Nirenberg}) and Lemma~\ref{lemmeDE} since $r>3$. 
\end{preuve}


\subsection{Changes of variable in velocity or in space}\label{SectionChangeVar}

We present in this section two changes of variables that will allow us to derive crucial estimates on the moments and the Brinkman force.

First, we adapt what \cite{hank-mou-moy} calls the \emph{straightening} change of variables. Under a smallness condition on $\norme{\na_x\uegd}_{\L^1(0,t;\L^\infty(\T^3))}$, we define a diffeomorphism in velocity as follows.

\begin{lemma}\label{changeVarV}
    Fix $c_*>0$ such that $c_* e^{c_*}<1/9$. Then, for any $t\in\R_+$ satisfying
    \begin{equation}\label{hyp-na-u}
        \norme{\na\uegd}_{\L^1(0,t;\L^\infty(\T^3))}\le c_*,
    \end{equation}
    and any $x\in\T^3$, the map
    \[
        \Gamma_{\egd}^{t,x}:v\mapsto V_{\egd}(0;t,x,v)
    \]\index{Gamma@$\Gamma_{\egd}^{t,x}(v)$: change of variables in velocity}
    is a $\CCC^1$-diffeomorphism from $\R^3$ to itself and satisfies
    \[
        \forall v\in\R^3,\qquad \det\D_v\Gamma_{\egd}^{t,x}(v)\ge\frac{e^{\frac{3t}{\eps}}}{2}.
    \]
\end{lemma}

\begin{preuve}
    Defining $Y_\egd=\frac{\sigma}{\eps} X_{\egd}$, the equations of characteristics in Definition~\ref{defCharac} can be written as follows: for any $t\ge0$ and $(y,v)\in\frac{\sigma}{\eps}\T^3\times\R^3$,
    \[
        \left\{\begin{aligned}
            \dot{Y}_\egd(s;t,y,v)&=\frac{1}{\eps}\,V_\egd(s;t,x,v)\\
            \dot{V}_\egd(s;t,y,v)&=\frac{1}{\eps}\left(\sigma\uegd\left(s,\frac{\eps}{\sigma} Y_\egd(s;t,y,v)\right)-V_\egd(s;t,y,v)\right)\\
            Y_\egd(t;t,x,v)&=y\\
            V_\egd(t;t,y,v)&=v.
            \end{aligned}\right.
    \]
    Following \cite{hank-mou-moy}, we define, for any $s\in\R_+$ and $z=(y,v)\in\frac{\sigma}{\eps}\T^3\times\R^3$,
    \[
        w_\egd(s,z)=\left(\frac{v}{\eps},\frac{1}{\eps}\left(\sigma\uegd\left(s,\frac{\eps}{\sigma} y\right)-v\right)\right)
    \]
    and $Z_\egd=(Y_\egd,V_\egd)$, so that
    \[
        \partial_s Z_\egd(s;t,z)=w_\egd(s,Z_\egd(s;t,y,v)),
    \]
    which after differentiation with respect to $z$ yields
    \[
        \partial_s\D_z Z_\egd(s;t,z)=\D_z w_\egd(s,Z(s;t,z))\D_z Z_\egd(s;t,z).
    \]
    Since for any $s\in\R_+$ and $z=(y,z)\in\T^3\times\R^3$,
    \[
        \norme{\D_z w_\egd(s,\cdot)}_{\L^\infty(\T^3\times\R^3)}\le\frac{1}{\eps}+\norme{\na\uegd}_{\L^\infty(\T^3)},
    \]
    Grönwall's lemma yields, for any $0\le s\le t$,
    \[
        \norme{\D_z Z_\egd(s;t,\cdot)}_{\L^\infty(\T^3\times\R^3)}\le e^{\frac{t-s}{\eps}}e^{\norme{\na_x\uegd}_{\L^1(0,t;\L^\infty(\T^3))}}.
    \]
    In particular,
    \[
        \norme{\D_{v}Y_\egd(s;t,\cdot)}_{\L^\infty(\T^3\times\R^3)}\le e^{\frac{t-s}{\eps}}e^{\norme{\nabla_{x}\uegd}_{\L^1(0,t;\L^\infty(\T^3))}}.
    \]
    Using
    \[
        V_\egd(0;t,y,v)=e^{\frac{t}{\eps}}v-\frac{\sigma}{\eps}\int_0^te^{\frac{s}{\eps}}\uegd(s,\eps Y_\egd(s;t,y,v)/\sigma)\dd s,
    \]
    we get
    \begin{align*}
        &\norme{e^{\frac{-t}{\eps}}\D_v V_\egd(0;t,\cdot)-\Id}_{\L^\infty(\T^3\times\R^3)}\\
        &\qquad \qquad\le\int_0^te^{\frac{s-t}{\eps}}\norme{\na_x\uegd(s)}_{\L^\infty(\T^3)}\norme{\D_v Y_\egd(s;t,\cdot)}_{\L^\infty(\T^3\times\R^3)}\dd s\\
        &\qquad \qquad\le\norme{\na_x\uegd}_{\L^1(0,t;\L^\infty(\T^3))}e^{\norme{\nabla_{x}\uegd}_{\L^1(0,t;\L^\infty(\T^3))}} \\
        &\qquad \qquad \le c_*e^{c_*}<\frac{1}{9}
    \end{align*}
    and we conclude using Theorem \ref{perturbationId}. 
\end{preuve}

\begin{remark}
We set in the following $c_*=1/30$ and the inequality $c_*e^{c_*}<1/9$ holds.
\end{remark}

This change of variables in velocity leads to the following crucial $\L^\infty_{\loc}(\R_+\times\T^3)$ bound.

\begin{corollary}\label{RhoLinfini}
    Under Assumption~\ref{hypGeneral}, if~\eqref{hyp-na-u} holds, then
    \[
        \norme{\rhoegd}_{\L^\infty((0,T)\times\T^3)}\lesssim \norme{\fegd^0}_{\L^1(\R^3;\L^\infty(\T^3))}\le M.
    \]
\end{corollary}

\begin{preuve}
    Following the method of characteristics, we write, for any $(t,x)\in[0,T]\times\T^3$,
    \[
        \rhoegd(t,x)=\int_{\R^3}\fegd(t,x,v)\dd v
        =e^{\frac{3t}{\eps}}\int_{\R^3}\feps^0(X_{\egd}(0;t,x,v),V_{\egd}(0;t,x,v))\dd v.
    \]
    We perform the change of variable $w=V_{\egd}(0;t,x,v)$ thanks to Lemma \ref{changeVarV} and get
    \begin{align*}
        |\rhoegd(t,x)|&\lesssim\int_{\R^3}\feps^0\left(X_{\egd}\left(0;t,x,[\Gamma_{\egd}^{t,x}]^{-1}(w)\right),w\right)\dd w\\
        &\lesssim\int_{\R^3}\norme{\fegd^0(\cdot,w)}_{\L^\infty(\T^3)}\dd w
        =\norme{\fegd^0}_{\L^1(\R^3;\L^\infty(\T^3))}.
    \end{align*}
\end{preuve}

\begin{remark}
    Thanks to this result, we can apply Lemma \ref{lemmeDE}. In particular, $\lambda_{\egd}$ is bounded from above and below by positive constants that are independent of $\eps$.
\end{remark}

In the proof of the previous corollary, the uniform bound $$|\feps^0\left(X_{\egd}\left(0;t,x,[\Gamma_{\egd}^{t,x}]^{-1}(w)\right),w\right)| \leq \norme{\fegd^0(\cdot,w)}_{\L^\infty(\T^3)} $$ 
was sufficient, but we will encounter cases in Sections~\ref{SectionLandLFParticle} and~\ref{SectionFineParticle} where we will need more precision. For this purpose, we will apply a change of variables relative to the position variable after the above change of velocity variables has been performed.

\begin{lemma}\label{changeVarX}
    For any $t\in\R_+$ satisfying
    \[
        \norme{\na_x\uegd}_{\L^1(0,t;\L^\infty(\T^3))}\le\frac{1}{30},
    \]
    for any $(x,v)\in\T^3\times\R^3$, we define
    \begin{equation}
    \label{eq-defXtilde}
        \index{X@$\tilde{X}_{\egd}^{t,x,v}$: spacial characteristics after an application of the change of variables in velocity}
        \tilde{X}_{\egd}^{t,x,v}:s\mapsto X_{\egd}(s;t,x,[\Gamma_{\egd}^{t,x}]^{-1}(v)).
    \end{equation}
    Then for any $0\le s\le t$ and $v\in\R^3$, the map
    \[
        \Phi_{\egd}^{s;t,v}:x\mapsto\tilde{X}_{\egd}^{t,x,v}(s)
        \index{P@$\Phi_{\egd}^{s;t,v}$: change of variables in space}
    \]
    is a $\CCC^1$-diffeomorphism from $\T^3$ to itself and satisfies
    \[
        \forall x\in\T^3,\qquad\det(\na_x\tilde{X}_{\egd}^{t,x,v}(s))\ge\frac{1}{2}.
    \]
\end{lemma}

\begin{preuve}
    The assumptions of Lemma \ref{changeVarV} are satisfied thus $[\Gamma_{\egd}^{t,x}]^{-1}(v)$ is well-defined for any $(x,v)\in\T^3\times\R^3$. The equations of characteristics yield, for any $0\le\tau\le t$,
    \[
        V_\egd(\tau;t,x,[\Gamma_{\egd}^{t,x}]^{-1}(v))=e^{\frac{t-\tau}{\eps}}[\Gamma_{\egd}^{t,x}]^{-1}(v)-\frac{\sigma}{\eps}\int_{\tau}^te^{\frac{\zeta-\tau}{\eps}}\uegd\left(\zeta,\tilde{X}_{\egd}^{t,x,v}(\zeta)\right)\dd\zeta,
    \]
    where
    \[
        [\Gamma_{\egd}^{t,x}]^{-1}(v)=e^{-\frac{t}{\eps}}v+\frac{\sigma}{\eps}\int_0^te^{\frac{\tau-t}{\eps}}\uegd\left(\tau,\tilde{X}_{\egd}^{t,x,v}(\tau)\right)\dd\tau.
    \]
    Therefore, for any $s\in[0,t]$,
    \begin{multline*}
        \tilde{X}_{\egd}^{t,x,v}(s)-x=-\frac{1}{\sigma}\int_s^tV_\egd\left(\tau;t,x,[\Gamma_{\egd}^{t,x}]^{-1}(v)\right)\dd\tau\\
        =\eps(e^{-\frac{t}{\eps}}-e^{-\frac{s}{\eps}})v
        +\int_0^\infty\left[e^{\frac{\tau-t}{\eps}}\1_{\tau\le t}-e^{\frac{\tau-s}{\eps}}\1_{\tau\le s}-\1_{s\le\tau\le t}\right]\uegd\left(\tau,\tilde{X}_{\egd}^{t,x,v}(\tau)\right)\dd\tau.
    \end{multline*}
    By differentiating, we get
    \begin{multline*}
        \norme{\D_{x}\tilde{X}_{\egd}^{t,\cdot,\cdot}-\Id}_{\L^\infty((0,t)\times\T^3\times\R^3)}\le3\norme{\na_{x}\uegd}_{\L^1(0,t;\L^\infty(\T^3)}\\
        +3\norme{\D_{x}\tilde{X}_{\egd}^{t,\cdot,\cdot}-\Id}_{\L^\infty((0,t)\times\T^3\times\R^3)}\norme{\na_{x}\uegd}_{\L^1(0,t;\L^\infty(\T^3))},
    \end{multline*}
    hence, thanks to the assumption on $\norme{\na_x\uegd}_{\L^1(0,t;\L^\infty(\T^3))}$,
    \[
        \norme{\D_{x}\tilde{X}_{\egd}^{t,\cdot,\cdot}-\Id}_{\L^\infty((0,t)\times\T^3\times\R^3)}\le\frac{1}{9}
    \]
    and we conclude with Theorem \ref{perturbationId}.
\end{preuve}

We now define \emph{strong existence times}, for which all the assumptions required in Lemma~\ref{EstimationNSsupGeneral},~\ref{changeVarV}, and~\ref{changeVarX} are satisfied.

\begin{definition}\label{StrongExistenceTime}
\index{Strong existence time}
    We say $T>0$ is a strong existence time if
    \[
        \norme{\na_x\uegd}_{\L^1(0,T;\L^\infty(\T^3))}\le\frac{1}{30}
    \]
    and
    \[
        \max\left(\int_0^T\norme{e^{t\Delta}\uegd^0}_{\dot\H^1(\T^3)}^4\dd t,\int_0^T\norme{\Fegd(t)}_{\L^2(\T^3)}^2\dd t\right)\le\frac{C^*}{2}.
    \]
\end{definition}


\subsection{Estimates on the convection term}
\label{sec-conv}

We will need several estimates on the convective term in the Navier-Stokes equation. First, we show that this term is in $\L^p((0,T)\times\T^3)$, where $p$ is the regularity index of Assumption~\ref{hypGeneral}, for appropriate values of $T>0$.  We provide a local in time and non uniform bound that we will use to initialize the bootstrap argument. We conclude this section by a uniform bound that is only valid for strong existence times.

\begin{lemma}\label{convLp}  Under Assumption~\ref{hypGeneral}, the following holds. Let $p$ be the regularity index appearing in the statement. 
  Let $T>0$ such that~\eqref{ENSH1}  holds on $[0,T]$ and 
 $\norme{\Fegd}_{\L^2((0,T)\times\T^3)}<+\infty$. Then $(\uegd\cdot\na_x)\uegd\in\L^p((0,T)\times\T^3)$.
    
    In particular there exists $T_{\egd}>0$ such that $(\uegd\cdot\na_x)\uegd\in\L^p((0,T_{\egd})\times\T^3)$.
\end{lemma}

\begin{preuve}
    Let $T>0$ such that~\eqref{ENSH1}  holds on $[0,T]$ and 
 $\norme{\Fegd}_{\L^2((0,T)\times\T^3)}<+\infty$. Following the proof of \cite[Lemma 6.2]{hank-mou-moy} and thanks to~\eqref{ENSH1}, we obtain
    \begin{equation}\label{interpolationConv}
    \begin{aligned}
        &\norme{\uegd\cdot\na_x\uegd}_{\L^a(0,T;\L^b(\T^3))}\\
        &\qquad\lesssim\norme{\uegd}_{\L^\infty(0,T;\L^6(\T^3))}^a\norme{\na_x\uegd}_{\L^{r_1}(0,T;\L^{r_2}(\T^3))}^{\frac{r_1}{a}}\norme{\na_x\uegd}_{\L^\infty(0,T;\L^2(\T^3))}^{1-\frac{r_1}{a}}\\
        &\qquad\lesssim_M\norme{\na_x\uegd}_{\L^{r_1}(0,T;\L^{r_2}(\T^3))}^{\frac{r_1}{a}},
        \end{aligned}
    \end{equation}
    for any $(a,b,r_1,r_2)\in(1,\infty)$ such that $r_1\le a$, $2\le b\le r_2$, and
    \begin{equation}\label{cond-interpolationConv}
        \frac{1}{b}=\frac{1}{6}+\frac{r_1}{a}\frac{1}{r_2}+\frac{1}{2}\left(1-\frac{r_1}{a}\right).
    \end{equation}
    
    We begin by taking $(a,b,r_1,r_2)=(2,3,2,6)$ in~\eqref{interpolationConv} and get, thanks to~\eqref{ENSH1},
    \[
        \norme{\uegd\cdot\na_x\uegd}_{\L^2(0,T;\L^3(\T^3))}\lesssim_M 1.
    \]
    Recall that thanks to Lemmas~\ref{EstimNormLp} and~\ref{EstimateRhoJLinfini}, we also have
    \[
        \norme{\jegd-\rhoegd\uegd}_{\L^2(0,T;\L^3(\T^3))}\lesssim_{M,T,\eps}1.
    \]
    Therefore, Theorem~\ref{ParabolicReg} yields
    \[
        \norme{\partial_t\uegd}_{\L^2(0,T;\L^3(\T^3))}+\norme{\Delta_x\uegd}_{\L^2(0,T;\L^3(\T^3))}\lesssim_{M,T,\eps}1.
    \]
    Then, by Sobolev embedding, we infer that, for any $r_2\in[2,+\infty)$,
    \[
        \norme{\na_x\uegd}_{\L^2(0,T;\L^{r_2}(\T^3))}\lesssim_{M,T,\eps}1.
    \]
    
    With $a=r\in[2,3)$, $b=3$ and $r_1=2$,~\eqref{interpolationConv} yields,
    \[
        \norme{\uegd\cdot\na_x\uegd}_{\L^r(0,T;\L^3(\T^3))}\lesssim_{M,T,\eps}1. 
    \]
    Again, we have
    \[
        \norme{\jegd-\rhoegd\uegd}_{\L^r(0,T;\L^3(\T^3))}\lesssim_{M,T,\eps}1
    \]
    and thus
    \[
        \norme{\partial_t\uegd}_{\L^r(0,T;\L^3(\T^3))}+\norme{\Delta_x\uegd}_{\L^r(0,T;\L^3(\T^3))}\lesssim_{M,T,\eps}1,
    \]
    which yields, by Sobolev embedding, for any $r\in[2,3)$ and $r_2\in[2,+\infty)$,
    \[
        \norme{\na_x\uegd}_{\L^r(0,T;\L^{r_2}(\T^3))}\lesssim_{M,T,\eps}1.
    \]
    To conclude, we aim to apply~\eqref{interpolationConv} with $a=b=p>3$ and $r_1\in[2,3)$. The condition~\eqref{cond-interpolationConv} implies
    \[
        r_2=\frac{6r_1}{-4p+3r_1+6}
    \]
    so that the inequality $p\le r_2$ becomes
    \[
        p<\frac{3(r_1+2)}{4}.
    \]
    Finally the second part of the statement follows thanks to Lemmas~\ref{EstimationNSsupGeneral} and~\ref{lem-FL2local}: there exists $T=T_{\egd}>0$ such that~\eqref{ENSH1} and  $\norme{\Fegd}_{\L^2((0,T)\times\T^3)}<+\infty$ hold.
\end{preuve}

%

\begin{lemma}\label{EstimConvectStrongT}
    Under Assumption~\ref{hypGeneral}, if $T>0$ is a strong existence time, then
    \begin{equation*}
        \norme{\left(\uegd\cdot\na_x\right)\uegd}_{\L^p((0,T)\times\T^3)}
        \lesssim \Psi_{\egd,0}^{\frac{1}{2}}\EEE_{\egd}(0)^{\frac{1-\beta_p}{2}}\norme{\Delta_x\uegd}_{\L^p((0,T)\times\T^3)}^{\beta_p},
    \end{equation*}
    where $p$ is the regularity index appearing in the assumption and $\beta_p$ is given by Corollary~\ref{Gagliardo-Nirenberg}.
\end{lemma}

\begin{preuve}
    Since $T$ is a strong existence time, thanks to Corollary~\ref{EstimNormLp} and Lemma~\ref{estimNabla}, we have
    \[
        \norme{\uegd}_{\L^\infty(0,T;\L^p(\T^3))}\lesssim\Psi_{\egd,0}^{\frac{1}{2}}
    \]
    and, for every $t\in(0,T)$,
    \[
        \norme{\na_x\uegd(t)}_{\L^\infty(\T^3)}
        \lesssim\EEE_{\egd}(0)^{\frac{1-\beta_p}{2}}e^{-\frac{(1-\beta_p)\lambda_{\egd}}{2}t}\norme{\Delta_x\uegd(t)}_{\L^p(\T^3)}.
    \]
    Therefore, using Hölder's inequality, we get
    \begin{align*}
        &\norme{(\uegd\cdot\na_x)\uegd}_{\L^p((0,T)\times\T^3)}^p\\
        &\qquad \qquad\lesssim\Psi_{\egd,0}^{\frac{p}{2}}\EEE_{\egd}(0)^{\frac{(1-\beta_p)p}{2}}\int_0^Te^{-\frac{(1-\beta_p)p\lambda_{\egd}}{2}t}\norme{\Delta_x\uegd(t)}_{\L^p(\T^3)}^{\beta_pp}\dd t\\
        &\qquad \qquad \lesssim \Psi_{\egd,0}^{\frac{p}{2}}\EEE_{\egd}(0)^{\frac{(1-\beta_p)p}{2}}\norme{\Delta_x\uegd}_{\L^p((0,T)\times\T^3)}^{\beta_pp}
    \end{align*}
    and the lemma follows.
\end{preuve}

\subsection{Initialization of the bootstrap argument}\label{SubsectionInitialization}

In this section, the last one that is common to all the regimes, we initiate the bootstrap argument that we will use to determine strong existence times, under several sets of assumptions.



First, we prove a non-uniform existence result.

\begin{lemma}\label{ExistsStrongTime}
    Under Assumption~\ref{hypGeneral}, for every $\eps>0$, there exists a strong existence time.
\end{lemma}

\begin{preuve}
Let $p$ be the regularity index of Assumption~\ref{hypGeneral}.
    We apply the parabolic regularity estimates from Theorem~\ref{ParabolicReg} to obtain, for any $T>0$,
    \begin{multline*}
        \norme{\partial_t\uegd}_{\L^p((0,T)\times\T^3)}+\norme{\Delta_x\uegd}_{\L^p((0,T)\times\T^3)}\\
        \lesssim \norme{\Fegd}_{\L^p((0,T)\times\T^3)}+\norme{(\uegd\cdot\na_x)\uegd}_{\L^p((0,T)\times\T^3)}+\norme{\uegd^0}_{\B^{s,p}_p(\T^3)}.
    \end{multline*}
    
    To control the first term, we remember that thanks to Lemma~\ref{EstimationNSsupGeneral}, there exists $T_{\egd}>0$ such that \eqref{ENSH1} and Lemma~\ref{lem-FL2local} hold. Then, thanks to Lemmas~\ref{EstimationNSsupGeneral} and~\ref{EstimateRhoJLinfini}, there exists a continuous function $\varphi_{\egd,0}^{\mathrm{force}}$ such that
    \begin{equation}
    \label{eqFLprough}
        \norme{\Fegd}_{\L^p((0,T_{\egd})\times\T^3)}\le\varphi_{\egd,0}^{\mathrm{force}}(T_{\egd}).
    \end{equation}
  Moreover, Lemma~\ref{convLp} shows that there exists a continuous function $\varphi_{\egd,0}^{\mathrm{conv}}$ that vanishes at 0 and such that
    \begin{equation*}
        \norme{(\uegd\cdot\na_x)\uegd}_{\L^p((0,T_{\egd})\times\T^3)}\\
        \le\varphi_{\egd,0}^{\mathrm{conv}}(T_{\egd}).
    \end{equation*}
    Therefore, we can take $T=T_{\egd}$ so that
    \begin{multline}\label{initBootstrapParabolic}
        \norme{\partial_t\uegd}_{\L^p((0,T)\times\T^3)}+\norme{\Delta_x\uegd}_{\L^p((0,T)\times\T^3)}\\
        \lesssim \varphi_{\egd,0}(T)+\norme{\uegd^0}_{\B^{s,p}_p(\T^3)},
    \end{multline}
    for some continuous function $\varphi_{\egd,0}$. 
    
    Using the energy--dissipation estimate, the Gagliardo-Nirenberg inequality (see Corollary~\ref{Gagliardo-Nirenberg}) and Hölder's inequality, we get
    \begin{align*}
        &\norme{\na_x\uegd}_{\L^1(0,T;\L^\infty(\T^3))}\\
        &\quad\lesssim\norme{\uegd}_{\L^\infty((0,T);\L^2(\T^3))}^{1-\beta_p}\int_0^{T}\norme{\Delta_x\uegd(t)}_{\L^p(\T^3)}^{\beta_p}\dd t
        +T\norme{\uegd}_{\L^\infty((0,T);\L^2(\T^3))}\\
        &\quad \lesssim T^{1-\frac{\beta_p}{p}}\E_{\egd}(0)^{\frac{1-\beta_p}{2}}\norme{\Delta_x\uegd}_{\L^p((0,T)\times\T^3)}^{\beta_p}+T\,\E_{\egd}(0)^{\frac{1}{2}}.
    \end{align*}
    Thanks to \eqref{initBootstrapParabolic}, we conclude that we can take $T$ small enough such that
    \[
        \norme{\na_x\uegd}_{\L^1(0,T;\L^{\infty}(\T^3))}\le \frac{1}{30}.
    \]
    
    Finally, as in the proof of Lemma~\ref{EstimationNSsupGeneral}, we show that we can reduce $T>0$ to ensure that
    \[
        \int_0^T\norme{e^{t\Delta}\uegd}_{\dot\H^1(\T^3)}^4\dd t\le\frac{C^*}{2}.
    \]
    Furthermore, thanks to~\eqref{eqFLprough} and Hölder's inequality, we have
    \[
     \norme{\Fegd}_{\L^2((0,T)\times\T^3)}^2\le T^{\frac{p}{(p-1)}}\varphi_{\egd,0}^{\mathrm{force}}(T),
    \]
    so that, up to reducing $T>0$, we have
    \[
    \norme{\Fegd}_{\L^2((0,T)\times\T^3)}^2 \le\frac{C^*}{2}.
    \]
    Gathering all the pieces together, this leads to the existence of a strong existence time.
\end{preuve}

Arguing as in the proof of Lemma~\ref{ExistsStrongTime}, we     also obtain the following result, which is useful to justify that the computations of the following sections make sense.
\begin{lemma}
Let $T>0$ be a strong existence time. Then
\[
   \norme{\partial_t\uegd}_{\L^p((0,T)\times\T^3)}+\norme{\Delta_x\uegd}_{\L^p((0,T)\times\T^3)} < +\infty.
\]
\end{lemma}

We can now define, for every $\eps>0$,
\[
\index{T@$T^*_\eps$: supremum of strong existence times}
    T^*_\eps=\sup\{T>0,\,T\text{ is a strong existence time}\}.
\]

Our goal is twofold. First, we aim to prove that under Assumption~\ref{hypGeneral} and additional smallness assumptions on the initial data depending on the regime, we have $T^*_\eps=+\infty$, at least for small values of $\eps$.

Then, we shall take a step back and prove that we can dispense with some of these additional smallness assumptions on the initial data at the cost of constraining the time horizon. The issue will be to prove that we can find a minimal time horizon that does not depend on $\eps$ (and therefore does not vanish as $\eps$ goes to 0).

To achieve these goals, we cannot apply the same strategy as in Lemma~\ref{ExistsStrongTime}, since the estimates we used do not provide a uniform bound with respect to $\eps$. We will thus need a new approach to control the Brinkman force that will provide an adequate desingularization, which we will achieve thanks to the changes of variables presented in Section~\ref{SectionChangeVar}. In the \emph{light} and \emph{light and fast} particle regimes, we can use the method of characteristics and suitable identities relying on integration by parts to obtain satisfactory uniform bounds.  
In the fine particle regime, though the previous argument could be refined and applied, we choose to present another strategy based on another functional, which we refer to as higher  dissipation, following \cite{hank}. These different approaches lead to different arguments to prove the expected convergences. In Section \ref{SectionLandLFParticle}, we present the strategy used for the \emph{light} and \emph{light and fast} particle regimes.  Section \ref{SectionFineParticle} is then dedicated to the fine particle regime.


 \section{The \emph{light} and \emph{light and fast} particle regimes}\label{SectionLandLFParticle}

In this section, we study in a unified manner the \emph{light} and \emph{light and fast} particle regimes, that correspond to the choices $(\gamma,\sigma)=(1,\eps^\alpha)$, for some $\alpha\in[0,1/2]$. For the sake of readability, we systematically drop the parameter $\gamma$: for instance, $u_{\eps,1,\sigma}$ will be referred to as $\ued$.
In short we study the behavior as $\eps \xrightarrow{} 0$ of solutions to the systems
\begin{equation*}
    \left\{
\begin{aligned}
&\partial_t u_{\eps,\sigma}+ (u_{\eps,\sigma}\cdot\na_{x})u_{\eps,\sigma}-\Delta_{x}u_{\eps,\sigma}+\na_{x}p_{\eps,\sigma}= j_{f_{\eps,\sigma}}-\rho_{f_{\eps,\sigma}}u_{\eps,\sigma},\\
&\div_{x}u=0,\\
  &\partial_tf_{\eps,\sigma}+ \frac{1}{\sigma}  v\cdot\na_{x}f_{\eps,\sigma}+ \frac{1}{\eps} \div_{v}\left[f_{\eps,\sigma}( \sigma u_{\eps,\sigma}-v)\right]=0, \\
 &\rho_{f_{\eps,\sigma}}(t,x)=\int_{\R^3}f_{\eps,\sigma}(t,x,v)\ddv,  \quad
j_{f_{\eps,\sigma}}(t,x)=\frac{1}{\sigma} \int_{\R^3}vf_{\eps,\sigma}(t,x,v)\ddv,
\end{aligned}
\right.
\end{equation*}
with $\sigma=1$ for the light particle regime and $\sigma= \eps^\alpha$ for the light and fast particle regime.

The main part of the proof is to complete the bootstrap argument set up in Section~\ref{sec-prelim}. 
To this purpose, we shall establish precise estimates on the Brinkman force
\[
F_{\eps,\sigma}= j_{f_{\eps,\sigma}}-\rho_{f_{\eps,\sigma}}u_{\eps,\sigma}
\]
thanks to a desingularization with respect to $\eps$, based on the lagrangian structure of the Vlasov equation.
We then study the limit $\eps\to0$ and prove convergence (and convergence rates) to the solutions of the Transport-Navier-Stokes equations~\eqref{TNS}. This eventually leads to the proof of Theorems~\ref{thm1} and~\ref{thm2}, but with an additional assumption on the initial data:
\begin{equation}
\label{eq-WP-LF}
        \norme{|v|^p\fed^0}_{\L^1(\T^3\times\R^3)}\lesssim\frac{\sigma^p}{\eps^{1-\kappa}},
\end{equation}
for some $\kappa >0$. Note that this is a true restriction only for $\alpha>1/3$. Treating the cases $\alpha \in (1/3,1/2]$ without \eqref{eq-WP-LF} requires a refined strategy and will be the focus of Section~\ref{sec-betterLF}.

This section is organized as follows:
\begin{itemize}
\item The desingularization of the Brinkman force is achieved in Section~\ref{SectionBrinkman-12}, leading to uniform $\L^p$ (in time and space) bounds.
 \item Thanks to these crucial bounds, in Section~\ref{sec-conclubootstrapLLF}, the bootstrap argument set up in Section~\ref{sec-prelim} is concluded under the assumption of a small initial modulated energy. We rely on an $\L^p$ parabolic maximal estimate combined with the exponential decay of the modulated energy through an interpolation argument.
 We therefore obtain that all positive times are \emph{strong existence times} in the sense of Definition~\ref{StrongExistenceTime}.
 \item In Section~\ref{SectionAsymptotic-12}, we perform the proof of convergence for the mildly well-prepared and well-prepared cases of Theorems~\ref{thm1} and~\ref{thm2}, mainly using an $\L^\infty\L^2$--$\L^2\dot\H^1$ energy estimate for the Navier-Stokes part and Wasserstein stability estimates for the Vlasov part.

 \item Finally the general case of Theorems~\ref{thm1} and~\ref{thm2}, in which we obtain a short-time convergence result, is discussed in Section~\ref{sec-smallT}.
\end{itemize}

\subsection{Uniform estimates on the Brinkman force in \texorpdfstring{$\L^p$}{L\string^p}}\label{SectionBrinkman-12}

Our aim is to use the method of characteristics and the changes of variables described in Subsection \ref{SectionChangeVar} to get a precise bound for the Brinkman force. For a strong existence time $T>0$, we shall rely on the exponential decay in time of the modulated energy on $[0,T]$ as proven in Lemma \ref{lemmeDE}.

Consider $\eps\in(0,1)$ and let $T>0$ be a strong existence time. For any $t\in[0,T]$ and $x\in\T^3$, we have, by the method of characteristics,
\[
    F_{\ed}(t,x)=e^{\frac{3t}{\eps}}\int_{\R^3}\fed^0(X_{\ed}(0;t,x,v),V_{\ed}(0;t,x,v))\left(\frac{v}{\sigma}-\ued(t,x)\right)\dd v,
\]
where we recall $(X_{\ed}(0;t,x,v),V_{\ed}(0;t,x,v))$ is defined in~\eqref{characV}.
Thanks to Lemma \ref{changeVarV}, we can apply the change of variables in velocity $v=\Gamma_{\ed}^{t,x}(w)$ so that
\begin{multline*}
    \Fed(t,x)=\\
e^{\frac{3t}{\eps}}\int_{\R^3}\feps^0\left(\tilde{X}_{\ed}^{t,x,w}(0),w\right)\left(\frac{1}{\sigma}[\Gamma_{\ed}^{t,x}]^{-1}(w)-\ued(t,x)\right)|\det\nabla_{w}[\Gamma_{\ed}^{t,x}]^{-1}(w)|\,\dd w,
\end{multline*}
where we recall $\tilde{X}_{\ed}^{t,x,w}(0)$ is defined in~\eqref{eq-defXtilde}
and therefore, still by Lemma~\ref{changeVarV}, this yields
\begin{equation}\label{BrinkmanDecompose-12}
    |\Fed(t,x)|\lesssim\int_{\R^3}\fed^0\left(\tilde{X}_{\ed}^{t,x,w}(0),w\right)\left|\frac{1}{\sigma}[\Gamma_{\ed}^{t,x}]^{-1}(w)-\ued(t,x)\right|\dd w.
\end{equation}
By integrating the equation of characteristics~\eqref{characV} relative to the velocity, we get
\[
    [\Gamma_{\ed}^{t,x}]^{-1}(w)=e^{-\frac{t}{\eps}}w+\frac{\sigma}{\eps}\int_0^te^{\frac{s-t}{\eps}}\ued\left(s,\tilde{X}_{\ed}^{t,x,w}(s)\right)\,\dd s.
\]
In order to desingularize this expression with respect to $\eps$, we perform an integration by parts in time and obtain
\begin{align*}
    [\Gamma_{\ed}^{t,x}]^{-1}(w)
    &=\sigma\ued(t,x)+e^{-\frac{t}{\eps}}\left(w-\sigma\ued^0\left(\tilde{X}_{\ed}^{t,x,w}(0)\right)\right)\\
    &-\sigma\int_0^te^{\frac{s-t}{\eps}}\partial_s\ued\left(s,\tilde{X}_{\ed}^{t,x,w}(s)\right)\dd s\\
    &-\sigma\int_0^te^{\frac{s-t}{\eps}}V_{\ed}\left(s;t,x,[\Gamma_{\ed}^{t,x}]^{-1}(w)\right)\cdot\na_{x}\ued\left(s,\tilde{X}_{\ed}^{t,x,w}(s)\right)\dd s.
\end{align*}
Therefore, \eqref{BrinkmanDecompose-12} can be written as follows.
\begin{lemma}
\label{lem-decompF}
    For every $\eps>0$, if $T>0$ is a strong existence time then for any $(t,x)\in[0,T]\times\T^3$,
    \[
        |\Fed(t,x)|\lesssim F_{\ed}^0+F_{\ed}^{dt}+F_{\ed}^{dx},
    \]
    where
   \begin{align*}
          \Fed^0(t,x)&=e^{\frac{-t}{\eps}}\int_{\R^3}\fed^0\left(\tilde{X}_{\ed}^{t,x,w}(0),w\right)\left|\frac{w}{\sigma}-\ued^0\left(\tilde{X}_{\ed}^{t,x,w}(0)\right)\right|\dd w, \\
          \Fed^{dt}(t,x)&=\int_{\R^3}\int_0^te^{\frac{s-t}{\eps}}\fed^0\left(\tilde{X}_{\ed}^{t,x,w}(0),w\right)\left|\partial_s\ued\left(s,\tilde{X}_{\ed}^{t,x,w}(s)\right)\right|\dd s\dd w, \\
          \Fed^{dx}(t,x)&=
        \int_{\R^3}\int_0^te^{\frac{s-t}{\eps}}\fed^0\left(\tilde{X}_{\ed}^{t,x,w}(0),w\right)\\
           &\qquad\qquad\qquad \times \left|V_{\ed}\left(s;t,x,[\Gamma_{\ed}^{t,x}]^{-1}(w)\right)\cdot\na_{x}\ued\left(s,\tilde{X}_{\ed}^{t,x,w}(s)\right)\right|\dd s\dd w.
    \end{align*}
    \index{$\Fed^0$, $\Fed^{dt}$, $ \Fed^{dx}$: decomposition of the Brinkman force for the \emph{light} and \emph{light and fast} particle regimes}
\end{lemma}

Let us study these three terms separately. In what follows, $p$ stands for the regularity index of Assumption~\ref{hypGeneral}.

\begin{lemma}\label{lemmaF0-12}
    Under Assumption~\ref{hypGeneral}, there exists $\mu_p>0$ such that for every $\eps>0$, if $T>0$ is a strong existence time and if there exists $\kappa\in(0,1)$ such that
    \begin{equation}\label{hyp-wellprep-2}
        \norme{|v|^p\fed^0}_{\L^1(\T^3\times\R^3)}\lesssim\frac{\sigma^p}{\eps^{1-\kappa}},
    \end{equation}
    then
    \[
        \norme{\Fed^0}_{\L^p((0,T)\times\T^3)}\lesssim M^{\mu_p}\eps^{\frac{\kappa}{p}}.
    \]
\end{lemma}

\begin{remark}This is the only term for which the extra assumption~\eqref{hyp-wellprep-2} is useful. As already mentioned, it is actually restrictive only for $\alpha>1/3$, since $p>3$ can be taken arbitrarily close to $3$.
It will be the purpose of Section~\ref{sec-betterLF} to remove this unnecessary assumption.
\end{remark}

\begin{preuve}
    Thanks to Hölder's inequality, we have
    \begin{align*}
        &\norme{\Fed^0}_{\L^{p}((0,T)\times\T^3)}^{p}\\
        &\qquad \lesssim\int_0^Te^{-\frac{pt}{\eps}}\int_{\T^3}\left|\int_{\R^3}\fed^0\left(\tilde{X}_{\ed}^{t,x,w}(0),w\right)\left|\frac{w}{\sigma}-\ued^0\left(\tilde{X}_{\ed}^{t,x,w}(0)\right)\right|\dd w\right|^{p}\dd x\dd t\\
        &\qquad\lesssim\int_0^Te^{-\frac{pt}{\eps}}\int_{\T^3}\left(\int_{\R^3}\fed^0\left(\tilde{X}_{\ed}^{t,x,w}(0),w\right)\dd w\right)^{p-1}\\
        &\qquad \qquad \times\left(\int_{\R^3}\fed^0\left(\tilde{X}_{\ed}^{t,x,w}(0),w\right)\left|\frac{w}{\sigma}-\ued^0\left(\tilde{X}_{\ed}^{t,x,w}(0)\right)\right|^{p}\dd w\right)\dd x.
    \end{align*}
    We apply Lemma \ref{changeVarX} and the change of variables in space $x'=\tilde{X}_{\ed}^{t,x,w}(0)$ so that
    \begin{align*}
        &\int_{\T^3\times\R^3}\fed^0\left(\tilde{X}_{\ed}^{t,x,w}(0),w\right)\left|\frac{w}{\sigma}-\ued^0\left(\tilde{X}_{\ed}^{t,x,w}(0)\right)\right|^p\dd x\dd w\\
        &\qquad \qquad \lesssim\sigma^{-p}\int_{\T^3\times\R^3}|w|^p\fed^0(x,w)\dd x\dd w+\norme{\fed^0}_{\L^1(\R^3;\L^\infty(\T^3))}\norme{\ued^0}_{\L^p(\T^3)}^p\\
        &\qquad \qquad \lesssim\eps^{\kappa-1}+M^{\mu_p},
    \end{align*}
    for some $\mu_p>0$, hence
    \[
        \norme{\Fed^0}_{\L^p((0,T)\times\T^3)}^p\lesssim\eps^{\kappa}M^{\mu_pp},
    \]
   up to modifying the value of $\mu_p>0$.
\end{preuve}

We deal with the second term in a similar fashion.

\begin{lemma}\label{lemmaFdt-12}
    Under Assumption~\ref{hypGeneral}, for every $\eps>0$, if $T>0$ is a strong existence time,
    \[
        \norme{\Fed^{dt}}_{\L^p((0,T)\times\T^3)}\lesssim \eps M\norme{\partial_t\ued}_{\L^p((0,T)\times\T^3)}.
    \]
\end{lemma}

\begin{preuve}
    Thanks to Hölder's inequality, for almost all $(t,x)\in(0,T)\times\T^3$,
    \begin{align*}
        |\Fed^{dt}(t,x)|^p
        &\le\left(\int_{\R^3}\int_0^te^{\frac{s-t}{\eps}}\fed^0\left(\tilde{X}_{\ed}^{t,x,w}(0),w\right)\dd s\dd w\right)^{p-1}\\
        &\quad    \times\left(\int_{\R^3}\int_0^te^{\frac{s-t}{\eps}}\fed^0\left(\tilde{X}_{\ed}^{t,x,w}(0),w\right)\left|\partial_s\ued\left(s,\tilde{X}_{\ed}^{t,x,w}(s)\right)\right|^{p}\dd s\dd w\right)\\
        &\lesssim\norme{\fed^0}_{\L^1(\R^3;\L^\infty(\T^3))}^{p-1}\eps^{p-1}\\
         &\quad   \times\int_0^te^{\frac{s-t}{\eps}}\int_{\R^3}\fed^0\left(\tilde{X}_{\ed}^{t,x,w}(0),w\right)\left|\partial_s\ued\left(s,\tilde{X}_{\ed}^{t,x,w}(s)\right)\right|^{p}\dd w\dd s.
    \end{align*}
    Applying Lemma \ref{changeVarX}, we can now use the change of variables in space $x'=\tilde{X}_{\ed}^{t,x,w}(s)=\Phi_{\ed}^{s;t,w}(x)$ in the following integral to obtain
    \begin{align*}
        &\int_{\T^3\times\R^3}\fed^0\left(\tilde{X}_{\ed}^{t,x,w}(0),w\right)\left|\partial_s\ued\left(s,\tilde{X}_{\ed}^{t,x,w}(s)\right)\right|^{p}\dd x\dd w\\
        &\qquad\lesssim\int_{\T^3\times\R^3}\fed^0\left(\tilde{X}_{\ed}^{t,[\Phi_{\ed}^{s;t,w}]^{-1}(x),w}(0),w\right)|\partial_s\ued(s,x)|^{p}\dd x\dd w\\
        &\qquad\lesssim\norme{\fed^0}_{\L^1(\R^3;\L^\infty(\T^3))}\norme{\partial_s\ued(s)}_{\L^{p}(\T^3)}^{p}.
    \end{align*}
    Therefore, we get
    \begin{align*}
        &\norme{\Fed^{dt}}_{\L^{p}((0,T)\times\T^3)}^{p} \\
        &\qquad\lesssim\norme{\fed^0}_{\L^1(\R^3;\L^\infty(\T^3))}^{p}\eps^{p-1}\int_0^T\norme{\partial_s\ued(s)}_{\L^p(\T^3)}^{p}\left(\int_s^Te^{\frac{s-t}{\eps}}\dd t\right)\dd s\\
          &\qquad\lesssim\eps^p\norme{\fed^0}_{\L^1(\R^3;\L^\infty(\T^3))}^p\norme{\partial_t\ued}_{\L^p((0,T)\times\T^3)}^p,
    \end{align*}
    which concludes the proof of the lemma.
\end{preuve}

The control of the last term in the Brinkman force is more intricate. We shall rely on the exponential decay of the modulated energy.

\begin{lemma}\label{lemmaFdx-12}
    Under Assumption~\ref{hypGeneral}, there exists $\mu_p\ge0$ such that for every $\eps>0$, if $T>0$ is a strong existence time,
    \[
        \norme{\Fed^{dx}}_{\L^p((0,T)\times\T^3)}\lesssim\eps\norme{\Delta_x\ued}_{\L^p((0,T)\times\T^3)}+\eps M^{\mu_p}\EEE_{\ed}(0)^{\frac{1}{2}}.
    \]
\end{lemma}

\begin{preuve}
    By integrating the equation of characteristics, we get, for all $0\le s\le t\le T$,
    \[
        V_{\ed}(s;t,x,[\Gamma_{\ed}^{t,x}]^{-1}(w))=e^{-\frac{s}{\eps}}w+\frac{\sigma}{\eps}\int_0^se^{\frac{\tau-s}{\eps}}\ued\left(\tau,\tilde{X}_{\ed}^{t,x,w}(\tau)\right)\dd\tau,
    \]
    so that
    \begin{equation*}
        \begin{aligned}
            |\Fed^{dx}(t,x)|&\le\int_{\R^3}e^{-\frac{t}{\eps}}\fed^0\left(\tilde{X}_{\ed}^{t,x,w}(0),w\right)|w|\int_0^t|\na_x\ued\left(s,\tilde{X}_{\ed}^{t,x,w}(s)\right)|\dd s\dd w\\
            &+\frac{\sigma}{\eps}\int_{\R^3}\fed^0\left(\tilde{X}_{\ed}^{t,x,w}(0),w\right) \times\\
            &\quad  \int_0^t\int_0^se^{\frac{s-t}{\eps}}e^{\frac{\tau-s}{\eps}}\left|\ued\left(\tau,\tilde{X}_{\ed}^{t,x,w}(\tau)\right)\right|\left|\nabla_{x}\ued\left(s,\tilde{X}_{\ed}^{t,x,w}(s)\right)\right|
            \dd\tau\dd s\dd w\\
            &=:I_1(t,x)+I_2(t,x).
        \end{aligned}
    \end{equation*}
    Thanks to Hölder's inequality, we have
    \begin{align*}
        |I_1(t,x)|^p&\le e^{-\frac{pt}{\eps}}\left[\int_0^t\left(\int_{\R^3}\fed^0\left(\tilde{X}_{\ed}^{t,x,w}(0),w\right)|w|^{\frac{p}{p-1}}\dd w\right)^{\frac{p-1}{p}}\right.\\
        &\quad\times\left.\left(\int_{\R^3}\fed^0\left(\tilde{X}_{\ed}^{t,x,w}(0),w\right)\left|\na_{x}\ued\left(s,\tilde{X}_{\ed}^{t,x,w}(s)\right)\right|^{p}\dd w\right)^{\frac{1}{p}}\dd s\right]^p\\
        &\le\norme{|v|^{\frac{p}{p-1}}\fed^0}_{\L^1(\R^3;\L^\infty(\T^3))}^{p-1}e^{-\frac{pt}{\eps}}\left[\int_0^te^{\frac{s}{\eps}}\eps(1-e^{-\frac{t}{\eps}})\right.\\
        &\quad\times\left.\left(\int_{\R^3}\fed^0(\tilde{X}_{\ed}^{t,x,w}(0),w)|\nabla_{x}\ued(s,\tilde{X}_{\ed}^{t,x,w}(s))|^{p}\dd w\right)^{\frac{1}{p}}\frac{e^{-\frac{s}{\eps}}\dd s}{\eps(1-e^{\frac{-t}{\eps}})}\right]^{p},
    \end{align*}
    and thanks to Jensen's inequality applied on the probability space 
    $$\left((0,t),\frac{e^{-\frac{s}{\eps}}\dd s}{\eps(1-e^{-\frac{t}{\eps}})}\right),$$ 
    we obtain
    \begin{multline*}
        |I_1(t,x)|^p\le\eps^{p-1}\norme{|v|^{\frac{p}{p-1}}\fed^0}_{\L^1(\R^3;\L^\infty(\T^3))}^{p-1}\\
        \times\int_0^te^{\frac{p(s-t)}{\eps}}\int_{\R^3}\fed^0\left(\tilde{X}_{\ed}^{t,x,w}(0),w\right)\left|\na_x\ued\left(s,\tilde{X}_{\ed}^{t,x,w}(s)\right)\right|^p\dd w\dd s.
    \end{multline*}
    We apply Lemma \ref{changeVarX} and use the change of variables $x'=\tilde{X}_{\ed}^{t,x,w}(s)$. Therefore, as in the proof of Lemma~ \ref{lemmaFdt-12}, we get
    \begin{multline*}
        \int_0^T\int_{\T^3}|I_1(t,x)|^p\dd x\dd t
        \le\eps^{p-1}\norme{|v|^{\frac{p}{p-1}}\fed^0}_{\L^1(\R^3;\L^\infty(\T^3))}^{p-1}\\
            \times\int_0^T\int_0^te^{\frac{p(s-t)}{\eps}}\left(\int_{\T^3\times\R^3}\fed^0\left(\tilde{X}_{\ed}^{t,x,w}(0),w\right)\left|\nabla_{x}\ued\left(s,\tilde{X}_{\ed}^{t,x,w}(s)\right)\right|^{p}\dd x\dd w\right)\dd s\dd t\\
        \lesssim\eps^{p-1}\norme{|v|^{\frac{p}{p-1}}\fed^0}_{\L^1(\R^3;\L^\infty(\T^3))}^{p-1}\norme{\fed^0}_{\L^1(\R^3;\L^\infty(\T^3))}\\
        \times\int_0^T\int_0^te^{\frac{p(s-t)}{\eps}}\norme{\nabla_{x}\ued(s)}_{\L^p(\T^3)}^{p}\dd s\dd t.
    \end{multline*}
    Since $T$ is a strong existence time, we can apply Lemmas~\ref{lemmeDE}--\ref{estimNabla}. Then, thanks to Hölder's and Young's inequality as well as the value of $\lambda_{\ed}$ given by Lemma~\ref{lemmeDE} and Corollary~\ref{RhoLinfini}, we conclude that
    \begin{align*}
        &\int_0^T\int_{\T^3}|I_1(t,x)|^p\dd x\dd t\\
        &\qquad \qquad\lesssim\eps^{p-1}\norme{|v|^{\frac{p}{p-1}}\fed^0}_{\L^1(\R^3;\L^\infty(\T^3))}^{p-1}\norme{\fed^0}_{\L^1(\R^3;\L^\infty(\T^3))}\EEE_{\ed}(0)^{\frac{p(1-\alpha_p)}{2}}\\
          &\qquad \qquad \qquad  \times\int_0^T\norme{\Delta_{x}\ued(s)}_{\L^p(\T^3)}^{p\alpha_p}e^{-\frac{p(1-\alpha_p)\lambda_\eps}{2}s}\left(\int_s^Te^{\frac{p(s-t)}{\eps}}\dd t\right)\dd s\\
        &\qquad \qquad\lesssim\eps^p\norme{|v|^{\frac{p}{p-1}}\fed^0}_{\L^1(\R^3;\L^\infty(\T^3))}^{p-1}\norme{\fed^0}_{\L^1(\R^3;\L^\infty(\T^3))}\EEE_{\ed}(0)^{\frac{p(1-\alpha_p)}{2}}\frac{1}{\lambda_{\ed}^{1-\alpha_p}}\\
        &\qquad \qquad \qquad  \times\norme{\Delta_{x}\ued}_{\L^p((0,T)\times\T^3)}^{p\alpha_p}\\
        &\qquad \qquad \lesssim\eps^{p}\norme{|v|^{\frac{p}{p-1}}\fed^0}_{\L^1(\R^3;\L^\infty(\T^3))}^{\frac{p-1}{1-\alpha_p}}\norme{\fed^0}_{\L^1(\R^3;\L^\infty(\T^3))}^{\frac{1}{1-\alpha_p}}\EEE_{\ed}(0)^{\frac{p}{2}}\\
        &\qquad \qquad \qquad    +\eps^{p}\norme{\Delta_{x}\ued}_{\L^p((0,T)\times\T^3)}^{p}.
    \end{align*}
    
    The control of $I_2$ is a little more technical. We apply Hölder's inequality twice and obtain
    \begin{multline*}
        |I_2(t,x)|^p
        \le\left[\eps^{-1}\int_{\R^3}\fed^0\left(\tilde{X}_{\ed}^{t,x,w}(0),w\right)\int_0^t\left|\na_x\ued\left(s,\tilde{X}_{\ed}^{t,x,w}(s)\right)\right|e^{\frac{s-t}{\eps}}\right.\\
            \times\left.\left(\int_0^s\left|\ued\left(\tau,\tilde{X}_{\ed}^{t,x,w}(\tau)\right)\right|^pe^{\frac{p(\tau-s)}{2\eps}}\dd\tau\right)^{\frac{1}{p}}\left(\int_0^s e^{\frac{p(\tau-s)}{2(p-1)\eps}}\dd\tau\right)^{\frac{p-1}{p}}\dd s\dd w\right]^p\\
        \lesssim\eps^{-1}\left[\int_{\R^3}\int_0^t\fed^0\left(\tilde{X}_{\ed}^{t,x,w}(0),w\right)\left|\na_x\ued\left(s,\tilde{X}_{\ed}^{t,x,w}(s)\right)\right|e^{\frac{s-t}{2\eps}}\right.\\
            \times\left. e^{\frac{s-t}{2\eps}}\left(\int_0^s\left|\ued\left(\tau,\tilde{X}_{\ed}^{t,x,w}(\tau)\right)\right|^pe^{\frac{p(\tau-s)}{2\eps}}\dd\tau\right)^{\frac{1}{p}}\dd s\dd w\right]^p\\
        \lesssim\eps^{-1}\left(\int_{\R^3}\int_0^t\fed^0\left(\tilde{X}_{\ed}^{t,x,w}(0),w\right)e^{\frac{p(s-t)}{2\eps}}\int_0^s\left|\ued\left(\tau,\tilde{X}_{\ed}^{t,x,w}(\tau)\right)\right|^pe^{\frac{p(\tau-s)}{2\eps}}\dd\tau\dd s\dd w\right)\\
            \times\left(\int_{\R^3}\int_0^t\fed^0\left(\tilde{X}_{\ed}^{t,x,w}(0),w\right)\left|\na_x\ued\left(s,\tilde{X}_{\ed}^{t,x,w}(s)\right)\right|^{\frac{p}{p-1}}e^{\frac{p(s-t)}{2(p-1)\eps}}\dd s\dd w\right)^{p-1}.
    \end{multline*}
    Therefore,
    \[
        \int_0^T\int_{\T^3}|I_2(t,x)|^p\dd x\dd t\lesssim\eps^{-1}\norme{\fed^0}_{\L^1(\R^3;\L^\infty(\T^3))}^{p-1}\left(\sup_{t\in(0,T)}I_3(t)\right)\times I_4,
    \]
    where
    \begin{multline*}
        I_3(t):=\int_0^te^{\frac{p(s-t)}{2\eps}}\int_0^se^{\frac{p(\tau-s)}{2\eps}}\\
            \int_{\T^3\times\R^3}\fed^0\left(\tilde{X}_{\ed}^{t,x,w}(0),w\right)\left|\ued\left(\tau,\tilde{X}_{\ed}^{t,x,w}(\tau)\right)\right|^p\dd x\dd w\dd\tau\dd s,
    \end{multline*}
    and
    \[
        I_4:=\int_0^T\left(\int_0^t\norme{\na_x\ued(s)}_{\L^\infty(\T^3)}^{\frac{p}{p-1}}e^{\frac{p(s-t)}{2(p-1)\eps}}\dd s\right)^{p-1}\dd t.
    \]
    
    On the one hand, we can apply Lemma~\ref{EstimNormLp} and Lemma~\ref{changeVarX} with the change of variable $x'=\tilde{X}_{\ed}^{t,x,w}(\tau)$ as before and obtain
    \[
        \sup_{t\in(0,T)}I_3(t)\lesssim\eps^2\norme{\fed^0}_{\L^1(\R^3;\L^\infty(\T^3))}\Psi_{\eps,0}^{\frac{p}{2}}.
    \]
    On the other hand, using Jensen's inequality, we have
    \[
        I_4\lesssim\eps^{p-2}\int_0^T\int_0^t\norme{\na_x\ued(s)}_{\L^\infty(\T^3)}^p e^{\frac{p(s-t)}{2(p-1)\eps}}\dd s\dd t.
    \]
    Thanks to Lemma~\ref{estimNabla} and Hölder's inequality, we find that
    \[
        I_4\lesssim\eps^{p-1}\EEE_{\ed}(0)^{\frac{(1-\beta_p)p}{2}}\norme{\Delta_x\ued}_{\L^p((0,T)\times\T^3)}^{\beta_pp}.
    \]
    
    We can conclude that
    \begin{align*}
        &\int_0^T\int_{\T^3}|I_2(t,x)|^p\dd x\dd t\\
        &\qquad\lesssim\eps^p\norme{\fed^0}_{\L^1(\R^3;\L^\infty(\T^3))}^p\Psi_{\ed,0}^{\frac{p}{2}}\EEE_{\ed}(0)^{\frac{(1-\beta_p)p}{2}}\norme{\Delta_x\ued}_{\L^p((0,T)\times\T^3))}^{\beta_pp}\\
        &\qquad \lesssim\eps^p\norme{\fed^0}_{\L^1(\R^3;\L^\infty(\T^3))}^{\frac{p}{1-\beta_p}}\Psi_{\ed,0}^{\frac{p}{2(1-\beta_p)}}\EEE_{\ed}(0)^{\frac{p}{2}}+\eps^p\norme{\Delta_x\ued}_{\L^p((0,T)\times\T^3))}^p.
    \end{align*}
    Gathering all pieces together, this allows to conclude the proof.
\end{preuve}


\subsection{Conclusion of the bootstrap argument in the mildly well-prepared case}
\label{sec-conclubootstrapLLF}

Recall that, as stated in Section~\ref{SubsectionInitialization}, we consider, under Assumption~\ref{hypGeneral}, for every $\eps>0$,
\[
    T^*_\eps=\sup\{T>0,\,T\text{ is a strong existence time}\}.
\]

The aim of this section is to prove that under the additional Assumption~\ref{hypSmallDataNS}, \eqref{hyp-wellprep-2} and the smallness of the initial modulated energy, we have $T^*_\eps=+\infty$, at least for $\eps$ close to 0. We only need to prove that, for every $T<T^*_\eps$,
\[
    \norme{\na_x\ued}_{\L^1(0,T;\L^\infty(\T^3))}\le\frac{1}{40}\qquad\text{and}\qquad\norme{F_{\ed}}_{\L^2((0,T)\times\T^3)}\le\frac{C^*}{4}.
\]
Indeed, as we have seen in the proof of Lemma~\ref{EstimationNSsupGeneral}, Assumption~\ref{hypSmallDataNS} already ensures that, for all $T\ge0$,
\[
    \int_0^T\norme{e^{t\Delta}\ued}_{\dot\H^1(\T^3)}^4\dd t\le\frac{C^*}{4}.
\]
Then, thanks to the timewise continuity of the norms that are involved in the definition of a strong existence time, we deduce that $T^*_\eps=+\infty$.

We begin by proving the following lemma.
\begin{lemma}\label{lemmePara-12}
    Under Assumption~\ref{hypGeneral}, there exists $\eps_0>0$ 
    such that for all $\eps\in(0,\eps_0)$ if~\eqref{hyp-wellprep-2} holds
    and $T<T^*_\eps$, then
    \[
        \norme{\partial_t\ued}_{\L^p((0,T)\times\T^3)}+\norme{\Delta_x\ued}_{\L^p((0,T)\times\T^3)}\lesssim M^{\omega_p},
    \]
    for some $\omega_p>0$.
\end{lemma}

\begin{preuve}
    Thanks to Theorem~\ref{ParabolicReg}, we have
    \begin{multline*}
        \norme{\partial_t\ued}_{\L^p((0,T)\times\T^3)}+\norme{\Delta_x\ued}_{\L^p((0,T)\times\T^3)}\\
        \lesssim\norme{\Fed}_{\L^p((0,T)\times\T^3)}+\norme{(\ued\cdot\na_x)\ued}_{\L^p((0,T)\times\T^3)}+\norme{\ued^0}_{\B_p^{s,p}(\T^3)}.
    \end{multline*}
    We can now apply the results from Lemmas \ref{EstimConvectStrongT} and \ref{lemmaF0-12}--\ref{lemmaFdx-12}. We find that
    \begin{align*}
        &\norme{\partial_t\ued}_{\L^p((0,T)\times\T^3)}+\norme{\Delta_x\ued}_{\L^p((0,T)\times\T^3)}\\
        &\qquad \qquad \qquad \lesssim \eps M\norme{\partial_t\ued}_{\L^p((0,T)\times\T^3)}+\eps\norme{\Delta_x\ued}_{\L^p((0,T)\times\T^3)}\\
        &\qquad \qquad \qquad+\eps^{\kappa}M^{\mu_p}+\eps M^{\mu_p}\EEE_{\ed}(0)^{\frac{1-\beta_p}{2}}\\
        &\qquad \qquad \qquad+\Psi_{\ed,0}^{\frac{1}{2}}\EEE_{\ed}(0)^{\frac{1-\beta_p}{2}}\norme{\Delta_x\ued}_{\L^p((0,T)\times\T^3)}^{\beta_p}+M
    \end{align*}
    and we conclude thanks to Young's inequality and taking $\eps$ small enough.
\end{preuve}

We have the following immediate consequence.

\begin{corollary}
\label{coro-nau-12}
    Under Assumption~\ref{hypGeneral}, there exist $\eps_0>0$ and $\eta>0$ 
    such that, for all $\eps\in(0,\eps_0)$, if~\eqref{hyp-wellprep-2} holds and
    \begin{equation}\label{hypEnergieModuleePetite-12}
        \EEE_{\ed}(0)\le\eta,
    \end{equation}
    then for any $T<T^*_\eps$,
    \[
        \norme{\na\ued}_{\L^1(0,T;\L^\infty(\T^3))}\le\frac{1}{40}.
    \]
\end{corollary}

\begin{preuve}
    Since $p>3$, we can apply Lemmas~\ref{estimNabla} and~\ref{lemmePara-12} and find that there exist $\eps_0$ 
    such that for all $\eps\in(0,\eps_0)$
    , for some $\omega_p>0$, we have
    \begin{align*}
        \int_0^T\norme{\na\ued(t)}_{\L^\infty(\T^3)}\dd t&\lesssim\EEE_{\ed}(0)^{\frac{1-\beta_p}{2}}\norme{\Delta_x\ued}_{\L^p((0,T)\times\T^3)}^{\beta_p}\\
        &\lesssim \EEE_{\ed}(0)^{\frac{1-\beta_p}{2}}M^{\omega_p\beta_p},
    \end{align*}
    which can be made as small as necessary by reducing the value of $\eta>0$ in~\eqref{hypEnergieModuleePetite-12}.
\end{preuve}

There remains to verify that
\begin{equation}
\label{Fedpetit}
    \norme{F_\ed}_{\L^2((0,T)\times\T^3)}\le\cfrac{C^*}{4}.
\end{equation}

We shall provide two proofs of this fact. The first one uses~\eqref{hypEnergieModuleePetite-12} and is straightforward, while the second one does not and is more technical, being based on the methods developed earlier to prove $\L^p\L^p$ estimates for the Brinkman force. It will however come in handy in order to 
\begin{itemize}
\item obtain convergence rates for the fluid velocity (Section~\ref{SectionAsymptotic-12});
\item treat the general case (Section~\ref{sec-smallT}).
\end{itemize}

\begin{lemma}
\label{firstproof}
   Under Assumption~\ref{hypGeneral} and~\ref{hypEnergieModuleePetite-12}, if $T$ is a strong existence time,
   \[
   \norme{F_\ed}_{\L^2(\R_+;\L^2(\T^3))} \leq \norme{f^0_\ed}_{\L^{1}(\R^3; \L^\infty(\T^3))}^{1/2} \mathscr{E}_{\ed}(0)^{1/2}.
   \]
\end{lemma}

\begin{preuve}
The argument is almost the same as in Remark~\ref{RemarkBorneUnifF}, except that we use the modulated energy--dissipation inequality of Lemma~\ref{ModulatedEnergyDissipation}. This yields
$$
  \norme{F_\ed}_{\L^2(\R_+;\L^2(\T^3))}^2 \leq \norme{\rho_\ed}_{\L^{\infty}(\R_+; \L^\infty(\T^3))} \mathscr{E}_{\ed}(0),
  $$
  and we conclude by Corollary~\ref{RhoLinfini}.
\end{preuve}

We deduce that~\eqref{Fedpetit} holds, taking $\eta$ small enough in~\eqref{hypEnergieModuleePetite-12}.
Let us now present the second approach leading to~\eqref{Fedpetit}. The idea is to mimic  Lemmas~\ref{lemmaF0-12}--\ref{lemmaFdx-12}, which results in the following lemma.
\begin{lemma}\label{lemmaFL2-12}
    Under Assumption~\ref{hypGeneral}, there exists $\eps_0>0$ and $\eta>0$ such that, for all $\eps\in(0,\eps_0)$, if~\eqref{hyp-wellprep-2}--\eqref{hypEnergieModuleePetite-12} hold, then for any $T<T^*_\eps$,
    \begin{equation}\label{FL2-bootstrap-12}
        \norme{\Fed}_{\L^2((0,T)\times\T^3)}\lesssim\eps M\norme{\partial_t\ued}_{\L^2((0,T)\times\T^3)}+\eps\norme{\Delta_x\ued}_{\L^2((0,T)\times\T^3)}+\eps^{\frac{\kappa}{2}}M^{\mu_2},
    \end{equation}
    for some $\mu_2>0$.
\end{lemma}

\begin{preuve}
Note that we can reproduce the proofs of Lemmas~\ref{lemmaF0-12}--\ref{lemmaFdt-12} and obtain
\begin{align*}
    \norme{\Fed^0}_{\L^2((0,T)\times\T^3)}&\lesssim M^{\mu_2}\eps^{\frac{\kappa}{2}}, \\
    \norme{\Fed^{dt}}_{\L^2((0,T)\times\T^3)}&\lesssim\eps M\norme{\partial_t\ued}_{\L^2((0,T)\times\T^3)}.
\end{align*}
Nevertheless,  the Gagliardo-Nirenberg inequality  cannot be applied exactly in the same way and some adaptation of the proof of Lemma~\ref{lemmaFdx-12} is required. We have
\[
    |\Fed^{dx}|\le I_1+I_2,
\]
where, for $(t,x)\in\R_+\times\T^3$,
\[
    I_1(t,x)=\int_{\R^3}e^{-\frac{t}{\eps}}\fed^0\left(\tilde{X}_{\ed}^{t,x,w}(0),w\right)|w|\int_0^t\left|\na_x\ued\left(s,\tilde{X}_{\ed}^{t,x,w}(s)\right)\right|\dd s\dd w,
\]
and
\begin{multline*}
    I_2(t,x)
    =\frac{\sigma}{\eps}\int_{\R^3}\fed^0\left(\tilde{X}_{\ed}^{t,x,w}(0),w\right)\\
    \times\int_0^t\int_0^se^{\frac{s-t}{\eps}}e^{\frac{\tau-s}{\eps}}\left|\ued\left(\tau,\tilde{X}_{\ed}^{t,x,w}(\tau)\right)\right|\left|\na_x\ued\left(s,\tilde{X}_{\ed}^{t,x,w}(s)\right)\right|\dd\tau\dd s\dd w.
\end{multline*}
We can still follow the proof of Lemma~\ref{lemmaFdx-12} and obtain
\[
    \norme{I_1}_{\L^2((0,T)\times\T^3)}\lesssim\eps M+\eps\norme{\Delta_x\ued}_{\L^2((0,T)\times\T^3)}.
\]
For $I_2$, we get
\[
    \int_0^T\int_{\T^3}|I_2(t,x)|^2\dd x\dd t\lesssim\eps^{-1}\norme{\fed^0}_{\L^1(\R^3;\L^\infty(\T^3))}\left(\sup_{t\in(0,T)}I_3(t)\right)\times I_4,
\]
where
\begin{multline*}
    I_3(t)=\int_0^te^{\frac{s-t}{\eps}}\int_0^se^{\frac{\tau-s}{\eps}}\\
        \int_{\T^3\times\R^3}\fed^0\left(\tilde{X}_{\ed}^{t,x,w}(0),w\right)\left|\ued\left(\tau,\tilde{X}_{\ed}^{t,x,w}(\tau)\right)\right|^2\dd x\dd w\dd\tau\dd s,
\end{multline*}
and
\[
    I_4(t)=\int_0^T\left(\int_0^t\norme{\na_x\ued(s)}_{\L^\infty(\T^3)}^2e^{\frac{s-t}{\eps}}\dd s\right)\dd t.
\]
We still have
\[
    \sup_{t\in(0,T)}I_3(t)\lesssim\eps^2\norme{\fed^0}_{\L^1(\R^3;\L^\infty(\T^3))}\Psi_{\eps,\sigma,0}\lesssim\eps^2 M^{\mu_2}.
\]
Finally, for what concerns $I_4$, thanks to Jensen's inequality and Lemmas~\ref{estimNabla} and~\ref{lemmePara-12}, for $\eps$ small enough, we have
\[
    I_4\lesssim\eps\EEE_{\ed}(0)^{\frac{1-\beta_p}{2}}\norme{\Delta_x\ued}_{\L^p((0,T)\times\T^3)}^{\beta_p}\lesssim\eps M^{\nu_2}.
\]
The proof of the result is complete.
\end{preuve}

Therefore, the $\L^2$ parabolic estimate can be written as follows.
\begin{lemma}\label{lemmePara2-12}
    Under Assumption~\ref{hypGeneral}, there exists $\eps_0>0$ and $\eta>0$ such that, for all $\eps\in(0,\eps_0)$, if~\eqref{hyp-wellprep-2}--\eqref{hypEnergieModuleePetite-12} hold, then for any $T<T^*_\eps$,
    \[
    \norme{\partial_t\ued}_{\L^2((0,T)\times\T^3)}+\norme{\Delta_x\ued}_{\L^2((0,T)\times\T^3)}\lesssim M^{\omega_2},
    \]
    for some $\omega_2>0$.
\end{lemma}

\begin{preuve}
    Theorem~\ref{ParabolicReg} ensures that
    \begin{multline*}
        \norme{\partial_t\ued}_{\L^2((0,T)\times\T^3)}+\norme{\Delta_x\ued}_{\L^2((0,T)\times\T^3)}\\
        \lesssim\norme{\Fed}_{\L^2((0,T)\times\T^3)}+\norme{(\ued\cdot\na_x)\ued}_{\L^2((0,T)\times\T^3)}+\norme{\ued^0}_{\H^1(\T^3)}.
    \end{multline*}
    Lemma~\ref{lemmaFL2-12} enables us to control the first right-hand side term. On the other hand, thanks to Hölder's inequality and Lemmas~\ref{EstimationNSsupGeneral},~\ref{lemmeDE}--\ref{estimNabla} and~\ref{lemmePara-12}, we have, for $\eps$ small enough,
    \begin{align*}
        \norme{(\ued\cdot\na)\ued}_{\L^2((0,T)\times\T^3)}&
        \lesssim\norme{\ued}_{\L^\infty(0,T;\L^2(\T^3))}^2\int_0^T\norme{\na_x\ued(t)}_{\L^\infty(\T^3)}^2\dd t\\
        &\lesssim\Psi_{\eps,\sigma,0}\EEE_\ed(0)^{1-\beta_p}\int_0^T\norme{\Delta_x\ued(t)}_{\L^p(\T^3)}e^{-(1-\beta_p)\lambda_\eps t}\dd t\\
        &\lesssim\Psi_{\eps,\sigma,0}\EEE_\ed(0)^{1-\beta_p}\norme{\Delta_x\ued}_{\L^p((0,T)\times\T^3)}^{2\beta_p}
        \lesssim M^{\omega_p}
    \end{align*}
    for some $\omega_p>0$, and the result follows.
\end{preuve}

Injecting this in~\eqref{FL2-bootstrap-12}, we obtain, by taking $\eps$ small enough, that
\[
    \norme{\Fed}_{\L^2((0,T)\times\T^3)}\le\frac{C^*}{4}.
\]

The bootstrap argument is finally complete.

\begin{remark}\label{RemarkEnergieModuleePetite}
    Let us make the small modulated energy assumption from~\eqref{hypEnergieModuleePetite-12} more explicit. This will justify the assumptions stated in Theorems~\ref{thm1}--\ref{thm2} (in the so-called mildly-well-prepared case).
    \begin{itemize}
        \item In the light particle regime, we have
            \begin{multline*}
                \EEE_{\eps,1,1}(0)=
                \frac{\eps}{2}\int_{\T^3\times\R^3}\left|v-\frac{\langle j_{\eps,1,1}^0\rangle}{\langle\rho_{\eps,1,1}^0\rangle}\right|^2f_{\eps,1,1}^0(x,v)\dd x\dd v\\
                +\frac{1}{2}\int_{\T^3}|u_{\eps,1,1}^0(x)-\langle u_{\eps,1,1}^0\rangle|^2\dd x
                +\frac{\eps\langle\rho_{\eps,1,1}^0\rangle}{2(1+\eps \langle\rho_{\eps,1,1}^0\rangle)}|\langle j_{\eps,1,1}^0\rangle-\langle u_{\eps,1,1}^0\rangle|^2.
            \end{multline*}
            Assumption~\ref{hypGeneral} ensures that the first and last term can be made as small as needed by reducing $\eps$. There only remains the middle term, and we can thus replace~\eqref{hypEnergieModuleePetite-12}, up to taking $\eps>0$ small enough, by the assumption
                \begin{equation}
            \label{eq-Erond0petit-L}
                \norme{u_{\eps,1,1}^0-\langle u_{\eps,1,1}^0\rangle}_{\L^2(\T^3)}\le\eta,
            \end{equation}
            up to reducing $\eta$.
        \item In the light and fast particle regime, we have
            \begin{multline*}
                \EEE_{\eps,1,\eps^\alpha}(0)=
                \frac{\eps}{2}\int_{\T^3\times\R^3}\left|\frac{v}{\eps^\alpha}-\frac{\langle j_{\eps,1,\eps^\alpha}^0\rangle}{\langle\rho_{\eps,1,\eps^\alpha}^0\rangle}\right|^2f_{\eps,1,\eps^\alpha}^0(x,v)\dd x\dd v\\
                +\frac{1}{2}\int_{\T^3}|u_{\eps,1,\eps^\alpha}^0(x)-\langle u_{\eps,1,\eps^\alpha}^0\rangle|^2\dd x
                +\frac{\eps \langle\rho_{\eps,1,\eps^\alpha}^0\rangle}{2(1+\eps \langle\rho_{\eps,1,\eps^\alpha}^0\rangle)}|\langle j_{\eps,1,\eps^\alpha}^0\rangle-\langle u_{\eps,1,\eps^\alpha}^0\rangle|^2,
            \end{multline*}
            so that, thanks to Assumption~\ref{hypGeneral} and~\eqref{hyp-wellprep-2}, to ensure that~\eqref{hypEnergieModuleePetite-12} holds, we need only assume
            \begin{equation}
            \label{eq-Erond0petit-LF}
                \norme{u_{\eps,1,\eps^\alpha}^0-\langle u_{\eps,1,\eps^\alpha}^0\rangle}_{\L^2(\T^3)}\le\eta.
            \end{equation}
            Without the extra assumption~\eqref{hyp-wellprep-2}, we need to impose that 
            \begin{equation*}
             {\eps^{1-2\alpha}}\int_{\T^3\times\R^3}|{v}|^2f_{\eps,1,\eps^\alpha}^0(x,v)\dd x\dd v \le \eta,
            \end{equation*}
            which is only restrictive for $\alpha=1/2$, and accounts for the statement of Theorem~\ref{thm2}.
    \end{itemize}
\end{remark}


\subsection{Convergence in the mildly well-prepared and  well-prepared cases}\label{SectionAsymptotic-12}

We study the convergence of the sequences $(\fed)_{\eps>0}$, $(\rhoed)_{\eps>0}$ and $(\ued)_{\eps>0}$ under several sets of assumptions for the initial data. We are able to provide rates of convergences in some of the cases. For the sake of readability, we do not keep track very precisely of the dependence on the initial data and express these rates in terms of the global initial bound $M>1$ defined in Assumption~\ref{hypGeneral}.

Our first result only requires the weak convergence of the fluid velocity and the particle density.

\begin{theorem}\label{thCV-12}
    Under Assumptions~\ref{hypGeneral}--\ref{hypSmallDataNS}, there exists $\eta>0$ such that if~\eqref{hyp-wellprep-2}--\eqref{hypEnergieModuleePetite-12} hold
    and if
    \[
        \ued\xrightharpoonup[\eps\to0]{} u^0\text{ in }w\text{-}\L^2(\T^3)\qquad\text{and}\qquad \rhoed^0\xrightharpoonup[\eps\to0]{}\rho^0\text{ in }w^*\text{-}\L^\infty(\T^3),
    \]
    then, for any $T>0$, $(\ueps)_{\eps>0}$ converges to $u$ in $\L^2((0,T)\times\T^3)$, $(\rhoeps)_{\eps>0}$ converges weakly-$*$ to $\rho$ in $\L^\infty((0,T)\times\T^3)$, where $(\rho,u)$ is a strong solution of
    \begin{equation}\label{limitRhoU-12}
        \left\{\begin{aligned}
            &\partial_t\rho+\div_x(\rho u)=0,\\
            &\rho|_{t=0}=\rho^0,\\
            &\partial_tu+(u\cdot\na_x)u-\Delta_xu+\na_xp=0,\\
            &\div_xu=0,\\
            &u|_{t=0}=u^0,
        \end{aligned}\right.
    \end{equation}
\end{theorem}

\begin{preuve}

    The  energy--dissipation estimate~\eqref{EstimationEnergieEq} shows that $(\ued)$ is uniformly bounded in $\L^2(0,T;\H^1(\T^3))$. Therefore, up to a subsequence that we shall not write explicitly, $(\ued)$ converges to some~$u$ in $w$-$\L^2(0,T;\H^1(\T^3))$. Furthermore, by the results of Section~\ref{sec-conclubootstrapLLF}, for $\eps$ small enough, all times are  strong existence times. Thanks to Lemma~\ref{lemmePara-12}, we can apply the Aubin-Lions lemma (see for example~\cite{boy-fab02}) and obtain the strong convergence of $(\ued)$ to $u$ in $\L^2((0,T)\times\T^3)$.
    
    Moreover, we can now use Corollary~\ref{RhoLinfini}, which ensures that $(\rhoed)$ is bounded in $\L^\infty((0,T)\times\T^3)$, so that it converges, up to a subsequence, to some~$\rho$ in $w^*$-$\L^\infty((0,T)\times\T^3)$. Therefore, $(\rhoed\ued)$ converges to $\rho u$ in $w$-$\L^2((0,T)\times\T^3)$.
    
    Furthermore, we can also apply Lemmas~\ref{BrinkmanDecompose-12}--\ref{lemmePara-12} to find that $(\Fed)$ converges to 0 in $\L^p((0,T)\times\T^3)$, which implies the same convergence in $\L^2((0,T)\times\T^3)$. As a consequence, $(\jed)$ converges to $\rho u$ in $w$-$\L^2((0,T)\times\T^3)$.
    
    We can finally take the limit $\eps\to0$ in the distributional sense in the Navier-Stokes equation
    \[
    \pa_t \ued + \ued \cdot \na_x \ued - \Delta \ued +\na_x \ued = \Fed
    \]
    and in the conservation law
    \[
    \partial_t\rhoegd+\div_x\jegd=0
    \]
    and obtain the derivation of~\eqref{limitRhoU-12}. 
    
    Let us note that thanks to the higher order energy estimate~\eqref{ENSH1} and another weak compactness argument, we  have $u \in \L^\infty(0,+\infty; \H^1(\T^3)) \cap \L^2 (0,+\infty; \H^2(\T^3))$.
    
        To conclude, we remark that the limit system has a unique such smooth solution by standard uniqueness results in the Fujita-Kato framework (see e.g. \cite[Theorem 3.3]{che-des-gal-gre}). Consequently the convergence of $(\ued)$ and $(\rhoed)$ holds without taking a subsequence.
\end{preuve}

If we assume a stronger form of convergence of the initial fluid velocity, then we have a stronger result for $(\ued)$ and $(\rhoed)$ as well. We can also state the convergence of $(\fed)$ and provide explicit rates in term of $\eps$. Recall that $W_1$ denotes the Wasserstein-1 distance (see Definition~\ref{def-W1}).

\begin{theorem}\label{thCVprecis-12}
    Under Assumptions~\ref{hypGeneral}--\ref{hypSmallDataNS}, there exists $\eps_0>0$ and $\eta>0$ such that for all $\eps \in (0,\eps_0)$,  if~\eqref{hyp-wellprep-2}--\eqref{hypEnergieModuleePetite-12} hold,
     there exists $\omega_p>0$ and $C>0$ such that, if $T>1$, for any $t\in[0,T]$,
    \[
        \norme{\ued(t)-u(t)}_{\L^2(\T^3)}\lesssim e^{CM^2T}M^{\omega_p}\left(\norme{\ued^0-u^0}_{\L^2(\T^3)}+\eps^{\frac{\kappa}{2}}\right),
    \]
    \[
        W_1(\rhoed(t),\rho(t))
        \lesssim W_1(\rhoed^0,\rho^0)
            +T^{1/2} e^{CM^2T}M^{\omega_p}\left(\norme{\ued^0-u^0}_{\L^2(\T^3)}+\eps^{\frac{\kappa}{2}}\right).
    \]
    and
    \begin{multline*}
        \int_0^TW_1(\fed(t),\rho(t)\otimes\delta_{v=\sigma u(t)})\dd t\\
        \lesssim T W_1(\rhoed^0,\rho^0)+ T^{\frac{3}{2}}e^{CM^2T}M^{\omega_p}\left(\norme{\ued^0-u^0}_{\L^2(\T^3)}+\eps^{\frac{\kappa}{2}}\right),
    \end{multline*}
    where $(\rho,u)$ is the unique solution of~\eqref{limitRhoU-12}.
\end{theorem}

\begin{remark}\label{rem-conv12-T}
    We only assume $T>1$ in order to simplify the time-dependence of the estimates above. In particular, the following result holds for any $T>0$ with the same convergence rate in term of $\eps$.
\end{remark}

As a straightforward consequence we deduce
\begin{corollary}
    Under Assumptions~\ref{hypGeneral}--\ref{hypSmallDataNS}, there exists $\eta>0$ such that if
    \begin{itemize}
        \item \eqref{hyp-wellprep-2}--\eqref{hypEnergieModuleePetite-12} hold,
        \item $(\ued^0)_{\eps>0}$ converges to $u^0$ in $\L^2(\T^3)$,
        \item $(\rhoed^0)_{\eps>0}$ weakly-$*$ converges to $\rho^0$ in $(\mathcal{P}_1(\T^3),W_1)$,
    \end{itemize}
    then, for any $T>0$,
    \begin{itemize}
        \item $(\ued)_{\eps>0}$ converges to $u$ in $\L^2((0,T)\times\T^3)$,
        \item $(\rhoed)_{\eps>0}$ converges to $\rho$ in $\L^\infty(0,T;\mathcal{P}_1(\T^3))$,
        \item $(\fed-\rho\otimes\delta_{v=\sigma u})_{\eps>0}$ weakly-$*$ converges to $0$ in $\L^1(0,T;\mathcal{P}_1(\T^3))$.
    \end{itemize}
\end{corollary}

Let us prove Theorem~\ref{thCVprecis-12}.
\begin{preuve}
    Let $(\rho,u)$ be the unique strong solution to the system \eqref{limitRhoU-12}. We set 
    $$\wed=\ued-u$$
which satisfies the equation
    \[
        \partial_t\wed+(u\cdot\na_x)\wed-\Delta_x\wed+\na_x(p_\eps-p)=\Fed-(\wed\cdot\na_x)\ued.
    \]
    We multiply this by $\wed$ and integrate over $(0,t)\times\T^3$ for any $t\in(0,T)$ to get
    \begin{multline*}
        \frac{1}{2}\norme{\wed(t)}_{\L^2(\T^3)}^2+\int_0^t\int_{\T^3}\wed(s,x)\cdot(u(s)\cdot\na_x)\wed(s,x)\dd s\dd x\\
        +\int_0^t\norme{\na_x \wed(s)}_{\L^2(\T^3)}\dd s\\
        =\cfrac{1}{2}\norme{\wed(0)}_{\L^2(\T^3)}^2+\int_0^t\int_{\T^3}\Fed(s,x)\cdot\wed(s,x)\dd s\dd x\\
        -\int_0^t\int_{\T^3}(\wed(s,x)\cdot\na_x)\ued(s,x)\cdot\wed(s,x)\dd s\dd x.
    \end{multline*}
    Since $\div_xu=0$, we have
    \[
        \int_0^t\int_{\T^3}\wed(s,x)\cdot(u(s,x)\cdot\na_x)\wed(s,x)=0.
    \]
    Furthermore, thanks to the Gagliardo-Nirenberg inequality (Theorem~\ref{th-Gagliardo-Nirenberg}), we can write
    \begin{align*}
        \left|\int_0^t\right.&\left.\int_{\T^3}(\wed(s,x)\cdot\na_x)\ued(s,x)\cdot\wed(s,x)\dd s\dd x\right|\\
        &\le\int_0^t\norme{\wed(s)}_{\L^4(\T^3)}^2\norme{\na_x\ued(s)}_{\L^2(\T^3)}\dd s\\
        &\lesssim\int_0^t\norme{\wed(s)}_{\L^2(\T^3)}^{\frac{1}{2}}\norme{\na_x\wed(s)}_{\L^2(\T^3)}^{\frac{3}{2}}\norme{\na_x\ued(s)}_{\L^2(\T^3)}\dd s\\
        &\qquad\qquad\qquad+\int_0^t\norme{\wed(s)}_{\L^2(\T^3)}^2\norme{\na_x\ued(s)}_{\L^2(\T^3)}\dd s.
    \end{align*}
    Since $T$ is a strong existence time, we can apply Lemma~\ref{EstimationNSsupGeneral} and Young's inequality to find that for any $a>0$,
    \begin{multline*}
        \left|\int_0^t\int_{\T^3}(\wed(s,x)\cdot\na_x)\ued(s,x)\cdot\wed(s,x)\dd s\dd x\right|
        \\
        \lesssim a^{\frac{4}{3}}\int_0^t\norme{\na_x\wed(s)}_{\L^2(\T^3)}^2\dd s
            +\frac{1}{a^4}\left(\Psi_{\eps,0}+\Psi_{\eps,0}^{\frac{1}{2}}\right)\int_0^t\norme{\wed(s)}_{\L^2(\T^3)}^2\dd s.
    \end{multline*}
    We can take $a>0$ small enough so that Grönwall's inequality and Lemma~\ref{EstimationNSsupGeneral} show that there exists a universal constant $C>0$ such that for any $t\in[0,T]$,
    \begin{equation*}
        \norme{\wed(t)}_{\L^2(\T^3)}^2\lesssim e^{CM^2T}\left(\norme{\wed(0)}_{\L^2(\T^3)}^2+\norme{\Fed}_{\L^2((0,T)\times\T^3)}^2\right),
    \end{equation*}
    hence
    by Lemmas~\ref{lemmaFL2-12} and~\ref{lemmePara2-12}, since $T$ is a strong existence time, 
    we obtain
    \begin{equation}\label{estimDiffU-12}
        \norme{\ued-u}_{\L^2((0,T)\times\T^3)}^2
        \lesssim e^{CM^2T}M^{\omega_p}\left(\norme{\ued^0-u(0)}_{\L^2(\T^3)}^2+\eps^{\kappa}\right),
    \end{equation}
    for some $\omega_p>0$.
    
    As we have seen previously, $(\rhoed)$ converges weakly-$*$ to $\rho$ in $\L^\infty((0,T)\times\T^3)$ such that
    \[
        \partial_t\rho+\div_x(\rho u)=0,
    \]
    and $\rho|_{t=0}=\rho^0$. Let us now compute a rate of convergence. As before, we study the equation satisfied by the difference $\chi_{\ed}=\rhoed-\rho$. We have
    \[
        \partial_t\chi_{\ed}+\div_x(\chi_{\ed}\ued)=-\div_xG_{\ed},
    \]
    where $G_{\ed}=\Fed+\rho(\ued-u)$. We consider the characteristic curves associated to this transport equation, defined by
    \[
       \frac{\dd}{\dd s} Y(s;t,x)=\ued(s,Y(s;t,x))\qquad\text{and}\qquad Y(t;t,x)=x.
    \]
    Since $\div_x\ued=0$, we have $\det\na_xY(s;t,x)=1$. Furthermore, thanks to Grönwall's inequality, since $T$ is a strong existence time,
    \[
        \norme{\na_xY(s;t,x)}_{\L^\infty(\T^3)}\le e^{\norme{\na_x\ued}_{\L^1(0,T;\L^\infty(\T^3))}}\lesssim 1.
    \]
    Therefore, for any $\psi\in\CCC^\infty(\T^3)$ such that $\norme{\na_x\psi}_{\L^\infty(\T^3)}\le1$, the method of characteristics and the change of variable $x'=Y(s;t,x)$ yield, for any $t\in[0,T]$,
    \begin{multline*}
        \int_{\T^3}\chi_{\ed}(t)\psi \dd x=\int_{\T^3}\chi_{\ed}(0)\psi \dd x-\int_0^t\int_{\T^3}(\div_xG_\eps)(s,Y(s;t,x))\psi(x)\dd x\dd s\\
        =\int_{\T^3}(\rhoed^0-\rho^0)\psi \dd x+\int_0^t\int_{\T^3}G_\eps(s,x)\na_x\left(\psi(Y(t;s,x))\right)\dd x\dd s.
    \end{multline*}
    Moreover,
    \[
        \norme{\na_x\left(\psi(Y(t;s,x))\right)}_{\L^2(\T^3)}\le\norme{\na_xY(t;s,x)\na_x\psi(Y(t;s,x))}_{\L^\infty(\T^3)}\lesssim1,
    \]
    so that, thanks to the previous estimate~\eqref{estimDiffU-12} for $\ued-u$, we have
    \begin{align*}
        \Big|\int_0^t\int_{\T^3}G_{\ed}&(s,x)\na_x\left(\psi(Y(t;s,x))\right)\dd x\dd s\Big|\\
        &\le\norme{G_{\ed}}_{\L^2((0,T)\times\T^3)}\norme{\na_x\left(\psi(Y(t;s,x))\right)}_{\L^2((0,T)\times\T^3)}\\
        &\lesssim T^{1/2} \left(\norme{\Fed}_{\L^2((0,T)\times\T^3)}+ M\norme{\ued-u}_{\L^2((0,T)\times\T^3)}\right)\\
        &\lesssim T^{1/2} e^{CM^2T}M^{\omega_p}\left(\norme{\ued^0-u^0}_{\L^2(\T^3)}+\eps^{\frac{\kappa}{2}}\right).
    \end{align*}
    for some $\omega_p>0$. Therefore, for every $t\in[0,T]$,
    \begin{multline}\label{estimDiffRho-12}
        W_1(\rhoed(t),\rho(t))
        \lesssim W_1(\rhoed^0,\rho^0)\\
            +T^{1/2} e^{CM^2T}M^{\omega_p}\left(\norme{\ued^0-u^0}_{\L^2(\T^3)}+\eps^{\frac{\kappa}{2}}\right).
    \end{multline}
    We can now prove the estimate on $(\fed)_{\eps>0}$. We write, for almost all $t\in(0,T)$,
    \begin{multline*}
        \fed(t)-\rho(t)\otimes\delta_{v=\sigma u(t)}
        =\left(\fed(t)-\rhoed(t)\otimes\delta_{v=\sigma\ued(t)}\right)\\
        +\left(\rhoed(t)\otimes\delta_{v=\sigma\ued(t)}-\rho(t)\otimes\delta_{v=\sigma u(t)}\right).
    \end{multline*}
    For $\psi\in\CCC^\infty_c(\T^3\times\R^3)$ such that $\norme{\na_{x,v}\psi}_{\L^\infty(\T^3\times\R^3)}\le1$, we have
    \begin{multline*}
        |\langle\rhoed(t)\otimes\delta_{v=\sigma\ued(t)}-\rho(t)\otimes\delta_{v=\sigma u(t)},\psi\rangle|\\
        \le\left|\int_{\T^3}(\rhoed-\rho)(t)\psi(x,\sigma\ued(t,x))\dd x\right|\\
            +\int_{\T^3}\rho(t,x)|\psi(x,\sigma\ued(t,x))-\psi(x,\sigma u)|\dd x.
    \end{multline*}
    On the one hand, for almost every $t\in(0,T)$,
    \begin{multline*}
        \left|\int_{\T^3}(\rhoed-\rho)(t)\psi(x,\sigma\ued(t,x))\dd x\right|\\
        \le W_1(\rhoed(t),\rho(t))\sup_{x\in\T^3}\left| \na_x\left(\psi(x,\sigma\ued(t,x))\right)\right|
    \end{multline*}
    and, for every $x\in\T^3$,
    \begin{align*}
        |\na_x\left(\psi(x,\sigma\ued(t,x))\right)|&\le\norme{\na_x\psi}_{\L^\infty(\T^3\times\R^3)}+\sigma\norme{\na_x\ued(t)}_{\L^\infty(\T^3)}\norme{\na_v\psi}_{\L^\infty(\T^3\times\R^3)} \\
        &\le 1 +\sigma\norme{\na_x\ued(t)}_{\L^\infty(\T^3)}.
    \end{align*}
    We use the fact that since $T$ is a strong existence time,
    \[
    \norme{\na_x\ued(t)}_{\L^1(0,T;\L^\infty(\T^3))} \lesssim 1,
    \]
  so that
    \begin{multline*}
        \int_0^T W_1(\rhoed(t),\rho(t))\sup_{x\in\T^3}\left| \na_x\left(\psi(x,\sigma\ued(t,x))\right)\right| \dd t\\
        \lesssim (T+\sigma)\left(\sup_{s\in(0,T)}W_1(\rhoed(s),\rho(s))\right)
        \lesssim T \left(\sup_{s\in(0,T)}W_1(\rhoed(s),\rho(s))\right).
    \end{multline*}
    On the other hand
    \begin{multline*}
      \int_{\T^3}\rho(t,x)|\psi(x,\sigma\ued(t,x))-\psi(x,\sigma u)|\dd x\\
        \le\sigma \norme{\rho}_{\L^\infty((0,T)\times\T^3)}\norme{\na_v\psi}_{\L^\infty(\T^3\times\R^3)}\norme{\ued(t)-u(t)}_{\L^2(\T^3)}.
    \end{multline*}
Hence, recalling that
    \begin{multline*}
      W_1(\rhoed(t)\otimes\delta_{v=\sigma\ued(t)},\rho(t)\otimes\delta_{v=\sigma u(t)}) \\
      = \sup_{\norme{\na_{x,v} \psi}=1} |\langle\rhoed(t)\otimes\delta_{v=\sigma\ued(t)}-\rho(t)\otimes\delta_{v=\sigma u(t)},\psi\rangle|,
    \end{multline*}
 we obtain, since $T$ is a strong existence time and thanks to the previous estimates~\eqref{estimDiffU-12}--\eqref{estimDiffRho-12},
    \begin{multline*}
        \int_0^TW_1(\rhoed(t)\otimes\delta_{v=\sigma\ued(t)},\rho(t)\otimes\delta_{v=\sigma u(t)})\dd t\\
        \lesssim T W_1(\rhoed^0,\rho^0)
            +T^{\frac{3}{2}}e^{CM^2T}M^{\omega_p}\left(\norme{\ued^0-u^0}_{\L^2(\T^3)}+\eps^{\frac{\kappa}{2}}\right),
    \end{multline*}
    for some $\omega_p>0$. 
    Furthermore, we also obtain
    \begin{align*}
        &\int_0^TW_1(\fed(t),\rhoed(t)\otimes\delta_{v=\sigma \ued(t)})\dd t\\
        &\qquad \qquad  \le\int_0^T\int_{\T^3\times\R^3}\fed(t,x,v)|v-\sigma\ued(t,x)|\dd x\dd v\dd t\\
         &\qquad \qquad\le \sigma T^{1-\frac{1}{p}}\left(\int_0^T\int_{\T^3}\left(\int_{\R^3}\fed(t,x,v)\left|\frac{v}{\sigma}-\ued(t,x)\right|\dd v\right)^p\dd x\dd t\right)^{\frac{1}{p}},
    \end{align*}
    and thanks to Lemmas \ref{BrinkmanDecompose-12}--\ref{lemmePara-12}, we get
    \[
        \int_0^TW_1(\fed(t),\rhoed(t)\otimes\delta_{v=\sigma\ued(t)})\dd t
        \lesssim \sigma T^{1-\frac{1}{p}}M^{\omega_p}\eps^{\frac{\kappa}{p}}.
    \]

    We can now conclude that
    \begin{multline*}
        \int_0^TW_1(\fed(t),\rho(t)\otimes\delta_{v=\sigma u(t)})\dd t\\
        \lesssim T W_1(\rhoed^0,\rho^0)+ T^{\frac{3}{2}}e^{CM^2T}M^{\omega_p}\left(\norme{\ued^0-u^0}_{\L^2(\T^3)}+\eps^{\frac{\kappa}{2}}\right).
    \end{multline*}
    All claimed estimates are finally proved.
\end{preuve}

Under an additional well-preparedness assumption, we can even prove pointwise convergence in time for the distribution function.

\begin{theorem}\label{corollaryCVf-12}
    Under Assumptions~\ref{hypGeneral}--\ref{hypSmallDataNS}, there exists $\eps_0>0$ and $\eta>0$ such that for all $\eps \in (0,\eps_0)$, if~\eqref{hyp-wellprep-2}--\eqref{hypEnergieModuleePetite-12} holds
    then, if $T>1$, for almost every $t\in(0,T)$,
    \begin{multline*}
        W_1(\fed(t),\rho(t)\otimes\delta_{v=\sigma u(t)})
        \lesssim\sigma e^{CM^2T}M^{\omega_p}\left(\norme{\ued^0-u^0}_{\L^2(\T^3)}+\eps^{\frac{\kappa}{2}}\right)\\
        +\sigma\norme{\left|\frac{v}{\sigma}-\ued^0\right|\fed}_{\L^1(\T^3\times\R^3)}
        +\sigma M^{\omega_p}\eps^{1-\frac{1}{p}}.
    \end{multline*}
\end{theorem}

As stated in Remark~\ref{rem-conv12-T}, we only work under the assumption $T>1$ in order to have a simpler time-dependence in the convergence rates. In particular, the following convergence results hold for any $T>0$ with the same rate in terms of the parameter $\eps$.

\begin{corollary}
    In the light particle regime, under Assumptions~\ref{hypGeneral}--\ref{hypSmallDataNS}, there exists $\eps_0>0$ and $\eta>0$, such that if
    \[
        \forall\eps\in(0,\eps_0),\qquad\norme{u_{\eps,1}^0-\langle u_{\eps,1}^0\rangle}_{\L^2(\T^3)}\le\eta,
    \]
    and if $(u_{\eps,1,1}^0)_{\eps>0}$ converges to $u^0$ in $\L^2(\T^3)$, and
    \[
        \int_{\T^3\times\R^3}|v-u_{\eps,1}^0(x)|f_{\eps,1}^0(x,v)\,\dd x\dd v\xrightarrow[\eps\to0]{}0,
    \]
    then $(f_{\eps,1,1})_{\eps>0}$ weakly-$*$ converges to $\rho\otimes\delta_{v=u}$ in $\L^\infty(0,T;\mathcal{P}_1(\T^3))$.
\end{corollary}

\begin{corollary}
    In the light and fast particle regime, under Assumption~\ref{hypGeneral}--\ref{hypSmallDataNS}, there exists $\eps_0>0$ and $\eta>0$, such that if there exists $\kappa\in(0,1)$ for which
    \[
        \forall\eps\in(0,\eps_0),\qquad\norme{|v|^p\fed^0}_{\L^1(\T^3\times\R^3)}\lesssim \sigma^p \eps^{\kappa-1},
    \]
    and
    \[
        \forall\eps\in(0,\eps_0),\qquad\norme{u_{\eps,\sigma}^0-\langle u_{\eps,\sigma}^0\rangle}_{\L^2(\T^3)}\le\eta
    \]
    then $(f_{\eps,\sigma})_{\eps>0}$ converges to $\rho\otimes\delta_{v=0}$ in $\L^\infty(0,T;\mathcal{P}_1(\T^3))$.
\end{corollary}

The proof of Theorem~\ref{corollaryCVf-12} is similar to the derivation of uniform estimates for the Brinkman force in Section~\ref{SectionBrinkman-12}.
\begin{preuve}
    For $\psi\in\CCC^\infty_c(\T^3\times\R^3)$ with $\norme{\na_{x,v}\psi}_{\L^\infty(\T^3\times\R^3)}\le1$, we have, thanks to the method of characteristics and the change of variable $w=\Gamma_{\ed}^{t,x}(v)$, for almost every $t\in(0,T)$,
    \begin{align*}
        &\langle\fed(t)-\rho(t)\otimes\delta_{v=\sigma u(t)},\psi\rangle \\
        &\qquad =\int_{\T^3\times\R^3}\fed(t,x,v)[\psi(x,v)-\psi(x,\sigma u(t,x))]\dd x\dd v\\
        &\qquad =\int_{\T^3\times\R^3}\fed^0\left(\tilde{X}_{\ed}^{t,x,w}(0),w\right)\\
          &\qquad \quad  \times\left[\psi\left(x,[\Gamma_{\ed}^{t,x}]^{-1}(w)\right)-\psi(x,\sigma u(t,x))\right]\left|\det\na_w[\Gamma_{\ed}^{t,x}]^{-1}(w)\right|e^{\frac{3t}{\eps}}\dd x\dd w,
    \end{align*}
    and Lemma \ref{changeVarV} yields
    \begin{multline*}
        |\langle\fed(t)-\rho(t)\otimes\delta_{v=\sigma u(t)},\psi\rangle|\\
        \lesssim\int_{\T^3\times\R^3}\fed^0\left(\tilde{X}_{\ed}^{t,x,w}(0),w\right)\left|[\Gamma_{\ed}^{t,x}]^{-1}(w)-\sigma u(t,x)\right|\dd x\dd w.
    \end{multline*}
    As in the proof of Lemma~\ref{BrinkmanDecompose-12} we integrate the characteristics equation and perform an integration by parts to get
    \begin{align*}
        [\Gamma_{\ed}^{t,x}]^{-1}(w)&=\sigma \ued(t,x)+e^{-\frac{t}{\eps}}\left(w-\sigma u^0\left(\tilde{X}_{\ed}^{t,x,w}(0)\right)\right)\\
        &-\sigma\int_0^te^{\frac{s-t}{\eps}}\partial_s\ued\left(s,\tilde{X}_{\ed}^{t,x,w}(s)\right)\\
        &-\sigma\int_0^te^{\frac{s-t}{\eps}}V_{\ed}\left(s;t,x,[\Gamma_{\ed}^{t,x}]^{-1}(w)\right)\cdot\na_x\ued\left(s,\tilde{X}_{\ed}^{t,x,w}(s)\right)\dd s.
    \end{align*}
    Then, we separate the integral into four parts:
    \begin{align*}
        &\frac{1}{\sigma}|\langle\fed(t)-\rho(t)\otimes\delta_{v=\sigma u(t)},\psi\rangle|\\
        &\qquad \lesssim\int_{\T^3\times\R^3}\fed^0\left(\tilde{X}_{\ed}^{t,x,w}(0),w\right)|\ued(t,x)-u(t,x)|\dd x\dd w\\
          &\qquad   +\left|\int_{\T^3\times\R^3}\fed^0\left(\tilde{X}_{\ed}^{t,x,w}(0),w\right)\left(\frac{w}{\sigma}-u^0\left(\tilde{X}_{\ed}^{t,x,w}(0)\right)\right)\dd x\dd w\right|\\
           &\qquad  +\int_{\T^3\times\R^3}\int_0^se^{\frac{s-t}{\eps}}\fed^0\left(\tilde{X}_{\ed}^{t,x,w}(0),w\right)|\partial_s\ued\left(s,\tilde{X}_{\ed}^{t,x,w}(s)\right)|\dd s\dd x\dd w\\
           &\qquad  +\int_{\T^3\times\R^3}\int_0^se^{\frac{s-t}{\eps}}\fed^0\left(\tilde{X}_{\ed}^{t,x,w}(0),w\right)\\
             &\qquad \qquad   \times\left|V\left(s;t,x,[\Gamma_{\ed}^{t,x}]^{-1}(w)\right)\cdot\na_x\ued\left(s,\tilde{X}_{\ed}^{t,x,w}(s)\right)\right|\dd s\dd x\dd w\\
        &\qquad \le I_1+I_2+I_3+I_4.
    \end{align*}
    Let us provide, for each of these integrals, a rate of convergence (to 0). Using the Cauchy-Schwarz inequality and the previously stated convergence of $(\ued)_{\eps>0}$, we get
    \begin{align*}
        |I_1|&\le\norme{\fed^0}_{\L^1(\R^3;\L^\infty(\T^3))}\norme{\ued(t)-u(t)}_{\L^2(\T^3)}\\
        &\lesssim e^{CM^2T}M^{\omega_p}\left(\norme{\ued^0-u^0}_{\L^2(\T^3)}+\eps^{\frac{\kappa}{2}}\right).
    \end{align*}
    Furthermore, thanks to Lemma \ref{changeVarX} and the corresponding change of variables,
    \[
        I_2\lesssim\int_{\T^3\times\R^3}\left|\frac{v}{\sigma}-\ued^0(x)\right|\fed^0(x,v)\,\dd x\dd v.
    \]
    We deal with $I_3$ and $I_4$ as in Section~\ref{SectionBrinkman-12}, but here the computations are simpler. Thanks to Hölder's inequality and Lemmas \ref{changeVarX} and \ref{lemmePara-12}, we have
    \begin{align*}
        I_3&\le\int_0^te^{\frac{s-t}{\eps}}\left(\int_{\T^3\times\R^3}\fed^0\left(\tilde{X}_{t,x,w}(0),w\right)\dd x\dd w\right)^{1-\frac{1}{p}}\\
            &\quad\times\left(\int_{\T^3\times\R^3}\fed^0\left(\tilde{X}_{t,x,w}(0),w\right)\left|\partial_s\ued\left(s,\tilde{X}_{t,x,w}(s)\right)\right|^p\dd x\dd w\right)^{\frac{1}{p}}\dd s\\
        &\le\norme{\fed^0}_{\L^1(\R^3;\L^\infty(\T^3))}\int_0^te^{\frac{s-t}{\eps}}\norme{\partial_s\ued(s)}_{\L^p(\T^3)}\dd t\\
        &\le\eps^{1-\frac{1}{p}}\norme{\fed^0}_{\L^1(\R^3;\L^\infty(\T^3))}\norme{\partial_s\ued}_{\L^p((0,T)\times\T^3)}
        \lesssim M^{\omega_p}\eps^{1-\frac{1}{p}},
    \end{align*}
    for some $\omega_p>0$. Finally, we also have
    \begin{align*}
        I_4&\le\int_0^t\int_{\T^3\times\R^3}e^{-\frac{t}{\eps}}\fed^0\left(\tilde{X}_{t,x,w}(0),w\right)|w|\left|\na_x\ued\left(s,\tilde{X}_{t,x,w}(s)\right)\right|\dd x\dd w\dd s\\
            &\quad+\frac{\sigma}{\eps}\int_0^t\int_{\T^3\times\R^3}e^{\frac{s-t}{\eps}}\fed^0\left(\tilde{X}_{t,x,w}(0),w\right)\int_0^se^{\frac{\tau-s}{\eps}}\left|\ued\left(\tau,\tilde{X}_{t,x,w}(\tau)\right)\right|\\
            &\qquad\qquad \times\left|\na_x\ued\left(s,\tilde{X}_{t,x,w}(s)\right)\right|\dd\tau\dd x\dd w\dd s\\
        &\le I_5+I_6.
    \end{align*}
    We conclude thanks to Hölder's inequality and Lemmas~\ref{estimNabla} and~\ref{lemmePara-12} that
    \begin{align*}
        I_5&\le\norme{\fed^0|v|^{\frac{p}{p-1}}}_{\L^1(\R^3;\L^\infty(\T^3))}^{1-\frac{1}{p}}\norme{\fed^0}_{\L^1(\R^3;\L^\infty(\T^3))}^{\frac{1}{p}}\int_0^te^{-\frac{t}{\eps}}\norme{\na_x\ued(s)}_{\L^p(\T^3)}\dd s\\
        &\lesssim M^{\omega_p}\eps^{1-\frac{\alpha_p}{p}},
    \end{align*}
    and
    \begin{align*}
        I_6&\le\frac{\sigma}{\eps}\int_0^t\int_0^se^{\frac{s-t}{\eps}}e^{\frac{\tau-s}{\eps}}\left(\int_{\T^3\times\R^3}\fed^0(\tilde{X}_{t,x,w}(0),w)|\ued(\tau,\tilde{X}_{t,x,w}(\tau))|^p\dd x\dd w\right)^{\frac{1}{p}}\\
            &\quad\times\left(\int_{\T^3\times\R^3}\fed^0(\tilde{X}_{t,x,w}(0),w)|\na_x\ued(s,\tilde{X}_{t,x,w}(s))|^{\frac{p}{p-1}}\dd x\dd w\right)^{1-\frac{1}{p}}\dd\tau\dd s\\
        &\le\frac{\sigma}{\eps}\norme{\fed^0}_{\L^1(\R^3;\L^\infty(\T^3))}\int_0^t\int_0^se^{\frac{s-t}{\eps}}e^{\frac{\tau-s}{\eps}}\norme{\ued(\tau)}_{\L^p(\T^3)}\norme{\na_x\ued(s)}_{\L^\infty(\T^3)}\dd\tau\dd s\\
        &\lesssim\frac{\sigma}{\eps}\norme{\fed^0}_{\L^1(\R^3;\L^\infty(\T^3))}\EEE_{\ed}(0)^{\frac{1-\beta_p}{2}}\Psi_{\eps,0}^{\frac{1}{2}}\\
          &\quad  \times\int_0^te^{\frac{\tau-s}{\eps}}\norme{\Delta_x\ued(s)}_{\L^p(\T^3)}^{\beta_p}\left(\int_\tau^te^{\frac{s-\tau}{\eps}}\dd s\right)\dd\tau\\
        &\lesssim M^{\omega_p}\sigma\eps^{1-\frac{\beta_p}{p}}.
    \end{align*}
    Therefore, there exists $\omega_p>0$ such that
    \begin{align*}
        |\langle\fed(t)-\rho(t)\otimes\delta_{v=\sigma u(t)},\psi\rangle|
        &\lesssim\sigma e^{CM^2T}M^{\omega_p}\left(\norme{\ued^0-u^0}_{\L^2(\T^3)}+\eps^{\frac{\kappa}{2}}\right)\\
        &+\sigma\norme{\left|\frac{v}{\sigma}-\ued^0\right|\fed}_{\L^1(\T^3\times\R^3)}
        +\sigma M^{\omega_p}\eps^{1-\frac{1}{p}},
    \end{align*}
    hence the theorem.
\end{preuve}

In Section~\ref{SectionRetour-12}, we shall apply the tools developed to deal with the fine particle regime to the \emph{light} and \emph{light and fast} particle regimes and obtain alternative proofs of the above results under slightly different sets of assumptions.


\subsection{Application to long time behavior}
\label{sec-long}
As a side consequence of the analysis, we obtain the description of the long time behavior of solutions to~\eqref{VNS-general} in the \emph{light} and \emph{light and fast} particle regime. What we specifically have proved are variants of the main result of 
\cite{hank-mou-moy} with assumptions which are somewhat tailored to the particular regime we consider.

Loosely speaking, \cite[Theorem 2.1]{hank-mou-moy} asserts, forgetting for a while about the small parameter $\eps$, that if the initial data $(u^0,f^0)$ is close to equilibrium,
then any solution remains close and the dynamics ultimately resembles that of a Dirac mass in velocity.

We observe that as corollaries of the previous analysis, we get  similar results but with assumptions bearing solely on the velocity field in the \emph{light} and \emph{light and fast} particle regime for $\alpha \in(0,1/3]$:

\begin{corollary}
\label{coro-longtime-LLF}
    In the \emph{light} and  \emph{light and fast} particle regime, for $\alpha\in [0,1/3]$, under Assumptions~\ref{hypGeneral} and~\ref{hypSmallDataNS}, for $\eps$ small enough, there is $\eta>0$ such that, if 
    $$
        \norme{u_{\eps,\eps^\alpha}^0-\langle u_{\eps,\eps^\alpha}^0\rangle}_{\L^2(\T^3)}\le\eta,
    $$
    then there exist $C, \lambda>0$ independent of $\eps$ such that, for all $t \geq 0$, 
    $$
        \EEE_{\eps,\eps^\alpha}(t) \leq C  \EEE_{\eps,\eps^\alpha}(0) e^{-\lambda t}. 
    $$
\end{corollary}

The fact that the decay of the modulated energy entails a monokinetic behavior as $t\to +\infty$  is for instance proved in \cite[Theorem 1.1]{hank-mou-moy}

\begin{preuve}
This is a consequence of Lemmas~\ref{ModulatedEnergyDissipation} and~\ref{lemmeDE}, of the bootstrap analysis led in Section~\ref{sec-conclubootstrapLLF}, of Remark~\ref{RemarkEnergieModuleePetite} and the observation that~\eqref{hyp-wellprep-2} is restrictive only for $\alpha>1/3$.
\end{preuve}


\subsection{Convergence in the general case}
\label{sec-smallT}

In this paragraph we aim to prove similar convergence results as the previous theorems 
under the sole Assumption~\ref{hypGeneral}.

In turn, the time horizon under which we are able to prove convergence results shall be constrained.


Recall that we consider, under Assumption~\ref{hypGeneral}, for every $\eps>0$,
\[
    T^*_\eps=\sup\{T>0,\,T\text{ is a strong existence time}\}.
\]
Our goal is to find a time $T_M$, independent of $\eps$, such that $T_M<T^*_\eps$. Let us first note that in view of the proof of Lemma~\ref{EstimationNSsupGeneral}, under Assumption~\ref{hypGeneral}, there exists a time $T_0$, independent of $\eps$, such that
\[
    \int_0^{T_0}\norme{e^{t\Delta}\ued}_{\dot\H^1(\T^3)}^4\dd t\le\frac{C^*}{4}.  
\]
Thus, it suffices to show that there exists $T_M\le\min(T^*_\eps,T_0)$ independent of $\eps$, such that
\[
    \norme{\na\ued}_{\L^p((0,T)\times\T^3)}\le\frac{1}{40}\qquad\text{and}\qquad\norme{\Fed}_{\L^2((0,T)\times\T^3)}\le\frac{C^*}{4}.
\]
Let $T\in(0,\min(T^*_\eps,T_M))$. 
We can apply Lemmas~\ref{lemmaFL2-12} and~\ref{lemmePara2-12} which ensure that for $\eps$ small enough,
\[
    \norme{\Fed}_{\L^2((0,T)\times\T^3)}\le\frac{C^*}{4}.
\]
To conclude, by a variant of the proof of Corollary~\ref{coro-nau-12}, we get
\[
 \norme{\na_x\ued}_{\L^1(0,T;\L^\infty(\T^3))} \lesssim T^{\omega_p}M^{\omega_p'}(1+\eps^{\frac{\kappa}{p}}),
 \]
 so that, imposing  $T$ small enough instead of relying on the smallness of the initial modulated energy $\mathscr{E}_\ed(0)$,
 this yields
 \[
  \norme{\na_x\ued}_{\L^1(0,T;\L^\infty(\T^3))} \le\frac{1}{40}.
  \]

Therefore, we have proved that there exists a time $T_M>0$, independent of $\eps$, that is a strong existence time.

The convergence results in Theorems~\ref{thCV-12}--\ref{thCVprecis-12} and Corollary~\ref{corollaryCVf-12} still hold when replacing the smallness assumptions on the initial condition by the smallness of the time horizon. Their proof are indeed almost identical. 

\section{The fine particle regime}\label{SectionFineParticle}

We study in this section the fine particle regime, which corresponds to $(\gamma,\sigma)=(\eps,1)$ in~\eqref{VNS-general}. For the sake of readability, we shall no longer write the parameters $(\gamma,\sigma)$: for instance, $u_{\eps,\eps,1}$ will be referred to as $\ueps$. 
For the reader's convenience, let us write down the system corresponding to this regime:
\begin{equation*}
    \left\{
\begin{aligned}
    &\partial_t \ueps+ (\ueps\cdot\na_{x})\ueps-\Delta_{x}\ueps+\na_{x}p= \frac{1}{\eps}\left(j_{\feps}-\rho_{\feps}\ueps\right),\\
  &\div_{x}\ueps=0,\\
   &\partial_t\feps+v\cdot\na_{x}\feps+ \frac{1}{\eps} \div_{v}\left[\feps( \ueps-v)\right]=0, \\
    &  \rho_{\feps}(t,x)=\int_{\R^3}\feps(t,x,v)\ddv,  \quad
  j_{\feps}(t,x)=\int_{\R^3}v\feps(t,x,v)\ddv.
\end{aligned}
\right.
\end{equation*}
As explained in the introduction, the fine particle regime is more singular than the \emph{light} and \emph{light and fast} ones. Yet, in the case of well-prepared initial data, the same desingularization (as in Section~\ref{SectionBrinkman-12}) of the Brinkman force
\[
    F_\eps=\frac{1}{\eps}\left(j_{\feps}-\rho_{\feps}\ueps\right)
\]
could be applied and would allow to prove the expected convergence. Nevertheless, we choose to present a different strategy to desingularize $F_\eps$ which yields stronger results.
For this purpose, we introduce and study what \cite{hank} calls \emph{higher} fluid-kinetic dissipation. This allows to prove qualitative convergence results for $\feps$ and $(\rhoeps,\ueps)$. However this approach is not sufficient to provide convergence rates. To that purpose, we shall rely on a relative entropy method.
Let us emphasize the fact that this method \emph{by itself does not} allow to conclude, but rather has to be used in conjunction with the previous analysis.

\begin{remark}
    The higher fluid-kinetic dissipation allows, in~\cite{hank}, to extend the study of the long-time-behavior of the Vlasov-Navier-Stokes system from $\R_+\times\T^3\times\R^3$ to the domain $\R_+\times\R^3\times\R^3$. It might also be the right approach to study high friction limits in the domain $\R_+\times\R^3\times\R^3$, which we do not consider here.
\end{remark}

This section is organized as follows:
\begin{itemize}
\item In Section~\ref{SubsectionHigherOrder}, we introduce the relevant higher dissipation functionals and prove a key identity, leading to the desingularization of the Brinkman force and eventually to uniform $\L^p$ (in time and space) bounds.
This process requires to consider mildly well-prepared initial data as defined in the statement of Theorem~\ref{thm3}.
\item Thanks to these crucial bounds, for well-prepared initial data, in Section~\ref{sec-bootstrap-fine}, we conclude the bootstrap argument, obtaining that all positive times are \emph{strong existence times}. Whereas the general strategy is the same as for the \emph{light} and \emph{light and fast} particle regime (recall Section~\ref{sec-conclubootstrapLLF}), a special care is required here to handle the $\L^2$ norm of the Brinkman force (recall Remark~\ref{RemarkBorneUnifF}).
 \item Section~\ref{sec-first-fine} concludes the convergence proof in the well-prepared case.
 \item In Section~\ref{sec-relat-fine}, we introduce a \emph{relative entropy} method, allowing to prove  quantitative convergence results in the well-prepared case of Theorem~\ref{thm3}, when one assumes that the solutions of the limiting system are smooth enough. This analysis  can be performed only because we have already proved that all positive times are strong existence times.
 \item Finally Section~\ref{sec-smallT-fine} explains how to get a small time convergence result in the mildly well-prepared case.
\end{itemize}


\subsection{Uniform estimates on the Brinkman force in \texorpdfstring{$\L^p$}{L\string^p}}\label{SubsectionHigherOrder}

 To conclude the bootstrap argument initiated in Section~\ref{SubsectionInitialization}, we need to control the Brinkman force in various  $\L^r$ spaces. We will rely on the exponential decay of the modulated energy on $[0,T]$, when $T$ is a strong existence time (recall Lemma~\ref{lemmeDE}). Notice, in this regime, the factor $\eps^{-1}$ in the expression of $F_\eps$ which is therefore more singular than in Section~\ref{SectionLandLFParticle}. To overcome this difficulty, we follow \cite{hank} and introduce the \emph{higher} fluid-kinetic dissipation. This section is an adaptation of \cite[Section 4]{hank} which tracks down the dependence with respect to $\eps$ and which is detailed for the sake of completeness.


\begin{definition}
    Let $r\ge2$. The higher fluid-kinetic dissipation (of order $r$) is defined, for almost every $t\ge0$, by
    \[
        \D_{\eps}^{(r)}(t)=\int_{\T^3\times\R^3}\feps\frac{|v-\ueps(t,x)|^r}{\eps^r}\dd x\dd v.
    \]
    \index{$\D_{\eps}^{(r)}(t)$: higher dissipation of order $r$}
\end{definition}

In the rest of this section, we work under Assumption~\ref{hypGeneral}, and consider a strong existence time $T>0$.

The following lemma ensures that a control  of $\D_\eps^{(r)}$ implies estimates for~$F_\eps$ in $\L^r$.

\begin{lemma}\label{LienFD-3}
    For $r\ge2$, we have
    \[
        \norme{F_\eps}_{\L^r((0,T)\times\T^3)}\le\norme{\feps^0}_{\L^1(\R^3;\L^\infty(\T^3))}^{1-\frac{1}{r}}\left(\int_0^T\D_{\eps}^{(r)}(t)\,\dd t\right)^{\frac{1}{r}}.
    \]
\end{lemma}

\begin{preuve}
    This is a direct consequence of Hölder's inequality and Corollary \ref{RhoLinfini}, since $T$ is a strong existence time:
    \begin{align*}
        \norme{F_\eps(t)}_{\L^p((0,T)\times\T^3)}^r
        &\le\int_0^T\int_{\T^3}\rhoeps(t,x)^{r-1}\left(\int_{\R^3}\frac{|v-\ueps|^r}{\eps^r}\feps(t,x,v)\dd v\right)\dd x\dd t\\
        &\lesssim\norme{\feps^0}_{\L^1(\R^3;\L^\infty(\T^3))}^{r-1}\int_0^T\D_{\eps}^{(r)}(t)\,\dd t,
    \end{align*}
        which corresponds to the desired estimate.
\end{preuve}

The key to desingularize $\D_{\eps}^{(r)}$is the following identity, obtained by using the method of characteristics and integration by parts.

\begin{lemma}
\label{fine-key}
    For all $r\ge2$, the following identity holds.
    \begin{multline*}
        \int_0^T\D_\eps^{(r)}(t)\dd t=-\frac{1}{r}\left[\int_{\T^3\times\R^3}\feps(t,x,v)\frac{|v-\ueps(t,x)|^r}{\eps^{r-1}}\dd x\dd v\right]_0^T\\
        -\int_0^T\int_{\T^3\times\R^3}\feps(t,x,v)(\partial_t\ueps+(\na_x\ueps)v)\cdot\frac{(v-\ueps(t,x))|v-\ueps(t,x)|^{r-2}}{\eps^{r-1}}\dd x\dd v\dd t.
    \end{multline*}
\end{lemma}

\begin{preuve}
    By the method of characteristics and a change of variables,
    \begin{align*}
        \D_\eps^{(r)}(t)&=\eps^{-r}e^{\frac{3t}{\eps}}\int_{\T^3\times\R^3}\feps^0(X_\eps(0;t,x,v),V_\eps(0;t,x,v))|v-\ueps(t,x)|^r\dd x\dd v\\
        &=\eps^{-r}\int_{\T^3\times\R^3}\feps^0(x,v)|V_\eps(t;0,x,v)-\ueps(t,X_\eps(t,0,x,v))|^r\dd x\dd v.
    \end{align*}
    Note that
    \begin{align*}
        &\frac{\dd}{\dd t}|V_\eps(t;0,x,v)-\ueps(t,X_\eps(t;0,x,v))|^r\\
        &\quad= r\frac{\dd}{\dd t}\left(V_\eps(t;0,x,v)-\ueps(t,X_\eps(t;0,x,v))\right)\cdot\left(V_\eps(t;0,x,v)-\ueps(t,X_\eps(t;0,x,v))\right)\\
        &\qquad \qquad    \times |V_\eps(t;0,x,v)-\ueps(t,X_\eps(t;0,x,v))|^{r-2}\\
        &\quad=r\left(\frac{\ueps(t,X_\eps(t;0,x,v))-V_\eps(t;0,x,v)}{\eps}-\frac{\dd}{\dd t}\left(\ueps(t,X_\eps(t;0,x,v))\right)\right)\\
        &\qquad \quad \cdot\left(V_\eps(t;0,x,v)-\ueps(t,X_\eps(t;0,x,v))\right)|V_\eps(t;0,x,v)-\ueps(t,X_\eps(t;0,x,v))|^{r-2}.
    \end{align*}
    This yields
    \begin{multline*}
        \frac{1}{\eps}|V_\eps(t;0,x,v)-\ueps(t,X_\eps(t;0,x,v))|^r\\
        =-\frac{1}{r}\frac{\dd}{\dd t}|V_\eps(t;0,x,v)-\ueps(t,X_\eps(t;0,x,v))|^r\\
        -\frac{\dd}{\dd t}\left(\ueps(t,X_\eps(t;0,x,v))\right)\cdot\left(V_\eps(t;0,x,v)-\ueps(t,X_\eps(t;0,x,v))\right)\\
        \times|V_\eps(t;0,x,v)-\ueps(t,X_\eps(t;0,x,v))|^{r-2}
    \end{multline*}
    and the result follows after integrating by parts.
\end{preuve}

Once in this form, the control of $\D_\eps^{(r)}$ is a consequence of the following estimates, which are straightforward applications of Hölder's inequality.

\begin{lemma}
    For all $r\ge2$, we have the following estimates
    \begin{align*}
        &\left|\int_0^T\int_{\T^3\times\R^3}\feps(t,x,v)\partial_t\ueps\cdot\frac{(v-\ueps(t,x))|v-\ueps(t,x)|^{r-2}}{\eps^{r-1}}\dd x\dd v\dd t\right|\\
        &\qquad \qquad \qquad \qquad \lesssim\norme{\feps^0}_{\L^1(\R^3;\L^\infty(\T^3))}^{\frac{1}{r}}\norme{\partial_t\ueps}_{\L^r((0,T)\times\T^3)}\left(\int_0^T\D_\eps^{(r)}(t)\dd t\right)^{1-\frac{1}{r}}, \\
        &\left|\int_0^T\int_{\T^3\times\R^3}\feps(t,x,v)(\na_x\ueps)v)\cdot\frac{(v-\ueps(t,x))|v-\ueps(t,x)|^{r-2}}{\eps^{r-1}}\dd x\dd v\dd t\right|\\
        &\qquad \qquad \qquad \qquad\lesssim\norme{|\na_x\ueps|m_r^{\frac{1}{r}}}_{\L^r((0,T)\times\T^3)}\left(\int_0^T\D_\eps^{(r)}(t)\dd t\right)^{1-\frac{1}{r}}.
    \end{align*}
\end{lemma}

A combination of the two results leads, thanks to Young's inequality, to the following lemma.
\begin{lemma}\label{lemmeExpressionD}
    The higher fluid-kinetic dissipation of order $r\ge2$ satisfies, for almost every $t\in(0,T)$,
    \begin{multline*}
        \int_0^t\D_\eps^{(r)}(s)\dd s+\int_{\T^3\times\R^3}\feps(t,x,v)\frac{|v-\ueps(t,x)|^r}{\eps^{r-1}}\dd x\dd v\\
        \lesssim\norme{\feps^0}_{\L^1(\R^3;\L^\infty(\T^3))}\norme{\partial_t\ueps}_{\L^r((0,T)\times\T^3)}^r\\
        +\norme{|\na_x\ueps|m_{r,\eps}^{\frac{1}{r}}}_{\L^r((0,T)\times\T^3)}^r
        +\int_{\T^3\times\R^3}\feps^0\frac{|v-\ueps^0|^r}{\eps^{r-1}}\dd x\dd v,
    \end{multline*}
    where
    $$
    m_{r,\eps} := \int_{\R^3} \feps |v|^r \dd v.
    $$
\end{lemma}

We now apply this estimate with $r=p$, the regularity index appearing in Assumption~\ref{hypGeneral}. The first term will be absorbed in the parabolic estimate, under a smallness condition on $\feps^0$ and the last term will lead us to assume well-preparedness of the initial data. We still have to deal with the middle one.

\begin{lemma}\label{lemmeMixed}
Let $T>0$ be a strong existence time.
    The following estimate holds
    \begin{align*}
        \norme{|\na_x\ueps|m_{p,\eps}^{\frac{1}{p}}}&_{\L^p((0,T)\times\T^3)}^p
        \lesssim \eps^{1-\alpha_p}\norme{\feps^0|v|^p}_{\L^1(\R^3;\L^\infty(\R^3))}^{\frac{1}{1-\alpha_p}}\EEE_\eps(0)^{\frac{p}{2}}\\
        &+\norme{\feps^0}_{\L^1(\R^3;\L^\infty(\T^3))}\EEE_\eps(0)^{\frac{(1-\beta_p)p}{2}}\Psi_{\eps,0}^{\frac{p}{2(1-\beta_p)}}\\
        &+\left(\eps^{1-\alpha_p}+\norme{\feps^0}_{\L^1(\R^3;\L^\infty(\T^3))}\EEE_\eps(0)^{\frac{(1-\beta_p)p}{2}}\right)\norme{\Delta_{x}\ueps}_{\L^p((0,T)\times\T^3)}^{p}.
    \end{align*}
\end{lemma}

\begin{preuve}
    The method of characteristics yields
    \begin{multline*}
        \norme{|\na_x\ueps|m_{p,\eps}^{\frac{1}{p}}}_{\L^p((0,T)\times\T^3)}^p\\
        =\int_0^T\int_{\T^3\times\R^3}|\na_x\ueps(t,x)|^pe^{\frac{3t}{\eps}}\feps^0(X_\eps(0;t,x,v),V_\eps(0;t,x,v))||v|^p\dd x\dd v\dd t
    \end{multline*}
    so that by the change of variables in velocity of Lemma \ref{changeVarV}, we get
    \begin{align*}
        &\norme{|\na_x\ueps|m_{p,\eps}^{\frac{1}{p}}}_{\L^p((0,T)\times\T^3)}^p\\
        &\qquad=\int_0^T\int_{\T^3\times\R^3}|\nabla_x\ueps(t,x)|^p e^{\frac{3t}{\eps}}\feps^0(\tilde{X}_{t,x,w}(0),w)|[\Gamma_{\eps}^{t,x}]^{-1}(w)|^p\\
        &\qquad \qquad \qquad \qquad \times|\det\nabla_{w}\Gamma_{\eps}^{t,x}(w)|^{-1}\dd x\dd w\dd t\\
        &\qquad \lesssim\int_0^T\int_{\T^3\times\R^3}|\na_x\ueps(t,x)|^p\feps^0(\tilde{X}^{t,x,w}_\eps(0),w)|[\Gamma_{\eps}^{t,x}]^{-1}(w)|^p\dd x\dd w\dd t.
    \end{align*}
    By integration of the characteristics equation on velocity~\eqref{characV}, we have
    \[
        [\Gamma_{\eps}^{t,x}]^{-1}(w)=e^{-\frac{t}{\eps}}w+\frac{1}{\eps}\int_0^te^{\frac{s-t}{\eps}}\ueps(s,\tilde{X}^{t,x,w}_\eps(s))\dd s.
    \]
    Therefore,
    \begin{align*}
        &\norme{|\na_x\ueps|m_{p,\eps}^{\frac{1}{p}}}_{\L^p((0,T)\times\T^3)}^p\\
        &\qquad \lesssim\int_0^T\int_{\T^3\times\R^3}|\na_x\ueps(t,x)|^p\feps^0(\tilde{X}^{t,x,w}_\eps(0),w)e^{-\frac{pt}{\eps}}|w|^p\dd x\dd v\dd t\\
        &\qquad +\int_0^T\int_{\T^3\times\R^3}|\na_{x}\ueps(t,x)|^p\feps^0(\tilde{X}^{t,x,w}_\eps(0),w)\\
        &\qquad \qquad \qquad \qquad \times\left(\frac{1}{\eps}\int_0^t e^{\frac{s-t}{\eps}}|\ueps(s,\tilde{X}^{t,x,w}_\eps(s))|\dd s\right)^p\dd x\dd w\dd t \\
        &\qquad=: I_1+I_2.
    \end{align*}
    On the one hand, thanks to Lemma~\ref{estimNabla} and Hölder's and Young's inequalities,
    \begin{align*}
        I_1&\le\norme{\feps^0|v|^p}_{\L^1(\R^3;\L^\infty(\T^3))}\int_0^Te^{-\frac{pt}{\eps}}\norme{\na_{x}\ueps(t)}_{\L^p(\T^3)}^p\dd t\\
        &\lesssim\norme{\feps^0|v|^p}_{\L^1(\R^3;\L^\infty(\T^3))}\EEE_\eps(0)^{\frac{(1-\alpha_p)p}{2}}\int_0^Te^{-\frac{pt}{\eps}}\norme{\Delta_{x}\ueps(t)}_{\L^p(\T^3)}^{\alpha_pp}\dd t\\
        &\lesssim\eps^{1-\alpha_p}\norme{\feps^0|v|^p}_{\L^1(\R^3;\L^\infty(\R^3))}^{\frac{1}{1-\alpha_p}}\EEE_\eps(0)^{\frac{p}{2}}+\eps^{1-\alpha_p}\norme{\Delta_{x}\ueps}_{\L^p((0,T)\times\T^3)}^p.
    \end{align*}

    On the other hand, applying Jensen's inequality and Lemma~\ref{changeVarX},
    \begin{align*}
        I_2&\lesssim\frac{1}{\eps}\int_0^T\int_{\T^3\times\R^3}\int_0^t|\na_x\ueps(t,x)|^p|\ueps(s,\tilde{X}^{t,x,w}_\eps(s))|^p\\
        &\qquad \qquad \times\feps^0(\tilde{X}^{t,x,w}_\eps(0),w)e^{\frac{s-t}{\eps}}\dd s\dd x\dd w\dd t\\
        &\lesssim\frac{1}{\eps}\norme{\feps^0}_{\L^1(\R^3;\L^\infty(\T^3))}\int_0^T\int_0^t\norme{\na_{x}\ueps(t)}_{\L^{\infty}(\T^3)}^p\norme{\ueps(s)}_{\L^p(\T^3)}^pe^{\frac{s-t}{\eps}}\dd s\dd t.
    \end{align*}
    Recall that Lemmas~\ref{EstimNormLp} and~\ref{estimNabla} ensure that for every $s\in(0,T)$,
    \[
        \norme{\ueps(s)}_{\L^p(\T^3)}\lesssim\Psi_{\eps,0}^{\frac{1}{2}},
    \]
    and every $t\in(0,T)$,
    \[
        \norme{\nabla_{x}\ueps(t)}_{\L^\infty(\T^3)}
        \lesssim\EEE_\eps(0)^{\frac{1-\beta_p}{2}}e^{-\frac{(1-\beta_p)\lambda_\eps}{2}t}\norme{\Delta_{x}\ueps(t)}_{\L^p(\T^3)}^{\beta_p}.
    \]
    Therefore, Hölder's and Young's inequalities yield
    \begin{align*}
        I_2&\lesssim\norme{\feps^0}_{\L^1(\R^3;\L^\infty(\T^3))}\Psi_{\eps,0}^{\frac{p}{2}}\EEE_\eps(0)^{\frac{(1-\beta_p)p}{2}}\norme{\Delta_{x}\ueps}_{\L^p((0,T)\times\T^3)}^{\beta_pp}\\
        &\lesssim\norme{\feps^0}_{\L^1(\R^3;\L^\infty(\T^3))}\EEE_\eps(0)^{\frac{(1-\beta_p)p}{2}}\Psi_{\eps,0}^{\frac{p}{2(1-\beta_p)}}\\
        &\quad+\norme{\feps^0}_{\L^1(\R^3\L^\infty(\T^3))}\EEE_\eps(0)^{\frac{(1-\beta_p)p}{2}}\norme{\Delta_{x}\ueps}_{\L^p((0,T)\times\T^3)}^{p},
    \end{align*}
    which concludes the proof.
\end{preuve}

In conclusion, we have obtained the following estimate.

\begin{lemma}\label{EstimateHighOrderD}
    Under Assumptions~\ref{hypGeneral}, if $T>0$ is a strong existence time, then
    \begin{align*}
        \int_0^T\D_\eps^{(p)}(t)&\dd t + \int_{\T^3\times\R^3}\feps(t,x,v)\frac{|v-\ueps(t,x)|^p}{\eps^{p-1}}\dd x\dd v\\
        &\lesssim
            \norme{\feps^0}_{\L^1(\R^3;\L^\infty(\T^3))}\norme{\partial_t\ueps}_{\L^p((0,T)\times\T^3)}^p\\
        &\quad+   
            \left(\eps^{1-\alpha_p}+\norme{\feps^0}_{\L^1(\R^3;\L^\infty(\R^3))}\EEE_\eps(0)^{\frac{(1-\beta_p)p}{2}}\right)\norme{\Delta_x\ueps}_{\L^p(0,T)\times\T^3)}^{p}\\
        &\quad+
            \eps^{1-\alpha_p}\norme{\feps^0|v|^p}_{\L^1(\R^3;\L^\infty(\T^3))}^{\frac{1}{1-\alpha_p}}\EEE_\eps(0)^{\frac{p}{2}}\\
        &\quad+
            \norme{\feps^0}_{\L^1(\R^3;\L^\infty(\T^3))}\EEE_\eps(0)^{\frac{(1-\beta_p)p}{2}}\Psi_{\eps,0}^{\frac{p}{2(1-\beta_p)}}\\
        &\quad+
            \int_{\T^3\times\R^3}\feps^0(x,v)\frac{|v-\ueps^0(x)|^p}{\eps^{p-1}}\dd x\dd v.
    \end{align*}
\end{lemma}

\subsection{Conclusion of the bootstrap argument in the well-prepared case}
\label{sec-bootstrap-fine}

This  section is the analogue of Section~\ref{sec-conclubootstrapLLF} for the \emph{light} and \emph{light and fast} particle regimes. Recall that, as is now usual, we consider, under Assumption~\ref{hypGeneral}, for every $\eps>0$,
\[
    T^*_\eps=\sup\{T>0,\,T\text{ is a strong existence time}\}.
\]

The aim of this Subsection is to prove that under Assumption~\ref{hypSmallDataNS} and additional well-preparedness assumptions on the initial data, we have $T^*_\eps=+\infty$, at least for $\eps$ small enough. We need only prove that, for every $T<T^*_\eps$,
\[
    \norme{\na_x\ueps}_{\L^1(0,T;\L^\infty(\T^3))}\le\frac{1}{40}\qquad\text{and}\qquad\norme{F_\eps}_{\L^2((0,T)\times\T^3)}\le\frac{C^*}{4}.
\]
Indeed,  as seen from the proof of Lemma~\ref{EstimationNSsupGeneral}, Assumption~\ref{hypSmallDataNS} ensures that for all $T\ge0$,
\[
    \int_0^T\norme{e^{t\Delta}\ueps}_{\dot\H^1(\T^3)}^4\dd t\le\frac{C^*}{4}.
\]
Then, thanks to the timewise continuity of the norms involved, we deduce that $T^*_\eps=+\infty$.

We begin by proving the following lemma.

\begin{lemma}\label{lemmePara-3}
    Under Assumption~\ref{hypGeneral}, there exist $\eps_0>0$ and $\eta>0$ such that, for all $\eps\in(0,\eps_0)$, if
    \begin{equation}\label{FineHypf0}
        \norme{\feps^0}_{\L^1(\R^3;\L^\infty(\T^3))}\le\eta
    \end{equation}
    and
    \begin{equation}\label{wellPrep-3}
        \int_{\T^3\times\R^3}\frac{|v-u_{\eps}^0(x)|^p}{\eps^{p-1}}f_{\eps}^0(x,v)\dd x\dd v\le M,
    \end{equation}
    then for any $T<T^*_\eps$,
    \[
        \norme{\partial_t\ueps}_{\L^p((0,T)\times\T^3)}+\norme{\Delta_x\ueps}_{\L^p((0,T)\times\T^3)}\lesssim M^{\omega_p}
    \]
    for some $\omega_p>0$.
\end{lemma}

\begin{preuve}
    Thanks to Theorem \ref{ParabolicReg}, we have
    \begin{multline*}
        \norme{\partial_t\ueps}_{\L^p((0,T)\times\T^3)}+\norme{\Delta_x\ueps}_{\L^p((0,T)\times\T^3)}\\
        \lesssim\norme{F_\eps}_{\L^p((0,T)\times\T^3)}+\norme{(\ueps\cdot\na_x)\ueps}_{\L^p((0,T)\times\T^3)}+\norme{\ueps^0}_{\B_p^{s,p}(\T^3)}.
    \end{multline*}
    Combining the estimates in Lemmas \ref{EstimConvectStrongT}, \ref{LienFD-3}, and \ref{EstimateHighOrderD}, we get
    \begin{align*}
        &\norme{\partial_t\ueps}_{\L^p((0,T)\times\T^3)}+\norme{\Delta_x\ueps}_{\L^p((0,T)\times\T^3)}\\
        &\qquad\qquad \lesssim
            \norme{F_\eps}_{\L^p((0,T)\times\T^3)}+\norme{(\ueps\cdot\na_x)\ueps}_{\L^p((0,T)\times\T^3)}+\norme{\ueps^0}_{\B_p^{s,p}(\T^3)}\\
        &\qquad\qquad\lesssim
            \norme{\feps^0}_{\L^1(\R^3;\L^\infty(\T^3))}\norme{\partial_t\ueps}_{\L^p((0,T)\times\T^3)}\\
        &\qquad\qquad+   
            \left(\norme{\feps^0}_{\L^1(\R^3;\L^\infty(\T^3))}^{1-\frac{1}{p}}\eps^{\frac{1-\alpha_p}{p}}+\norme{\feps^0}_{\L^1(\R^3;\L^\infty(\R^3))}\EEE_\eps(0)^{\frac{1-\beta_p}{2}}\right)\\
        &\qquad\qquad \qquad \qquad    \times\norme{\Delta_x\ueps}_{\L^p(0,T)\times\T^3)}\\
        &\qquad\qquad+
            \norme{\feps^0}_{\L^1(\R^3;\L^\infty(\T^3))}^{1-\frac{1}{p}}\eps^{\frac{1-\alpha_p}{p}}\norme{\feps^0|v|^p}_{\L^1(\R^3;\L^\infty(\T^3))}^{\frac{1}{p(1-\alpha_p)}}\EEE_\eps(0)^{\frac{1}{2}}\\
        &\qquad\qquad+
            \norme{\feps^0}_{\L^1(\R^3;\L^\infty(\T^3))}\EEE_\eps(0)^{\frac{1-\beta_p}{2}}\Psi_{\eps,0}^{\frac{1}{2(1-\beta_p)}}\\
        &\qquad\qquad+
            \norme{\feps^0}_{\L^1(\R^3;\L^\infty(\T^3))}^{1-\frac{1}{p}}\norme{\feps^0(x,v)\frac{|v-\ueps^0(x)|^p}{\eps^{p-1}}}_{\L^1(\T^3\times\R^3)}^{\frac{1}{p}}\\
        &\qquad\qquad+   
            \Psi_{\eps,0}^{\frac{1}{2}}\EEE_\eps(0)^{\frac{(1-\beta_p)}{2}}\norme{\Delta_x\ueps}_{\L^p((0,T)\times\T^3)}^{\beta_p}
        +
            \norme{\ueps^0}_{\B^{s,p}_p(\T^3)}.
    \end{align*}
    We can use the smallness condition~\eqref{FineHypf0}, apply Young's inequality and, thanks to Assumption~\ref{hypGeneral}, we find that there exists $\eps_0>0$ and $\eta>0$ such that for every $\eps\in(0,\eps_0)$, if~\eqref{FineHypf0} holds, then
    \begin{align*}
    &\norme{\partial_t\ueps}_{\L^p((0,T)\times\T^3)}+\norme{\Delta_x\ueps}_{\L^p((0,T)\times\T^3)}\\
        &\qquad\qquad\lesssim
            \norme{\feps^0}_{\L^1(\R^3;\L^\infty(\T^3))}^{1-\frac{1}{p}}\eps^{\frac{1-\alpha_p}{p}}M^{\omega_p}
        +
            \norme{\feps^0}_{\L^1(\R^3;\L^\infty(\T^3))}M^{\omega_p}\\
        &\qquad\qquad+
            \norme{\feps^0}_{\L^1(\R^3;\L^\infty(\T^3))}^{1-\frac{1}{p}}M^{\omega_p}
        +   
            M^{\omega_p},
    \end{align*}
    hence the result.
\end{preuve}

\begin{remark}
     The well-preparedness assumption~\eqref{wellPrep-3} implies the following convergence
     \[
         W_1(f_{\eps}^0,\rho_{\eps}^0\otimes\delta_{v=u_{\eps}^0})\lesssim M\eps^{1-\frac{1}{p}}.
     \]
 \end{remark}

We have the following immediate consequence.
\begin{corollary}
    Under Assumption~\ref{hypGeneral}, there exist $\eps_0>0$ and $\eta>0$ such that, for all $\eps\in(0,\eps_0)$, if~\eqref{FineHypf0}--\eqref{wellPrep-3} hold and if
    \begin{equation}\label{hypSmallModulatedEnergy-3}
        \EEE_\eps(0)\leq \eta,
    \end{equation}
    then for any $T<T^*_\eps$,
    \[
        \norme{\na_x\ueps}_{\L^1(0,T;\L^\infty(\T^3))}\le\frac{1}{40}.
    \]
\end{corollary}

\begin{preuve}
    Since $p>3$, we can apply Lemmas~\ref{estimNabla} and~\ref{lemmePara-3} and find that there exist $\eps_0$ and $\eta>0$ such that for all $\eps\in(0,\eps_0)$, if~\eqref{FineHypf0}--\eqref{wellPrep-3} hold
    then, for some $\omega_p>0$,
    \[
        \int_0^T\norme{\na_x\ueps(t)}_{\L^\infty(\T^3)}\dd t\lesssim\EEE_{\eps}(0)^{\frac{1-\beta_p}{2}}\norme{\Delta_x\ueps}_{\L^p((0,T)\times\T^3)}^{\beta_p}
        \lesssim \EEE_{\eps}(0)^{\frac{1-\beta_p}{2}}M^{\omega_p},
    \]
    which can be made as small as necessary by reducing the value of $\eta>0$ in~\eqref{hypSmallModulatedEnergy-3}.
\end{preuve}

There remains to verify that
\[
    \norme{F_\eps}_{\L^2((0,T)\times\T^3)}\le\cfrac{C^*}{4}.
\]

We cannot use the exact same strategy as in Lemma \ref{EstimateHighOrderD} because the Gagliardo-Nirenberg inequality does not apply exactly similarly, but we are able to prove the following result.

\begin{lemma}\label{HighOrderD2}
    Under Assumptions~\ref{hypGeneral}, for every $T<T^*_\eps$,
    \begin{align*}
    \int_0^T\D_{\eps}^{(2)}(t)\dd t&\lesssim\norme{\feps^0}_{\L^1(\R^3;\L^\infty(\T^3))}\norme{\partial_t\ueps}_{\L^2((0,T)\times\T^3)}^2+\eps^{1-\alpha_2}\norme{\Delta_x\ueps}_{\L^2((0,T)\times\T^3)}^2\\
    &+\eps^{1-\alpha_2}\norme{\feps^0|v|^2}_{\L^1(\R^3;\L^\infty(\T^3))}^{\frac{1}{1-\alpha_2}}\EEE_\eps(0)\\
    &+\norme{\feps^0}_{\L^1(\R^3;\L^\infty(\T^3))}\Psi_{\eps,0}\EEE_\eps(0)^{1-\beta_p}\norme{\Delta_x\ueps}_{\L^p((0,T)\times\T^3)}^{2\beta_p}\\
    &+\int_{\T^3\times\R^3}\feps^0(x,v)\frac{|v-\ueps^0(x)|^2}{\eps}\dd x\dd v.
    \end{align*}
\end{lemma}

\begin{preuve}
    Applying Lemma \ref{lemmeExpressionD} with $r=2$, we obtain
    \begin{multline*}
        \int_0^T\D_\eps^{(2)}(t)\dd t
        \lesssim\norme{\feps^0}_{\L^1(\R^3;\L^\infty(\T^3))}\norme{\partial_t\ueps}_{\L^2((0,T)\times\T^3)}^2\\
        +\norme{|\na_x\ueps|m_2^{\frac{1}{2}}}_{\L^2((0,T)\times\R^3)}^2
        +\norme{\feps^0\frac{|v-\ueps^0|^2}{\eps}}_{\L^1(\T^3\times\R^3)}.
    \end{multline*}
    Following the proof of Lemma \ref{lemmeMixed}, we have
    \begin{align*}
        &\norme{|\na_x\ueps|m_2^{\frac{1}{2}}}_{\L^2((0,T)\times\T^3)}^2\\
       &\qquad \lesssim
            \int_0^T\int_{\T^3\times\R^3}|\na_x\ueps(t,x)|^2\feps^0(\tilde{X}^{t,x,w}_\eps(0),w)e^{-\frac{t}{\eps}}|w|^2\dd x\dd w\dd t\\
      &\qquad  +
            \int_0^T\int_{\T^3\times\R^3}|\na_x\ueps(t,x)|^2\feps^0(\tilde{X}^{t,x,w}_\eps(0),w)\\
    &\qquad \quad \qquad \qquad \times\left(\frac{1}{\eps}\int_0^te^{\frac{s-t}{\eps}}|\ueps(s,\tilde{X}^{t,x,w}_\eps(s))|\dd s\right)^2\dd x\dd w\dd t\\
     &\qquad \lesssim I_1+I_2,
    \end{align*}
    with
    \[
        I_1\lesssim\eps^{1-\alpha_2}\norme{\feps^0|v|^2}_{\L^1(\R^3;\L^\infty(\T^3))}^{\frac{1}{1-\alpha_2}}\EEE_\eps(0)+\eps^{1-\alpha_2}\norme{\Delta_x\ueps}_{\L^2((0,T)\times\T^3)}^2.
    \]
    To control $I_2$, we need to rely on Lemma~\ref{lemmePara-3}. Thanks to Jensen's inequality the Gagliardo-Nirenberg theorem, we have
        \begin{align*}
        I_2
        &\lesssim
            \eps^{-1}\int_0^T\int_{\T^3\times\R^3}\int_0^te^{\frac{s-t}{\eps}}|\na_x\ueps(t,x)|^2\feps^0(\tilde{X}^{t,x,w}_\eps(0),w)\\
        &\qquad \qquad \qquad \qquad \qquad  \times|\ueps(s,\tilde{X}_{x,t,w}(s))|^2\dd s\dd x\dd w\dd t\\
        &\lesssim
            \eps^{-1}\norme{\feps^0}_{\L^1(\R^3;\L^\infty(\T^3))}\int_0^T\int_0^te^{\frac{s-t}{\eps}}\norme{\na_x\ueps(t)}_{\L^\infty(\T^3)}^2\norme{\ueps(s)}_{\L^2(\T^3)}^2\dd s\dd t\\
        &\lesssim
            \eps^{-1}\norme{\feps^0}_{\L^1(\R^3;\L^\infty(\T^3))}\Psi_{\eps,0}\,\EEE_\eps(0)^{1-\beta_p}\\
        &\qquad \qquad \qquad \times\int_0^T\int_0^te^{\frac{s-t}{\eps}} e^{(1-\beta_p)\lambda_\eps t} \norme{\Delta_x\ueps(t)}_{\L^p(\T^3)}^{2\beta_p}\dd s\dd t\\
        &\lesssim
            \norme{\feps^0}_{\L^1(\R^3;\L^\infty(\T^3))}\Psi_{\eps,0}\,\EEE_\eps(0)^{1-\beta_p}\norme{\Delta_x\ueps}_{\L^p((0,T)\times\T^3)}^{2\beta_p}.
    \end{align*}
    Gathering all pieces together, we obtain the claimed estimate.
\end{preuve}

Therefore, the $\L^2$ parabolic estimate can be written as follows.

\begin{lemma}\label{lemmePara2-3}
    Under Assumption~\ref{hypGeneral}, there exist $\eps_0>0$ and $\eta>0$ such that, for all $\eps\in(0,\eps_0)$, if~\eqref{FineHypf0}--\eqref{wellPrep-3} are satisfied,
    then for any $T<T^*$,
    \[
        \norme{\partial_t\ueps}_{\L^2((0,T)\times\T^3)}+\norme{\Delta_x\ueps}_{\L^2((0,T)\times\T^3)}\lesssim M^{\omega_2},
    \]
    for some $\omega_2>0$.
\end{lemma}

\begin{preuve}
    Thanks to Theorem~\ref{ParabolicReg}, we have
    \begin{multline}
        \norme{\partial_t\ueps}_{\L^2((0,T)\times\T^3)}+\norme{\Delta_x\ueps}_{\L^2((0,T)\times\T^3)}\\
        \lesssim\norme{F_\eps}_{\L^2((0,T)\times\T^3)}+\norme{(\ueps\cdot\na_x)\ueps}_{\L^2((0,T)\times\T^3)}+\norme{\ueps^0}_{\H^1(\T^3)}.
    \end{multline}
    On the one hand, combining Lemmas~\ref{LienFD-3} and~\ref{HighOrderD2}, we get
    \begin{align*}
        \norme{F_\eps}_{\L^2((0,T)\times\T^3)}
        &\lesssim\norme{\feps^0}_{\L^1(\R^3;\L^\infty(\T^3))}\norme{\partial_t\ueps}_{\L^2((0,T)\times\T^3)}\\
            &+\eps^{\frac{1-\alpha_2}{2}}\norme{\feps^0}_{\L^1(\R^3;\L^\infty(\T^3))}\norme{\Delta_x\ueps}_{\L^2((0,T)\times\R^3)}\\
            &+\eps^{\frac{1-\alpha_2}{2}}\norme{\feps^0}_{\L^1(\R^3;\L^\infty(\T^3))}^{\frac{1}{2}}\norme{\feps^0|v|^2}_{\L^1(\R^3;\L^\infty(\T^3))}^{\frac{1}{2(1-\alpha_2)}}\EEE_\eps(0)^{\frac{1}{2}}\\
            &+\norme{\feps^0}_{\L^1(\R^3;\L^\infty(\T^3))}\Psi_{\eps,0}^{\frac{1}{2}}\EEE_\eps(0)^{\frac{1-\beta_p}{2}}\norme{\Delta_x\ueps}_{\L^p((0,T)\times\T^3)}^{\beta_p}\\
            &+\norme{\feps^0}_{\L^1(\R^3;\L^\infty(\T^3))}^{\frac{1}{2}}\left(\int_{\T^3\times\R^3}\feps^0(x,v)\frac{|v-\ueps^0(x)|^2}{\eps}\dd x\dd v\right)^{\frac{1}{2}}.
    \end{align*}
    Since
    \begin{align*}
        \norme{\feps^0|v|^2}_{\L^1(\R^3;\L^\infty(\T^3))}&\le\int_{|v|\le1}\norme{\feps^0(v)}_{\L^\infty(\T^3)}\dd v+\int_{|v|\ge1}\norme{\feps^0(v)}_{\L^\infty(\T^3)}|v|^p\dd v\\
        &\lesssim\eta+M,
    \end{align*}
    and
    \begin{multline*}
        \int_{\T^3\times\R^3}\feps^0(x,v)\frac{|v-\ueps^0(x,v)|^2}{\eps}\dd x\dd v\\
        \le
            \left(\int_{\T^3\times\R^3}\feps^0(x,v)\dd x\dd v\right)^{1-\frac{2}{p}}\left(\int_{\T^3\times\R^3}\feps^0\frac{|v-\ueps^0|^p}{\eps^{\frac{p}{2}}}\dd x\dd v\right)^{\frac{2}{p}}\\
        \le
            \eta^{1-\frac{2}{p}} \eps^{1-\frac{2}{p}}\left(\int_{\T^3\times\R^3}\feps^0\frac{|v-\ueps^0|^p}{\eps^{\frac{p}{2}}}\dd x\dd v\right)^{\frac{2}{p}}
        \lesssim \eta^{1-\frac{2}{p}}M^{\frac{2}{p}} \eps^{1-\frac{2}{p}},
    \end{multline*}
    we have, thanks to Lemma~\ref{lemmePara-3}, for $\eps$ and $\eta$ small enough,
    \begin{equation}\label{estimFL2Preuve-3}
        \norme{F_\eps}_{\L^2((0,T)\times\T^3)}\lesssim\eta\norme{\partial_t\ueps}_{\L^2((0,T)\times\T^3)}+\eta\eps^{\frac{1-\alpha_2}{2}}\norme{\Delta_x\ueps}_{\L^2((0,T)\times\T^3)}
        +\eps^{1-\frac{2}{p}}.
    \end{equation}
    
    On the other hand, thanks to Hölder's inequality and Lemmas~\ref{EstimationNSsupGeneral},~\ref{lemmeDE}--\ref{estimNabla} and~\ref{lemmePara-3}, we have, for $\eps$ and $\eta$ small enough,
    \begin{align*}
        \norme{(\ueps\cdot\na)\ueps}_{\L^2((0,T)\times\T^3)}
        &\lesssim\norme{\ueps}_{\L^\infty(0,T;\L^2(\T^3))}^2\int_0^T\norme{\na_x\ueps(t)}_{\L^\infty(\T^3)}^2\dd t\\
        &\lesssim\Psi_{\eps,0}\EEE_\eps(0)^{1-\beta_p}\int_0^T\norme{\Delta_x\ueps(t)}_{\L^p(\T^3)}e^{-(1-\beta_p)\lambda_\eps t}\dd t\\
        &\lesssim\Psi_{\eps,0}\EEE_\eps(0)^{1-\beta_p}\norme{\Delta_x\ueps}_{\L^p((0,T)\times\T^3)}^{2\beta_p} \\
        &\lesssim \eta^{1-\beta_p}M^{\omega_p},
    \end{align*}
    for some $\omega_p>0$. Hence the result.
\end{preuve}

We inject the estimate of Lemma~\ref{lemmePara2-3} into the estimate~\eqref{estimFL2Preuve-3} that was obtained in the course the proof to find that if $\eta$ and $\eps$ are small enough, then
\[
    \norme{F_\eps}_{\L^2((0,T)\times\T^3)}\le\frac{C^*}{4},
\]
and the bootstrap argument is finally complete.


\subsection{Non quantitative convergence in the well-prepared case}
\label{sec-first-fine}

In the following paragraphs, we will consider a fixed time horizon $T>0$ and study the convergence of the sequences $(\feps)_{\eps>0}$, $(\rhoeps)_{\eps>0}$ and $(\ueps)_{\eps>0}$ under several sets of convergence assumptions for the initial data. 

Our first result only requires the weak convergence of the fluid velocity and the particle density, but the higher singularity of the Brinkman force prevents us from finding explicit convergence rates as in Section~\ref{SectionLandLFParticle}.

\begin{theorem}\label{FineFirstConv}
    Under Assumptions \ref{hypGeneral}--\ref{hypSmallDataNS}, there exist $\eps_0>0$ and $\eta>0$ such that if, for all $\eps\in(0,\eps_0)$, \eqref{FineHypf0}--\eqref{hypSmallModulatedEnergy-3} hold
    and if
    \[
        \ueps\xrightharpoonup[\eps\to0]{} u^0\text{ in }w\text{-}\L^2(\T^3)\qquad\text{and}\qquad\rhoeps^0\xrightharpoonup[\eps\to0]{}\rho^0\text{ in }w^*\text{-}\L^\infty(\T^3),
    \]
    then $(\ueps)_{\eps>0}$ converges to $u$ in $\L^2((0,T)\times\T^3)$, $(\rhoeps)_{\eps>0}$ converges weakly-$*$ to $\rho$ in $\L^\infty((0,T)\times\T^3)$ where $(\rho,u)$ satisfies, with $\tilde{\rho}=1+\rho$,
    \begin{equation*}
        \left\{
\begin{aligned}
        &\partial_t\tilde{\rho}+\div_x(\tilde{\rho}u)=0,\\
        &\tilde{\rho}|_{t=0}= 1+ \rho^0, \\
        &\partial_t(\tilde{\rho}u)+\div_x(\tilde{\rho}u\otimes u)-\Delta_x u+\nabla_xp=0,\\
        &\div_xu=0, \\
          &u|_{t=0}= u^0.
    \end{aligned}
\right.
    \end{equation*}
    Furthermore, for almost all $t\in (0,T)$,
    \begin{equation}
    \label{W1rhofine}
        W_1(\rhoeps(t),\rho(t))\lesssim W_1(\rhoeps^0,\rho^0)+\norme{\ueps-u}_{\L^2((0,t)\times\T^3)}+\eps\xrightarrow[\eps\to0]{}0,
    \end{equation}
    and $(\feps)_{\eps>0}$ converges to $\rho\otimes\delta_{v=u}$ in the sense that 
    \[
     \int_0^TW_1(\feps(t),\rho(t)\otimes\delta_{v=u(t)})\dd t\xrightarrow[\eps\to0]{}0.
    \]
\end{theorem}

\begin{preuve}
    The  energy dissipation estimate~\eqref{EstimationEnergieEq} first shows that $(\ueps)$ is uniformly bounded in $\L^2(0,T;\H^1(\T^3))$. Therefore, up to a subsequence that we shall not write, $(\ueps)$ converges to some $u$ in $w$-$\L^2(0,T;\H^1(\T^3))$. Furthermore, since we have proven in the previous section that all times are strong existence times and thanks to Lemma~\ref{lemmePara2-3}, we can apply the Aubin-Lions lemma and obtain the strong convergence of $(\ueps)$ to $u$ in $\L^2((0,T)\times\T^3)$.
    
    Moreover, the fact that $T$ is a strong existence time enables us to apply Corollary~\ref{RhoLinfini}, which ensures that $(\rhoeps)$ is bounded in $\L^\infty((0,T)\times\T^3)$, so that it converges, up to a subsequence, to some $\rho$ in $w^*$-$\L^\infty((0,T)\times\T^3)$. Therefore, $(\rhoeps\ueps)$ converges to $\rho u$ in $w$-$\L^2((0,T)\times\T^3)$.
    
    Contrary to the \emph{light} and \emph{light and fast} regimes, we cannot state that the Brinkman force $F_\eps$ converges to 0. Yet, thanks to Lemmas~\ref{LienFD-3},~\ref{EstimateHighOrderD} and~\ref{lemmePara-3}, $(F_\eps)$ is bounded and therefore converges, up to a subsequence, to some $F$ in $w$-$\L^2((0,T)\times\T^3)$.
    
    Finally, $(\jeps)$ also converges to $\rho u$ in $w$-$\L^2((0,T)\times\T^3)$. Indeed, for any $\varphi\in\L^2((0,T)\times\T^3)$, thanks to Hölder's inequality and the modulated energy--dissipation estimate from Lemma~\ref{lemmeDE}, we have, for $\eps\in(0,\eps_0)$,
     \begin{align*}
        &\left|\int_0^T\int_{\T^3}(\jeps(t,x)-\rhoeps(t,x)\ueps(t,x))\varphi(t,x)\dd x\dd t\right|\\
        &\le\int_0^T\int_{\T^3\times\R^3}\feps(t,x,v)|v-\ueps(t,x)||\varphi(t,x)|\dd x\dd v\dd t\\
        &\le\norme{\feps^0}_{\L^1(\R^3;\L^\infty(\T^3)}^{\frac{1}{2}}\norme{\varphi}_{\L^2((0,T)\times\T^3)} \left(\int_0^T\int_{\T^3\times\R^3}\feps(t,x,v)|v-\ueps(t,x)|^2\dd x\dd v\dd t\right)^{\frac{1}{2}}\\
        &\le \norme{\feps^0}_{\L^1(\R^3;\L^\infty(\T^3)}^{\frac{1}{2}}\norme{\varphi}_{\L^2((0,T)\times\T^3)}\eps^{\frac{1}{2}}\EEE_\eps(0)^{\frac{1}{2}}\\
        &\lesssim_M \sqrt{\eps}\,\eta\norme{\varphi}_{\L^2((0,T)\times\T^3)}.
    \end{align*}
    
    We can take the limit $\eps\to0$ in the conservation of mass~\eqref{conservationMass} and obtain
    \[
        \partial_t\rho+\div_x(\rho u)=0.
    \]
    Let us now take the limit $\eps\to0$ in the conservation of momentum~\eqref{conservationMomentum}
    \[
        \partial_t\jeps+\div_x\left(\int_{\R^3}\feps\,v\otimes v\,\dd v\right)=-F_\eps.
    \]
    We have
    \begin{align*}
         \int_{\R^3}\feps v\otimes f\dd v
         &=\int_{\R^3}\feps(v-\ueps)\otimes(v-\ueps)\dd v\\
        &\quad+\int_{\R^3}\feps(v\otimes\ueps+\ueps\otimes v)\dd v-\int_{\R^3}\feps\ueps\otimes\ueps\dd v\\
        &:= I_{\eps,1}+I_{\eps,2}+I_{\eps,3}.
    \end{align*}
    Thanks to the modulated energy estimate from Lemma~\ref{lemmeDE},
    \[
        \int_0^T\int_{\T^3}|I_{\eps,1}|\lesssim\int_0^T\int_{\T^3\R^3}\feps|v-\ueps|^2\dd v\xrightarrow[\eps\to0]{}0.
    \]
    Furthermore, we can write
    \begin{multline*}
        I_{\eps,2}-2\rho u\otimes u
        =\jeps\otimes \ueps-\rho u\otimes u+\ueps\otimes\jeps-\rho u\otimes u\\
        =(\jeps-\rho u)\otimes u+\jeps\otimes(\ueps-u)+\ueps\otimes(\jeps-\rho u)+\rho(\ueps-u)\otimes u,
    \end{multline*}
    so that we can prove that $(I_{\eps,2})_{\eps>0}$ converges to $2\rho\,u\otimes u$ in the distribution sense thanks to the previous bounds, weak convergences and strong convergence of $(\ueps)_{\eps>0}$. Finally
    \[
        I_{\eps,3}+\rho u\otimes u=(\rho-\rhoeps)\, u\otimes u+\rhoeps(\ueps-u)\otimes u,
    \]
    so $(I_{\eps,3})_{\eps>0}$ converges to $-\rho\,u\otimes u$ in the distribution sense. Therefore, we have proven, in the distribution sense,
    \begin{equation}\label{EqLimRho-3}
        \partial_t(\rho u)+\div_x(\rho\,u\otimes u)=-F.
    \end{equation}
    Thanks to the strong convergence of $(\ueps)_{\eps>0}$ and the weak convergence of $(F_\eps)_{\eps>0}$ to $F$, we can also take the limit in the Navier-Stokes equations to obtain that $u$ satisfies
    \begin{equation}\label{EqLimU-3}
        \partial_tu+\div_x(u\otimes u)-\Delta_xu+\na_xp=F.
    \end{equation}
    Summing \eqref{EqLimRho-3} and \eqref{EqLimU-3} leads to
    \[
        \partial_t(\tilde{\rho}u)+\div_x(\tilde{\rho}u\otimes u)-\Delta_xu+\na_xp=0,
    \]
    where $\tilde{\rho}=1+\rho$.

       We conclude that no subsequence is required by the same weak compactness and uniqueness argument as at the end of of the proof Theorem~\ref{thCV-12}.
       The fact that the limit system has a unique solution follows from~\cite{pai-zha-zha}.
    
    
    Let us now prove the desired convergence for $(\rho_\eps)$. Consider
    $$\chieps=\rho-\rhoeps.$$ 
    We have
    \[
        \partial_t\chieps=-\div_x(\rho u)+\div_x\jeps
        =-\div_x((\rho-\rhoeps)u)-\div_x(\rhoeps u)+\div_x\jeps
    \]
    so that $\chieps$ is a solution to the transport equation with source
    \[
        \partial\chieps+\div_x(\chieps u)=-\div_xJ_\eps,
    \]
    where $J_\eps=\rhoeps(u-\ueps)-\eps F_\eps$. Following the proof of Theorem~\ref{thCVprecis-12}, we obtain that, for almost every $t\in(0,T)$,
    \begin{align*}
        W_1(\rhoeps(t),\rho(t))&\lesssim W_1(\rhoeps^0,\rho^0)+\norme{\ueps-u}_{\L^2((0,t)\times\T^3)}+\eps\norme{F_\eps}_{\L^2((0,t)\times\T^3)}\\
        &\lesssim W_1(\rhoeps^0,\rho^0)+\norme{\ueps-u}_{\L^2((0,t)\times\T^3)}+\eps,
    \end{align*}
    since $T$ is a strong existence time.

    
    To prove the convergence of $(\feps)_{\eps>0}$, we first note that, for any $\varphi\in\CCC^\infty(\T^3\times\R^3)$ such that $\norme{\na_{x,v}\varphi}_{\L^\infty(\T^3\times\R^3)}\le1$, thanks to Lemma~\ref{lemmeDE},
    \begin{align*}
       &|\langle\feps(t)-\rhoeps(t)\otimes\delta_{v=\ueps(t)},\varphi\rangle| \\
        &\qquad \qquad \le\norme{\na_{x,v}\varphi}_{\L^\infty(\T^3\times\R^3)}
        \int_{\T^3\times\R^3}\feps(t,x,v)|v-\ueps(t,x)|\dd x
        \\ 
        &\qquad \qquad \leq \int_{\T^3\times\R^3}\feps(t,x,v)|v-\ueps(t,x)|\dd x \dd x
    \end{align*}
    so that
    $$
    W_1(\feps(t), \rhoeps(t)\otimes\delta_{v=\ueps(t)}) \leq \int_{\T^3\times\R^3}\feps(t,x,v)|v-\ueps(t,x)|\dd x,
    $$
    and thus by the modulated energy--dissipation inequality this yields
    $$
    \int_0^T      W_1(\feps(t), \rhoeps(t)\otimes\delta_{v=\ueps(t)}) \dd t \le\sqrt{\eps}\EEE_\eps(0)^{\frac{1}{2}}\sqrt{T}.
    $$
    Furthermore,
    \begin{multline*}
       |\langle\rhoeps(t)\otimes\delta_{v=\ueps(t)}-\rho(t) \otimes \delta_{v=u(t)},\varphi\rangle| \\
        \lesssim\left|\int_{\T^3}(\rhoeps(t,x)-\rho(t,x))\varphi(x,u(t,x))\dd x\right|+\norme{\rhoeps(t)}_{\L^\infty(\T^3)}\norme{\ueps(t)-u(t)}_{\L^2(\T^3)}.
    \end{multline*}
    We have the bound
    $$
    \left|\int_{\T^3}(\rhoeps(t,x)-\rho(t,x))\varphi(x,u(t,x))\dd x\right| \leq \left(1+ \| \na_x u (t) \|_{\L^\infty(\T^3)}\right) W_1(\rhoeps(t), \rho(t)).
    $$
    and thus it will be treated thanks to~\eqref{W1rhofine}. We 
    eventually obtain
    \begin{multline*}
    \int_0^T W_1(\rhoeps(t)\otimes\delta_{v=\ueps(t)},\rho(t) \otimes \delta_{v=u(t)}) \dd t \\
    \lesssim_M \int_0^T W_1(\rhoeps(t), \rho(t)) \dd t + \norme{\ueps(t)-u(t)}_{\L^2((0,T) \times \T^3)}^2.
    \end{multline*}
        This leads to the expected convergence for $(\feps)$.



\end{preuve}

Assume now that 
\begin{equation}
\label{hyp-u0-strong}
 \ueps\xrightarrow[\eps\to0]{} u^0\text{ in }\L^2(\T^3).
\end{equation}
(As already explained in the introduction, thanks to Assumption~\ref{hypGeneral}, this can always be ensured up to taking a subsequence.)
By the Aubin-Lions lemma, up to a subsequence, we then have
$$
\sup_{[0,T]} \norme{u_\eps-u}_{\L^2(\T^3)} \xrightarrow[\eps\to0]{} 0.
$$
As opposed to the \emph{light} and \emph{light and fast} regimes, it is not straightforward to ensure that this holds without taking a subsequence, and if it is the case, to provide a rate of convergence. This will be the object of the next section. We can still provide quantitative pointwise convergence  for $(\feps)$, that will depend on this yet unclear convergence.

\begin{corollary}\label{CorollaryFineFirstConv}
    Under the assumptions of Theorem~\ref{FineFirstConv} and~\eqref{hyp-u0-strong}, if
    \[
        \norme{\na_xu}_{\L^\infty((0,T)\times\T^3)}<+\infty,
    \]
    then, for almost every $t\in(0,T)$,
    \begin{multline*}
        W_1(\feps(t),\rho(t)\otimes\delta_{v=u(t)})\\
        \lesssim W_1(\rhoeps^0,\rho^0)+\norme{\ueps(t)-u(t)}_{\L^2(\T^3)}+\norme{\ueps-u}_{\L^2((0,t)\times\T^3)}+\eps.
    \end{multline*}
\end{corollary}

\begin{preuve}

    To prove the pointwise convergence of $(\feps)_{\eps>0}$, we first note that, for any $\psi\in\CCC^\infty(\T^3\times\R^3)$ such that $\norme{\na_{x,v}\psi}_{\L^\infty(\T^3\times\R^3)}\le1$, thanks Hölder's inequality and Lemmas~\ref{EstimateHighOrderD} and~\ref{lemmePara-3}, for almost every $t\in(0,T)$,
    \begin{align*}
        &\left|\langle\feps(t)-\rhoeps(t)\otimes\delta_{v=\ueps(t)},\psi\rangle\right|\\
        &\qquad \le\norme{\na_{x,v}\psi}_{\L^\infty(\T^3\times\R^3)}\int_{\T^3\times\R^3}\feps(t,x,v)|v-\ueps(t,x)|\dd x\dd v\\
        &\qquad \lesssim M^{\omega_p}\eps^{p-1}\int_{\T^3\times\R^3}\feps(t,x,v)\frac{|v-\ueps(t,x)|^p}{\eps^{p-1}}\dd x\dd v
        \lesssim M^{\omega_p'}\eps^{p-1},
    \end{align*}
    for some $\omega_p,\omega_p'>0$. Furthermore,
    \begin{multline*}
        |\langle\rhoeps(t)\otimes\delta_{v=\ueps(t)}-\rho(t)\otimes\delta_{v=u(t)},\psi\rangle|
        \le \left|\int_{\T^3}(\rhoeps(t,x)-\rho(t,x))\psi(t,x,u(t,x))\dd x\right|\\
        +\left|\int_{\T^3}\rhoeps(t,x)(\psi(t,x,\ueps(t,x)-\psi(t,x,u(t,x)))\dd x\right|.
    \end{multline*}
    On the one hand,
    \begin{align*}
        \left|\int_{\T^3}(\rhoeps(t,x)-\rho(t,x))\psi(t,x,u(t,x))\dd x\right|
        &\le W_1(\rhoeps(t),\rho(t))\sup_{x\in\T^3}|\na_x\left(\psi(x,u(t,x))\right)|\\
        &\le W_1(\rhoeps(t),\rho(t))\left(1+\norme{\na_xu(t)}_{\L^\infty(\T^3)}\right) \\
        &\lesssim W_1(\rhoeps(t),\rho(t)).
    \end{align*}
    On the other hand, since $T$ is a strong existence time, thanks to Corollary~\ref{RhoLinfini},
    \begin{multline*}
        \left|\int_{\T^3}\rhoeps(t,x)(\psi(t,x,\ueps(t,x)-\psi(t,x,u(t,x)))\dd x\right|\\
        \lesssim\int_{\T^3}|\ueps(t,x)-u(t,x)|\dd x\le\norme{\ueps(t)-u(t)}_{\L^2(\T^3)},
    \end{multline*}
    which completes the proof of the theorem.
   
\end{preuve}

\begin{remark} As for the \emph{light} and \emph{light and fast} particle regime, we obtain as a straightforward side consequence of the analysis, the description of the long time behavior of solutions to the Vlasov-Navier-Stokes system in the fine particle regime. In this regime, it turns out that the required assumptions are essentially the same as in \cite{hank-mou-moy}, as opposed to the results of Section~\ref{sec-long}.
More interestingly,  by uniformity with respect to $\eps$ of the exponential decay of the modulated energy, we also recover as a corollary a description of the long time behavior of solutions to the Inhomogeneous Navier-Stokes equations \eqref{inNS}. This is of course a well-known result, see e.g. \cite[Theorem 1.4]{pou}.
\end{remark}


\subsection{Relative entropy estimates}
\label{sec-relat-fine}

We introduce in this section a relative entropy method  for the fine particle regime with the aim to 
\begin{itemize}
\item justify the pointwise in time convergence of $(u_\eps)$ towards $u$ (without taking a subsequence),
\item obtain full quantitative (with respect to $\eps$) estimates. 
\end{itemize}
This will be achieved after assuming that the solution to the limit system is smooth enough.
A side objective is to justify why the related entropy method by itself  does not allow to prove any convergence result but rather has to be used in conjunction with the previous bootstrap analysis. Namely, the relative entropy method requires exactly the same uniform bound for $\rhoeps$ in $\L^\infty$  as for the other parts of this work.





We introduce the relative entropy
\begin{equation*}
\index{H@$\He(t)$: relative entropy}
\begin{aligned}
\He(t) &:=\frac{1}{2}\int_{\T^3\times\R^3}\feps(t,x,v)|v-u(t,x)|^2\dd x\dd v
                +   \frac{1}{2}\int_{\T^3}|\ueps(t,x)-u(t,x)|^2\dd x,
\end{aligned}
\end{equation*}
where $(\rho,u)$ is a smooth solution to the Inhomogeneous incompressible Navier-Stokes equations (note that $\rho$ does not appear explicitly in $\He$, though). The precise required smoothness shall be clarified in the upcoming statements.

The following is the analogue of \cite[Lemma 3]{gou-jab-vas04b} in the case without diffusion for the Vlasov equation.
\begin{lemma}
\label{lem-evolution-HH}
Under Assumption~\ref{hypGeneral}, the relative entropy satisfies for all $t \geq 0$
\begin{equation}
\label{eq-derivH}
\begin{aligned}
\He(t) &+ \int_0^t \int_{\T^3}|\na_x(\ueps - u )|^2\dd x \dd s + \frac{1}{\eps} \int_0^t \int_{\T^3 \times \R^3}|v-\ueps|^2\feps \dd x\dd v \dd s \\
&\leq \He(0) + \int_0^t  \sum_{j=1}^{4} I_j(s)  \dd s.
\end{aligned}
\end{equation}
with
\begin{align*}
&I_1 := - \int_{\T^3 \times \R^3}  f_\eps (v-u) \otimes (v-u) : \na_x u \dd x \dd v, \\
&I_2 := - \int_{\R^3}  (\ueps - u ) \otimes (\ueps -u) : \na_x u \dd x, \\
&I_3 :=  \int_{\T^3 \times \R^3}  \feps (v-\ueps) \cdot G  \dd x \dd v , \\
&I_4 :=  \int_{\T^3} (\rhoeps- \rho) (\ueps -u ) \cdot G \dd x,
\end{align*}
where $G = \frac{\na_x p - \Delta_x u}{1+\rho}$.

\end{lemma}

\begin{preuve}
Let us argue as if $f_\eps$ were a smooth function.
We first compute
\begin{align*}
\frac{\dd}{\dd t}  \frac{1}{2}&\int_{\T^3\times\R^3}\feps |v-u|^2 \dd x\dd v \\
&= \frac{1}{2}\int_{\T^3\times\R^3}|v-u|^2 \pa_t \feps\dd x\dd v  -\int_{\T^3\times\R^3}\feps (v-u) \cdot \pa_t u \dd x\dd v \\
&= - \int_{\T^3\times\R^3} \feps ((v-u) \cdot \na_x u)\cdot v  \dd x \dd v +  \frac{1}{\eps} \int_{\T^3\times\R^3}  \feps (v-u) \cdot (\ueps -v)  \dd x \dd v \\
&\quad +\int_{\T^3\times\R^3} \feps (v-u) \cdot (u\cdot \na_x u) \dd x \dd v  \int_{\T^3\times\R^3} \feps (v-u)  \cdot G \dd x \dd v \\
&= I_1 +   \frac{1}{\eps} \int_{\T^3\times\R^3}  \feps (v-u) \cdot (\ueps -v)  \dd x \dd v  + \int_{\T^3\times\R^3} \feps (v-u)  \cdot  G \dd x \dd v . 
\end{align*}
Furthermore,
\begin{align*}
\frac{\dd}{\dd t} \frac{1}{2}&\int_{\T^3}|\ueps-u|^2\dd x = \int_{\T^3}(\ueps-u) \cdot \pa_t (\ueps-u) \dd x \\
&= \int_{\T^3}  (\ueps-u) \cdot (u \cdot \na_x u - \ueps \cdot \na_x \ueps)  \dd x + \frac{1}{\eps} \int_{\T^3\times\R^3}  \feps (\ueps-u)\cdot (v-\ueps)   \dd x \dd v  \\
&+ \int_{\T^3}  (\ueps-u) \cdot G   \dd x + \int_{\T^3} (\ueps-u) \cdot \Delta_x \ueps  \dd x.
\end{align*}
Since $\div \ueps =0$, we have
\begin{align*}
 \int_{\T^3}  &(\ueps-u) \cdot (u \cdot \na_x u - \ueps \cdot \na_x \ueps)  \dd x  
 \\
 &= -  \int_{\T^3}   (\ueps - u ) \otimes (\ueps -u) : \na_x u  \dd x +  \frac{1}{2} \int_{\T^3}     \ueps \cdot \na_x |\ueps-u|^2  \dd x\\
& = I_2.
\end{align*}
Since $\div u =0$ and by definition of $G$,
$$
\Delta_x u =  - (1+ \rho) G + \na_x p,
$$
we can also rewrite
\begin{align*}
\int_{\T^3} (\ueps-u) \cdot \Delta_x \ueps  \dd x &=  \int_{\T^3} (\ueps-u) \cdot \Delta_x (\ueps-u)  \dd x +   \int_{\T^3} (\ueps-u) \cdot \Delta_x u  \dd x \\
&= -\int_{\T^3} |\na_x (\ueps-u)|^2  \dd x  - \int_{\T^3} (1+\rho) (\ueps-u) \cdot G   \dd x
\end{align*}
All in all, we have obtained
\begin{align*}
\frac{\dd}{\dd t} \frac{1}{2}\int_{\T^3}|\ueps-u|^2\dd x &= I_2 -\int_{\T^3} |\na_x (\ueps-u)|^2  \dd x \\
 &+ \frac{1}{\eps} \int_{\T^3\times\R^3}  \feps (\ueps-u)\cdot (v-\ueps)   \dd x \dd v
- \int_{\T^3} \rho (\ueps-u) \cdot G   \dd x.
\end{align*}
We finally write
\begin{align*}
&\frac{1}{\eps} \int_{\T^3\times\R^3}  \feps (\ueps-u)\cdot (v-\ueps)   \dd x \dd v + \frac{1}{\eps} \int_{\T^3\times\R^3}  \feps (\ueps-u)\cdot (v-\ueps)   \dd x \dd v \\
&= - \frac{1}{\eps} \int_{\T^3 \times \R^3}|v-\ueps|^2\feps \dd x\dd v
\end{align*}
and 
\begin{align*}
\int_{\T^3\times\R^3} &\feps (v-u)  \cdot  G\dd x \dd v - \int_{\T^3} \rho (\ueps-u) \cdot G   \dd x \\
&= \int_{\T^3\times\R^3} \feps (v-\ueps)  \cdot G\dd x \dd v + \int_{\T^3} (\rhoeps-\rho) (\ueps-u)  \cdot  G\dd x \\
&= I_3 + I_4.
\end{align*}
To conclude, adding everything up, this yields
\begin{equation*}
\frac{\dd}{\dd t}  \He + \int_{\T^3}|\na_x(\ueps - u )|^2\dd x  + \frac{1}{\eps} \int_{\T^3 \times \R^3}|v-\ueps|^2\feps \dd x\dd v = \sum_{j=1}^4 I_j,
\end{equation*}
which concludes the computation. To obtain~\eqref{eq-derivH} in the general case, we notice that we can take  $|v-u(t,x)|^2$ as a test function in the weak formulation of the Vlasov equation (by a standard density argument).
\end{preuve}

The control of $I_4$ requires a pointwise estimate on $(\rhoeps-\rho)$ in $\dot\H^{-1}(\T^3)$, as opposed to the Wasserstein control obtained in Theorem~\ref{FineFirstConv}. Nevertheless, the same strategy leads to the following result. 

\begin{lemma}
\label{lem-rho-H-1}
Let $T>0$. Assume there exist $C>0$ and $\eps_0\in(0,1)$ such that
\begin{equation}
\label{eq-condi-rho}
    \forall\eps\in(0,\eps_0),\qquad\| \rhoeps \|_{\L^\infty(0,T \times \T^3)} \leq C.
\end{equation}
Then for all $t \in (0,T)$,
    \begin{align*}
        \norme{\rhoeps(t)-\rho(t)}_{\dot{\H}^{-1}(\T^3)}
        &\lesssim\norme{\rhoeps^0-\rho^0}_{\dot{\H}^{-1}(\T^3)}\\
            &+\eps\int_0^t\norme{F_\eps(s)}_{\L^2(\T^3)}+\int_0^t\norme{\ueps(s)-u(s)}_{\L^2(\T^3)}.
    \end{align*}
\end{lemma}

\begin{preuve}
    As in the proof of Theorem~\ref{FineFirstConv}, $\chieps=\rho-\rhoeps$ satisfies the following transport equation
    \[
        \partial\chieps+\div_x(\chieps u)=-\div_xJ_\eps,
    \]
    where $J_\eps=\rhoeps(u-\ueps)-\eps F_\eps$. Recall the characteristics of the equation: for $t\ge0$ and $x\in\T^3$, we consider the solution $Y(\cdot;t,x)$ to the Cauchy problem
    \[
        \frac{\dd }{\dd s}{Y}(s;t,x)=u(s,Y(s;t,x))\qquad\text{and}\qquad Y(t;t,x)=x.
    \]
    Since $\div_xu=0$, we have $\det\na_xY(s;t,x)=1$ for every $s,t\ge0$ and $x\in\T^3$. 
    Therefore, for every $\varphi\in\CCC^\infty(\T^3)$ such that $\norme{\na_x\varphi}_{\L^2(\T^3)}\le1$, for every $t\in[0,T]$, thanks to the method of characteristics,
    \begin{align*}
        \int_{\T^3}\chieps(t,x)\varphi(x)\dd x
        &=\int_{\T^3}\chieps(0,x)\varphi(x)\dd x-\int_0^t\int_{\T^3}(\div_xJ_\eps)(s,Y(s,t,x))\varphi(x)\dd x\dd s\\
        &=\int_{\T^3}(\rhoeps^0-\rho^0)\varphi+\int_0^t\int_{\T^3}J_\eps(s,Y(s,t,x))\cdot\na_x\varphi(x)\dd x\dd s.
    \end{align*}
    Applying the change of variable $x'=Y(s;t,x)$ yields
    \[
        \left(\int_{\T^3}|J_\eps(s,Y(s,t,x))|^2\dd x\right)^{\frac{1}{2}}
        \lesssim\eps\norme{F_\eps}_{\L^2(\T^3)}+\norme{\rhoeps}_{\L^\infty((0,T)\times\T^3)}\norme{\ueps-u}_{\L^2(\T^3)}
    \]
    so that, thanks to the Cauchy-Schwarz inequality,
    \begin{multline*}
        \left|\int_{\T^3}\chieps(t,x)\varphi(x)\dd x\right|
        \lesssim\norme{\rhoeps^0-\rho^0}_{\dot{\H}^{-1}(\T^3)}\\
            +\eps\int_0^t\norme{F_\eps(s)}_{\L^2(\T^3)}+\int_0^t\norme{\ueps(s)-u(s)}_{\L^2(\T^3)},
    \end{multline*}
    which leads to the expected result.
\end{preuve}

Let us now state the estimate on $\He$ we have obtained thanks to the previous results.

\begin{proposition}\label{prop-relativeentropy}
    Let $T>0$. Assume that
    \begin{equation}\label{hypFineSol}
        \| G \|_{\L^\infty((0,T)\times \T^3)} + \| \na_x G \|_{\L^\infty((0,T)\times \T^3)} \lesssim 1.
    \end{equation}
    There exist $\eps_0>0$ and $C_T>0$ such that, if~\eqref{eq-condi-rho} holds, then for all $\eps\in(0,\eps_0)$ and $t \in (0,T)$,
    \begin{equation*}
        \begin{aligned}
            \He(t) 
            \leq C_T \left(\He(0) + \| \rhoeps(0) - \rho(0)\|^2_{\dot\H^{-1}(\T^3)} +   \int_0^T |I_3| \dd t \right).
        \end{aligned}
    \end{equation*}
\end{proposition}


\begin{preuve}Clearly, we have
\begin{align*}
|I_1(t)| &\lesssim \|\na_x u(t) \|_{\L^\infty( \T^3)} \int_{\T^3\times\R^3}|v-u(t,x)| ^2\feps(t,x,v)\dd x\dd v \\
&\lesssim \|\na_x u(t) \|_{\L^\infty( \T^3)} \He(t), \\
|I_2(t)| &\lesssim \|\na_x u(t) \|_{\L^\infty( \T^3)} \int_{\T^3}|\ueps-u(t,x)| ^2 \dd x \\
&\lesssim \|\na_x u(t) \|_{\L^\infty( \T^3)} \He(t)
\end{align*}
On the other hand, applying Lemma~\ref{lem-rho-H-1} and Young's inequality, we obtain for all $t \in (0,T)$, for any $a>0$,
\begin{align*}
    \int_0^t |I_4|(s)  \dd s
    &\le\int_0^t\norme{\rhoeps(s)-\rho(s)}_{\dot\H^{-1}(\T^3)}\norme{\na_x((\ueps(s)-u(s))\cdot G(s))}_{\L^2(\T^3)}\dd s\\
    &\le 2aT\| \rhoeps(0) - \rho(0)\|^2_{\dot\H^{-1}(\T^3)}
    +2aT^2\eps^2\int_0^t\norme{F_\eps(s)}_{\L^2(\T^3)}^2\dd s\\
    &+\left( 2aT^2+
    +\frac{\norme{\na_x G}_{\L^\infty(\T^3\times\R^3)}^2}{4a}\right)\int_0^t\norme{\ueps(s)-u(s)}_{\L^2(\T^3)}^2\dd s \\
    &+  \frac{\norme{G}_{\L^\infty(\T^3\times\R^3)}^2}{4a}\int_0^t\norme{\na_x(\ueps(s)-u(s))}_{\L^2(\T^3)}^2\dd s.
\end{align*}
Thanks to the Cauchy-Schwarz inequality, we have
\begin{align*}
    \eps^2\int_0^t\norme{F_\eps(s)}_{\L^2(\T^3)}^2&\dd s \\
    &\le\norme{\rhoeps}_{\L^\infty((0,T)\times\T^3)}\int_0^t\int_{\T^3\times\R^3}\feps(s,x,v)|v-\ueps(s,x)|^2\dd x\dd v\dd s\\
    &\le C\int_0^t\int_{\T^3\times\R^3}\feps(s,x,v)|v-\ueps(s,x)|^2\dd x\dd v\dd s.
\end{align*}
The conclusion thus follows from Lemma~\ref{lem-evolution-HH} and Gronwall's lemma, by taking $a=\norme{G}_{\L^\infty(\T^3\times\R^3)}^2$ and $\eps_0$ sufficiently small.

\end{preuve}

We are now ready to state quantitative convergence results. We need only to verify that the bound~\eqref{eq-condi-rho} holds. To do so, we rely on the analysis of the previous sections.

\begin{theorem}\label{thCVprecis-3}
    Let $T>0$. Under Assumptions~\ref{hypGeneral}--\ref{hypSmallDataNS}, if we assume that
    \[
        \norme{G}_{\L^\infty((0,T)\times\T^3)}+\norme{\na_xG}_{\L^\infty((0,T)\times\T^3)}\lesssim1,
    \]
    there exists $\eps_0>0$ and $\eta>0$ such that if, for all $\eps\in(0,\eps_0)$, \eqref{FineHypf0}--\eqref{hypSmallModulatedEnergy-3} hold
    and if
    \[
        \ueps\xrightarrow[\eps\to0]{} u^0\text{ in }\L^2(\T^3)\qquad\text{and}\qquad\rhoeps^0\xrightarrow[\eps\to0]{}\rho^0\text{ in }\dot{\H}^{-1}(\T^3)
    \]
    then $(\ueps)_{\eps>0}$ converges to $u$ in $\L^\infty(0,T;\L^2(\T^3))$ and $(\rhoeps)_{\eps>0}$ converges to $\rho$ in $\L^\infty(0,T;\dot{\H}^{-1}(\T^3))$. Furthermore, for every $t\in(0,T)$
    \begin{multline}
    \label{urhoconvfine}
        \norme{\ueps(t)-u(t)}_{\L^2(\T^3)}+\norme{\rhoeps(t)-\rho(t)}_{\dot{\H}^{-1}(\T^3)}\\
        \lesssim_{T,M,G}\norme{\ueps^0-u^0}_{\L^2(\T^3)}+\norme{\rhoeps^0-\rho^0}_{\dot\H^{-1}(\T^3)}+\eps^{\frac{1}{2}}.
    \end{multline}
\end{theorem}

\begin{preuve}
    Under this set of assumptions, we have proven all times are strong existence times, see Section~\ref{sec-bootstrap-fine}. Therefore, Lemma~\ref{RhoLinfini} applies and ensures that the bound \eqref{eq-condi-rho} holds.
    
    We use the Cauchy-Schwarz inequality and Lemmas~\ref{RhoLinfini},~\ref{HighOrderD2} and~\ref{lemmePara2-3} to find that
    \begin{equation*}
        \int_0^T|I_3(t)|\dd t\le\norme{\rhoeps}_{\L^\infty((0,T)\times\T^3)}^{\frac{1}{2}}\left(\eps^2\int_0^T\D_\eps^{(2)}(t)\dd t\right)^{\frac{1}{2}}\norme{G}_{\L^2((0,T)\times\T^3)}
        \lesssim_{T,M,G} \eps.
    \end{equation*}
    Furthermore, thanks to Hölder's inequality
    \[
        \int_{\T^3\times\R^3}\feps^0|v-\ueps^0|^2\lesssim_M\left(\int_{\T^3\times\R^3}\feps^0|v-\ueps^0|^p\right)^{\frac{2}{p}}\lesssim_M\eps^{2-\frac{2}{p}}.
    \]
    Thus, applying Proposition~\ref{prop-relativeentropy}, for every $t\in[0,T]$,
    \[
        \He(t)\lesssim_{T,M,G}\norme{\ueps^0-u^0}_{\L^2(\T^3)}^2
            +\norme{\rhoeps^0-\rho^0}_{\dot\H^{-1}(\T^3)}^2+\eps,
    \]
    hence the result for $\norme{\ueps(t)-u(t)}_{\L^2(\T^3)}$.
    
    Recall that since $T$ is a strong existence time,
    \[
        \norme{F_\eps}_{\L^2((0,T)\times\T^3)}\le\frac{C^*}{2}.
    \]
    Therefore, Lemma~\ref{lem-rho-H-1} yields
    \begin{align*}
        \norme{\rhoeps(t)-\rho(t)}_{\dot\H^{-1}(\T^3)}
        &\lesssim_T\norme{\rhoeps^0-\rho^0}_{\dot\H^{-1}(\T^3)}
            +\eps\norme{F_\eps}_{\L^2((0,T)\times\T^3)}\\
            &\qquad +\norme{\ueps-u}_{\L^1(0,T;\L^2(\T^3))}\\
        &\lesssim_{T,M,G}\norme{\rhoeps^0-\rho^0}_{\dot\H^{-1}(\T^3)}+\norme{\ueps^0-u^0}_{\L^2(\T^3)}+\eps^{\frac{1}{2}},
    \end{align*}
    hence the theorem.
\end{preuve}

To conclude, Theorem~\ref{thCVprecis-3} combined with Corollary~\ref{CorollaryFineFirstConv} allows us to derive the quantitative convergence of $(f_\eps)$.

\begin{remark}
Another benefit of the introduction of higher dissipation functionals appears here. To estimate the contribution of $I_3$, one can also use the energy--dissipation estimate. However this would yield $\eps^{\frac{1}{4}}$ instead of $\eps^{\frac{1}{2}}$ as in~\eqref{urhoconvfine}.
\end{remark}

\begin{remark}
    Note that we could have kept track of the time dependence and provided  explicit quantitative results, as we have done in Section~\ref{SectionLandLFParticle}.
\end{remark}


\subsection{Convergence in the mildly well-prepared case}
\label{sec-smallT-fine}

In the previous subsections, we derived convergence results for any time horizon under 
\begin{itemize}
   \item a smallness assumption for the initial fluid velocity in the sense of Assumption~\ref{hypSmallDataNS},
       \item a smallness assumption for the initial modulated energy in the sense of~\eqref{hypSmallModulatedEnergy-3},
    \item a smallness assumption for the initial distribution function in the sense of~\eqref{FineHypf0},
    \item well-preparedness for the initial data in the sense of~\eqref{wellPrep-3}.
\end{itemize}
In this section, we show how to obtain short-time convergence results when dispensing with Assumptions~\ref{hypSmallDataNS} and~\eqref{hypSmallModulatedEnergy-3}, which corresponds to the mildly well-prepared case of Theorem~\ref{thm3}.

According to Lemma~\ref{ExistsStrongTime}, there exists a strong existence time $T_{M,\eps}>0$. Let us prove that we can find a strong existence time that does not depend on $\eps$. We will do so using a bootstrap argument similar to the one used in the previous sections. We define
\[
    T_{\eps}^*=\sup\left\{T>0,\,T\text{ is a strong existence time}\right\}.
\]
Our goal is to find a time $T_M$, independent of $\eps$, such that $T_M<T_\eps^*$. As we have proven in Section~\ref{sec-smallT} in the \emph{light} and \emph{light and fast} particle regimes, there exists $T_M>0$, independent of $\eps$, such that for every $\eps\in(0,1)$,
\[
    \int_0^{T_M}\norme{e^{t\Delta}\ueps^0}_{\dot\H^1(\T^3)}^4\dd t\le\frac{C^*}{4}.
\]
For what concerns $\norme{F_\eps}_{\L^2((0,T)\times \T^3)}$, we can  apply Lemma~\ref{lemmePara2-3} and its proof to ensure that there exists $\eps_0>0,\,\eta>0$ and $T_M>0$ (independent of $\eps$) such that if~\eqref{FineHypf0}--\eqref{wellPrep-3} hold, then for every $T<\min(T_\eps^*,T_M)$,
\[
    \norme{F_\eps}_{\L^2((0,T)\times\T^3)}\le\frac{C^*}{4}.
\]
There only remains to obtain a satisfactory bound for $\norme{\na_x\ueps}_{\L^1(0,T;\L^\infty(\T^3))}$. To achieve this, we apply Lemmas~\ref{estimNabla} and~\ref{lemmePara-3}, which yield
\[
    \norme{\na_x\ueps}_{\L^1(0,T;\L^\infty(T^3))}\lesssim T^{1-\beta_p}\EEE_{\eps}(0)^{\frac{1-\beta_p}{2}}\norme{\Delta_x\ueps}_{\L^p((0,T)\times T^3)}^{\beta_p}\lesssim T^{1-\beta_p}M^{\omega_p},
\]
for some $\omega_p>0$. The right-hand side can be made as mall as necessary (e.g. <1/40) by reducing the value of $T_M>0$, still independently from $\eps$.

This concludes the bootstrap argument, and we can now perform the same analysis as in Subsections~\ref{sec-first-fine} and~\ref{sec-relat-fine} and obtain the same convergence results when we replace Assumptions~\ref{hypSmallDataNS} and~\eqref{hypSmallModulatedEnergy-3} by the time horizon constraint $T<T_M$.


\section{Back to the \emph{light} and \emph{light and fast} particle regimes}
\label{sec-betterLF}

In this section, we come back to the \emph{light} and \emph{light and fast} particle regimes and the improve the results of Section~\ref{SectionLandLFParticle}, precisely by removing the well-preparedness assumption~\eqref{eq-WP-LF}.

 To this purpose we modify the main strategy by now studying estimates in $\L^{\frac{p}{p-1}}(0,T; \L^p(\T^3))$ (as is now usual, $p$ is the regularity index of Assumption~\ref{hypGeneral}).
It turns out that this refined approach is also efficient to lower down the required regularity for the initial fluid velocity $u^0$. Namely we may require an uniform bound in $\B^{s-\frac{2 (p-1)}{p},p}_p$ instead of $\B^{s-\frac{2}{p},p}_p$ in Assumption~\ref{hypGeneral}.

\begin{itemize}
\item In Section~\ref{sec-LF-nohyp}, we provide the key $\L^{\frac{p}{p-1}} \L^p$ estimates for the convection and Brinkman force, which pave the way for the improved $\L^{\frac{p}{p-1}} \L^p$ parabolic estimates of 
Section~\ref{sec-conv12-re}. This new strategy then allows for the claimed convergence results of Theorems~\ref{thm1} and~\ref{thm2} under optimal assumptions on the initial conditions.

 
\item In  Section~\ref{SectionRetour-12}, we explain how the objects introduced and the methods developed to study the fine particle regime in Section~\ref{SectionLandLFParticle} can be used as well for the \emph{light} and \emph{light and fast} particle regimes. This yields alternative proofs of the results of Section~\ref{SectionLandLFParticle}.

\end{itemize}

\subsection{Uniform \texorpdfstring{$\L^{\frac{p}{p-1}} \L^p$}{L\string^\{p/(p-1)\}L\string^p} estimates}\label{sec-LF-nohyp}

In this section, we shall establish
$\L^{\frac{p}{p-1}} \L^p$ estimates for the convective term and the Brinkman force in the Navier-Stokes equation. 

The  first lemma, which concerns the convective term $\ued\cdot\na_x\ued$ is a straightforward adaptation of Lemma~\ref{EstimConvectStrongT}. 

\begin{lemma}\label{EstimConvectStrongT-re}
    Under Assumption~\ref{hypGeneral}, there exists $\nu_p>0$ such that if $T>0$ is a strong existence time, then
    \begin{equation*}
        \norme{\left(\ued\cdot\na_x\right)\ued}_{\L^{\frac{p}{p-1}}(0,T;\L^p(\T^3))}
        \lesssim \Psi_{\ed,0}^{\frac{1}{2}}\EEE_{\egd}(0)^{\frac{1-\beta_p}{2}}\norme{\Delta_x\ued}_{\L^{\frac{p}{p-1}}(0,T;\L^p(\T^3))}^{\beta_p}.
    \end{equation*}
\end{lemma}

The three following lemmas are adaptations of Lemmas~\ref{lemmaF0-12}, \ref{lemmaFdt-12} and \ref{lemmaFdx-12} in the $\L^{\frac{p}{p-1}}\L^p$ framework.

The main benefit of looking at $\L^{\frac{p}{p-1}}\L^p$ estimates appears in the next lemma. Note that  the right-hand side in~\eqref{eq-better} is always vanishing as $\eps \to 0$, contrary to Lemma~\ref{lemmaF0-12} in which the extra assumption~\eqref{hyp-wellprep-2} was needed.

\begin{lemma}\label{lemmaF0-12-re}
    Under Assumption~\ref{hypGeneral}, there exists $\mu_p>0$ such that for every $\eps>0$, if $T>0$ is a strong existence time,
    \begin{equation}
    \label{eq-better}
        \norme{\Fed^0}_{\L^{\frac{p}{p-1}}(0,T;\L^p(\T^3))}\lesssim M^{\mu_p}\frac{\eps^{1-\frac{1}{p}}}{\sigma}.
    \end{equation}
\end{lemma}

\begin{preuve}
This is a straightforward adaptation of the proof of Lemma~\ref{lemmaF0-12}
\end{preuve}

\begin{lemma}\label{lemmaFdt-12-re}
    Under Assumption~\ref{hypGeneral}, for every $\eps>0$, if $T>0$ is a strong existence time,
    \[
        \norme{\Fed^{dt}}_{\L^{\frac{p}{p-1}}(0,T;\L^p(\T^3))}\lesssim \eps M\norme{\partial_t\ued}_{\L^{\frac{p}{p-1}}(0,T;\L^p(\T^3))}
    \]
\end{lemma}

\begin{preuve}The beginning of the proof is similar to that of Lemma~\ref{lemmaFdt-12}, except that we use Minkowski's integral inequality~\cite[Theorem 202]{hardy-littlewood-polya} instead of H\"older's inequality. We obtain
  \begin{align*}
        &\norme{\Fed^{dt}}_{\L^{\frac{p}{p-1}}(0,T;\L^p(\T^3))}^{\frac{p}{p-1}}\\
        &\quad \lesssim\norme{\fed^0}_{\L^1(\R^3;\L^\infty(\T^3))}^{\frac{p}{p-1}}\int_0^T \left(\int_0^t\norme{\partial_s\ued(s)}_{\L^p(\T^3)} e^{\frac{s-t}{\eps}}\dd s\right)^{\frac{p}{p-1}} \dd t\\
         &\quad\lesssim \norme{\fed^0}_{\L^1(\R^3;\L^\infty(\T^3))}^{\frac{p}{p-1}} \eps^{\frac{1}{p-1}}  \int_0^T\norme{\partial_s\ued(s)}_{\L^p(\T^3)}^{\frac{p}{p-1}} \left(\int_s^T e^{\frac{s-t}{\eps}} \dd t\right) \dd s\\
         &\quad \lesssim \norme{\fed^0}_{\L^1(\R^3;\L^\infty(\T^3))}^{\frac{p}{p-1}} \eps^{\frac{p}{p-1}}  \norme{\partial_t\ued}_{\L^{\frac{p}{p-1}}(0,T;\L^p(\T^3))}^{\frac{p}{p-1}},
    \end{align*}
    where we have also applied Jensen's inequality 
    in the antepenultimate line.
\end{preuve}

\begin{lemma}\label{lemmaFdx-12-re}
    Under Assumption~\ref{hypGeneral}, there exists $\mu_p\ge0$ such that for every $\eps>0$, if $T>0$ is a strong existence time,
    \[
        \norme{\Fed^{dx}}_{\L^{\frac{p}{p-1}}(0,T;\L^p(\T^3))}\lesssim\eps\norme{\Delta_x\ued}_{\L^{\frac{p}{p-1}}(0,T;\L^p(\T^3))}+\eps M^{\mu_p}\EEE_{\ed}(0)^{\frac{1}{2}}.
    \]
\end{lemma}

\begin{preuve}
As in the beginning of the proof of Lemma~\ref{lemmaFdx-12}, we have
    \begin{equation*}
            |\Fed^{dx}(t,x)|\leq I_1(t,x)+I_2(t,x),
    \end{equation*}
    with
    \begin{align*}
    I_1 &=\int_{\R^3}e^{-\frac{t}{\eps}}\fed^0\left(\tilde{X}_{\ed}^{t,x,w}(0),w\right)|w|\int_0^t|\na_x\ued\left(s,\tilde{X}_{\ed}^{t,x,w}(s)\right)|\dd s\dd w ,\\
    I_2 &= \frac{\sigma}{\eps}\int_{\R^3}\fed^0\left(\tilde{X}_{\ed}^{t,x,w}(0),w\right)\\
    &\qquad \qquad \int_0^t\int_0^se^{\frac{s-t}{\eps}}e^{\frac{\tau-s}{\eps}}\left|\ued\left(\tau,\tilde{X}_{\ed}^{t,x,w}(\tau)\right)\right|\left|\nabla_{x}\ued\left(s,\tilde{X}_{\ed}^{t,x,w}(s)\right)\right|
            \dd\tau\dd s\dd w.
    \end{align*}
    The term $I_1$ is treated as in the proof of Lemma~\ref{lemmaFdx-12}, except that we apply here  Minkowski's integral inequality as in the proof of Lemma~\ref{lemmaFdt-12-re}. The outcome is the estimate
    \begin{align*}
    \norme{I_1}_{\L^{\frac{p}{p-1}}(0,T;\L^p(\T^3))}\lesssim \eps\norme{|v|\fed^0}_{\L^1(\R^3;\L^\infty(\T^3))}^{\frac{1}{1-\alpha_p}}\EEE_{\ed}(0)^{\frac{1}{2}}\\
            +\eps\norme{\Delta_{x}\ued}_{\L^{\frac{p}{p-1}}(0,T;\L^p(\T^3))}.
    \end{align*}
    The term $I_2$ is also treated as in Lemma~\ref{lemmaFdx-12}. We obtain
       \[
        \int_0^T\left(\int_{\T^3}|I_2(t,x)|^p\dd x\right)^{\frac{1}{p-1}}\dd t\lesssim\eps^{-1}\norme{\fed^0}_{\L^1(\R^3;\L^\infty(\T^3))}\left(\sup_{t\in(0,T)}I_3(t)\right)^{\frac{1}{p-1}}\times \tilde{I}_4,
    \]
    where
    \begin{multline*}
        I_3(t)=\int_0^te^{\frac{p(s-t)}{2\eps}}\int_0^se^{\frac{p(\tau-s)}{2\eps}}\\
            \int_{\T^3\times\R^3}\fed^0\left(\tilde{X}_{\ed}^{t,x,w}(0),w\right)\left|\ued\left(\tau,\tilde{X}_{\ed}^{t,x,w}(\tau)\right)\right|^p\dd x\dd w\dd\tau\dd s,
    \end{multline*}
    and
    \[
        \tilde{I}_4:=\int_0^T\left(\int_0^t\norme{\na_x\ued(s)}_{\L^\infty(\T^3)}^{\frac{p}{p-1}}e^{\frac{p(s-t)}{2(p-1)\eps}}\dd s\right)\dd t.
    \]
    We finally conclude that
    \begin{align*}
        &\int^T\left(\int_{\T^3}|I_2(t,x)|^p\dd x\right)^{\frac{1}{p-1}}\dd t\\
        &\qquad\lesssim\eps^{\frac{p}{p-1}}\norme{\fed^0}_{\L^1(\R^3;\L^\infty(\T^3))}^{\frac{p}{p-1}}\Psi_{\ed,0}^{\frac{1}{2}\frac{p}{p-1}}\EEE_{\ed}(0)^{\frac{(1-\beta_p)}{2}\frac{p}{p-1}}\norme{\Delta_x\ued}_{\L^{\frac{p}{p-1}}(0,T;\L^p(\T^3))}^{\beta_p \frac{p}{p-1}}\\
        &\qquad \lesssim\eps^{\frac{p}{p-1}}\norme{\fed^0}_{\L^1(\R^3;\L^\infty(\T^3))}^{\frac{p}{p-1}\frac{1}{1-\beta_p}}\Psi_{\ed,0}^{\frac{p}{p-1}\frac{1}{2(1-\beta_p)}}\EEE_{\ed}(0)^{\frac{1}{2}\frac{p}{p-1}} \\
        &\qquad \qquad  +\eps^{\frac{p}{p-1}}\norme{\Delta_x\ued}_{\L^{\frac{p}{p-1}}(0,T;\L^p(\T^3))}^{\frac{p}{p-1}}.
    \end{align*}
    The proof of the lemma is finally complete.
\end{preuve}

\subsection{Conclusion of the bootstrap and proof of convergence under optimal assumptions}
\label{sec-conv12-re}
Equipped with the uniform $\L^{\frac{p}{p-1}} \L^p$ estimates of the previous section, we are in position to conclude the bootstrap argument and eventually obtain the convergence towards smooth solutions of the Transport-Navier-Stokes equations~\eqref{TNS}. As the concluding arguments are very close to those of Sections~\ref{sec-conclubootstrapLLF}, \ref{SectionAsymptotic-12} and~\ref{sec-smallT}, we will provide very few details.

\subsubsection{Mildly well-prepared and well-prepared cases}

\begin{lemma}\label{lemmePara-12-re}
    Under Assumption~\ref{hypGeneral}, there exists $\eps_0>0$ 
    such that for all $\eps\in(0,\eps_0)$, if $T>0$ is a strong existence time then
    \[
        \norme{\partial_t\ued}_{\L^{\frac{p}{p-1}}(0,T;\L^p(\T^3))}+\norme{\Delta_x\ued}_{\L^{\frac{p}{p-1}}(0,T;\L^p(\T^3))}\lesssim M^{\omega_p},
    \]
    for some $\omega_p>0$.
\end{lemma}

\begin{preuve}
This is a consequence of Lemmas~\ref{lem-decompF}, \ref{lemmaF0-12-re}, \ref{lemmaFdt-12-re}, \ref{lemmaFdx-12-re} and \ref{EstimConvectStrongT-re},
and of the parabolic estimates of Theorem~\ref{ParabolicReg}.
\end{preuve}

\begin{corollary}
    Under Assumption~\ref{hypGeneral}, there exist $\eps_0>0$ and $\eta>0$ 
    such that, for all $\eps\in(0,\eps_0)$,
    if
    \begin{equation}\label{hypEnergieModuleePetite-12-re}
        \EEE_{\ed}(0)\le\eta,
    \end{equation}
    then for any $T<T^*_\eps$,
    \[
        \norme{\na\ued}_{\L^1(0,T;\L^\infty(\T^3))}\le\frac{1}{40}.
    \]
\end{corollary}

\begin{preuve}
This is a small variant of the proof of Corollary~\ref{coro-nau-12}.
    Since $p>3$, we can apply Lemmas~\ref{estimNabla} and~\ref{lemmePara-12-re} and find that there exist $\eps_0$ 
    such that for all $\eps\in(0,\eps_0)$
    , for some $\omega_p>0$, we have
    \begin{align*}
        \int_0^T\norme{\na\ued(t)}_{\L^\infty(\T^3)}\dd t&\lesssim\EEE_{\ed}(0)^{\frac{1-\beta_p}{2}}\norme{\Delta_x\ued}_{\L^{\frac{p}{p-1}}(0,T;\L^p(\T^3))}^{\beta_p }\\
        &\lesssim \EEE_{\ed}(0)^{\frac{1-\beta_p}{2}}M^{\omega_p\beta_p},
    \end{align*}
    which can be made as small as necessary by reducing the value of $\eta>0$.
\end{preuve}

Combining with Lemma~\ref{firstproof}, we conclude as in Section~\ref{sec-conclubootstrapLLF} that $T^\eps = +\infty$.

We have the following variant of Lemma~\ref{lemmaF0-12-re}.
\begin{lemma}\label{lemmaF0-12-L2}
Let $q \geq 1$.
    Under Assumption~\ref{hypGeneral}, there exists $\mu>0$ such that for every $\eps>0$, if $T>0$ is a strong existence time,
    \[
        \norme{\Fed^0}_{\L^{q}(0,T;\L^2(\T^3))}\lesssim M^{\mu} \eps^{\frac{1}{q}} \left(\int_{\T^3 \times \R^3} f^0_\ed \frac{|v|}{\sigma} \dd x \dd v +1\right).
    \]
\end{lemma}

We finally get the following improvement of Theorems~\ref{thCV-12} and~\ref{thCVprecis-12}.

\begin{theorem}
\label{improLF}
\begin{itemize}
\item Under the sole Assumptions~\ref{hypGeneral} and~\eqref{hypEnergieModuleePetite-12-re},
  if
    \[
        \ued\xrightharpoonup[\eps\to0]{} u^0\text{ in }w\text{-}\L^2(\T^3)\qquad\text{and}\qquad \rhoed^0\xrightharpoonup[\eps\to0]{} \rho^0\text{ in }w^*\text{-}\L^\infty(\T^3),
    \]
    then, for any $T>0$, $(\ueps)_{\eps>0}$ converges to $u$ in $\L^2((0,T) \times \T^3)$, $(\rhoeps)_{\eps>0}$ converges weakly-$*$ to $\rho$ in $\L^\infty((0,T)\times\T^3)$, where $(\rho,u)$ satisfies~\eqref{limitRhoU-12}.

\item  Additionally, if
\begin{equation}
\label{ass-L2}
\eps^{\frac{1}{2}}\int_{\T^3 \times \R^3} f^0_\ed \frac{|v|}{\sigma} \dd x \dd v \xrightarrow[\eps\to0]{} 0,
\end{equation}
for all $T>1$, for all $t\in[0,T]$,
        \begin{multline*}
 \norme{\ued(t)-u(t)}_{\L^2(\T^3)}\\ 
        \lesssim e^{CM^2T}M^{\omega_p}\left(\norme{\ued^0-u^0}_{\L^2(\T^3)}+\eps^{\frac{1}{2}}\int_{\T^3 \times \R^3} f^0_\ed \frac{|v|}{\sigma} \dd x \dd v\right),
        \end{multline*}
        \begin{multline*}
W_1(\rhoed(t),\rho(t)) \\
        \lesssim W_1(\rhoed^0,\rho^0)
            +T^{1/2}e^{CM^2T}M^{\omega_p}\left(\norme{\ued^0-u^0}_{\L^2(\T^3)}+\eps^{\frac{1}{2}}\int_{\T^3 \times \R^3} f^0_\ed \frac{|v|}{\sigma} \dd x \dd v\right),
        \end{multline*}
and
    \begin{multline*}
        \int_0^TW_1(\fed(t),\rho(t)\otimes\delta_{v=\sigma u(t)})\dd t\\
        \lesssim T W_1(\rhoed^0,\rho^0)+ T^{\frac{3}{2}}e^{CM^2T}M^{\omega_p}\left(\norme{\ued^0-u^0}_{\L^2(\T^3)}+\eps^{\frac{1}{2}}\int_{\T^3 \times \R^3} f^0_\ed \frac{|v|}{\sigma} \dd x \dd v\right),
    \end{multline*}
    where $(\rho,u)$ satisfies~\eqref{limitRhoU-12}.
\end{itemize}
\end{theorem}

\begin{remark}
For Theorem~\ref{improLF} to be significant, we need to ensure that 
$$
\eps^{\frac{1}{2}}\int_{\T^3 \times \R^3} f^0_\ed \frac{|v|}{\sigma} \dd x \dd v \xrightarrow[\eps\to 0]{} 0,
$$
which is restrictive only for $\alpha=1/2$. This accounts for the statement of Theorem~\ref{thm2}.
\end{remark}

\begin{remark}
As in Section~\ref{SectionLandLFParticle}, we assume $T>1$ to simplify the expression of the convergence rates, though we would obtain the same results, in terms of $\eps$, for any $T>0$.
\end{remark}

\begin{preuve}
For the first item, we use Lemma~\ref{lemmaF0-12-L2} for $q<2$ and Lemmas~\ref{lemmaFdt-12-re}, \ref{lemmaFdx-12-re} and~\ref{lemmePara-12-re} to ensure that
$F_\ed \xrightarrow[\eps\to 0]{} 0$ in ${\L^{q}(0,T;\L^2(\T^3))}$. The arguments of the proof of Theorem~\ref{thCV-12} then apply without significant modification.

For the second item, we use Lemmas~\ref{lemmaFL2-12} and~\ref{lemmePara2-12},
together with Lemma~\ref{lemmaF0-12-L2} for $q=2$ and~\eqref{ass-L2}. We get the bound
\[
\norme{\Fed}_{\L^2((0,T)\times\T^3)}\lesssim\eps M^\mu_p
    +\eps^{1/2- \alpha } \int_{\T^3\times\R^3}|v|^2f_{\eps,\sigma}^0(x,v),
\]
for some $\mu_p>0$. 
Again, the arguments of the proof of Theorem~\ref{thCV-12} apply \emph{mutatis mutandis}.
\end{preuve}

We finally state the following extension of Corollary~\ref{coro-longtime-LLF}.

\begin{corollary}
Corollary~\ref{coro-longtime-LLF} holds in the light and fast particle regime for all $\alpha \in (0,1/2)$.
\end{corollary}

\subsubsection{General case}

The arguments of Section~\ref{sec-smallT} apply \emph{mutantis mutandis}. The only noticeable difference is the following.
As in the proof of Theorem~\ref{improLF}, if $T>0$ is a strong existence time, we obtain the estimate
\[
\norme{\Fed}_{\L^2((0,T)\times\T^3)}\lesssim\eps M^\mu_p
    +\eps^{1/2- \alpha } \int_{\T^3\times\R^3}|v|^2f_{\eps,\sigma}^0(x,v),
\]
Therefore, to ensure $\norme{\Fed}_{\L^2((0,T)\times\T^3)} \leq C^\star/4$, in the critical case $\alpha=1/2$,
the assumption 
$$
\int_{\T^3\times\R^3}|v|^2f_{\eps,\eps^{1/2}}^0(x,v)\dd x\dd v\leq \eta,
$$
is required, as stated in Theorem~\ref{thm2}.


\subsection{Higher dissipation and relative entropy in the \emph{light} and \emph{light and fast} particle regimes}\label{SectionRetour-12}

In this section, we apply the higher dissipation and relative entropy functions (introduced for the fine particle regime in Section~\ref{SectionLandLFParticle}) to obtain alternative proofs for the \emph{light} and \emph{light and fast} particle regimes. 

\subsubsection{Higher dissipation}

Let us start by introducing the higher dissipation functionals which are analogous to those presented in Section \ref{SubsectionHigherOrder}; we then show how to recover the convergence results obtained in Section~\ref{SectionAsymptotic-12}.

\begin{definition}
    Let $r\ge2$. The higher (fluid-kinetic) dissipation (of order $q$) is defined, for almost every $t\ge0$, by
    \[
        \D_{\ed}(t)=\int_{\T^3\times\R^3}\fed\frac{|\frac{v}{\sigma}-\ued(t,x)|^r}{\eps^r}\dd x\dd v.
    \]
\end{definition}

Thanks to Hölder's inequality, we can relate $\D_{\ed}^{(r)}$ to $\Fed$ as follows.

\begin{lemma}
    For $r\ge2$, we have
    \[
        \norme{\Fed}_{\L^r((0,T)\times\T^3)}\le\eps\norme{\fed^0}_{\L^1(\R^3;\L^\infty(\T^3))}^{1-\frac{1}{r}}\left(\int_0^T\D_{\ed}^{(r)}(t)\dd t\right)^{\frac{1}{r}}.
    \]
\end{lemma}

We proceed as in Section~\ref{SubsectionHigherOrder} to desingularize $\D_{\ed}^{(r)}$ and obtain the following estimate.

\begin{lemma}
\label{lem-higher12}
    The higher fluid-kinetic dissipation of order $r\ge2$ satisfies
    \begin{multline*}
        \int_{\T^3\times\R^3}\fed(T)\frac{|\frac{v}{\sigma}-\ued(T)|^r}{\eps^{r-1}}\dd x\dd v+ \int_0^T\D_{\ed}^{(r)}(t)\dd t\\
        \lesssim\norme{\fed^0}_{\L^1(\R^3;\L^\infty(\T^3))}\norme{\partial_t\ued}_{\L^r((0,T)\times\T^3)}\\
            +\frac{1}{\sigma^r}\norme{|\na_x\ued|m_r^{\frac{1}{r}}}_{\L^r((0,T)\times\T^3)}^r+\int_{\T^3\times\R^3}\fed^0\frac{|\frac{v}{\sigma}-\ued^0|^r}{\eps^{r-1}}\dd x\dd v,
    \end{multline*}
    where $m_r(t,x)=\int_{\R^3}|v|^r\fed(t,x,v)\dd v$.
\end{lemma}

As in Section~\ref{SubsectionHigherOrder}, we apply this identity with $r=p$, the regularity index from Assumption~\ref{hypGeneral} and can prove the following estimate on the middle term.

\begin{lemma}
    Under Assumptions~\ref{hypGeneral}--\ref{hypSmallDataNS}, if $T>0$ is a strong existence time, we have
    \begin{multline*}
        \norme{|\na_x\ued|m_{p,\eps}^{\frac{1}{p}}}_{\L^p((0,T)\times\T^3)}^p\lesssim\eps^{1-\alpha_p}\norme{\fed^0|v|^p}_{\L^1(\R^3;\L^\infty(\T^3))}^{\frac{1}{1-\alpha_p}}\EEE_{\ed}(0)^{\frac{p}{2}}\\
        +\sigma^p\norme{\fed^0}_{\L^1(\R^3;\L^\infty(\T^3))}\EEE_{\ed}(0)^{\frac{(1-\beta_p)p}{2}}\Psi_{\ed,0}^{\frac{p}{2(1-\beta_p)}}\\
        +\left(\eps^{1-\alpha_p}+\sigma^p\norme{\fed^0}_{\L^1(\R^3;\L^\infty(\T^3))}\EEE_{\ed}(0)^{\frac{(1-\beta_p)p}{2}}\right)\norme{\Delta_x\ued}_{\L^p((0,T)\times\T^3)}^p.
    \end{multline*}
\end{lemma}

Therefore, combining the previous statements, we obtain the following estimate on the Brinkman force.

\begin{lemma}
    Under Assumptions~\ref{hypGeneral}--\ref{hypSmallDataNS}, if $T>0$ is a strong existence time the Brinkman force satisfies
    \begin{align*}
        \norme{\Fed}_{\L^p((0,T)\times\T^3)}&\lesssim\eps M\norme{\partial_t\ued}_{\L^p((0,T)\times\T^3)}
            +\eps M^{\omega_p}\norme{\Delta_x\ued}_{\L^p((0,T)\times\T^3)}\\
            &+\eps M^{\omega_p}
            +\eps^{\frac{1}{p}}M^{\omega_p}\left(\int_{\T^3\times\R^3}\fed^0(x,v)\left|\frac{v}{\sigma}-\ued^0(x)\right|^p\dd x\dd v\right)^{\frac{1}{p}}.
    \end{align*}
\end{lemma}

Equipped with this estimate, we can apply the same strategy as in Section~\ref{SectionLandLFParticle} and obtain quantitative convergence results similar to those of Theorems~\ref{thCVprecis-12} and~\ref{corollaryCVf-12}. 
This approach straightforwardly yields the following result. Thanks to Lemma~\ref{lem-higher12}, we directly obtain that for all $T>0$,
\[
 \int_{\T^3\times\R^3}\fed(T)\frac{|\frac{v}{\sigma}-\ued(T)|^p}{\eps^{p-1}}\dd x\dd v \xrightarrow[\eps \to 0]{} 0
 \]
 with a rate of convergence in $\eps$, assuming only that
 \[
    \int_{\T^3\times\R^3}\fed^0\frac{|\frac{v}{\sigma}-\ued^0|^p}{\eps^{p-1}}\dd x\dd v \xrightarrow[\eps \to 0]{} 0,
\]
and without assuming any convergence for the initial fluid velocity $\ued^0$.

\subsubsection{Relative entropy}

In this paragraph, we use the relative entropy method introduced in Section~\ref{sec-relat-fine} to prove the quantitative convergence obtained in Theorem~\ref{thCVprecis-12} with another point of view. 

We introduce the relative entropy
\begin{multline*}
    \Hed(t) :=\frac{\eps}{2}\int_{\T^3\times\R^3}\fed(t,x,v)\left|\frac{v}{\sigma}-u(t,x)\right|^2\dd x\dd v\\
                +   \frac{1}{2}\int_{\T^3}|\ued(t,x)-u(t,x)|^2\dd x,
\end{multline*}
where $u$ is a smooth solution to the  incompressible Navier-Stokes equations.

Following the proof of Lemma~\ref{lem-evolution-HH}, we obtain the following result.
\begin{lemma}\label{lem-evolution-Hj}
    The relative entropy satisfies for all $t \geq 0$
    \begin{equation*}
    \begin{aligned}
        \Hed(t) &+ \int_0^t \int_{\T^3}|\na_x(\ued - u )|^2\dd x \dd s +  \int_0^t \int_{\T^3 \times \R^3}\left|\frac{v}{\sigma}-\ued\right|^2\fed \dd x\dd v \dd s \\
                &\leq \Hed(0) + \int_0^t  \sum_{i=1}^{4} J_i(s)  \dd s.
    \end{aligned}
    \end{equation*}
    with
    \begin{align*}
        &J_1 := - {\eps}\int_{\T^3 \times \R^3}   \fed \left(\frac{v}{\sigma}-u\right) \otimes \left(\frac{v}{\sigma}-u\right) : \na_x u \dd x \dd v, \\
        &J_2 := - \int_{\R^3}  (\ued - u ) \otimes (\ued -u) : \na_x u \dd x, \\
        &J_3 := {\eps} \int_{\T^3 \times \R^3}  \fed \left(\frac{v}{\sigma}-\ued\right) \cdot G  \dd x \dd v , \\
        &J_4 := {\eps} \int_{\T^3}  \rhoed (\ued -u ) \cdot G \dd x,
    \end{align*}
where $G = {\na_x p - \Delta_x u}$.
\end{lemma}

Following the approach leading to Theorem~\ref{thCVprecis-3}, we also obtain quantitative convergence estimates similar to those of Theorems~\ref{thCVprecis-12} and~\ref{corollaryCVf-12}.


\section{Further developments}
\label{sec-further}

The purpose of this section is to describe applications of the methods developed in this paper for other high friction limits of Vlasov-Navier-Stokes type systems. Namely, we discuss
\begin{itemize}
\item the case of dimension two (Section~\ref{sec-dim2});
\item a generalization for mixtures of high friction limits, for models in which fragmentation and coalescence effects are taken into account for the dispersed phase (Section~\ref{sec-mix});
\item a (short time) derivation of a Boussinesq-Navier-Stokes system from the Vlasov-Navier-Stokes system on the whole space where the gravity force is taken into account (Section~\ref{sec-BNS}).
\end{itemize}
We conclude with a short discussion of some selected open problems related to these questions.

\subsection{The case of dimension two}
\label{sec-dim2}

The methods developed for $\T^3 \times \R^3$ can be readily extended to handle the two dimensional case, that is for $\T^2 \times \R^2$, but with some simplifications. As the propagation of higher regularity for Navier-Stokes system is much more favorable in dimension two (see e.g. \cite[Theorem 3.7]{che-des-gal-gre}), 
no smallness condition is required for the initial fluid velocity $u^0$, that is Assumption~\ref{hypSmallDataNS} can be dispensed with for Theorems~\ref{thm1}, \ref{thm2} and~\ref{thm3} to hold in dimension two.

\subsection{High friction limits for mixtures}
\label{sec-mix}
In this section, we consider a two-phase Vlasov-Navier-Stokes system which models a bidispersed aerosol, in which two types of particles with different radii $r_1 >r_2$ coexist. Particles with radius $r_1$, described by a distribution function $f_1$, can fragment, resulting after the breakup in particles with smaller radius $r_2$, whose density is described by another distribution function denoted by $f_2$. 
    We consider the regime where the radii $r_1$ and $r_2$ are very small compared to the observation length, which corresponds to a \emph{fine particle} regime as introduced in Section~\ref{sec-scalings}. 
This results in the following system:

\begin{equation}
\label{VNS-twofluid}
\left\{
\begin{aligned}
&\partial_t\ueps+ (\ueps\cdot\na_{x})\ueps-\Delta_{x}\ueps+\na_{x}p=
    \frac{1}{\kappa}\left(j_{f^{(1)}_{\eps}}-\rho_{f^{(1)}_{\eps}}\ueps\right)+ \frac{1}{\eps}\left(j_{f^{(2)}_{\eps}}-\rho_{f^{(2)}_{\eps}}\ueps\right),\\
   &\div_{x}\ueps=0,\\
  &\partial_t f^{(1)}_{\eps}+v\cdot\na_{x}f^{(1)}_{\eps}+  \frac{1}{\kappa}\div_{v}\left[f^{(1)}_{\eps}( \ueps-v)\right]=-f^{(1)}_{\eps}, \\
   &\partial_t f^{(2)}_{\eps}+v\cdot\na_{x}f^{(2)}_{\eps}+ \frac{1}{\eps} \div_{v}\left[f^{(2)}_{\eps}( \ueps-v)\right]={f^{(1)}_{\eps}}.
\end{aligned}
\right.
\end{equation}
We consider the regime
$\kappa= \eps^\alpha$, for  $\alpha> 2/3$. 
Arguing as in Section~\ref{sec-formal}, the formal limit turns out to be the following 
 two-fluid Inhomogeneous Navier-Stokes system
\begin{equation}
\label{twofluidINS}
\left\{
\begin{aligned}
&\partial_t \rho^{(1)} + \div_x (\rho^{(1)} u) =-\rho^{(1)}, \\
&\partial_t \rho^{(2)} + \div_x (\rho^{(2)} u) =\rho^{(1)}, \\
&\partial_t ( (1+\rho^{(1)} + \rho^{(2)}) u) + \na_x ((1+\rho^{(1)} + \rho^{(2)}) u\otimes u) -\Delta_x u +  \nabla_x p = 0, \\
&\div_{x}u=0
\end{aligned}
\right.
\end{equation}

\begin{remark}
The model~\eqref{VNS-twofluid} is inspired from the augmented fluid-kinetic system formally studied by Benjelloun, Desvillettes and  Moussa \cite{ben-des-mou}:
\begin{equation}
\label{BDM}
\left\{
\begin{aligned}
&\partial_t f +v\cdot\na_{x}f +  \div_{v}\left[f( u-v)\right]=-f, \\
&\partial_t \rho + \div_x (\rho u) =\int_{\R^3} f \dd v, \\
&\partial_t ( (1+\rho) u) + \na_x ((1+\rho) u\otimes u) -\Delta_x u +  \nabla_x p = 2 \int_{\R^3}(v-u) f \dd v, \\
&\div_{x}u=0
\end{aligned}
\right.
\end{equation}
This system is formally derived in \cite{ben-des-mou} from~\eqref{VNS-twofluid} for the parameter $\kappa=1$.
\end{remark}

The two-phase Vlasov-Navier-Stokes system~\eqref{twofluidINS} is also endowed with an energy--dissipation structure. Defining in this case the energy and dissipation functionals as
\begin{align}
\label{energy-mix}
 \E_\eps(t)
                &:=   \frac{1}{2}\int_{\T^3\times\R^3}|v|^2 f^{(1)}_{\eps} (t,x,v)\dd x\dd v
                + \frac{1}{2}\int_{\T^3\times\R^3}|v|^2 f^{(2)}_{\eps} (t,x,v)\dd x\dd v
                \\
           \nonumber     &+   \frac{1}{2}\int_{\T^3}| \ueps(t,x)|^2\dd x,
                \end{align}
              \begin{align}
\label{dissipation-mix}
                \D_\eps(t)&:=\frac{1}{2} \int_{\T^3}|\na_x \ueps(t,x)|^2\dd x +  \frac{1}{\kappa}\int_{\T^3 \times \R^3}|v-\ueps |^2 f^{(1)}_{\eps} (t,x,v)\dd x\dd v \\
                \nonumber &+ \frac{1}{\eps}  \int_{\T^3 \times \R^3}|v-\ueps |^2 f^{(2)}_{\eps} (t,x,v)\dd x\dd v.
\end{align}
we have (formally) the identity
\[
\frac{\dd}{\dd t} \E_\eps + \D_\eps =0.
\]


We obtain the following convergence result. Let us emphasize that we need to restrict to the range of parameters $\alpha>2/3$ for the parameter $\kappa = \eps^\alpha$.

\begin{theorem}
\label{thm3-mix}
Let $\alpha>2/3$.
Let $(u_{\eps},f^{(1)}_{\eps}, f^{(2)}_{\eps})$ be a global weak solution associated to the initial condition $(u_{0,\eps},f^{(1)}_{0,\eps}, f^{(2)}_{0,\eps})$. We have the following convergence results.

  \noindent {\bf 1. Mildly well-prepared case.}   Under Assumption\footnote{We shall not write down the two-phase analogues of the various previous  assumptions which prevail here.}~\ref{hypGeneral}, there exist $\eps_0$, $\eta>0$, $M'>0$,  such that, if for all $\eps \in (0,\eps_0)$,
  \[
   \norme{f^{(1)}_{0,\eps}}_{\L^1(\R^3;\L^\infty(\T^3))} +  \norme{ f^{(2)}_{0,\eps} }_{\L^1(\R^3;\L^\infty(\T^3))} \leq \eta
  \]
  and 
  \begin{align*}
   \int_{\T^3\times\R^3}\frac{|v-u_{\eps,\eps,1}^0(x)|^p}{\kappa^{p-1}}f^{(1)}_{0,\eps}(x,v)\dd x\dd v + \int_{\T^3\times\R^3}\frac{|v-u_{\eps,\eps,1}^0(x)|^p}{\eps^{p-1}}f^{(2)}_{0,\eps}(x,v)\dd x\dd v\le M',
  \end{align*}
  then
  there exists $T>0$ such that, up to a subsequence, for all $t \in [0,T]$,
  \begin{equation}
  \label{conv-1-mix}
  W_1(f^{(1)}_{\eps}(t),\rho^{(1)}(t)\otimes\delta_{v= u(t)})\xrightarrow[\eps\to0]{}0, \quad   W_1(f^{(2)}_{\eps}(t),\rho^{(2)}(t)\otimes\delta_{v= u(t)})\xrightarrow[\eps\to0]{}0,
  \end{equation}
    and for all $t \in [0,T]$,
    \begin{equation}
    \label{conv-2-mix}
     \norme{u_{\eps}(t)-u(t)}_{\L^2(\T^3)} \xrightarrow[\eps\to0]{}0,
    \end{equation}
   where $(\rho^{(1)}, \rho^{(2)},u)$ satisfies the two-fluid Inhomogeneous Navier-Stokes system \eqref{twofluidINS}.
    
     \noindent {\bf 2. Well-prepared case.} Under the additional Assumption~\ref{hypSmallDataNS}  there exists $\eta'>0$ small enough such that if
    \begin{multline*}
\int_{\T^3\times\R^3}\left|v- \frac{\langle   j_{f^{(1)}_{0,\eps}}+j_{f^{(2)}_{0,\eps}}\rangle}{\langle \rho_{f^{(1)}_{0,\eps}}+ \rho_{f^{(2)}_{0,\eps}}\rangle}\right|^2 f^{(1)}_{0,\eps}\dd x\dd v 
                + 
                \int_{\T^3\times\R^3}\left|v- \frac{\langle   j_{f^{(1)}_{0,\eps}}+j_{f^{(2)}_{0,\eps}}\rangle}{\langle \rho_{f^{(1)}_{0,\eps}}+ \rho_{f^{(2)}_{0,\eps}}\rangle}\right|^2 f^{(2)}_{0,\eps}\dd x\dd v \\
            +\left|\langle\ueps\rangle -  \frac{\langle   j_{f^{(1)}_{0,\eps}}+j_{f^{(2)}_{0,\eps}}\rangle}{\langle \rho_{f^{(1)}_{0,\eps}}+ \rho_{f^{(2)}_{0,\eps}}\rangle}\right|^2            +   \int_{\T^3}|u_{0,\eps} -\langle u_{0,\eps} \rangle|^2\dd x
        \leq \eta,
    \end{multline*}
     the convergences~\eqref{conv-1-mix} and~\eqref{conv-2-mix} hold for any $T>0$.

\end{theorem}

\begin{preuve}

The arguments leading to the existence of a global weak solution are for instance provided in in~\cite[Annexe C]{ben-phd}.
We shall not provide the complete proof of Theorem~\ref{thm3-mix}, as the analysis follows the arguments developed in the previous Sections. We still use a bootstrap argument.

Let us only highlight the main new objects and ingredients which are needed.
First, as the system involves two Vlasov equations, two systems of characteristics curves have to be introduced.
\begin{definition}
The equations of characteristics read as
\begin{equation}
\label{charac-1}
\left\{
\begin{aligned}
            &\dot{X}_\eps^{(1)}(s;t,x,v)=V_\eps^{(1)}(s;t,x,v),\\
            &\dot{V}_\eps^{(1)}(s;t,x,v)=\frac{\ueps(s,X_\eps^{(1)}(s;t,v,x))-V_\eps^{(1)}(s;t,x,v)}{\kappa},\\
            &X_\eps^{(1)}(t;t,x,v)=x,\\
            &V_\eps^{(1)}(t;t,x,v)=v.
\end{aligned}
\right.
\end{equation}
and
\begin{equation}
\label{charac-2}
\left\{
\begin{aligned}
            &\dot{X}_\eps^{(2)}(s;t,x,v)=V_\eps^{(2)}(s;t,x,v),\\
            &\dot{V}_\eps^{(2)}(s;t,x,v)=\frac{\ueps(s,X_\eps^{(2)}(s;t,v,x))-V_\eps^{(2)}(s;t,x,v)}{\eps},\\
            &X_\eps^{(2)}(t;t,x,v)=x,\\
            &V_\eps^{(2)}(t;t,x,v)=v.
\end{aligned}
\right.
\end{equation}
\end{definition}
Thanks to the method of characteristics, we deduce the representation formulas
    \begin{equation}
      \label{eq-f1}
f^{(1)}_{\eps} (t,x,v) =e^{-t} e^{\frac{3t}{\kappa}} f^{(1)}_{0,\eps} ({X}_\eps^{(1)}(0;t,x,v), {V}_\eps^{(1)}(0;t,x,v)).
\end{equation}
and
    \begin{equation}
    \label{eq-f2}
    \begin{aligned}
f^{(2)}_{\eps} (t,x,v) &= e^{\frac{3t}{\eps}} f^{(2)}_{0,\eps} ({X}_\eps^{(2)}(0;t,x,v), {V}_\eps^{(2)}(0;t,x,v)) \\
&+ \int_0^t e^{\frac{3(t-s)}{\eps}}f^{(1)}_{\eps}(s,{X}_\eps^{(2)}(s;t,x,v), {V}_\eps^{(2)}(s;t,x,v)) \dd s.
\end{aligned}
\end{equation}

Furthermore, the conservation of total momentum reads here as
$$
\frac{\dd}{\dd t} \langle \rho_{f^{(1)}_{\eps}}+  \rho_{f^{(2)}_{\eps}} \rangle = 0, \quad \frac{\dd}{\dd t} \langle    j_{f^{(1)}_{\eps}}+  j_{f^{(2)}_{\eps}} + \ueps \rangle = 0,
$$
which leads to the following definition of the modulated energy:
 \begin{multline*}
                \mathscr{E}_\eps(t)
                =
                \frac{1}{2  }\int_{\T^3\times\R^3}\left|v- \frac{\langle   j_{f^{(1)}_{\eps}}+j_{f^{(2)}_{\eps}}\rangle}{\langle \rho_{f^{(1)}_{\eps}}+ \rho_{f^{(2)}_{\eps}}\rangle}\right|^2 f^{(1)}_{\eps}\dd x\dd v \\
                + 
                \frac{1}{2}\int_{\T^3\times\R^3}\left|v- \frac{\langle   j_{f^{(1)}_{\eps}}+j_{f^{(2)}_{\eps}}\rangle}{\langle \rho_{f^{(1)}_{\eps}}+ \rho_{f^{(2)}_{\eps}}\rangle}\right|^2 f^{(2)}_{\eps}\dd x\dd v \\
            +\frac{\langle \rho_{f^{(1)}_{\eps}}+ \rho_{f^{(2)}_{\eps}}\rangle}{2(1+ \langle \rho_{f^{(1)}_{\eps}}+ \rho_{f^{(2)}_{\eps}}\rangle)}\left|\langle\ueps\rangle -  \frac{\langle   j_{f^{(1)}_{\eps}}+j_{f^{(2)}_{\eps}}\rangle}{\langle \rho_{f^{(1)}_{\eps}}+ \rho_{f^{(2)}_{\eps}}\rangle}\right|^2            +   \frac{1}{2}\int_{\T^3}|\ueps -\langle\ueps\rangle|^2\dd x.
            \end{multline*}

We have the following two-phase variants of Lemmas~\ref{ModulatedEnergyDissipation} and~\ref{lemmeDE}.
\begin{lemma}\label{ModulatedEnergyDissipation-mix}
    Under Assumption~\ref{hypGeneral}, for every $\eps>0$, for almost every $0\le s\le t$ (including $s=0$),
    \[
        \mathscr{E}_\eps(t)+\int_s^t\D_\eps(\tau)\dd\tau\le\EEE_\eps(s).
    \]
\end{lemma}

\begin{lemma}\label{lemmeDE-mix}
    Under Assumption \ref{hypGeneral}, if $T>0$ is such that $\rho_\eps^{(i)}\in\L^\infty((0,T)\times\T^3)$ for $i=1,2$, then there is $C>0$ (independent of $\eps$) such that
    \begin{equation}\label{DE-mix}
        \forall t\in[0,T],\qquad  \mathscr{E}_\eps(t)\le C e^{-\lambda_{\egd} t} \mathscr{E}_\eps(0),
    \end{equation}
    with
    \[
        \lambda_\eps=\min\left(\frac{c_P}{\kappa\left(c_P+4\left(\norme{\rho^{(1)}_\eps}_{\L^\infty((0,T)\times\T^3)}+\norme{\rho^{(2)}_\eps}_{\L^\infty((0,T)\times\T^3)}\right)\right)},\frac{c_P}{2}\right),
    \]
    where $c_P$ is a universal constant.
\end{lemma}

The following is a slight generalization of Lemma~\ref{changeVarV}.
\begin{lemma}\label{changeVarV-mix}
    Fix $c_*>0$ such that $c_* e^{c_*}<1/9$. Then, for any $t\in\R_+$ satisfying
    \begin{equation}\label{hyp-na-u-mix}
        \norme{\na_x u_\eps}_{\L^1(0,t;\L^\infty(\T^3))}\le c_*,
    \end{equation}
    for $i=1,2$, for any $s \in [0,t]$, $x\in\T^3$, the map
    \[
        \Gamma^{(i)}_{\eps, s;t,x}:v\mapsto V_\eps^{(i)}(s;t,x,v)
    \]
    is a $\CCC^1$-diffeomorphism from $\R^3$ to itself and satisfies
    \[
        \forall v\in\R^3,\qquad \det\D_v\Gamma^{(i)}_{\eps, t,x}(v)\ge\frac{e^{\frac{3t}{\eps}}}{2}.
    \]
\end{lemma}

We deduce the two-phase analogue of Lemma~\ref{RhoLinfini}.

\begin{lemma} Under the assumptions of Lemma~\ref{changeVarV-mix}, we have
\begin{align}
\label{rho1inf} \| \rho^{(1)}_\eps(t) \|_{\L^\infty(  \T^3)} &\lesssim e^{-t} \| f^{(1)}_{0,\eps} \|_{\L^1( \R^3; \L^\infty(\T^3))}, \\
\label{rho2inf} \| \rho^{(2)}_\eps \|_{\L^\infty((0,T) \times \T^3)} &\lesssim  \| f^{(1)}_{0,\eps} \|_{\L^1( \R^3; \L^\infty(\T^3))} +  \| f^{(2)}_{0,\eps} \|_{\L^1( \R^3; \L^\infty(\T^3))}. 
\end{align}

\end{lemma}

\begin{preuve}
The first estimate~\eqref{rho1inf} is straightforward, using~\eqref{eq-f1} and the change of variables of Lemma~\ref{changeVarV-mix} for $i=1$.
For the second one, using~\eqref{eq-f2}, thanks to Lemma~\ref{changeVarV-mix}, we write
\begin{align*}
\Big| \int_{\R^3} &\int_0^t  e^{\frac{3(t-s)}{\eps}} f^{(1)}_{\eps}(s,{X}_\eps^{(2)}(s;t,x,v), {V}_\eps^{(2)}(s;t,x,v)) \dd s \dd v  \Big| \\
&\lesssim \int_0^t  \int_{\R^3} f^{(1)}_{\eps}\left(s,{X}_\eps^{(2)}\left(s;t,x,[\Gamma^{(2)}_{\eps, s;t,x}]^{-1}(v)\right),v\right) \dd v \dd s \\
&\lesssim \int_0^t \int_{\R^3}  e^{-s} e^{\frac{3s}{\eps}} f^{(1)}_{0,\eps}\Big({X}_\eps^{(1)}\left(0;s,{X}_\eps^{(2)}(s;t,x,[\Gamma^{(2)}_{\eps,s;t,x}]^{-1}(v)),v\right), \\ &\qquad \qquad \qquad \qquad {V}_\eps^{(1)}\left(0;s,{X}_\eps^{(2)}(s;t,x,[\Gamma^{(2)}_{\eps,s;t,x}]^{-1}(v)),v\right)\Big)  \dd v \dd s.
\end{align*}
By a  slight variant of Lemma~\ref{changeVarV-mix}, under the same assumptions, the map 
  \[
        \widetilde\Gamma_{\eps, s, t,x}:v\mapsto {V}_\eps^{(1)}\left(0;s,{X}_\eps^{(2)}(s;t,x,[\Gamma^{(2)}_{\eps,s;t,x}]^{-1}(v)),v\right)
    \]
    is a $\CCC^1$-diffeomorphism of $\R^3$ and satisfies
    \[
        \forall v\in\R^3,\qquad \det\D_v\widetilde\Gamma_{\eps, s,x}(v)\ge\frac{e^{\frac{3s}{\eps}}}{2}.
    \]
We deduce that
\begin{align*}
\Big| \int_{\R^3} &\int_0^t  e^{\frac{3(t-s)}{\eps}} f^{(1)}_{\eps}(s,{X}_\eps^{(2)}(s;t,x,v), {V}_\eps^{(2)}(s;t,x,v)) \dd s \dd v  \Big| \\
&\lesssim \int_0^t \int_{\R^3}  e^{-s} f^{(1)}_{0,\eps}\left(Y_\eps(s,t,x,v),v\right)  \dd v \dd s \\ &\lesssim \| f^{(1)}_{0,\eps} \|_{\L^1( \R^3; \L^\infty(\T^3))}.
\end{align*}
We finally obtain~\eqref{rho2inf}, which completes the proof of the lemma.

\end{preuve}

  Let $p$ be the regularity index of Assumption~\ref{hypGeneral}. We now introduce the higher dissipation functionals
\begin{equation*}
\D_{1,\eps}^{(p)}(t) := \int_{\T^3\times\R^3}f^{(1)}_{\eps}(t,x,v)\frac{|v-\ueps(t)|^p}{\kappa^{p}}\dd x\dd v,
\end{equation*}
\begin{equation*}
\D_{2,\eps}^{(p)}(t) :=  \int_{\T^3\times\R^3}f^{(2)}_{\eps}(t,x,v)\frac{|v-\ueps(t)|^p}{\eps^{p}}\dd x\dd v.
\end{equation*}

We then obtain the two-phase analogue of Lemma~\ref{fine-key}.

\begin{lemma}
For all $T>0$,
       \begin{multline*}
        \int_0^T\D_{1,\eps}^{(p)}(t)\dd t=-\frac{1}{p}\left[\int_{\T^3\times\R^3}f^{(1)}_{\eps}(t,x,v)\frac{|v-\ueps(t)|^p}{\kappa^{p-1}}\dd x\dd v\right]_0^T\\
        -\int_0^T\int_{\T^3\times\R^3}f^{(1)}_{\eps} (t,x,v)(\partial_t\ueps+(\na_x\ueps)v)\cdot\frac{(v-\ueps(t,x))|v-\ueps(t,x)|^{p-2}}{\kappa^{p-1}}\dd x\dd v\dd t \\
        - \frac{1}{p} \frac{1}{\kappa^{p-1}} \int_0^T f^{(1)}_{\eps} |v-u(t,x)|^p \dd t
           \end{multline*}
      \begin{multline*}
        \int_0^T\D_{2,\eps}^{(p)}(t)\dd t=-\frac{1}{p}\left[\int_{\T^3\times\R^3}f^{(2)}_{\eps}(t,x,v)\frac{|v-\ueps(t)|^p}{\eps^{p-1}}\dd x\dd v\right]_0^T\\
        -\int_0^T\int_{\T^3\times\R^3}f^{(2)}_{\eps} (t,x,v)(\partial_t\ueps+(\na_x\ueps)v)\cdot\frac{(v-\ueps(t,x))|v-\ueps(t,x)|^{p-2}}{\eps^{p-1}}\dd x\dd v\dd t \\
        +  \frac{1}{p} \frac{1}{\eps^{p-1}} \int_0^T f^{(1)}_{\eps} |v-u(t,x)|^p \dd t.
    \end{multline*}
\end{lemma}
We emphasize that the assumption $\alpha>2/3$ is used in order to absorb the term 
$$
 \frac{1}{p} \frac{1}{\eps^{p-1}} \int_0^T f^{(1)}_{\eps} |v-u(t,x)|^p \dd t
$$
in the right-hand side of the above inequality.  Since $\alpha > 2/3$, we can always take $p>3$ close enough to $3$ and find $\gamma \in (0, 1)$ so that 
$$
 \frac{1}{\eps^{p-1}} \int_0^T f^{(1)}_{\eps} |v-u(t,x)|^p \dd t \leq \eps^\gamma  \int_0^T\D_{1,\eps}^{(p)}(t)\dd t.
 $$
We deduce the key estimate
\begin{align*}
 &\int_0^T\D_{1,\eps}^{(p)}(t)\dd t +   \int_0^T\D_{2,\eps}^{(p)}(t)\dd t 
 + \int_{\T^3\times\R^3}f^{(1)}_{\eps}(t,x,v)\frac{|v-\ueps(t)|^p}{\kappa^{p-1}}\dd x\dd v \\
 &\qquad \qquad+ \int_{\T^3\times\R^3}f^{(2)}_{\eps}(t,x,v)\frac{|v-\ueps(t)|^p}{\eps^{p-1}}\dd x\dd v 
 \lesssim  1 \\
&+\left| \int_0^T\int_{\T^3\times\R^3}f^{(1)}_{\eps} (t,x,v)(\partial_t\ueps+(\na_x\ueps)v)\cdot\frac{(v-\ueps(t,x))|v-\ueps(t,x)|^{p-2}}{\kappa^{p-1}}\dd x\dd v\dd t \right| \\
&+ \left| \int_0^T\int_{\T^3\times\R^3}f^{(2)}_{\eps} (t,x,v)(\partial_t\ueps+(\na_x\ueps)v)\cdot\frac{(v-\ueps(t,x))|v-\ueps(t,x)|^{p-2}}{\eps^{p-1}}\dd x\dd v\dd t\right|. 
\end{align*}
  At this stage of the proof, arguing as for the fine particle regime, 
 we deduce straightforward variants of Lemmas~\ref{lemmeExpressionD}, \ref{lemmeMixed}, \ref{EstimateHighOrderD}, \ref{lemmePara-3}, \ref{HighOrderD2}, \ref{lemmePara-3}. Using Lemma~\ref{LienFD-3}, we may conduct a bootstrap argument, yielding the required uniform estimates to eventually pass to the limit, arguing as in the proof of Theorem~\ref{FineFirstConv}.
 
 \end{preuve}
 

  \subsection{Derivation of a Boussinesq-Navier-Stokes system}
 \label{sec-BNS}
 
Consider the following Vlasov-Navier-Stokes system, set on $\R^3 \times \R^3$, in which the gravity force $g= -\overline{g} e_3$, where $\overline{g}>0$ is the gravitational constant, is taken into account:
\begin{equation}
    \label{VNS-g}
    \left\{
\begin{aligned}
    &\partial_t u_\eps + (u_\eps  \cdot\na_{x})u_\eps -\Delta_{x}u_\eps +\na_{x}p_\eps = j_{f_\eps }-\rho_{f_\eps} u_\eps ,\\
  &\div_{x}u_\eps =0,\\
    &\partial_t f_\eps +  v\cdot\na_{x}f_\eps + \frac{1}{\eps} \div_{v}\left[f_\eps (u_\eps -v+g)\right]=0, \\
   &\rho_{f_\eps} (t,x)=\int_{\R^3}f_\eps (t,x,v)\ddv,  \quad
     j_{f_\eps }(t,x)= \int_{\R^3}v f_\eps (t,x,v)\ddv.
\end{aligned}
\right.
\end{equation}
This is the exact analogue of the regime studied in \cite{hofe} for the Vlasov-Stokes case. The addition of gravity is important in view of the modeling of \emph{sedimentation} phenomena.

Arguing as in Section~\ref{sec-formal}, the formal limit as $\eps \to 0$ is the following Boussinesq-Navier-Stokes type system:
\begin{equation}
    \label{BNS}
    \left\{
\begin{aligned}
&\partial_t \rho + \div_x (\rho (u+g))=0, \\
&\partial_tu+ u\cdot\na_{x} u-\Delta_{x}u+\na_{x}p=\rho g, \\
&\div_x u =0,
\end{aligned}
\right.
\end{equation}
These equations slightly differ from the most standard form of the Boussinesq-Navier-Stokes equations (also referred to as Boussinesq without (thermal) diffusivity), in which the velocity field in the equation for $\rho$ is merely $u$ instead of $u+g$ as above. They form a classical geophysical model, see \cite{ped}.
Recently, in particular because of a formal resemblance with the incompressible Euler equations, it has attracted a lot of attention in the mathematical community, see e.g. \cite{hou-li,cha,dan-pai,hmi-rou,mas-sai-zha} and references therein. 

There are two main differences if one wants to apply the methods of this work:
\begin{itemize}
\item the equations are set in $\R^3$ instead of $\T^3$;
\item the addition of the gravity force results in a continuous injection of momentum and energy in the system.
\end{itemize}
Especially because of the gravity, we cannot expect the decay of the energy of the system (as obtain in the gravity-less case studied in \cite{hank}): therefore our methods do not allow us to obtain convergence results for arbitrarily large times.

However, all the other aspects of the proofs (maximal parabolic estimates, desingularizations...) still apply with straightforward modifications. Therefore, we can adapt the analysis of Section~\ref{sec-smallT} to obtain a \emph{short time} derivation of the Boussinesq-Navier-Stokes equations~\eqref{BNS}.  Details of the proof are left to the reader. In the following statement, one has to replace $\T^3$ by $\R^3$ when necessary.

\begin{theorem}
\label{thm-gravity} Let $(u_{\eps},f_{\eps})$ be a global weak solution associated to an initial condition $(u^0_{\eps},f^0_{\eps})$. 
 Under Assumption~\ref{hypGeneral}, there exists $T>0$ such that 
  \begin{equation*}
    \int_0^T W_1(f_{\eps}(t),\rho(t)\otimes\delta_{v= u(t)})\dd t \xrightarrow[\eps\to0]{}0
  \end{equation*}
    and for all $t \in [0,T]$,
    \begin{equation*}
     \norme{u_{\eps}(t)-u(t)}_{\L^2(\R^3)} \xrightarrow[\eps\to0]{}0,
    \end{equation*}
   where $(\rho,u)$ satisfies the Boussinesq-Navier-Stokes system \eqref{BNS}.
    
    If we further assume that
    \[
        \int_{\R^3\times\R^3}|v-u_{\eps}^0(x)|f_{\eps}^0(x,v)\,\dd x\dd v\xrightarrow[\eps\to0]{}0,
    \]
    then, for all $t \in [0,T]$,
    \begin{equation*}
        W_1(f_{\eps}(t),\rho(t)\otimes\delta_{v=u(t)}) \xrightarrow[\eps\to0]{}0.
    \end{equation*}
    
\end{theorem}

\subsection{Open problems}
\label{sec-open}
We end this section with a discussion of some open problems that are directly related to the high friction limits studied in this paper.


\noindent {\bf 1.} \emph{The case of rough data.} The present analysis requires to consider Fujita-Kato type data for the fluid part. One may wonder if it is possible to prove the limit from Leray solutions of~\eqref{VNS-general}  to Leray solutions of the limit Navier-Stokes equations~\eqref{TNS} or~\eqref{inNS}.


\noindent {\bf 2.} \emph{The case of bounded domains with absorption boundary conditions.} Let $\Omega$ be a smooth bounded domain of $\R^3$ and consider  the Vlasov-Navier-Stokes system~\eqref{VNS-general} set on $\Omega \times \R^3$ with the following absorbing boundary conditions:
\begin{equation*}
u(t,\cdot)|_{\pa\Omega} = 0, \quad f(t,\cdot)|_{\Sigma^-}=0
\end{equation*}
where, denoting by $n(x)$ the exterior normal at a point $x \in \pa\Omega$,
$$
\Sigma^-=\{ (x,v) \in \pa\Omega \times \R^3, \quad v\cdot n(x)<0\}.
$$
Global weak solutions to this system with these boundary conditions were built in \cite{bou-gra-mou}, while large time behavior was very recently studied in
\cite{ert-han-mou}, elaborating on the techniques of \cite{hank-mou-moy}. In this case, the right decaying functional in view of large time behavior is the kinetic energy itself (as opposed to the modulated energy for the torus case).
Building on these contributions, we conjecture that the results obtained in the torus extend to this setting, leading to the derivation of~\eqref{TNS} or \eqref{inNS} with the following boundary conditions:
$$
u(t,\cdot)|_{\pa\Omega} = 0, \quad \rho(t,\cdot)|_{\widetilde\Sigma^-(t)}=0
$$
where
$$
\widetilde\Sigma^-(t)-=\{ x \in \pa\Omega , \quad u(t,x)\cdot n(x)<0\}.
$$

\noindent {\bf 3.} \emph{The case of unbounded domains.}
Another natural extension of the results of this work concerns the case of smooth unbounded domains $\Omega$, still with absorbing boundary conditions on $\pa \Omega$. A first step would be the case $\Omega =\R^3$, for which large time behavior of the Vlasov-Navier-Stokes system is studied in \cite{hank}. It seems possible to extend the weighted in time estimates of \cite{hank} to handle as well high friction regimes for arbitrarily large times.

\noindent {\bf 4.} \emph{Global derivation of the Boussinesq-Navier-Stokes equations.}
The derivation of the Boussinesq-Navier-Stokes equations obtained in Section~\ref{sec-BNS} is a short time result. It is an open problem to extend it to arbitrarily large times. We remark that all global existence results for  Boussinesq-Navier-Stokes such as \cite{dan-pai,hmi-rou} do not come with decay estimates  for the velocity field $u$, which may even grow at most exponentially fast. This seems to indicate that the strategy of this work, which is based on a uniform control the velocity field in $\L^1_t \L^\infty$, is not adapted, and that other techniques should be devised.
On the other hand, for particular initial data, the recent work
\cite{mas-sai-zha} did prove better large time decay estimates, relying on additional mixing mechanisms. It would be interesting to investigate if this can be extended to  the study of~\eqref{VNS-g}.

\noindent {\bf 5.} \emph{General scalings for mixtures.} 
The study of the systems~\eqref{VNS-twofluid} is restricted to the range of parameters $\alpha>2/3$. We do not know if and how much this threshold can be improved.
We recall that the system that allows in \cite{ben-des-mou} to derive~\eqref{BDM} reads as
\begin{equation*}
\left\{
\begin{aligned}
   &\partial_tu+ (u\cdot\na_{x})u-\Delta_{x}u+\na_{x}p=
    \left(j_{f_1}-\rho_{f_1}u\right)+ \frac{1}{\eps}\left(j_{f_2}-\rho_{f_2}u\right),\\
    &\div_{x}u=0,\\
  &\partial_tf_1+v\cdot\na_{x}f_1+  \div_{v}\left[f_1( u-v)\right]=-f_1, \\
    &\partial_tf_2+v\cdot\na_{x}f_2+ \frac{1}{\eps} \div_{v}\left[f_2( u-v)\right]={f_1},
\end{aligned}
\right.
\end{equation*}
and thus corresponds to $\alpha =0$. We actually conjecture that the formal derivation towards~\eqref{BDM} is wrong in general, because the part involving $f^{(1)}_\eps$ in the expression~\eqref{eq-f2} for $f^{(2)}_\eps$ does not seem to concentrate in velocity, so that the term $\frac{1}{\eps}\left(j_{f_2}-\rho_{f_2}u\right)$ may be unbounded.


\appendix

\section{Sobolev and Besov (semi-)norms and Wasserstein-1 distance}

Let $k \in \Z^3$. The $k$-th Fourier coefficient of $f \in \L^1(\T^3)$ is given by
$$
c_k(f) := \int_{[0,2\pi]^3} f(x) e^{- i k\cdot x} \dd x.
$$
As usual, this definition is extended by duality to tempered distributions $f\in \mathcal{S}'(\T^3)$.
Let us first recall the very classical definition of inhomogeneous and homogeneous Sobolev spaces.
\begin{definition}
\label{def-sobolev}
Let $s \in \R$. The (inhomogeneous) Sobolev space $\H^s(\T^3)$ is defined by
$$
 \H^s (\T^3) := \left\{ f \in \mathcal{S}'(\T^3), \|u \|_{ \H^s (\T^3)}^2 :=\frac{1}{(2\pi)^3} \sum_{k \in \Z^3} (1+|k|^2)^{s} |c_k(f)|^2 < +\infty   \right\}.
$$
The homogeneous Sobolev space $\dot \H^s(\T^3)$ is defined by
$$
\dot \H^s (\T^3) := \left\{ f \in \mathcal{S}'(\T^3), \|u \|_{\dot \H^s (\T^3)}^2 := \frac{1}{(2\pi)^3} \sum_{k \in \Z^3} |k|^{2s} |c_k(f)|^2 < +\infty   \right\}.
$$
\end{definition}

We now provide a definition of  Besov spaces on the torus $\T^3$. We refer to \cite[Chapter 2]{bah-che-dan} for a reference on the topic.
The existence of the functions $\chi$ and $\varphi$ below is provided by \cite[Proposition 2.10]{bah-che-dan}
\begin{definition}
\label{def-besov}
Let $\mathcal{C} = \{ \xi \in \R^3, \, 3/4 \leq |\xi|\leq 8/3\}$. Let $\chi \in \mathcal{D}(B(0,4/3)), \varphi \in \mathcal{D}(\mathcal{C})$ be radial functions, with values in $[0,1]$ such that 
\begin{align*}
\forall \xi \in \R^3, \, \chi(\xi) + \sum_{j \geq 0 } \varphi(2^{-j} \xi) =1, \\
|j-j'|\geq 2 \implies \text{Supp } \varphi(2^{-j} \cdot) \cap \text{Supp } \varphi(2^{-j'} \cdot)  = \emptyset, \\
j\geq 1 \implies \text{Supp } \chi \cap \text{Supp } \varphi(2^{-j} \cdot)  = \emptyset.
\end{align*}
Consider the operators $\Delta_j$, $j \in \Z$, acting on $\mathcal{S}'(\T^3)$ defined by 
\begin{itemize}
\item $\Delta_j := 0$ for $j \leq -2$,
\item $\Delta_{-1} f := \chi(D)f$,
with
$
c_k(\chi(D)f) = \chi(k) c_k(f),
$
\item $\Delta_j f := \varphi(2^{-j} D) f$, with
$
c_k(\varphi(2^{-j} D) f)= \varphi(2^{-j} k) c_k(f).
$
\end{itemize}
Let $s \in \R, q,r \in [1,+\infty]$. The Besov space $\B^{s,q}_r(\T^3)$ is defined as
$$
\B^{s,q}_r(\T^3) := \left\{ f \in \mathcal{S}'(\T^3), \, \|u \|_{\B^{s,q}_r(\T^3)} := \left(\sum_{j \in \Z} 2^{r j s} \| \Delta_j f \|_{\L^q(\T^3)}^r \right)^{1/r}  \right\}.
$$
\index{$\B^{s,q}_r(\T^3)$: Besov space}
\end{definition}

To conclude this short section of functional analysis reminders, let us give the definition of the Wasserstein-1 distance $W_1$ that is intensively used in this work.

\begin{definition}
\label{def-W1}
Let $X$ be either $\T^3$ or $\T^3 \times \R^3$.
Let $\mu, \nu$ be two probability measures on $X$. The Wasserstein-1 distance between $\mu$ and $\nu$ is defined as 
$$
W_1(\mu,\nu) = \sup_{\| \nabla \psi \|_\infty \leq 1}  \int_{X} \psi \left( \dd\mu - \dd \nu\right). 
$$
\index{W@$W_1$: Wasserstein(-1) distance}
\end{definition}

 We recall that $W_1$ allows to metricize the weak-$\star$ convergence on $\mathcal{P}_1(X)$, the set of probability measures on $X$ with finite first moment.

\section{Gagliardo-Nirenberg interpolation estimates}

\begin{theorem}\label{th-Gagliardo-Nirenberg}
    Consider $1\le p,q,r\le\infty$ and $m\in\N$. Assume that $j\in\N$ and $\alpha\in\R$ satisfy
    \[
        \frac{1}{p}=\frac{j}{3}+\left(\frac{1}{r}-\frac{m}{3}\right)\alpha+\frac{1-\alpha}{q},
    \]
    \[
        \frac{j}{m}\le\alpha\le1,
    \]
    with the exception $\alpha<1$ if $m-j-3/r\in\N$. Then, the following holds. For any $g\in\L^q(\T^3)$, if $\D^mg\in\L^r(\T^3)$, then $\D^jg\in\L^p(\T^3)$ and we have the estimate 
    \[
        \norme{\D^jg}_{\L^p(\T^3)}\lesssim\norme{\D^mg}_{\L^r(\T^3)}^{\alpha}\norme{g}_{\L^q(\T^3)}^{1-\alpha}+\norme{g}_{\L^q(\T^3)},
    \]
    where the constant behind $\lesssim$ does not depend on $g$. If $\langle\D^jg\rangle=0$, then the term $\norme{g}_{\L^q(\T^3)}$ in the right-hand side can be dispensed with.
\end{theorem}

We introduce in the following corollary some notations we use in this paper.

\begin{corollary}\label{Gagliardo-Nirenberg}
    For any $p\ge2$, we set
    \[
        \alpha_p=\frac{5p-6}{7p-6},\qquad \beta_p=\frac{5p}{7p-6}.
    \]
    Let $g\in L^2(\T^3)$. We have
    \[
        \begin{array}{rcl}
            \forall p\ge2,\D^2g\in\L^p(\T^3)&\Longrightarrow&\norme{\na g}_{\L^p(\T^3)}\lesssim\norme{\Delta g}_{\L^p(\T^3)}^{\alpha_p}\norme{g-\langle g\rangle}_{\L^2(\T^3)}^{1-\alpha_p},\\
    
            \forall p>3,\D^2g\in\L^p(\T^3)&\Longrightarrow& \norme{\na g}_{\L^\infty(\T^3)}\lesssim\norme{\Delta g}_{\L^p(\T^3)}^{\beta_p}\norme{g-\langle g\rangle}_{\L^2(\T^3)}^{1-\beta_p}.
        \end{array}
    \]
\end{corollary}


\section{Estimates on the incompressible Navier-Stokes equations}\label{RappelsNS}

This brief section is dedicated to the presentation of some higher order energy estimates for the Navier-Stokes system with a source term:
\begin{equation}
\label{NS}
\left\{
\begin{aligned}
   &\partial_tu+(u\cdot\na_{x})u-\Delta_{x}u+\na_{x}p=F,\\
  &\div_{x}u=0,\\
     &u|_{t=0}=u^0.
\end{aligned}
\right.
\end{equation}

Let us first state an estimate on the solution of the heat equation.

\begin{lemma}\label{lemmaHeat}
    Let $u^0\in\dot\H^{\frac{1}{2}}(\T^3)$. For every $t\ge0$,
    \[
        \frac{1}{2}\norme{e^{t\Delta}u^0}_{\dot\H^{\frac{1}{2}}(\T^3)}^2+\int_0^t\norme{e^{s\Delta}u^0}_{\dot\H^{\frac{3}{2}}(\T^3)}^2\dd s\le\frac{1}{2}\norme{u^0}_{\dot\H^{\frac{1}{2}}(\T^3)}^2.
    \]
\end{lemma}


This result is one of the main ingredients in proving the following property of smooth solutions to the Navier-Stokes equations.

\begin{theorem}\label{propTpsLong}
    There exists a universal constant $C^*\in(0,1)$ such that the following holds. Consider $u^0\in \H^{1}_{\div}(\T^3)$, $F\in \L^2_{\loc}(\R_+;\H^{-\frac{1}{2}}(\T^3))$ and $T>0$ such that
    \begin{equation}\label{hypEE3}
        \int_0^T\norme{e^{t\Delta}u^0}_{\dot{\H}^1(\T^3)}^4\dd t+\int_0^T\norme{F(t)}_{\dot{\H}^{-\frac{1}{2}}(\T^3)}^2\dd t\le C^*.
    \end{equation}
    Then there exists on $(0,T)$ a unique Leray solution of \eqref{NS}. This solution belongs to $L^\infty(0,T;\dot{\H}^{\frac{1}{2}}(\T^3))\cap\L^2(0,T;\dot\H^{\frac{3}{2}}(\T^3))$ and satisfies, for a.e. $0\le t\le T$,
    \[
        \norme{u(t)}_{\dot\H^{\frac{1}{2}}(\T^3)}^2+\int_0^t\norme{u(s)}_{\dot\H^{\frac{3}{2}}(\T^3)}^2\dd s\lesssim\norme{u^0}_{\dot\H^{\frac{1}{2}}}^2+\int_0^t\norme{F(s)}_{\dot\H^{-\frac{1}{2}}(\T^3)}^2\dd s.
    \]
    Furthermore, if $F\in\L^2_{\loc}(\R_+;\L^2(\T^3))$, the solution belongs to $\L^\infty(0,T;\H^1(\T^3))\cap \L^2(0,T;\H^2(\T^3))$ and satisfies for a.e. $0\le t\le T$
    \begin{equation}\label{EEH1}
        \norme{u(t)}_{\H^1(\T^3)}^2+\int_0^t\norme{\Delta_{x}u(s)}_{\L^2(\T^3)}^2\dd s\lesssim\norme{u^0}_{\H^1(\T^3)}^2+\int_0^t\norme{F(s)}_{\L^2(\T^3)}^2\dd s.
    \end{equation}
\end{theorem}

\begin{preuve}
    We refer the reader to \cite[Theorem 10.1]{rob-rod-sad} for what concerns the estimate for $u$ in $\L^\infty(0,T;\dot{\H}^{\frac{1}{2}}(\T^3))\cap \L^2(0,T;\dot{\H}^{\frac{3}{2}}(\T^3))$. The additional assumption $u^0\in \H^1(\T^3)$ allows to follow the proof of \cite[Proposition 9.10]{hank-mou-moy} on $[0,T]$ and this yields the second part of the statement.
\end{preuve}




\section{Perturbation of the identity map}

We make use of the following version of the inverse function theorem.

\begin{theorem}\label{perturbationId}
    For $\Omega=\T^3$ or $\Omega=\R^3$, if $\phi:\Omega\rightarrow\Omega$ is $\CCC^1$ and satisfies $\norme{\nabla\phi}_\infty<1$, then $f=\Id+\phi$ is a $\CCC^1$-diffeomorphism of $\Omega$ to itself satisfying $\norme{\nabla f}_\infty\le(1-\norme{\nabla\phi}_\infty)^{-1}$. If furthermore $\norme{\nabla\phi}_\infty\le1/9$, then $\det\nabla f\ge1/2$.
\end{theorem}


\section{Maximal parabolic regularity}






Let us state the maximal parabolic regularity estimates for the Stokes equation on which we rely heavily in this work. The following follows from the classical \cite[Theorem 2.7]{gig-soh} (see also the expository notes \cite{sala}) and the transference argument for instance described in \cite[Corollary 9.8]{hank-mou-moy} .

\begin{theorem}\label{thReg}
    Let $q,r>1$. Let  $u^0\in\CCC^\infty(\T^3)$. Let $T>0$ and  $F\in{\L^q(\R_+^*, \L^r(\T^3))}$, and suppose $u$ is the unique tempered solution of the inhomogeneous Stokes equation
   \[
        \left\{\begin{aligned}
            &\partial_tu-\Delta_xu=F,\\
            &\div_x u=0, \\
            &u|_{t=0}=u^0.
       \end{aligned}\right.
   \]
   Then
   \[
        \norme{\partial_tu}_{\L^q(\R_+^*, \L^r(\T^3))}+  \norme{\Delta_xu}_{\L^q(\R_+^*, \L^r(\T^3))}\lesssim_{q,r}\norme{F}_{\L^q(\R_+^*, \L^r(\T^3))}+ \norme{u^0}_{\B^{s,r}_q(\T^3)},
    \]
    with $s= 2-\frac{2}{q}$.
\end{theorem}

We can now state the resulting maximal parabolic estimate for the Navier-Stokes equation, which follows from  Theorem~\ref{thReg}.

\begin{theorem}\label{ParabolicReg}
    Let $q,r>1$. Let $u^0\in\L^2_{\div}(\T^3)\cap\B^{s,r}_q(\T^3)$ for $s=2-2/q$ and $F\in{\L^q(\R_+^*, \L^r(\T^3))}$. Suppose $u$ is a Leray solution of the Navier-Stokes equations
    \[
        \left\{\begin{aligned}
            &\partial_tu+(u\cdot\na_x)u-\Delta_xu+\na_xp=F,\\
            &\div_xu=0,\\
            &u|_{t=0}=u^0,
        \end{aligned}\right.
    \]
    then
    \begin{multline*}
        \norme{\partial_tu}_{\L^q(\R_+^*, \L^r(\T^3))}+\norme{\Delta_xu}_{\L^q(\R_+^*, \L^r(\T^3))}\\
        \lesssim_{q,r} \norme{(u\cdot\na_x)u}_{\L^q(\R_+^*, \L^r(\T^3))}+\norme{F}_{\L^q(\R_+^*, \L^r(\T^3))}+\norme{u^0}_{\B^{s,r}_q(\T^3)},
    \end{multline*}
     with $s= 2-\frac{2}{q}$.
\end{theorem}





 \printindex

\bibliographystyle{abbrv}
\bibliography{references}

\end{document}